\documentclass{article}
\usepackage[english]{babel}

\usepackage[letterpaper,top=2cm,bottom=2cm,left=2cm,right=2cm,marginparwidth=1.75cm]{geometry}

\usepackage{float}
\usepackage{calrsfs}
\usepackage{amssymb}
\usepackage{amsmath}
\usepackage{amsthm}
\usepackage{subcaption}
\usepackage{multirow}
\usepackage{graphicx}
\usepackage{amsfonts}
\usepackage{mathtools}
\usepackage{mathrsfs}
\usepackage{placeins}
\usepackage{comment}
\usepackage{bbm}
\usepackage{dsfont}
\usepackage[export]{adjustbox}
\usepackage{overpic}

\usepackage{algorithm}
\usepackage{algpseudocode}

\usepackage{booktabs}
\newtheorem{remark}{Remark}

\usepackage[mathscr]{euscript}
\usepackage[dvipsnames,x11names]{xcolor}
\usepackage[colorlinks=true, allcolors=OliveGreen]{hyperref}
\usepackage{amsmath}   
\usepackage{tikz}
\usepackage{pgfplots}
\usetikzlibrary{external}
\usepackage{authblk}
\usepackage{array}
\newcolumntype{C}{>{\centering\arraybackslash}p{2.5cm}}
\usepackage{nicematrix}
\usepackage{caption}
\usepackage{subcaption}
\usepackage{tikz}

\usetikzlibrary{patterns,tikzmark}
\usetikzlibrary{matrix,decorations.pathreplacing,calc}
\usepackage{hf-tikz}

\usepackage{csvsimple}
\usepackage{siunitx}
\usepackage{tikz} 
\usepackage{pgfplots}
\pgfplotsset{compat=newest} 
\usepgfplotslibrary{units} 
\usepackage{tikz}
\usetikzlibrary{matrix,decorations.pathreplacing}
\pgfkeys{/pgfplots/legend entry/.code=\addlegendentry{#1}}

\usepackage{pdfpages}

\newcommand{\averagel}{\{\!\!\{}
\newcommand{\averager}{\}\!\!\}}

\newcommand{\jumpl}{[\![}
\newcommand{\jumpr}{]\!]}

\newcommand{\partition}{\mathcal{T}_h}

\usepackage[export]{adjustbox}
\usepackage{lipsum}                     
\usepackage{xargs}                      
\usepackage[colorinlistoftodos,prependcaption]{todonotes}
\newcommandx{\unsure}[2][1=]{\todo[linecolor=red,backgroundcolor=red!25,bordercolor=red,#1]{#2}}
\newcommandx{\change}[2][1=]{\todo[linecolor=blue,backgroundcolor=blue!25,bordercolor=blue,#1]{#2}}
\newcommandx{\info}[2][1=]{\todo[linecolor=OliveGreen,backgroundcolor=OliveGreen!25,bordercolor=OliveGreen,#1]{#2}}
\newcommandx{\improvement}[2][1=]{\todo[linecolor=Plum,backgroundcolor=Plum!25,bordercolor=Plum,#1]{#2}}
\newcommandx{\thiswillnotshow}[2][1=]{\todo[disable,#1]{#2}}

\DeclareMathAlphabet{\mathcalligra}{T1}{calligra}{m}{n}

\graphicspath{{Images/}}

\tikzset{  font={\fontsize{15pt}{12}\selectfont}}

\title{A $p$-adaptive polytopal discontinuous Galerkin method for high-order approximation of brain electrophysiology\footnote{\textbf{Funding}:
PFA is partially funded by the European Union (ERC SyG, NEMESIS, project number 101115663). Views and opinions expressed are however those of the authors only and do not necessarily reflect those of the European Union or the European Research Council Executive
Agency. 
PFA is partially supported by ICSC—Centro Nazionale di Ricerca in High Performance Computing, BigData, and Quantum Computing funded by European Union—NextGeneration EU. 
CBLS's has been funded by the National Recovery and Resilience Plan (NRRP), Mission 4, Component 1 – Investment 3.4, and Investment 4.1, funded by the European Union.
The present research is part of the activities of "Dipartimento di Eccellenza 2023-2027”. The authors are members of INdAM-GNCS. }}

\author[1]{Caterina B. Leimer Saglio}

\affil[1]{MOX-Dipartimento di Matematica, Politecnico di Milano, Piazza Leonardo da Vinci 32, Milan, 20133, Italy}

\author[1]{Stefano Pagani}

\author[1]{Paola F. Antonietti}

\begin{document}
\maketitle

\begin{abstract}
Multiscale mathematical models have shown great promise in computational brain electrophysiology but are still hindered by high computational costs due to fast dynamics and complex brain geometries, requiring very fine spatio-temporal resolution. This paper introduces a novel $p$-adaptive discontinuous Galerkin method on polytopal grids (PolyDG) coupled with Crank–Nicolson time integration to approximate such models efficiently. The $p$-adaptive method enhances local accuracy via dynamic, element-wise polynomial refinement/de-refinement guided by a-posteriori error estimators. A novel clustering algorithm automatizes the selection of elements for adaptive updates, further improving efficiency. A wide set of numerical tests, including epileptic seizure simulations in a sagittal section of a human brain stem, demonstrate the method’s ability to reduce computational load while maintaining the accuracy of the numerical solution in capturing the dynamics of multiple wavefronts. 
\end{abstract} 

\section{Introduction}\label{sec:1}
 Biophysical modeling and numerical simulation are pivotal in improving our understanding of biological processes. In silico approaches now offer tools to accelerate medical research and translate it into the clinic with novel personalized therapies \cite{viceconti2023position}.  Brain modeling constitutes a fundamental pillar of computational neuroscience, where mathematical laws, particularly partial differential equations, serve as essential tools for describing the spatiotemporal evolution of interconnected biophysical processes and the intricate interplay with the brain’s complex structures, \cite{goriely2015,mardal2022mathematical}.\\
 
 The field of computational brain electrophysiology is still advancing toward this objective: despite the availability of mathematical multiscale models and advanced numerical methods \cite{kager2007seizure,somjen2008computer,stefanescu2012computational,mardal2022mathematical,schreiner2022simulating,saetra2023neural,benedusi2024,berre2024cut,saglio2024high}, computational efficiency presents a significant barrier to broader application. Three factors significantly increase the computational burden: the spatial and temporal constraints associated with the cellular scale of the problem, the non-linearity of the coupled multiscale models, and the geometric complexity. Consequently, the number of degrees of freedom required to accurately discretize the corresponding mathematical problem when considering tissue- or organ-scale simulations increases dramatically. More precisely, the wavefronts of neuronal action potentials are characterized by steep and rapid profiles evolving on the millisecond time scale and localized at the synaptic level \cite{izhikevich2000neural,wei2014unification}, i.e., the connections between different neurons. The analysis of physiological and pathological phenomena must consider various scales that characterize neuronal activity, ranging from the membrane scale \cite{oyehaug2012dependence}, where chemical interactions take place across the cell membrane, to the entire brain. Furthermore, epileptic seizures can last for hours \cite{fisher2005epileptic}, making the problem multiscale in both space and time.\\ 

From a modeling perspective, computational neuronal electrophysiology relies on nonlinear ionic models derived from the Hodgkin–Huxley framework \cite{hodgkin1952propagation,cressman2009influence,barreto2011ion}, which describe the dynamics of transmembrane potential and the associated ionic concentrations. At the tissue and organ level, these ionic models are coupled with the monodomain or bidomain equations \cite{mardal2022mathematical,jaeger2022deriving} to model the propagation of transmembrane electric potential. Furthermore, the brain's geometry is highly complex, and its material properties are strongly anisotropic, which adds further challenges in the discretization step. The numerical approximation of the resulting coupled system of partial and ordinary differential equations typically leads to a large (monolithic or segregated) algebraic system, featuring a prohibitively large number of degrees of freedom.
These computational challenges are also common in cardiac electrophysiology, where a wide strand of discretization approaches are available in the literature, see, e.g., \cite{colli1990wavefront,vigmond2002,franzone2005simulating,goktepe2009computational,seemann2010,nagaiah2011,arthurs2012efficient,pezzuto2016space,quarteroni2017integrated,jaeger2022deriving}; we also refer to the comprehensive monographs \cite{franzone2014mathematical, quarteroni2019} and the reference therein.
In the specific context of high-order discretizations of electrophysiology models, we mention, e.g., spectral element methods
\cite{bordas2009simulation,cantwell2014high,africa2023matrix}, isogeometric analysis \cite{patelli2017isogeometric,bucelli2021multipatch,torre2022efficient,antonietti2025space}, Domain Decomposition schemes \cite{coudiere2019domain}, and high-order parallel-in-time approaches \cite{Pezzuto_2024}
High-order discontinuous Galerkin discretizations of electrophysiology models \cite{Hoermann_et_al_2018,Botta_et_al_2024} offer a highly flexible framework for implementing $h$-, $p$-, and $hp$-adaptivity strategies, thanks to the "local" nature of the discretization spaces. These methods enable enhanced accuracy while reducing the number of degrees of freedom and computational costs. In particular, discontinuous Galerkin methods formulated on polytopal meshes (PolyDG) \cite{Antonietti2013A1417,Bassi201245,cangiani2014hp,cangiani2022hp} naturally accommodate complex geometries, benefiting from the inherent flexibility of agglomerated grid structures. For brain tissue electrophysiology, in \cite{saglio2024high} we proposed a high-order discontinuous PolyDG method to handle inherent geometrical complexities. Nevertheless, the computational costs for these simulations remain substantial, as fine meshes and small time steps are required to currently capture the wavefronts. 
Given the inherited complexities of the problem at hand, $p$- \cite{hoermann2018adaptive} or $hp$-adaptive methods can provide an effective approach to lessen simulation costs while preserving the accuracy, see, e.g., \cite{hoermann2018adaptive,cherry2000space,schotzau2010posteriori} for triangular and  \cite{cangiani2016hp,Antonietti2013A1417} for polytopal meshes.\\

In this work, we propose a novel $p$-adaptive PolyDG method coupled with Crank–Nicolson time stepping for the numerical discretization of the monodomain equation coupled with the Barreto–Cressman neuronal ionic model. The proposed $p$-adaptive algorithm is driven by suitable a posteriori estimators and allows local and dynamic high-order accuracy in the numerical approximation of the solution while efficiently reducing the total number of degrees of freedom. 
To further enhance computational efficiency we also propose a novel technique based on $k$-means centroids for the automatic and dynamic detection of the subset of mesh elements candidate for refinement or coarsening of the local polynomial degree. We demonstrate the practical performance of the proposed $p$-adaptive method with extensive numerical experiments, including benchmark test cases and simulations of epileptic seizure propagation in a sagittal section of the human brainstem. Our numerical results confirm the capability of the proposed method to substantially reduce the total number of degrees of freedom, while preserving the accuracy required to capture such complex wavefront dynamics.
We point out that, while the proposed approach has been developed and validated within the framework of the monodomain equation coupled with the Barreto-Cressman model, its design is sufficiently general to be extended to other models of (brain) electrophysiology, including different ionic models or alternative macroscopic descriptions.  \\

The remainder of the manuscript is organized as follows:
In Section~\ref{sec:mathematicalmodel}, we introduce the mathematical model of brain electrophysiology under investigation. Section~\ref{sec:polyDG} presents the semi discrete formulation within the PolyDG framework and the fully discrete method based on employing the Crank–Nicolson time marching scheme. Section~\ref{sec::Padaptive} introduces and discusses in detail the proposed $p$-adaptive method, including a detailed discussion on implementation details. In Section~\ref{sec:NumericalResults}, we demonstrate the practical performance of the proposed $p$-adaptive method through a series of benchmark test cases.
Finally, in Section~\ref{sec:application} we illustrate the effectiveness of the proposed $p-$adaptive method in simulating physio-pathological scenarios in brain electrophysiology.

\section{The mathematical model}
\label{sec:mathematicalmodel}
\noindent In this section, we present the mathematical formulation of the brain electrophysiology model. Specifically, we consider the monodomain equation \cite{potse2006comparison,sundnes2006computational} coupled with suitable neuronal ionic models \cite{gutkin1998dynamics,barreto2011ion}.  

Given an open, bounded, polygonal domain $\Omega \in \mathbb{R}^d$, $(d=2,3)$ and a final time $T>0$, we introduce the transmembrane potential $u = u(\boldsymbol{x},t)$ with $u: \Omega \times [0,T] \rightarrow \mathbb{R}$, and the vector $\boldsymbol{y} = \boldsymbol{y}(\boldsymbol{x},t)$ with $\boldsymbol{y}: \Omega \times [0,T] \rightarrow \mathbb{R}^n, n\ge1,$ containing the ion concentrations and gating variables of the ionic model. The coupled multiscale problem reads as follows: 
\\
For any time $ t \in (0,T]$, find $u=u(\boldsymbol{x},t)$ and $\boldsymbol{y}=\boldsymbol{y}(\boldsymbol{x},t)$ such that:
\begin{equation}
    \label{eq:monodomain}
    \begin{dcases}
        \chi_m C_m  \frac{\partial u}{\partial t} - \nabla \cdot (\mathbf{\Sigma} \nabla u) +  \chi_m f(u,\boldsymbol{y}) = I_\mathrm{ext} & \mathrm{in} \; \Omega \times (0,T], \\
        \frac{\partial \boldsymbol{y}}{\partial t} + \boldsymbol{m}(u,\boldsymbol{y}) = \boldsymbol{0} &\mathrm{in} \; \Omega \times (0,T], \\
        \mathbf{\Sigma} \nabla u \cdot \boldsymbol{n} = 0  & \mathrm{on}\; \partial \Omega  \times (0,T],\\
        u(0) = u^0, \; \boldsymbol{y}(0) = \boldsymbol{y}^0 &\mathrm{in}\; \Omega.\\
    \end{dcases}
\end{equation}
Here, $\boldsymbol{\Sigma}$ represents the conductivity tensor, assumed to be constant in time and piecewise in space. Under the isotropic assumption, it is defined as $\boldsymbol{\Sigma}= \sigma \mathbbm{1}$, where $\sigma$ is the parameter defining the conductivity magnitude. %
In the first two test cases, we employ a fully isotropic conductivity tensor to focus only on the $p$-adaptive algorithm. Nonetheless, in the context of brain electrophysiology (test case 3), the 2D conductivity tensor is usually defined as:
\begin{equation}
\label{eq:sigma}
      \boldsymbol{\Sigma} = \sigma_{a} \boldsymbol{l}\otimes\boldsymbol{l} + \sigma_{n} \boldsymbol{n}\otimes\boldsymbol{n} = \sigma_{l} \mathbbm{1} + (\sigma_{n} - \sigma_{l})\boldsymbol{n}\otimes\boldsymbol{n},
\end{equation}
where $\sigma_{l}$ represents the conductivity along the principal axonal direction, $\sigma_{n}$ the conductivity in the orthogonal direction, $\boldsymbol{l} = \boldsymbol{l}(x)$ the axonal fibers direction in the point $x \in \Omega$ and $\boldsymbol{n} = \boldsymbol{n}(x)$ its normal.
We applied our $p$-adaptive algorithm to this context in the last test case.
\noindent In all the test cases, we assume homogeneous Neumann boundary conditions and, finally, we introduce suitable initial conditions $u^0$ and $\boldsymbol{y}^0$. 

For the modeling of ionic currents, we employ the Barreto-Cressman ionic model, a conductance-based model that describes a neuron's membrane potential and the dynamic interactions between intra- and extracellular ion concentrations~\cite{barreto2011ion}. Through parameters' variations, it can effectively represent neuronal bursting, a phase of epilepsy marked by rapid transmembrane potential spiking, as well as the classical action potential. A detailed derivation of this model can be found in \cite{cressman2009influence,barreto2011ion}. The variables of the models are three different ionic concentrations, namely intracellular sodium $s(t) =[\mathrm{Na}]_i(t)$, extracellular potassium $k(t) = [\mathrm{K}]_o(t)$, and intracellular calcium $c(t) = [\mathrm{Ca}]_i(t)$, and three gating variables $g^s,g^k,g^c$ that control the opening and closing of ion channels:
\begin{align}
    \label{eq:yBC}
      \boldsymbol{y}(\boldsymbol{x},t) &= 
      \left[\:
       s(\boldsymbol{x},t), \:
       k(\boldsymbol{x},t),\:
       c(\boldsymbol{x},t), \:
       g^s(\boldsymbol{x},t), \:
       g^k(\boldsymbol{x},t), \:
       g^c(\boldsymbol{x},t) \:
    \right]^\top,
\end{align}
with $\boldsymbol{y}(\boldsymbol{x},t): \Omega\times [0,T] \rightarrow \mathbb{R}^6$. The ionic reaction term of the model in Equation \eqref{eq:monodomain} is defined as follows:
\begin{equation}
    \label{eq:FBC}
        f(u,\boldsymbol{y}) = I_\mathrm{Na}(u,\boldsymbol{y}) + I_\mathrm{K}(u,\boldsymbol{y}) + I_\mathrm{Cl}(u,\boldsymbol{y}), 
\end{equation}
\noindent where the currents of sodium ($I_\mathrm{Na}$), potassium ($I_\mathrm{K}$), and chloride ($I_\mathrm{Cl}$) ions are:
\begin{equation*}
    \begin{aligned}
        I_\mathrm{Na}(u,\boldsymbol{y}) &= \left( \mathrm{G_{NaL}} + \mathrm{G_{Na}} (g^s(t))^3g^k(t) \right) (u(t) - \mathrm{E_{Na}}), \\
        I_\mathrm{K}(u,\boldsymbol{y}) &= \left( \mathrm{G_{K}}(g^c(t))^4 + \mathrm{G_{AHP}}\frac{c(t)}{1 + c(t)} + \mathrm{G_{KL}} \right) (u(t) - \mathrm{E_{K}}),\\
        I_\mathrm{Cl}(u,\boldsymbol{y}) &=  \mathrm{G_{CIL}}(u(t) - \mathrm{E_{Cl}}). 
    \end{aligned}
\end{equation*}
 Finally, the right-hand side of the ionic model is given by:
\begin{align*}
    \begin{aligned}
      {\boldsymbol{m}(u,\boldsymbol{y})} = 
      &\left[
    \begin{aligned}
       & \frac{c}{80} + \mathrm{G_{Ca}}\frac{0.002(u-\mathrm{E_{Ca}})}{1 + \exp\left(-\frac{25 + u}{2.5} \right)}  \\
       &   \frac{1}{\tau} \left( I_\mathrm{diff} - 14 I_\mathrm{pump} - I_\mathrm{glia} + 7 \gamma I_\mathrm{K} \right)  \\
       & \frac{1}{\tau} \left( \gamma I_\mathrm{Na} - 3I_\mathrm{pump} \right)  \\
       & 3 \frac{\boldsymbol g - \boldsymbol g_{\infty}}{ \tau_{g}}
    \end{aligned}
    \right].
\end{aligned}
\end{align*}

\par
For the weak formulation of the problem, we introduce the Sobolev space $V= H^1(\Omega)$, and we employ a standard definition of the scalar product in $L^2(\Omega)$, denoted by $(\cdot,\cdot)_{\Omega}$. The induced norm is denoted by $\|\cdot\|$. We remind that for vector-valued and tensor-valued functions, the definition extends component-wise \cite{salsa2022partial}.
\noindent Starting from Equation \eqref{eq:monodomain}, we introduce the ionic component and the dynamics component of the ionic model as:
\begin{equation}
    \begin{aligned}
       a(u,v) &= (\boldsymbol{\Sigma} \nabla u,\nabla v)_{\Omega}, & 
       r_{\text{ion}}(u,\boldsymbol{y},v) &= (f(u,\boldsymbol{y}),v)_{\Omega} &\forall \:  v \in V, \\\quad
       r_{\text{m}}(u,\boldsymbol{y},\boldsymbol{w}) &=  (\boldsymbol{m}(u,\boldsymbol{y}),\boldsymbol{w})_{\Omega} &\forall \:  \boldsymbol{w} \in [V]^n.
    \end{aligned}
    \label{eq:barreto_weak}
\end{equation}
We assume that the forcing terms, physical parameters, and initial conditions in Equation \eqref{eq:monodomain} are sufficiently regular, $i.e.$: $f(u,\boldsymbol{y}) \in L^2(0,T;L^2(\Omega))$, $I^\mathrm{ext}(u,\boldsymbol{y}) \in L^2(0,T;L^2(\Omega))$, $\boldsymbol{m}(u,\boldsymbol{y}) \in [L^2(0,T;L^2(\Omega))]^n$, $\chi_m$ and $C_m \in L_+^\infty(\Omega)$ where $L^{\infty}_+(\Omega):=\{v \in L^\infty(\Omega): v \ge 0 \text{ a.e. in } \Omega\}$, $u^0\in L^2(\Omega)$ and $\boldsymbol{y}^0\in [L^2(\Omega)]^n$.
\noindent The weak formulation of the problem \eqref{eq:monodomain} reads:
$ \forall \: t \in (0,T] $ find $u(t) \in V ,\boldsymbol{y}(t) \in V^n$ such that:
\begin{equation}
 \left\{ 
    \begin{aligned}
        &\chi_m C_m  \left(\frac{\partial u(t)}{\partial t},v\right)_{\Omega} + a(u(t),v) + \mathrm{\chi_m} r_{\text{ion}}(u(t),v) = (I^\text{ext},v)_\Omega  &\forall \:  v \in V, \\
        & \left(\frac{\partial \boldsymbol{y}(t)}{\partial t},\boldsymbol{w}\right)_{\Omega} + r_{\text{m}}(u(t),\boldsymbol{y}(t),\boldsymbol{w}) = 0  &\forall \:  \boldsymbol{w} \in V^n, \\
        &u(0) = u^0, \; \boldsymbol{y}(0) = \boldsymbol{y}^0 &\text{in } \Omega .
    \end{aligned}
    \right.
    \label{eq::weak_general}
\end{equation}

\section{PolyDG semi-discretization and fully-discrete formulation}
\label{sec:polyDG}
We now present the PolyDG semi-discrete formulation of the problem described in Equation \eqref{eq::weak_general}. Let \(\mathcal{T}_h\) represent a polytopal mesh partition of the domain \(\Omega\), consisting of disjoint elements \(K\). For each element \(K\), we define its diameter as \(h_K\) and set \(h = \max_{K \in \mathcal{F}_h} h_K < 1\). The interfaces are defined as the intersections of the \((d-1)\)-dimensional facets of neighboring elements. For \(d = 2\), the interfaces of a given element \(K\) are always segments. We denote by \(\mathcal{F}_h^I\) the union of all interior faces contained within \(\Omega\) and by \(\mathcal{F}_h^N\) those lying on the boundary \(\partial \Omega\).
In the following we assume that the underlying grid is  \emph{uniformly polytopal shape regular} in the sense of \cite{pietro2020hybrid}.\\
\noindent We define \(\mathbb{P}^p(K)\) as the space of polynomials of degree at most \(p \geq 1\) over the element \(K\) and the discontinuous finite element space as: 
\begin{equation*}
V_h^{\text{DG}} = \{v_h \in L^2(\Omega) : v_h|_{K} \in \mathbb{P}^{p_K}(K) \quad \forall \: K \in \mathcal{T}_h\},
\end{equation*}
with local degree $p_K\geq 1$ for any $K \in \mathcal{T}_h$.
Let \(F \in \mathcal{F}_h^I\) be the face shared by the elements \(K^{\pm}\), and let \(\boldsymbol{n}^{\pm}\) denote the normal unit vectors on the face \(F\) pointing outward to $K^\pm$. For  regular enough scalar-valued function \(v\) and a vector-valued function \(\boldsymbol{q}\), the trace operators are defined as follows \cite{arnold2002unified}:
\begin{equation*}
    \begin{aligned}
      &\averagel v \averager = \frac{1}{2} (v^+ + v^-) , \quad & \jumpl v  \jumpr& = v^+ \boldsymbol{n}^+ + v^- \boldsymbol{n}^-, \quad &\text{on $F \in \mathcal{F}_h^I$},&\\
      &\averagel \boldsymbol{q} \averager =  \frac{1}{2} (\boldsymbol{q}^+ + \boldsymbol{q}^-) , \quad & \jumpl \boldsymbol{q} \jumpr & = \boldsymbol{q}^+\cdot \boldsymbol{n}^+ + \boldsymbol{q}^-\cdot \boldsymbol{n}^-,  \quad &\text{on $F \in \mathcal{F}_h^{I}$}.&\\
    \end{aligned}
\end{equation*}
Here, the superscripts \(\pm\) indicate the traces of these functions on \(F\) taken in the interiors of \(K^{\pm}\), respectively. 
Moreover, we assume for the derivation of the formulation that the local polynomial degree vector $\boldsymbol{p} = \{p_K\}_{K \in \mathcal{T}_h}$ is given. 
This assumption allows us to give the following consistent definition of the penalization parameter \(\eta : \mathcal{F}_h^I \cup \mathcal{F}_h^D \rightarrow \mathbb{R}_+\):
\begin{equation}
    \eta = \eta(\boldsymbol{p},h,\boldsymbol{\Sigma}) = \eta_0 
    \begin{cases}
        \{\boldsymbol{\Sigma}_K\}_A \dfrac{\{p^2_K\}_A}{\{h_K\}_H} \quad \text{on } F \in \mathcal{F}_h^I, \\
        \boldsymbol{\Sigma}_K \dfrac{p^2_K}{h_K} \quad \text{on } F \in \mathcal{F}_h^D, \\
    \end{cases}
    \label{eq:eta}
\end{equation}
which depends explicitly on both the local degrees and the mesh size. 
In Equation \eqref{eq:eta} we set $\boldsymbol{\Sigma}_K = \|\sqrt{\boldsymbol{\Sigma}}|_K\|^2_{L^2(K)}$ and we consider both the harmonic average operator $\{\cdot\}_H$, that is $\{h\}_H = \frac{2h_+h_-}{h_+ + h_-}$, and the arithmetic average operator $\{\cdot\}_A$, that is $\{p\}_H = \frac{p_+ + p_-}{2}$  on $F \in \mathcal{F}_h^I$.  When applying our $p$-adaptive method (described in the following section), the polynomial degree vector $\boldsymbol{p} = \{p_K\}_{K \in \mathcal{T}_h}$ is updated at each time step. The penalization parameter $\eta$ is adjusted accordingly, while $\eta_0$ is chosen large enough to ensure stability. 
This setting enable us to introduce the following bilinear form $\mathcal{A}(\cdot,\cdot): V_h^{\text{DG}}\times V_h^{\text{DG}} \rightarrow \mathbb{R}$:
\begin{equation}
    \begin{aligned}
       \mathcal{A}(u,v) = \int_{\Omega} \mathbf{\Sigma}\nabla_h u \cdot  \nabla_h v \;dx+ \sum_{F \in \mathcal{F}_h^I} \int_F (\eta  \jumpl u \jumpr \cdot  \jumpl v  \jumpr - \averagel \boldsymbol{\Sigma} \nabla u \averager \cdot  \jumpl v\jumpr   - \jumpl u \jumpr \cdot \averagel \boldsymbol{\Sigma} \nabla v\averager) d\sigma \quad \forall \: u,v \in V^{\text{DG}}_h ,
    \end{aligned}
    \label{eq:coer}
\end{equation}
where $\nabla_h$ is the element-wise gradient. The semi-discrete formulation of problem in Equation \eqref{eq:monodomain} reads:
\vspace{1mm}\\
\noindent For any $t \in (0,T]$, find $(u_h(t),\boldsymbol{y}_h(t)) \in V^{\text{DG}}_h \times \left[V^{\text{DG}}_h\right]^n \text{ such that:}$
\begin{equation}
\label{eq:weak_barreto_2}
 \left\{ 
    \begin{aligned}
        \chi_m C_m  \left(\frac{\partial u_h(t)}{\partial t},v_h\right)_{\Omega}& +  \mathcal{A}(u_h(t),v_h)  + \chi_m\left( f(u_h(t),\boldsymbol{y}_h(t)),v_h\right)_{\Omega} =  (I_h^{\mathrm{ext}},v_h)_{\Omega}  &\forall \:  v_h \in V_h^{\text{DG}}, \\
         \left(\frac{d\boldsymbol{y}_h(t)}{d t},\boldsymbol{w}_h\right)_{\Omega} &+ \left( \boldsymbol{m}(u_h(t),\boldsymbol{y}_h(t)),\boldsymbol{w}_h\right)_{\Omega} =  0  &\forall \:  \boldsymbol{w}_h \in [V_h^{\text{DG}}]^n, \\
        u_h(0) = u_h^0,\;& \boldsymbol{y}_h(0) = \boldsymbol{y}_h^0 &\text{in } \Omega.
    \end{aligned}
    \right.
\end{equation}

As a consequence of our assumptions, we denote the dimension of the discrete space as $N_h(\boldsymbol{p})$, to make explicit its dependence on the vector of local polynomial distribution.
Let $N_h(\boldsymbol{p})$ be the dimension of $V_h^\mathrm{DG}$ and let $(\varphi_j)^{N_h(\boldsymbol{p})}_{j=0}$ be a suitable basis for $V_{h}^\mathrm{DG}$, then $u_h(t) = \sum_{j=0}^{N_h(\boldsymbol{p})} U_j(t)\varphi_j$ and $y_l(t) = \sum_{j=0}^{N_h(\boldsymbol{p})} Y^l_j(t)\varphi_j$ for all $l=1,...,n$. We denote $\mathbf{U} \in \mathbb{R}^{N_h(\boldsymbol{p})}$, $\mathbf{Y}_l \in \mathbb{R}^{N_h(\boldsymbol{p})}$ for all $l=1,...,n$ and $\mathbf{Y} = [\mathbf{Y}_1,...,\mathbf{Y}_n]^\top$. We define the matrices:
\begin{equation}
\small
    \begin{aligned}
      [\mathbf{M}]_{ij} &= (\varphi_i,\varphi_j)_{\Omega}, \; &\text{(Mass matrix),}& \quad  i,j = 1,...,N_h(\boldsymbol{p}) \\
      [\mathbf{F}]_{j} &= (I^\text{ext},\varphi_j)_{\Omega},  \; &\text{(Forcing term),}& \quad j = 1,...,N_h(\boldsymbol{p})\\
      [\mathbf{I}(u,\boldsymbol{y})]_{j} &= (f(u,\boldsymbol{y}),\varphi_j)_{\Omega},  \; &\text{(Non-linear ionic forcing term),}& \quad  j = 1,...,N_h(\boldsymbol{p})\\
      [\mathbf{G}_l(u,\boldsymbol{y})]_{j} &= (\boldsymbol{m}_l(u,\boldsymbol{y}),\boldsymbol{\varphi}_j)_{\Omega},  \; &\text{(Dynamics of the ionic model),}& \quad  j = 1,...,N_h(\boldsymbol{p}), \;l=1,...,n  \\ 
      [\mathbf{A}]_{ij} &= \mathcal{A}(\varphi_i,\varphi_j) &\text{(Stiffness matrix),}& \quad  i,j = 1,...,N_h(\boldsymbol{p}).
      \label{eq::matrixFull}
    \end{aligned}
\end{equation}
We are now ready to introduce the fully-discrete formulation.
We partition the interval \([0, T]\) into \(N\) sub-intervals \((t^{(k)}, t^{(k+1)}]\), each of length \(\Delta t\), such that \(t^{(k)} = k\Delta t\) for \(k = 0, \dots, N-1\). For temporal discretization, we adopt the second-order Crank-Nicolson scheme \cite{sundnes2007computing}, which handles the linear part (diffusion term) in a semi-implicit way, while the ion term in a fully explicit way. Given the initial conditions \(\mathbf{U}_0\) and \(\mathbf{Y}_0\), the discrete scheme is:
\\
Find $\mathbf{U}^{(k+1)} = \mathbf{U}(t^{(k+1)})$ and  $\mathbf{Y}^{(k+1)} = \mathbf{Y}(t^{(k+1)})$ for $k=0,...,N-1$, such that:
\begin{equation*}
\small
\left\{
\label{eq:Full_discrete_complete}
    \begin{aligned}
      & \left(\chi_m C_m \mathbf{M} + \frac{\Delta t}{2}\mathbf{A} \right) \mathbf{U}^{(k+1)}  =  \left(\chi_m C_m \mathbf{M} - \frac{\Delta t}{2}\mathbf{A} \right) \mathbf{U}^{(k)}  - \mathrm{\chi_m} \Delta t \mathbf{I}^{(k+1)} + \Delta t \mathbf{F}^{(k+1)} , \\
      &\mathbf{Y}^{(k+1)} = \mathbf{Y}^{(k)} - \Delta t \mathbf{G}^{(k)},\\
      &\mathbf{U}^0 = \mathbf{U}_0, \; \mathbf{Y}^0 = \mathbf{Y}_0,
    \end{aligned}
\right.
\end{equation*}
\noindent where the ionic current at step $k+1$ is computed as follows:
\begin{equation}
\begin{aligned}
  &\mathbf{I}^{(k+1)} = \frac{3}{2}\mathbf{I}^{(k)} - \frac{1}{2} \mathbf{I}^{(k-1)}.
\end{aligned}
\label{eq:Istim}
\end{equation}

\section{$p$-adative PolyDG framework}
\label{sec::Padaptive}
In practical applications, the solution of the problem in Equation \eqref{eq:monodomain} is characterized by action potentials with extremely steep and fast wavefronts that necessitate very fine spatial and temporal discretizations for their accurate numerical discretization. However, at a specific time instant, these traveling action potentials are localized in a specific portion of the domain, favoring adaptive numerical discretizations. These latter can accurately track the wavefront while removing unnecessary degrees of freedom, retaining the accuracy of the solution while reducing computational efforts.  
This section presents a spatial $p$-adaptive strategy based on the PolyDG method that locally adapts the polynomial order, using high-order discretization only in specific domain regions where wavefront tracking is needed. 
Unlike spatial adaptivity methods that rely on re-meshing, which demands significant computational resources and integration with a mesh generator, $p$-adaptivity avoids re-meshing during computation. An element-wise error indicator guides this flexibility in local polynomial degree selection \cite{cangiani2023posteriori} (described in Section \ref{API}).
The modification of the polynomial degree is then automatically driven by this local indicator and a suitable threshold numerically determined at the beginning of the simulation (see Section \ref{pformula}).
In Section \ref{palgo}, we present the complete algorithm that automatically uses higher-order discretizations only where needed to reduce the computational costs and reduce the total number of degrees of freedom while maintaining the accuracy required for accurate wavefront tracking. 

\subsection{A-posteriori indicator}\label{API}
In this Section, we discuss the \emph{a-posteriori} error indicator to determine/update the local polynomial order. We employ the comprehensive \emph{a-posteriori} indicator proposed in \cite{cangiani2014hp} and defined for each time $t$ and on each element $K \in \mathcal{T}_h$ as: 
\begin{align}
    \tau^2_{{K}}(\cdot, \cdot) = \tau^2_{{K},r}(\cdot, \cdot) + \tau^2_{{K},n}(\cdot)  + \tau^2_{{K},j}(\cdot)  + \tau^2_{{K},t}(\cdot) + \mathcal{O}^2_{{K},r},
    \label{eq::indicator}
\end{align}
where the components represent the norm of the residual ($\tau^2_{{K},r}$), the jump component of the normal gradient ($\tau^2_{{K},n}$), the jump component  ($\tau^2_{{K},j}$), and the tangential component of the gradient ($\tau^2_{{K},t}$) of the numerical solution, respectively. The local data oscillations on the element $K$, defined for the forcing term ($\mathcal{O}^2_{{K},r}$) is instead given by the error generated by the $L_2$-projection of the forcing term. In the following, we define all the terms appearing in Equation \eqref{eq::indicator} and refer to \cite{cangiani2020posteriori} for more details.

Let $u_h^{(k)} \in V_h^{DG}$ and $\boldsymbol{y}^{(k)}_h \in [V_h^{DG}]^n$ be the approximated solution at the time $t^{(k)}$, $k=1,2,\ldots $. Then, the first four components of the indicator \eqref{eq::indicator} are defined as
\begin{align}
  \tau_{{K},r}(u_h^{(k)}, \boldsymbol{y}^{(k)}_h) =& \|h_K \mathcal{R}(u_h^{(k)}, \boldsymbol{y}^{(k)}_h)\|_{L^2(K)},  \label{eq:indicator_res}\\
    \tau_{{K},n}({u_h^{(k)}}) =& \left (\|\sqrt{h_K}\jumpl (\boldsymbol{\Sigma} \nabla_h u_h^{(k)})  \cdot \boldsymbol{n}_K \jumpr \|^2_{L^2(\partial {K}) \cap \mathcal{F}_h^I} + \|\sqrt{h_K} (\boldsymbol{\Sigma} \nabla_h u_h^{(k)}  \cdot \boldsymbol{n}_K )\|^2_{L^2(\partial {K}) \cap \mathcal{F}_h^N}\right)^\frac{1}{2}, \label{fig:indicator_components_normal}\\
    \tau_{{K},j}({u_h^{(k)}})=& \|\sqrt{\eta}\jumpl u_h^{(k)} \jumpr \|^2_{L^2(\partial {K}) \cap \mathcal{F}_h^I},  \label{fig:indicator_components_jump}\\ 
    \tau_{{K},t}({u_h^{(k)}}) =& \|\sqrt{h_K}\jumpl (\boldsymbol{\Sigma} \nabla_h u_h^{(k)})  \cdot \boldsymbol{t}_K \jumpr \|_{L^2(\partial {K}) \cup \mathcal{F}_h^I},\label{fig:indicator_components_tang}
\end{align}
where $h_K$ is the diameter of $K\in \mathcal{T}_h$, $\Pi$ is the $L^2$ orthogonal projector onto the space of polynomials of degree $p_K$ on $K$, $\eta$ is the penalty coefficient defined in Equation \eqref{eq:eta}, and $\boldsymbol{n}_K$ and $\boldsymbol{t}_K$ are the unit normal and tangential vectors to $\partial K$. Moreover, in Equation \eqref{eq:indicator_res} the residual $\mathcal{R}(\cdot,\cdot)$ of the monodomain problem is given by
\begin{align*}
    & \mathcal{R}(u_h^{(k)}, \boldsymbol{y}^{(k)}_h) = \chi_m  f(u^{(k)}_h,\boldsymbol{y}^{(k)}_h) - \nabla_h \cdot (\boldsymbol{\Sigma} \nabla_h u^{(k)}_h) + \chi_m C_m \frac{u^{(k)}_h - u^{(k-1)}_h}{\Delta t}  - \Pi I_h^\text{ext}.
\end{align*}
Finally, for any $t^{(k)}$ and any $K\in \mathcal{T}_K$, the terms incorporating the oscillations with respect to the data are defined as follows: $\mathcal{O}_{{K},r} = \|h_K (\Pi I_h^\text{ext} - I_h^\text{ext}) \|_{L^2(K)}$. We remark that in the case of non-homogeneous Neumann boundary conditions, the indicator has to be suitably modified so as to take into account the oscillation term on the boundary datum. 

To assess the effectiveness of the \emph{a-posteriori} error indicator, we considered a benchmark test case for electrophysiology based solely on the monodomain equation with a nonlinear reaction term \cite{pezzuto2016space} defined by
$$
f(u) = a(u - V_{\text{rest}})(u - V_{\text{thres}})(u - V_{\text{depol}})
$$
where $V_\text{rest} \le  V_\text{thres} \le V_\text{depol}$ and $a>0$.
Under this assumption, we can find an analytical solution, as shown later in Section \ref{sec:test_case_1}. In this setting, we evaluate the error exploiting the energy norm that we derived in \cite{saglio2024high} and defined as follows: for any time $t$ :
\begin{equation}
\begin{aligned}
  \|v\|^2_{E} = \|v\|^2 + \int_{0}^t \frac{2\mu}{C_m \chi_m}\|v\|^2_{\mathrm{DG}}\mathrm{d}s + \int_0^t \frac{a}{C_m}\|v\|^4_{L^4(\Omega)}\mathrm{d}s \quad \forall v\in H^1(\partition).
\end{aligned}
\label{eq::energy_norm}
\end{equation}
In Equation \eqref{eq::energy_norm}, $\| \cdot \|_\mathrm{DG}^2$ is given by
\begin{equation}
    \begin{aligned}
       &\| {v} \|_\mathrm{DG}^2 =  \|\nabla_h {v}\|^2 +\|\eta^{\frac{1}{2}}\jumpl{v}\jumpr\|^2_{\mathcal{F}_h^\mathrm{I}} \qquad \forall v\in H^1(\mathcal{T}_h), 
    \end{aligned}
    \label{eq::DG_norm}
\end{equation}
where \(\|\cdot\| = \|\cdot\|_{L^2(\Omega)}\) and \(\|\cdot\|_{\mathcal{F}} = \left(\sum_{F \in \mathcal{F}} \|\cdot\|^2_{L^2(F)} \right)^{\frac{1}{2}}\). Here, $\eta$ is defined in Equation \eqref{eq:eta} and $H^s(\mathcal{T}_h)$ is the space of piecewise $H^s$ functions with $s\ge0$. 
For further details see \cite{saglio2024high}, where we provided the \emph{a-priori} stability, and the convergence analysis with respect to the energy norm defined in \eqref{eq::energy_norm} for the semi-discrete formulation of the problem in Equation \eqref{eq:weak_barreto_2}.

\subsubsection{Update of the polynomial degree} \label{pformula}
\begin{figure}[h]
    \centering
\begin{subfigure}[b]{0.3\textwidth}
\includegraphics[scale=0.8]{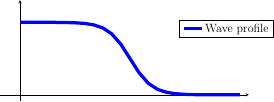}
\includegraphics[scale=0.8]{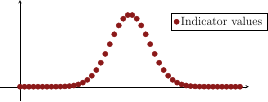}
\subcaption{\label{fig:cluster1}}
\end{subfigure}
\begin{subfigure}[b]{0.3\textwidth}
\includegraphics[scale=0.8]{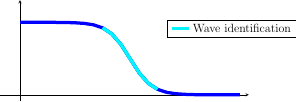}
\includegraphics[scale=0.8]{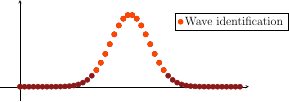}
\subcaption{\label{fig:cluster2}}
\end{subfigure}
\begin{subfigure}[b]{0.3\textwidth}
\includegraphics[scale=0.8]{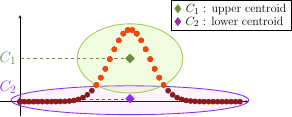}
\vspace{2em}
\subcaption{\label{fig:cluster3}}
\end{subfigure}
    \caption{Example of the clustering process for identifying the initial wave in a one-dimensional configuration. 
    \eqref{fig:cluster1} represents an example of a travelling wave and the corresponding distribution of the indicator $\tau_K$ defined as in Equation \eqref{eq::indicator}. \eqref{fig:cluster2} clustering of the local indicator values $\tau_K$, based on employing the $k$-means algorithm for the identification of the wavefront and the quiescence regions. \eqref{fig:cluster3} definition of the upper centroid ($C_1$) and lower centroid ($C_2$) resulting from the $k$-means clustering.}
    \label{fig:clustering}
\end{figure}

Let $u_h^{(k)} \in V_h^{DG}$ and $\boldsymbol{y}^{(k)}_h \in [V_h^{DG}]^n$ be the numerically approximated solution at time $t^{(k)}$, $k=1,2,\ldots$ and let $\tau_{{K}} =  \tau_{{K}}(u_h^{(k)}, \boldsymbol{y}^{(k)}_h)$ the local error indicator on $K \in \mathcal{T}_h$ defined as in \eqref{eq::indicator}. We can update the polynomial degree on each element $K \in \mathcal{T}_h$ exploiting the following function: 
\begin{align} 
    p_K^\text{arctan} = \left \lceil p^\text{max}\frac{2}{\pi}\arctan \left(\frac{\tau_{K}}{\tau_\text{threshold}} \right)\right \rceil,
    \label{eq:tanh}
\end{align}
where $\lceil \cdot \rceil$ is the ceiling function, $p^\text{max}$ is the maximum polynomial degree to be used in the simulation, and $\tau_\text{threshold}$ corresponds to a suitable threshold to be defined. 
To avoid large element-wise variations in the local polynomial degree, we opt for a smooth scaling, taking into account not only the information obtained from the evaluation of the adaptive function in Equation \eqref{eq:tanh} but also the polynomial degree $p_K^\text{old}$ at the previous time step. Therefore, for each element, we compute the new element-wise polynomial degree
$p_\text{new}$ as:
{
\begin{align} 
    p_K^\text{new} = \mathds{1}_{\{p_K^\text{old} \le p_K^\text{arctan}\}} \; \text{min}(p_K^\text{old} + 1 , p_K^\text{arctan}) +  \mathds{1}_{\{p_K^\text{old} > p_K^\text{arctan} \} } \; \text{max}(p_K^\text{old} - 1 , p_K^\text{arctan}).
    \label{eq:finalpformula}
\end{align}}
We remark that in the $p$-adaptive function defined in Equation \eqref{eq:tanh}, $\tau_\text{threshold}$ must be identified to select the subset of elements for which the polynomial degree should be updated (either increased or decreased). The selection of this threshold is typically subject to problem-specific tuning according to the mesh size of the underlying partition $\mathcal{T}_h$. 
To allow an automatic problem-specific selection of $\tau_\text{threshold}$, we propose a procedure relying once more on the physics of the problem. Since the action potential wavefront is localized in space, the error indicators should be subdivided into two clusters: respectively, the subsets of elements in close proximity to the action potential wavefront should be associated with higher error indicators, while the elements distant from the wavefront should correspond to lower error indicators.
For this reason, we employ a k-means algorithm to cluster the local error indicators into two subsets at each time step (see Figure~\ref{fig:clustering}). 

The clustering is performed only at the beginning of the simulation according to the value of the local indicator ($\tau_K$), obtained with an initial condition approximated using the maximum degree. Here $p_\text{max}$ is assigned to all the elements so that it is possible to straightforwardly reduce the polynomial degree far from the wavefront. 
More in detail, we consider two possible thresholds defined by exploiting the clustering centroids $C = [C_1 \; , \; C_2]$ as follows:
\begin{equation}
    \tau^\text{min}_\text{threshold} = \min (C_{1},C_{2}) \; , \; \qquad \tau^\text{mean}_\text{threshold} = \frac{C_{1} + C_{2}}{2}.
    \label{eq::centroid}
\end{equation}
Note that, in this way, the threshold for the $p$-adaptive strategy is fixed at the beginning of the simulation and it is not updated during time evolution. 

\subsection{$p$-adaptive algorithm} \label{palgo}
The polynomial degree is adapted locally and independently for each element using \eqref{eq:finalpformula}. 
To avoid an excessive computational overhead generated by the calculation of the local error indicator $\tau_{{K}}$ on each element $K\in\mathcal{T}_h$ at each time step, we design an updating strategy specifically tailored for the problem which involves the evolution of the action potential wavefront. The structure of the adaptive algorithm is summarized in the pseudo-code reported in Algorithm \ref{alg::algorithm} and detailed below. Specifically, for every time step $t^{(k)}$, the indicator $\tau_K(u_h^{(k)}, \boldsymbol{y}_h^{(k)})$ is not updated on all the mesh elements but only on a specific subset $\mathcal{S}^{(k)}$ of mesh elements defined as $\mathcal{S}^{(k)} = \mathcal{S}_K^{(k-1)} \cup \mathcal{S}_{\text{neigh}}^{(k-1)}$. 

\begin{figure}[h]
\centering
$\mathcal{S}^{(k-1)}$ \hspace{3cm} $\mathcal{S}^{(k-1)}_\text{neigh}$   \hspace{3cm} $\mathcal{S}^{(k)}$\\
\begin{subfigure}[t]{0.15\textwidth}
\vspace*{-5em}
\end{subfigure}
\begin{subfigure}[t]{0.8\textwidth}
\centering
\includegraphics[trim={15cm 0cm 15cm 0cm},clip,scale=0.08]{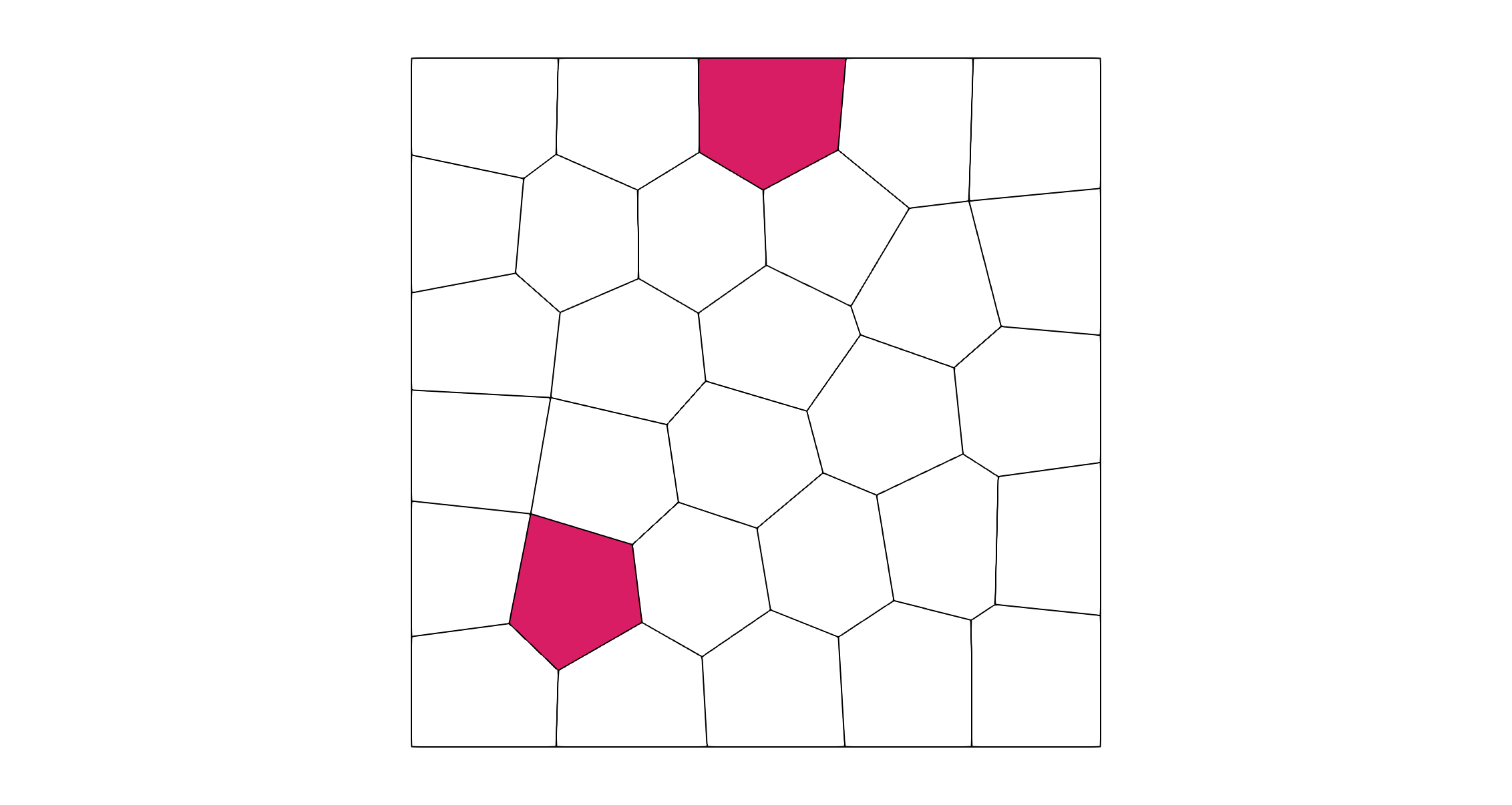}
\includegraphics[trim={15cm 0cm 15cm 0cm},clip,scale=0.08]{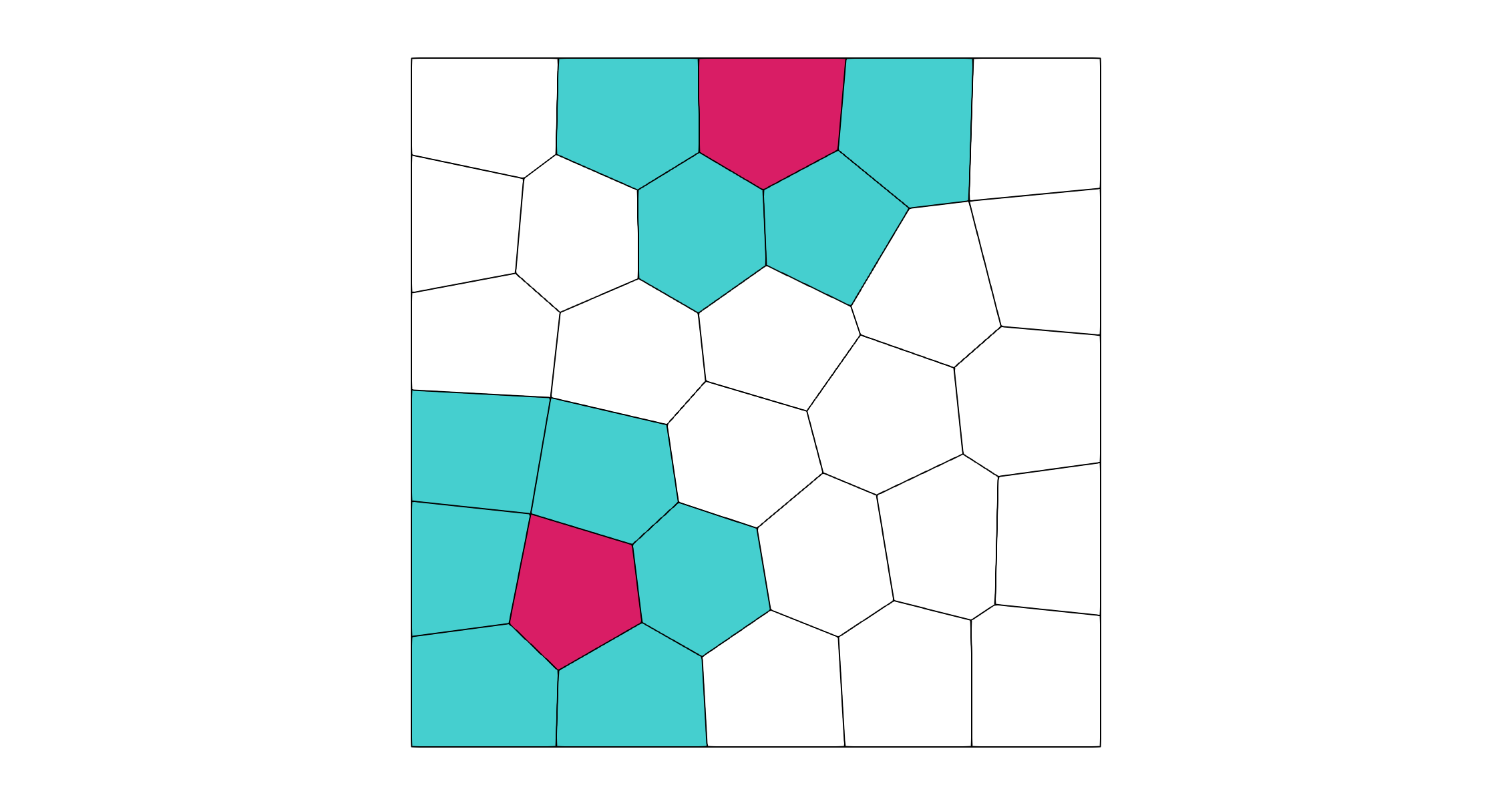}
\includegraphics[trim={15cm 0cm 15cm 0cm},clip,scale=0.08]{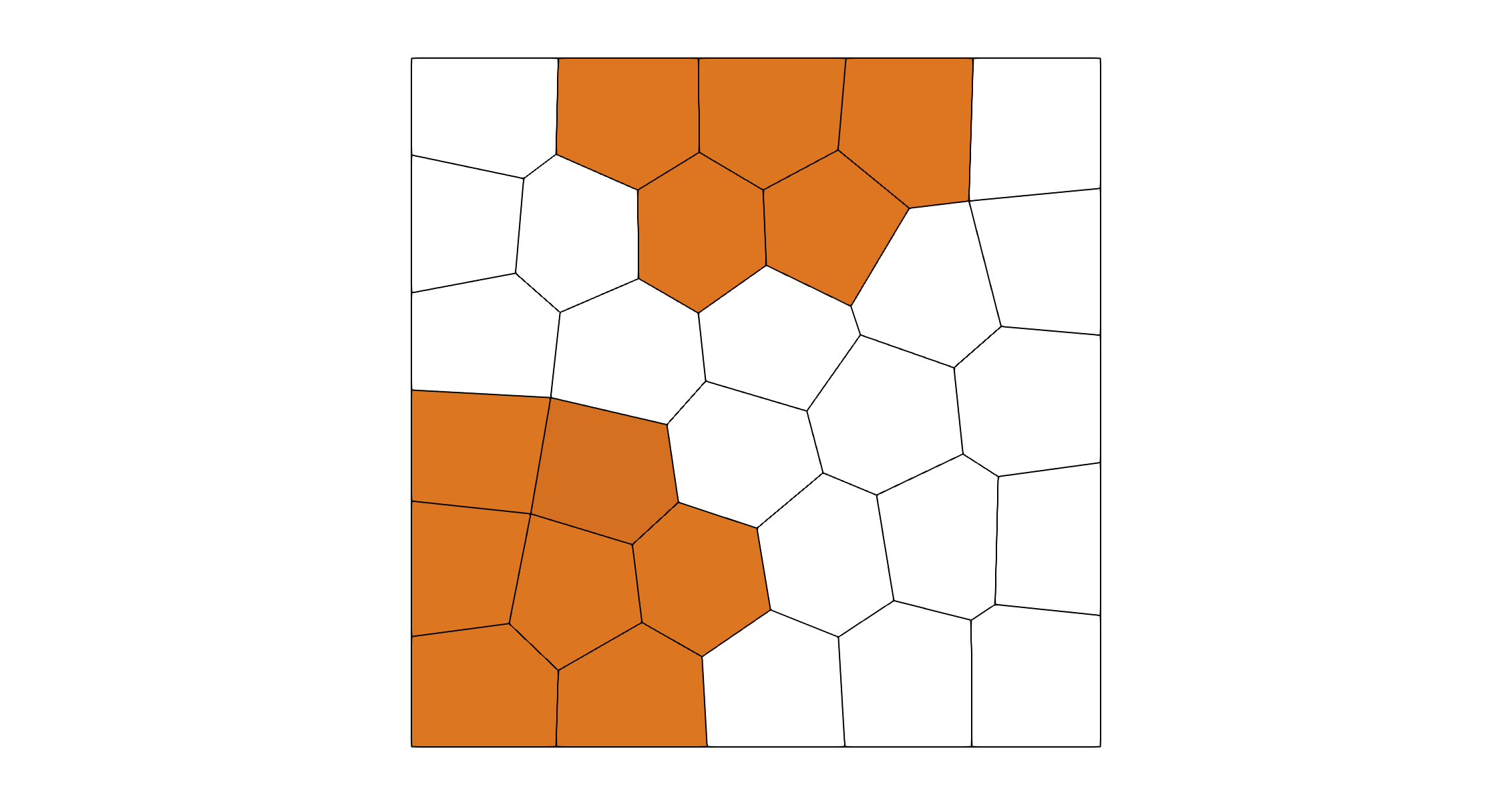}
\end{subfigure}
\caption{Example of subset of elements $\mathcal{S}^{(k)}$ where the indicator defined in \eqref{eq::indicator} is updated at the $k\!-\!th$ time step. In red are depicted the elements whose polynomial degree was either updated at the $(k-1)\!-\!th$ time step or was at the maximum allowed value; in blue are depicted the neighbors of such elements. Finally, the complete set of elements for computing the indicator is in orange.
} 
\label{fig::p-adapt}
\end{figure}

In Figure~\ref{fig::p-adapt}, we show how to construct $\mathcal{S}^{(k)}$. 
Specifically, we start from $\mathcal{S}^{(k-1)}$, which is defined as the union of the elements $K \in \mathcal{T}_h$ with a modified polynomial degree at the previous time step and the elements with maximum polynomial degree, i.e.,
\begin{align*}
\mathcal{S}^{(k-1)} = \{ K \in \mathcal{T}_h\, \:  |& \: p_K^{(k-1)} \textrm{ has been updated at the previous time step  } t^{(k-1)} \textrm{  or } p_K^{(k-1)}=p_\text{max}\}.
\end{align*}
Then, $\mathcal{S}^{(k-1)}_\text{neigh}$ is defined as the subset of neighboring elements to any $K \in \mathcal{S}^{(k-1)}$, i.e.,
\begin{equation*}
\mathcal{S}^{(k-1)}_\text{neigh}= \{ K \in \mathcal{T}_h \: | \:  K \: \textrm{is a neighboring element to } K^{\prime}\in \mathcal{S}^{(k-1)} \}.
\end{equation*}
At the beginning of the simulation, the initial configuration is $\mathcal{S}^{(0)} = \mathcal{T}_h$, that is all elements of the mesh are analyzed at the first iteration. Once the subset $\mathcal{S}^{(k)}$ has been defined, we update the indicator $\tau_K(u_h^{(k)}, \boldsymbol{y}_h^{(k)})$ for all the elements in $\mathcal{S}^{(k)}$, while for the other elements, we consider the indicator and the local polynomial degree at the previous time step.
\begin{algorithm}
  \caption{$p$-adaptive algorithm}\label{alg::algorithm}
  \begin{algorithmic}[1]
  \State \textbf{Input:} Initial conditions \{$u_h^{(0)}$, $\boldsymbol{y}_h^{(0)}$\}, Matrices: $\mathbf{A}, \mathbf{M}$ 
    
  \For {$k={1,\ldots,N}$}
       \State Compute $u_h^{(k)},\boldsymbol{y}_h^{(k)}$ the with the Crank-Nicholson scheme
        \If{$k=1$} 
            \State Exploit the $k$-means clustering and the compute threshold $\tau_\text{threshold}$.
        \EndIf 
       \State Compute $\mathcal{S}^{(k)}$
       \State Update the indicator $\tau_K(u_h^{(k)}, \boldsymbol{y}_h^{(k)})$ for all $K \in \mathcal{S}^{(k)}$, for all $K \not\in \mathcal{S}^{(k)}$ consider the indicator at the previous time-step
       \State Update the polynomial degree distribution $\boldsymbol{p}^k$ based on employing the adaptive function in Eq. \eqref{eq:finalpformula}
       \State Update the matrices $\mathbf{A}, \mathbf{M}$ with respect to the new degree distribution
       \State Update the forcing and ionic terms $\boldsymbol{I}_\text{ext}^{(k)}$, $\boldsymbol{I}_\text{ion}^{(k)}$ with respect to the new degree distribution
       \State Update the solution $u_h^{(k)},\boldsymbol{y}_h^{(k)}$ with respect to the new degree distribution
  \EndFor
\end{algorithmic}
\end{algorithm}
\begin{remark}
To further reduce the computational costs, we remark that the indicator $\tau_K(u_h^{(k)}, \boldsymbol{y}_h^{(k)})$ could also be updated not at every time step but only every $\bar{k}>1$ temporal iterations (here $\bar{k}$ should depend on the velocity of the wavefront and the choice of $\Delta t$ and $h$). 
\end{remark}

\subsection{Implementation details}
\begin{figure}[h!]
\centering
\vspace{-0.3cm}
    \begin{subfigure}[t]{\textwidth}
    \centering
    \includegraphics[scale=0.7]{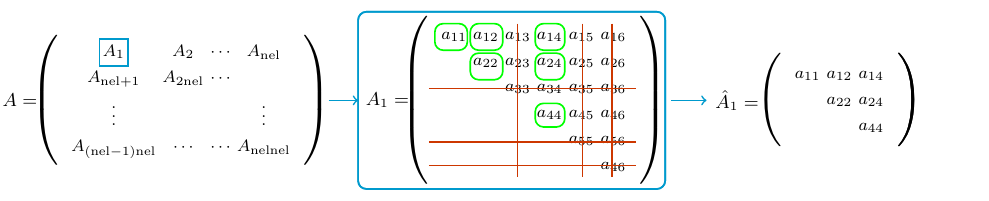}
    \caption{Dynamic update of the block-diagonal matrix $\mathbf{A}_1 \in \mathbf{R}^{6x6}$ with $\hat{\mathbf{A}}_1 \in \mathbf{R}^{3x3}$.}
    \end{subfigure}
    \begin{subfigure}[t]{\textwidth}
    \centering
\hspace{-1em}\includegraphics[scale=0.7]{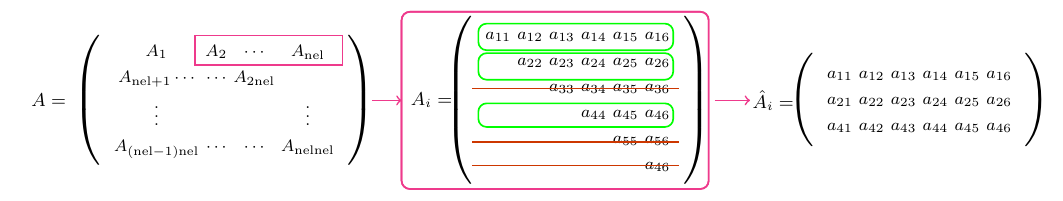}
    \caption{Dynamic update of the off-diagonal blocks  $\mathbf{A}_i \in \mathbf{R}^{6x6}$ with $\hat{\mathbf{A}}_i \in \mathbf{R}^{3x6}$, $i=1,\ldots,n_\text{el}n_\text{el}$.}
    \end{subfigure}\hfill
    \caption{Example of the dynamic update of the stiffness matrix $\mathbf{A}$ in two dimensions, supposing that the element $K_1$ is selected to decrease the local polynomial approximation degree from $p_{K_1}=2$ to $p_{K_1}=1$ whereas all the other mesh elements $p_K=2$.}
    \label{fig::examplematrix}
\end{figure}

From a practical implementation perspective, selecting an appropriate set of basis functions is essential for dynamically adapting the polynomial degree without reassembling the mass and stiffness matrices at each time step. In this work, we consider a set of hierarchical basis functions to avoid the reassembling phase at each time step due to the local adaptation algorithm.
The initial condition is constructed by employing $p_K=p^\text{max}$ for any $p_K \in \mathcal{T}_h$.
At the first step, we assemble $\mathbf{A}, \mathbf{M}$ by choosing $p_K=p^\text{max}$ for any $p_K \in \mathcal{T}_h$; we refer to \eqref{eq::matrixFull} and Section \ref{sec:polyDG} for their definition.
Then, at each time step, we can efficiently and locally update $\mathbf{A}, \mathbf{M}$ according to the local adaptation of the polynomial degree driven by the $p$-adaptive algorithm. The update is trivial and massively parallel for the mass matrix because, in the DG setting, the mass matrix is block diagonal. In Figure~\ref{fig::examplematrix}, we show an example of the dynamic update of the stiffness matrix $\mathbf{A}$ for a two-dimensional example. We suppose to have a polytopal mesh with $n_\text{el}$ elements, and we assume, for simplicity and clarity of the example, to have $p_\text{max}=2$. If the element $K_1$ is selected by our adaptive algorithm to decrease the local polynomial approximation degree from $p_{K_1}=2$ to $p_{K_1}=1$, the global stiffness matrix has to be locally modified as described in Figure~\ref{fig::examplematrix}.
Let $\mathbf{A}_1\in \mathbb{R}^{6x6}$ be the block diagonal part of $\mathbf{A}$ corresponding to the (six) degrees of freedom associated with $K_1$. The update of $\mathbf{A}_1$ into $\hat{\mathbf{A}}_1 \in \mathbb{R}^{3x3}$ is depicted in Figure~\ref{fig::examplematrix} (top panel).
Let now $\mathbf{A}_i\in \mathbb{R}^{6x6}$, $i=2, \ldots, n_\text{el}$, the off-diagonal blocks corresponding to the degrees of freedom associated with $K_1$.
The update of $\mathbf{A}_j$ into $\hat{\mathbf{A}}_j \in \mathbb{R}^{3x6}$ is depicted in Figure~\ref{fig::examplematrix} (bottom panel).
We remark that particular attention has to be given to the jump-jump term of the stiffness matrix, i.e., 
$\sum_{F \in \mathcal{F}_h^I} \int_F \eta  \jumpl u_h \jumpr \cdot  \jumpl v_h  \jumpr d\sigma \quad \forall \: u_h,v_h \in V^{\text{DG}}_h$.
Since the penalization parameter $\eta$ defined in \eqref{eq:eta} depends on the element-wise polynomial degree, this matrix has to be reassembled at each adaptive step. We can therefore write, exploiting the definition in \eqref{eq:coer}: $\mathbf{A}=\widetilde{\mathbf{A}}+\mathbf{S}$ where $\mathbf{S}$ is the matrix representation of the jump-jump term (cf. \eqref{eq::matrixFull} and \eqref{eq:eta}, last term) and  $\widetilde{\mathbf{A}}$ is the matrix representation of the first three terms in the bilinear form defined as in \eqref{eq:eta}. In the $p$-adaptive algorithm $\widetilde{\mathbf{A}}$ is updated as described above and illustrated in Figure~\ref{fig::examplematrix} while the matrix $\mathbf{S}$ is reassembled.

\section{Numerical results}
\label{sec:NumericalResults}
In this section, we present various numerical tests to evaluate the performance of the proposed $p$-adaptive method.
The numerical results have been obtained with a numerical code based the open source Lymph library \cite{antonietti2025lymph}, which implements
the PolyDG method for multiphysics differential problems.
\subsection{Test case 1: traveling wave solution}
\label{sec:test_case_1}
The first test case consists of a traveling wavefront benchmark.
For numerical verification, we consider a simplified model that describes an analytical traveling-wave, where the non-linear ionic term is defined by the following non-linear function: $f(u) = a(u - V_{\text{rest}})(u - V_{\text{thres}})(u - V_{\text{depol}})$ where $V_\text{rest} \le  V_\text{thres} \le V_\text{depol}$ and $a>0$. Here, $V_\text{rest}$ is the resting potential, $V_\text{thres}$ represents the threshold potential, and $V_\text{depol}$ represents the depolarization potential. Starting from \cite{pezzuto2016space}, where the convergence for the one-dimensional problem is analyzed, the solution is extended to the two-dimensional case by introducing $u_\mathrm{ex}$ as the exact solution defined as follows:
\begin{equation*}
    \begin{aligned}
      u_\mathrm{ex}(\boldsymbol{x},t) = \frac{V_{\text{dep}}-V_{\text{rest}}}{2}\left[1- \tanh\left(\frac{\boldsymbol{x}-\boldsymbol{c}t}{\epsilon}\right)\right] + V_{\text{rest}},
\end{aligned}
\label{eq:one-dimensional}
\end{equation*}
where $\epsilon$ characterizes the thickness of the wavefront, while $\boldsymbol{c}$ is the propagation speed of the wave.
\begin{table}[h]
\footnotesize
\centering
\caption{Test case 1. Values of the model parameters.}%
\label{table::analitical}
\begin{tabular}{lcl|lcl}  
\toprule
 \textbf{Parameter} & \textbf{Value} &  \textbf{Unit} &  \textbf{Parameter} & \textbf{Value} &  \textbf{Unit} \\
\midrule
    $\sigma_n$ &  $0.17$ & $\mathrm{mS\cdot mm^{-1}}$  & $a$ &  $1.4e-5$ & $\mathrm{mS \cdot mm^{-2} \cdot mV^{-2}}$ \\
    $\sigma_t$ &  $0.62$ & $\mathrm{mS\cdot mm^{-1}}$  & $\mathrm{\chi_m}$ &  $140$ & $\mathrm{mm^{-1}}$ \\
    $V_{\text{depol}}$ &  $30$ & $\mathrm{mV}$ &  $\mathrm{C_m}$ &  $0.01$ & $\mathrm{\mu F \cdot mm^{-2}}$  \\
    $V_{\text{rest}}$ &  -$85$ & $\mathrm{mV}$  & $c$ &  $0.5$ &  $\mathrm{mm \cdot ms^{-1}}$ \\
    $V_{\text{thres}}$ &  -$57.6$ & $\mathrm{mV}$  &$\epsilon$ &  $0.2$ & $\mathrm{mm}$\\
\bottomrule
\end{tabular}
\end{table}
\subsubsection{Test case 1a: convergence analysis}
\label{sec:convergence}
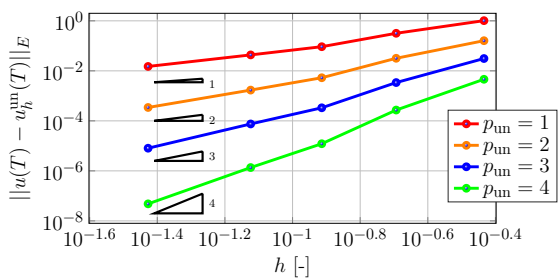
\begin{figure}[h]
\centering
\begin{tikzpicture}[scale=0.55,transform shape]{
\begin{axis}[%
width=3.875in,
height=2in,
at={(1.733in,0.687in)},
scale only axis,
xmode=log,
xmin=0.025,
xmax=0.4,
xminorticks=true,
xlabel = { $h$ [-]},
ylabel = { $||u(T)-u_h^{\text{un}}(T)||_E$},
ymode=log,
ymin=8e-9,
ymax=2,
yminorticks=true,
axis background/.style={fill=white},
xmajorgrids,
xminorgrids,
ymajorgrids,
yminorgrids,
legend style={at={(1.15,0.54)},legend cell align=left, align=left, draw=white!15!black}
]      	
\addplot [color=red, line width=2.0pt, mark=ball]
  table[row sep=crcr]{%
    0.3682  1.0163 \\
    0.2025  0.3178 \\
    0.1220  0.0924\\
    0.0752  0.0434\\
    0.0374  0.0150\\
};
\addlegendentry{$p_\text{un}=1$}			
\addplot [color=orange, line width=2.0pt, mark=ball]
  table[row sep=crcr]{%
    0.3682  0.1606 \\
    0.2025  0.0318 \\
    0.1220  0.0053 \\
    0.0752  0.0017 \\
    0.0374  3.4308e-4 \\
};
\addlegendentry{$p_\text{un}=2$}
\addplot [color=blue, line width=2.0pt, mark=ball]
  table[row sep=crcr]{%
    0.3682  0.0311 \\
    0.2025  0.0034 \\
    0.1220  3.3314e-4 \\
    0.0752  7.5196e-5 \\
    0.0374  8.0577e-6 \\
};
\addlegendentry{$p_\text{un}=3$}
\addplot [color=green, line width=2.0pt, mark=ball]
  table[row sep=crcr]{%
    0.3682  0.0046 \\
    0.2025  2.7075e-4 \\
    0.1220  1.2175e-5 \\
    0.0752  1.3578e-6 \\
    0.0374  4.7985e-8 \\
};
\addlegendentry{$p_\text{un}=4$}
\node[right, align=left, text=black, font=\footnotesize]
at (axis cs:0.0547,0.003) {$1$};
\addplot [color=black, line width=1.5pt]
  table[row sep=crcr]{%
0.0542   0.00485\\
0.0390   0.0035\\
0.0542   0.0035\\
0.0542   0.00485\\
};
\node[right, align=left, text=black, font=\footnotesize]
at (axis cs:0.0547,0.00011) {$2$};
\addplot [color=black, line width=1.5pt]
  table[row sep=crcr]{%
0.0542   0.0001748\\
0.0390   0.0001000\\
0.0542   0.0001000\\
0.0542   0.0001748\\
};
\node[right, align=left, text=black, font=\footnotesize]
at (axis cs:0.0547,4e-6) {$3$};
\addplot [color=black, line width=1.5pt]
  table[row sep=crcr]{%
0.0542   0.0000060\\
0.0390   0.0000025\\
0.0542   0.0000025\\
0.0542   0.0000060\\
};
\node[right, align=left, text=black, font=\footnotesize]
at (axis cs:0.0547,5e-08) {$4$};
\addplot [color=black, line width=1.5pt]
  table[row sep=crcr]{%
0.0542   1.23e-07\\
0.0390   2e-08\\
0.0542   2e-08\\
0.0542   1.23e-07\\
};
\end{axis}}
\end{tikzpicture}
\caption{Test case 1a. Computed errors in the energy norm defined in Equation \eqref{eq::energy_norm} at the final time $T=1e-4 \, \text{[s]}$ as a function of the mesh size (loglog scale) for different choices of the polynomial degree $p_\text{un}=1,2,3,4$.}
\label{fig:convergence_analysis}
\end{figure}

For the first test case, we consider the domain $\Omega = (-1,2) \times (-0.5,0.5)$, discretized with a mesh constructed by PolyMesher \cite{talischi2012polymesher}. We define the conductivity tensor as
$\boldsymbol{\Sigma} = 0.1336 \mathbbm{1} [\text{ mS} \cdot \text{mm}^{-1}$], $\mathds{1}$ being the identity tensor. The model's parameters are chosen as reported in Table~\ref{table::analitical}, cf. from \cite{pezzuto2016space}. In Figure~\ref{fig:convergence_analysis}, we report the computed errors in the energy norm defined in \eqref{eq::energy_norm} at the final time $T=1e-4 \, \text{[s]}$ with $\Delta t = 1e-6$. These results have been obtained considering a uniform polynomial degree for every element of the mesh, i.e., $p_K=p_\text{un}$ for any $K \in \mathcal{T}_h$, with $p_\text{un}=1,2,3,4$ and considering a sequence of uniformly refined meshes with mesh-size $h$ (the corresponding number of elements is equals to $70,250,800,3000,8000$).  
We observe that the theoretical convergence rates are achieved, consistently with the theoretical analysis proved in \cite{saglio2024high}. 

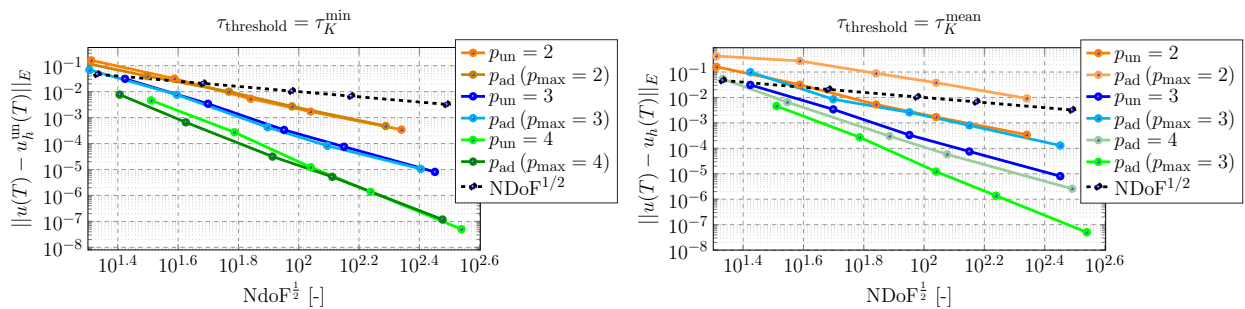
\begin{figure}[h!]
      \centering
\begin{subfigure}{0.45\textwidth}
\begin{tikzpicture}[scale = .53,transform shape]{
\begin{axis}[%
width=3.875in,
height=2in,
at={(1.733in,0.687in)},
scale only axis,
xmode=log,
xmin=20,
xmax=400,
xminorticks=true,
xlabel = { $\text{NdoF}^\frac{1}{2}$ [-]},
ylabel = { $||u(T)-u_h^\text{un}(T)||_E$},
ymode=log,
ymin=8e-9,
ymax=0.5,
yminorticks=true,
axis background/.style={fill=white},
title={\color{black} $\tau_\text{threshold}=\tau_K^\text{min}$},
xmajorgrids,
xminorgrids,
ymajorgrids,
yminorgrids,
major grid style={dashed, color=black!40},
minor grid style={dotted, color=black!30},
legend style={
at={(1.35,1.04)}, 
legend cell align=left, 
legend columns=1,
align=left, 
draw=white!15!black}
]

\addplot [color=orange, line width=2.0pt, mark=ball]
  table[row sep=crcr]{%
    20.49  0.1606 \\
    38.7  0.0318 \\
    69.28  0.0053 \\
    109.54  0.0017 \\
    219.08  3.4308e-4 \\
};
\addlegendentry{$p_\text{un}=2$}

\addplot [color=Orange3, line width=2.0pt, mark=ball]
  table[row sep=crcr]{%
    16.9  0.1792 \\
    31.4  0.0398 \\
    58.7  0.01 \\
    94.97  0.0027 \\
    194.2  4.7308e-4 \\
};
\addlegendentry{$p_\text{ad}\,(p_\text{max}=2)$}

\addplot [color=blue, line width=2.0pt, mark=ball]
  table[row sep=crcr]{%
    26.54  0.0311 \\
    50  0.0034 \\
    89.44  3.3314e-4 \\
    141.42  7.5196e-5 \\
    282.84  8.0577e-6 \\
};
\addlegendentry{$p_\text{un}=3$}

\addplot [color=DeepSkyBlue2, line width=2.0pt, mark=ball]
  table[row sep=crcr]{%
    20.22  0.0689 \\
    39.45  0.0076 \\
    78.8  4.2386e-4 \\
    124.42  8.165e-5 \\
    253.81  1.0577e-5 \\
};
\addlegendentry{$p_\text{ad}\,(p_\text{max}=3)$}

\addplot [color=green, line width=2.0pt,  mark=ball]
  table[row sep=crcr]{%
    32.40    0.0046 \\ 
    61.23   2.714e-4 \\ 
    109.544   1.2261e-5 \\ 
    173.205  1.3805e-6 \\ 
    346.41  5.035e-8 \\ 
};
\addlegendentry{$p_\text{un}=4$}

\addplot [color=Green4, line width=2.0pt,  mark=ball]
  table[row sep=crcr]{%
    25.49    0.00767 \\ 
    42.23    6.56e-4 \\ 
    81.86    3.151242e-5 \\ 
    129.205  5.21404e-6 \\ 
    299.71   1.181e-7\\ 
};
\addlegendentry{$p_\text{ad}\,(p_\text{max}=4)$}

\addplot [color=black, dashed, line width=2.0pt,  mark=ball]
  table[row sep=crcr]{%
    21.49    0.0465 \\ 
    48.23   0.0207 \\ 
    94.86   0.0105 \\ 
    149.205  0.0067 \\ 
    309.71  0.00322\\ 
};
\addlegendentry{$\mathrm{NDoF}^{1/2}$}
\end{axis}}
\end{tikzpicture}
\end{subfigure}\hspace*{2ex}
      \begin{subfigure}{0.45\textwidth}
\begin{tikzpicture}[scale = 0.53,transform shape]{
\begin{axis}[%
width=3.875in,
height=2in,
at={(1.733in,0.687in)},
scale only axis,
xmode=log,
xmin=20,
xmax=400,
xminorticks=true,
xlabel = { $\mathrm{NDoF}^\frac{1}{2}$ [-]},
ylabel = { $||u(T)-u_h(T)||_E$},
ymode=log,
ymin=1e-8,
ymax=0.9,
yminorticks=true,
axis background/.style={fill=white},
title={\color{black} $\tau_\text{threshold}=\tau_K^\text{mean}$},
xmajorgrids,
xminorgrids,
ymajorgrids,
yminorgrids,
major grid style={dashed, color=black!40},
minor grid style={dotted, color=black!30},
legend style={
at={(1.35,1.04)}, 
legend cell align=left, legend columns=1, align=left, draw=white!15!black}
]

\addplot [color=orange, line width=2.0pt, mark=ball]
  table[row sep=crcr]{%
    20.49  0.1606 \\
    38.7  0.0318 \\
    69.28  0.0053 \\
    109.54  0.0017 \\
    219.08  3.4308e-4 \\
};
\addlegendentry{$p_\text{un}=2$}

\addplot [color=Tan1, line width=2.0pt, mark=ball]
  table[row sep=crcr]{%
    20.49  0.4241 \\
    38.7  0.2731 \\
    69.28  0.0891 \\
    109.54  0.0380 \\
    219.08  0.0092 \\
};
\addlegendentry{$p_\text{ad}\,(p_\text{max}=2)$}

\addplot [color=blue, line width=2.0pt, mark=ball]
  table[row sep=crcr]{%
    26.54  0.0311 \\
    50  0.0034 \\
    89.44  3.3314e-4 \\
    141.42  7.5196e-5 \\
    282.84  8.0577e-6 \\
};
\addlegendentry{$p_\text{un}=3$}

\addplot [color=DeepSkyBlue2, line width=2.0pt, mark=ball]
  table[row sep=crcr]{%
    26.54  0.0989 \\
    50  0.00849 \\
    89.44  0.0026 \\
    141.42  7.9853e-4 \\
    282.84  1.29e-4 \\
};
\addlegendentry{$p_\text{ad}\,(p_\text{max}=3)$}

\addplot [color=DarkSeaGreen3, line width=2.0pt,  mark=ball]
  table[row sep=crcr]{%
    21.49    0.053762 \\ 
    35.23   0.006512 \\ 
    76.86   3.01242e-4 \\ 
    119.205  5.81404e-5 \\ 
    309.71  2.581e-6\\ 
};
\addlegendentry{$p_\text{ad}=4$}
\addplot [color=green, line width=2.0pt,  mark=ball]
  table[row sep=crcr]{%
    32.40    0.0046 \\ 
    61.23   2.714e-4 \\ 
    109.544   1.2261e-5 \\ 
    173.205  1.3805e-6 \\ 
    346.41  5.035e-8 \\ 
};
\addlegendentry{$p_\text{ad}\,(p_\text{max}=3)$}
\addplot [color=black, dashed, line width=2.0pt,  mark=ball]
  table[row sep=crcr]{%
    21.49    0.0465 \\ 
    48.23   0.0207 \\ 
    94.86   0.0105 \\ 
    149.205  0.0067 \\ 
    309.71  0.00322\\ 
};
\addlegendentry{$\mathrm{NDoF}^{1/2}$}
\end{axis}}
\end{tikzpicture}
\end{subfigure}
\caption{Test case 1a. Computed errors in the energy norm defined in \eqref{eq::energy_norm} at the final time $T=1e-4 \, \mathrm{[s]}$ as a function of total number of degrees of freedom $\mathrm{NDoF}$   (loglog scale) for uniform ($p_\text{un}$) and adaptive ($p_\text{ad}$) polynomial degree and different choices of $\tau_{\text{threshold}}=\tau_\text{threshold}^\text{min}, \tau_\text{threshold}^\text{mean}$,  cf. \eqref{eq::centroid}. Here $p_\text{max}$ represents the maximum polynomial degree in the adaptive algorithm.}
\label{fig:comparison_padapt}
\end{figure}

We next present the numerical results obtained based on employing our $p$-adaptive algorithm and compare them with the corresponding ones obtained  under uniform $p$-refinement.  Figure~\ref{fig:comparison_padapt} reports the computed errors in the energy norm defined in \eqref{eq::energy_norm} at the final time $T=1e-4$ as a function of total number of degrees of freedom ($\mathrm{NDoF}$) for adaptive and uniform polynomial distributions (loglog scale).  More precisely, in Figure~\ref{fig:comparison_padapt} we show the computed errors obtained with a uniform polynomial degree $p_\text{un}=2,3,4$ and compare them with the corresponding results obtained by employing the proposed $p$-adaptive method (denoted as $p_\text{ad}$).

These results are based on two different threshold parameters: Figure~\ref{fig:comparison_padapt} (left) show the results obtained setting $\tau_\text{threshold}=\tau_\text{threshold}^\text{min}$, i.e., the threshold is defined as the minimum of the centroids of the $k$-means clustering (cf. also \eqref{eq::centroid}), while Figure~\ref{fig:comparison_padapt} (right) show the results obtained setting  $\tau_\text{threshold}=\tau_\text{threshold}^\text{mean}$, i.e., the mean of the centroids obtained from the $k$-means clustering (cf. also \eqref{eq::centroid}).
In this set of numerical experiments, we impose the maximum polynomial degree in our $p$-adaptive algorithm as $p_\text{max}=2,3,4$.
From the results reported in Figure~\ref{fig:comparison_padapt}, we can conclude that choosing as threshold the mean of the centroids seems to lead to larger errors,  while exploiting a threshold based on the minimum of the centroids seems to provide better results and lead to approximations that are comparable to those obtained with a uniform polynomial degree. 
Given these numerical results, the simulation of the following sections will be carried out using the minimum threshold.

\subsubsection{Test case 1b: degrees of freedom analysis}
\label{dof}
\begin{figure}[!htbp]
\centering
\begin{subfigure}[t]{0.45\textwidth}
\centering
\color{white}{test}
\end{subfigure}\hfill
\begin{subfigure}[t]{0.45\textwidth}
\vspace*{-3ex}
\centering
\hspace*{0em}\includegraphics[trim={0 0cm 0 0cm},clip,scale=0.18]{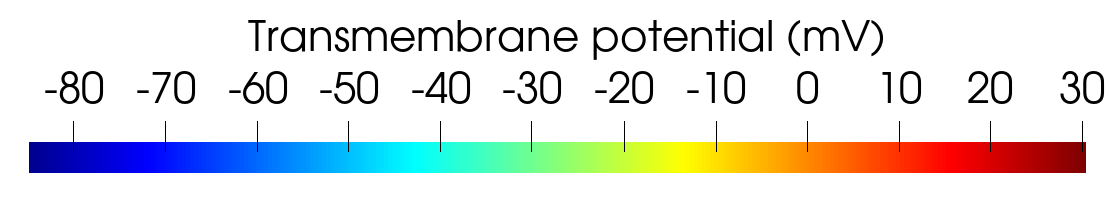}
\end{subfigure}
\begin{subfigure}[t]{0.45\textwidth}
\centering
\vspace*{-1ex}
\includegraphics[trim={0 11cm 0 13cm},clip, scale=0.09]{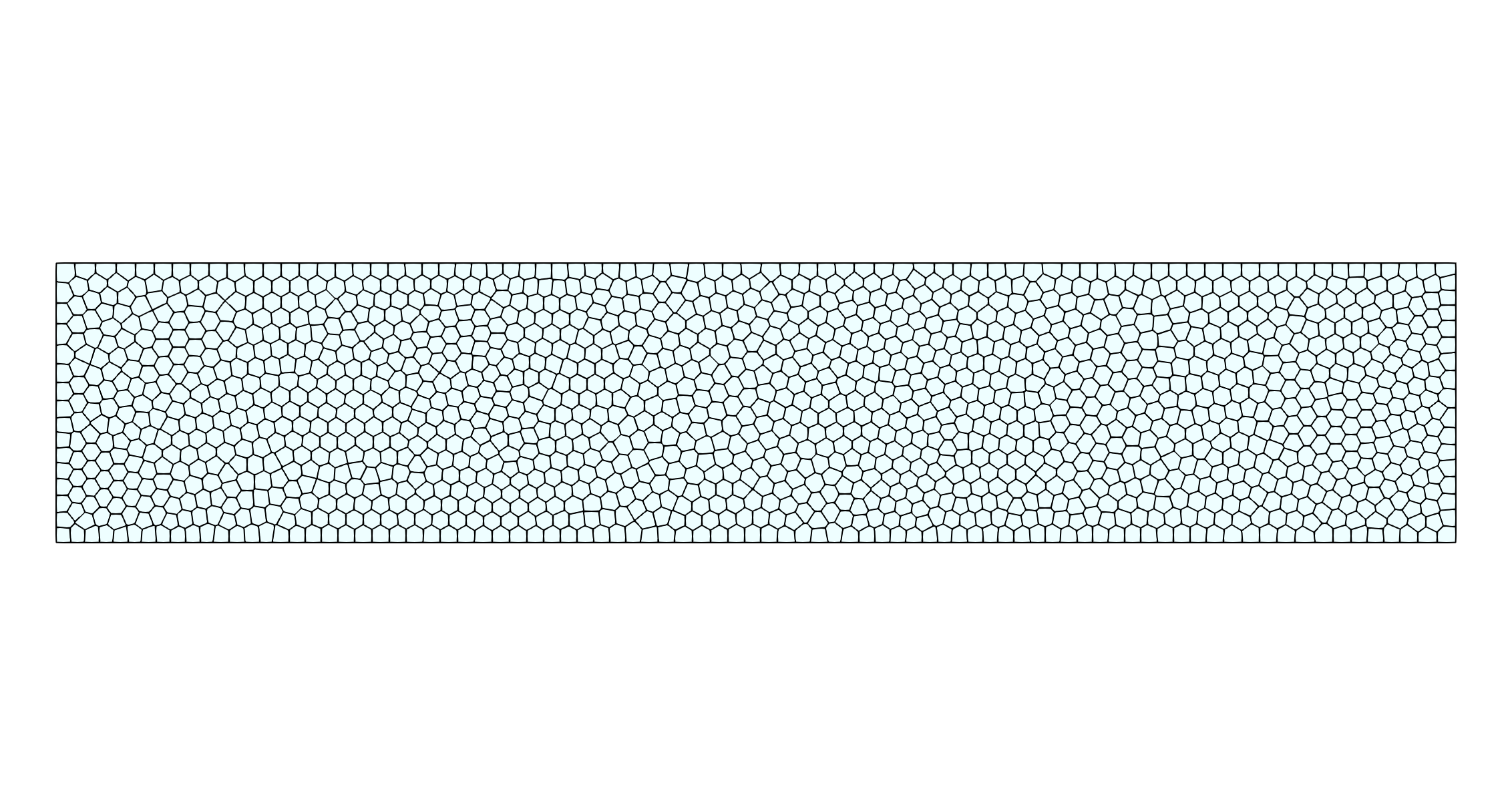}
\caption{\label{fig:comput_domain_square2}} 
\end{subfigure}\hfill
\begin{subfigure}[t]{0.45\textwidth}
\centering
\vspace*{-1ex}
\includegraphics[trim={0 11cm 0 13cm},clip, scale=0.09]{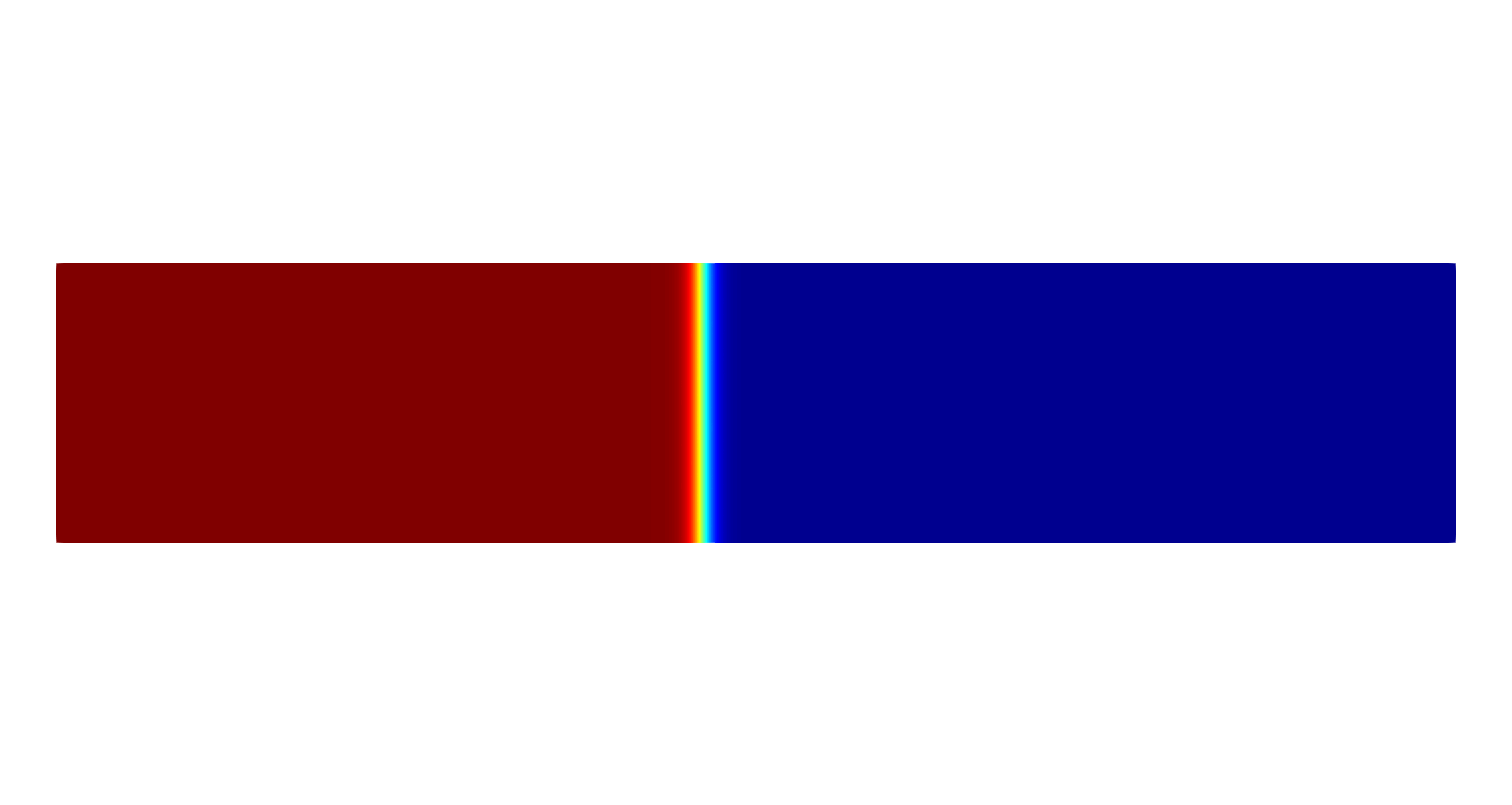} \caption{\label{fig:comput_domain_IC2}} 
\end{subfigure}
\caption{Test case 1b.  \eqref{fig:comput_domain_square2} Polytopal grid with 1500 elements ($h = 0.1170$); (\ref{fig:comput_domain_IC2}) initial condition of the transmembrane potential.}
\label{fig:tvinitial}
\end{figure}

We next analyze the performance of our $p$-adaptive scheme to correctly track the evolution of the wavefront. We set  $\Omega = (-2,3)\times (-0.5,0.5)$, and consider a polytopal mesh with 1500 elements, cf. Figure~\ref{fig:comput_domain_IC2}. 
Here, the  conductivity tensor is defined as
$\boldsymbol{\Sigma} = 0.0081 \mathbbm{1} \, \mathrm{[mS \cdot mm^{-1}]}$, and all the other parameters are defined as in Table~\ref{table::analitical}. For the $p$-adaptive scheme, the initial polynomial degree is uniform and equal to $p_{un}=5$, which corresponds to a total number of degrees of freedom equal to $\mathrm{NDoF} = 3.15e+4$.
The initial condition for the transmembrane potential is reported in Figure~\ref{fig:comput_domain_IC2}.

\begin{figure}[h!]
\centering
\begin{subfigure}[b]{0.45\textwidth}
\centering
\hspace*{35ex}\includegraphics[scale=0.17]{singlewave_stationary/scale_.png}
\includegraphics[trim={0 13cm 0 13cm},clip, scale=0.075]{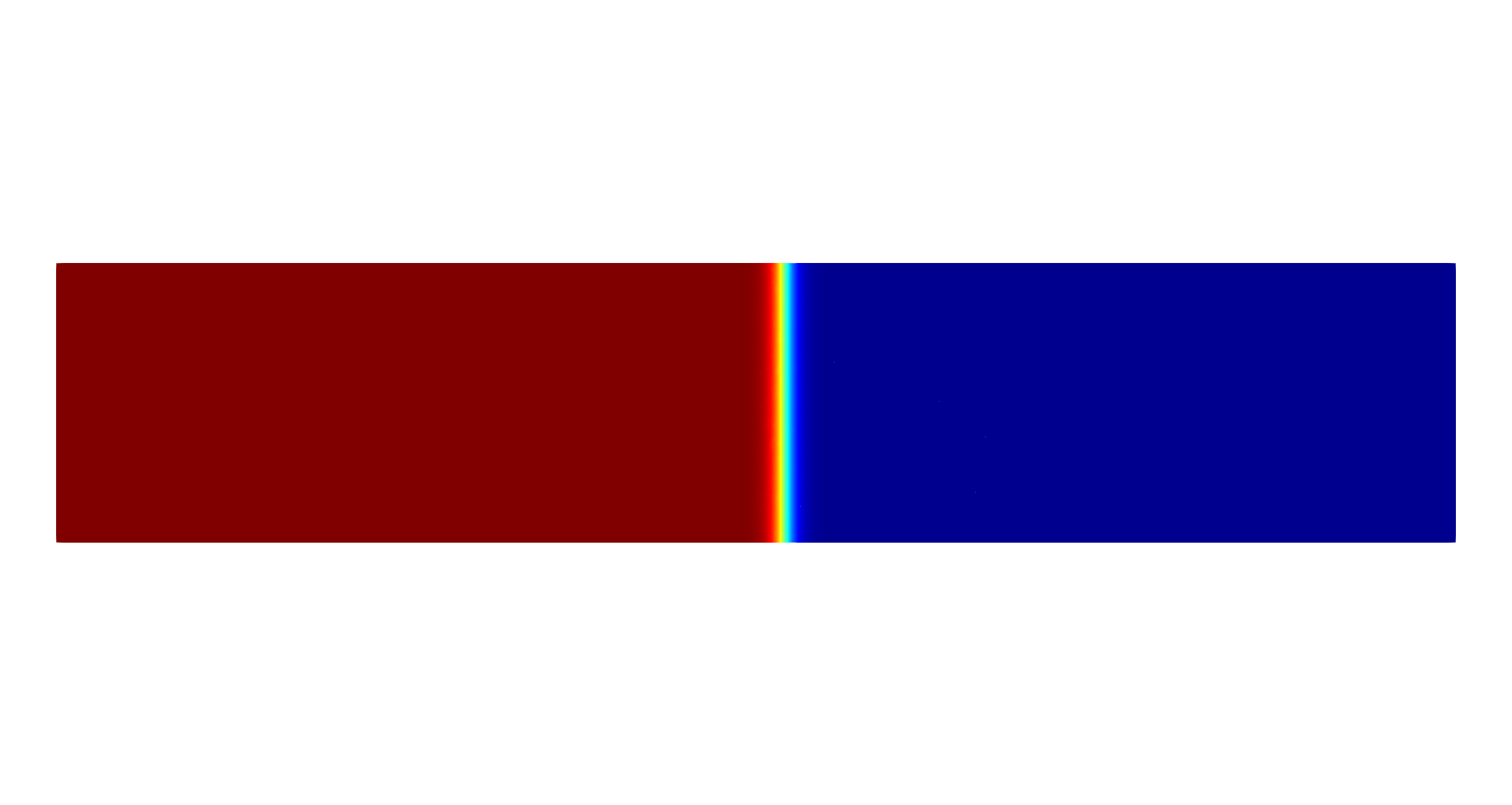} 
\includegraphics[trim={0 13cm 0 13cm},clip, scale=0.075]{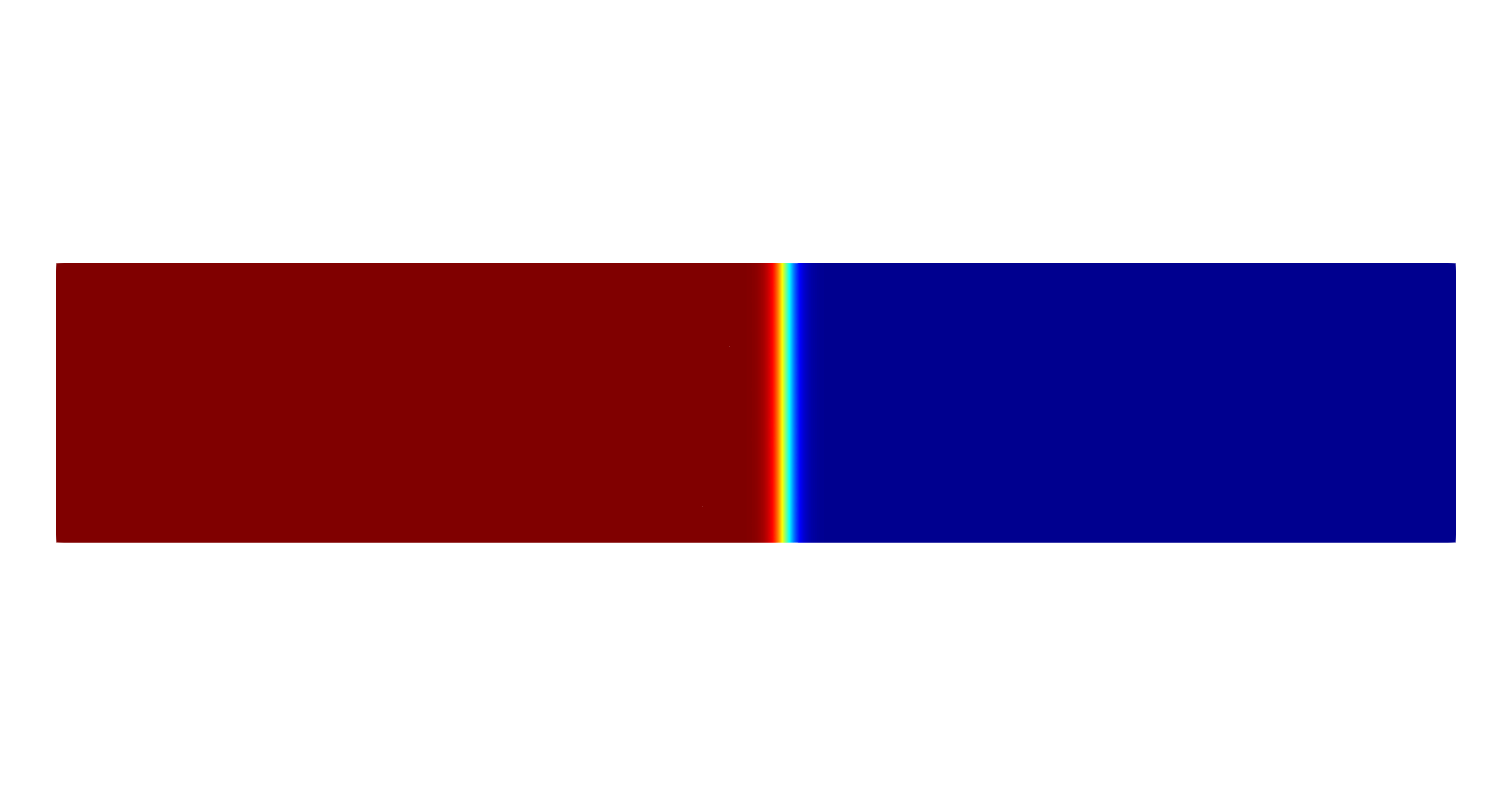}
\hspace*{30ex}\includegraphics[trim={0 0cm 0 0.3cm},clip, scale=0.17]{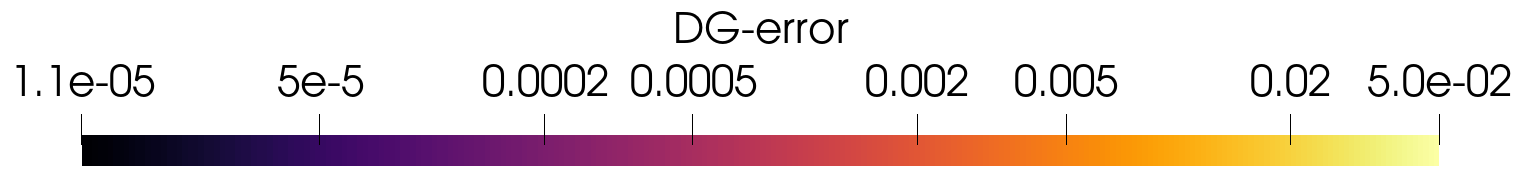}
\includegraphics[trim={0 13cm 0 13cm},clip, scale=0.075]{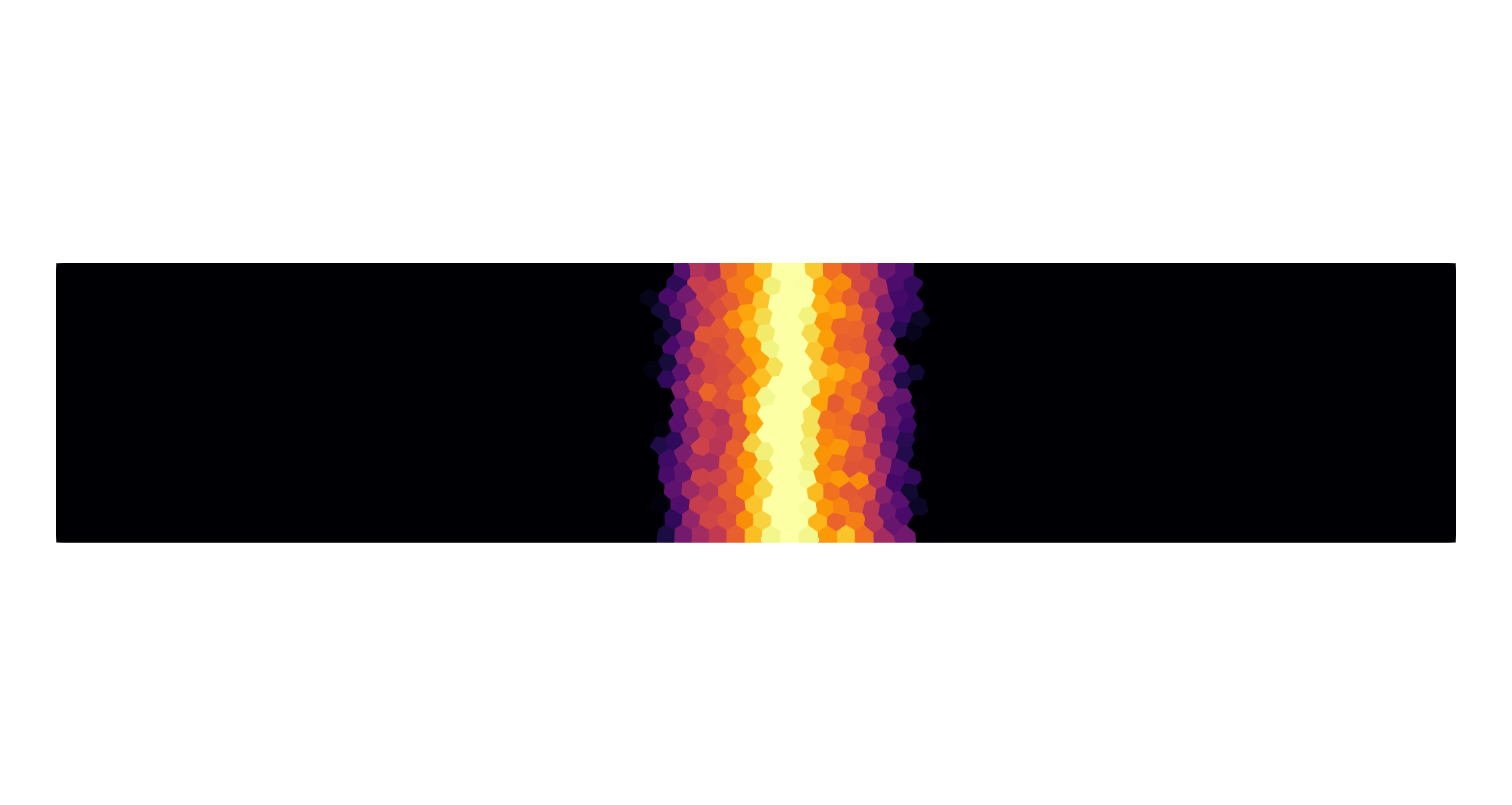}
\hspace*{30ex}\includegraphics[trim={0 0cm 0 0.3cm},clip, scale=0.17]{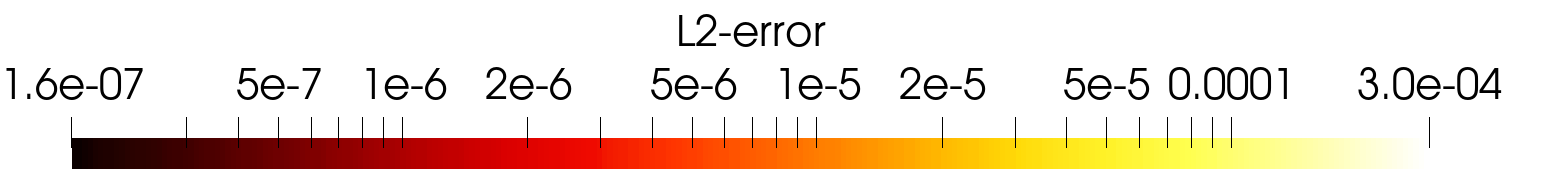}
\includegraphics[trim={0 13cm 0 13cm},clip, scale=0.075]{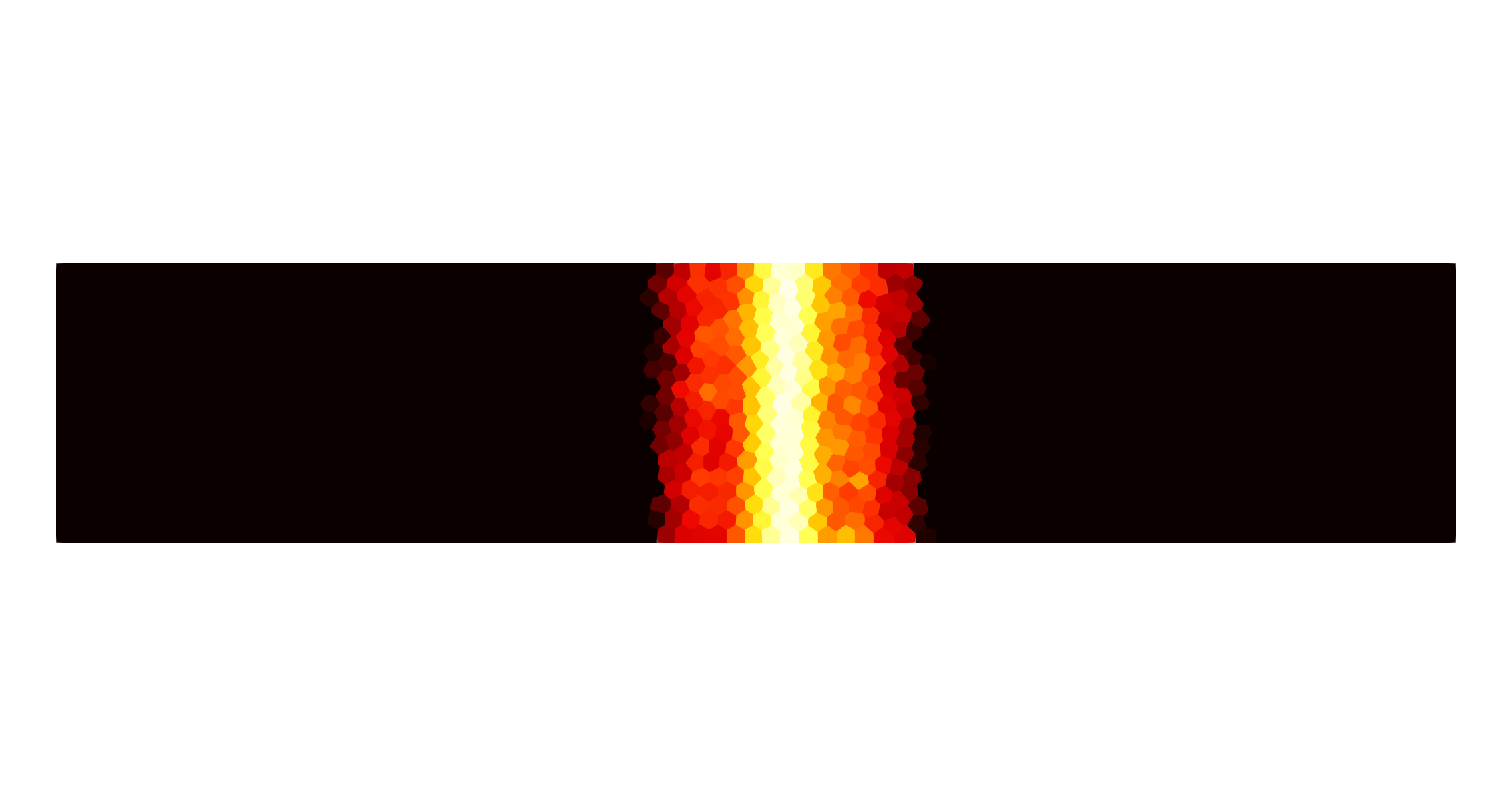}
\hspace*{40ex}\includegraphics[trim={0 0cm 0 0.3cm},clip, scale=0.17]{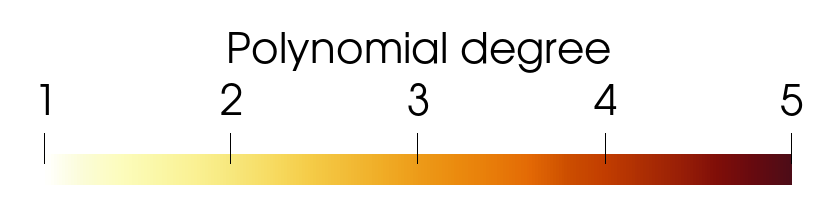}
\includegraphics[trim={0 13cm 0 13cm},clip, scale=0.075]{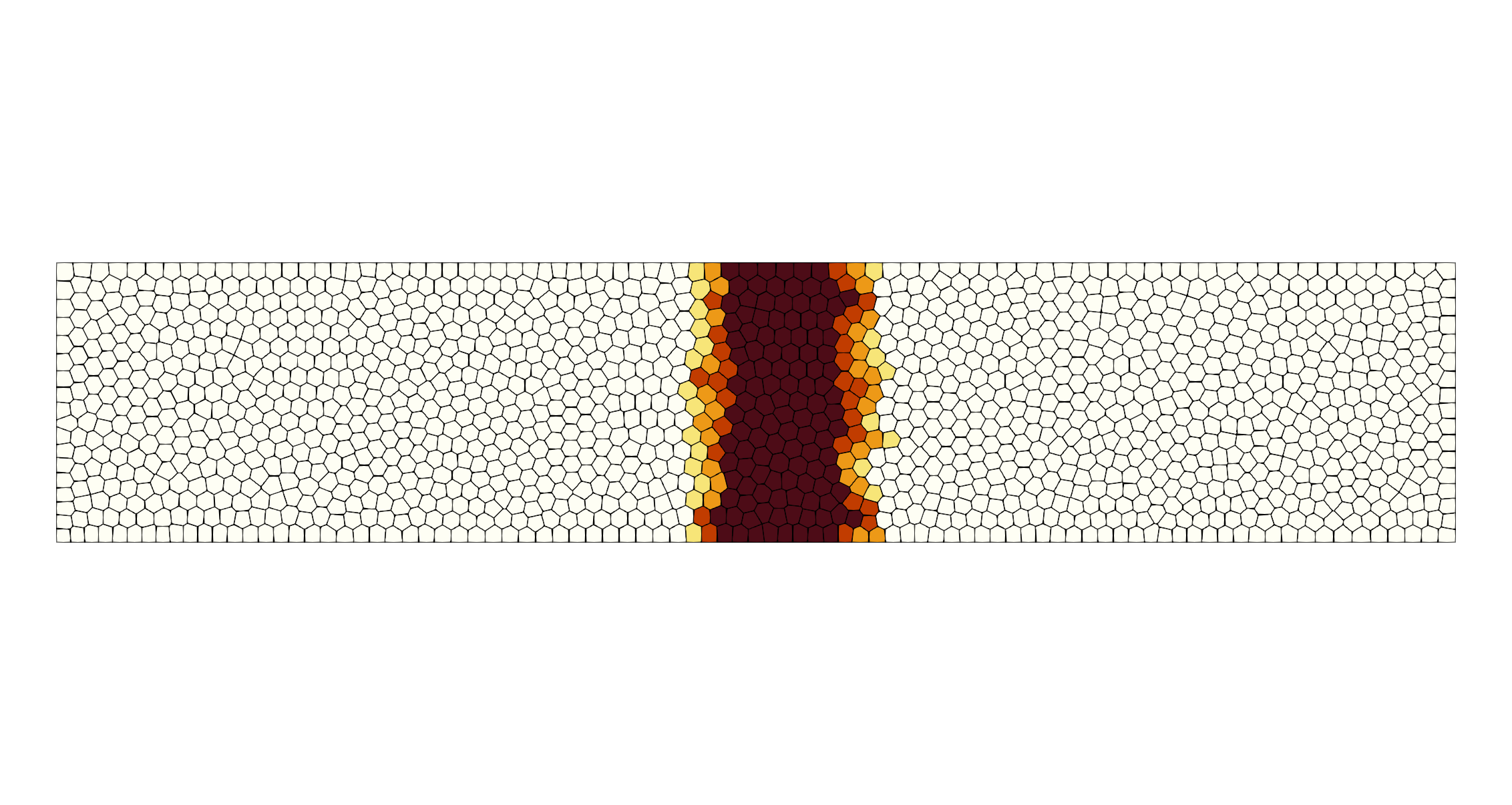}
\hspace*{30ex}\includegraphics[trim={0 0.1cm 0 0.1cm},clip, scale=0.17]{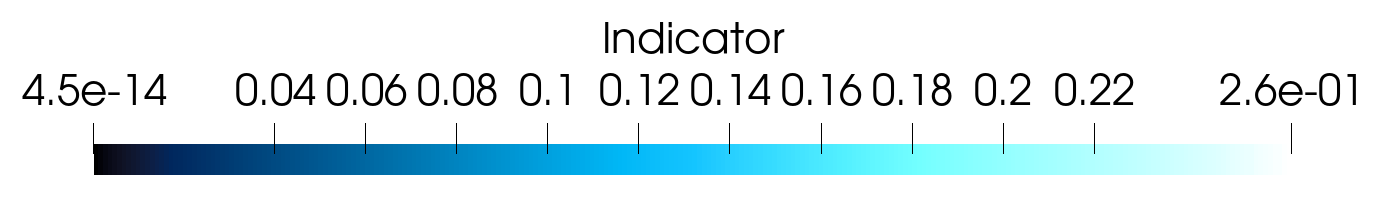}
\includegraphics[trim={0 13cm 0 13cm},clip, scale=0.075]{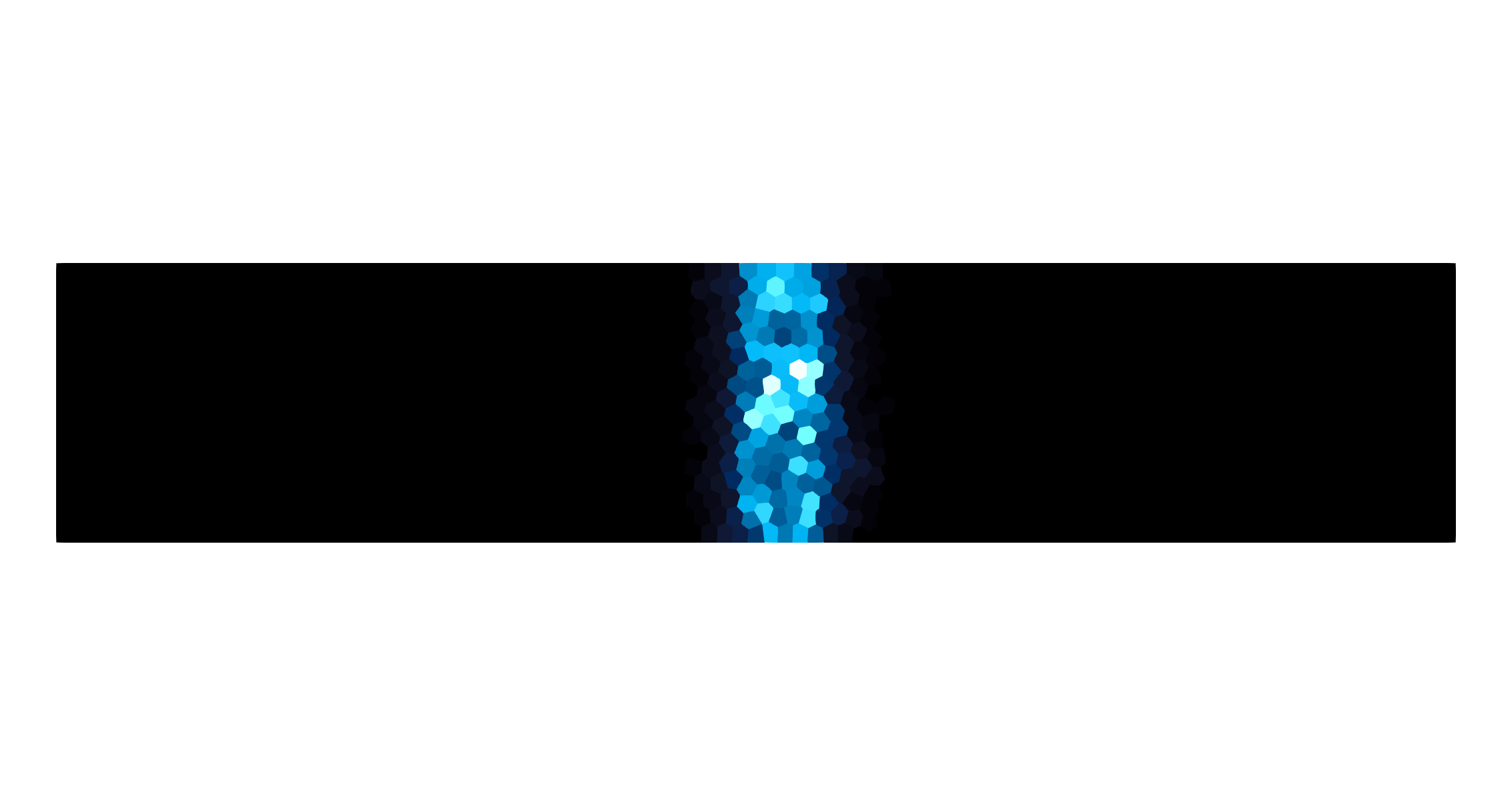} \caption{$t = 5 \, [\mathrm{ms}]$}  
\end{subfigure}
\begin{subfigure}[b]{0.45\textwidth}
    \centering
   \phantom{\includegraphics[scale=0.17]{singlewave_stationary/scale_.png}}
\includegraphics[trim={0 13cm 0 13cm},clip, scale=0.075]{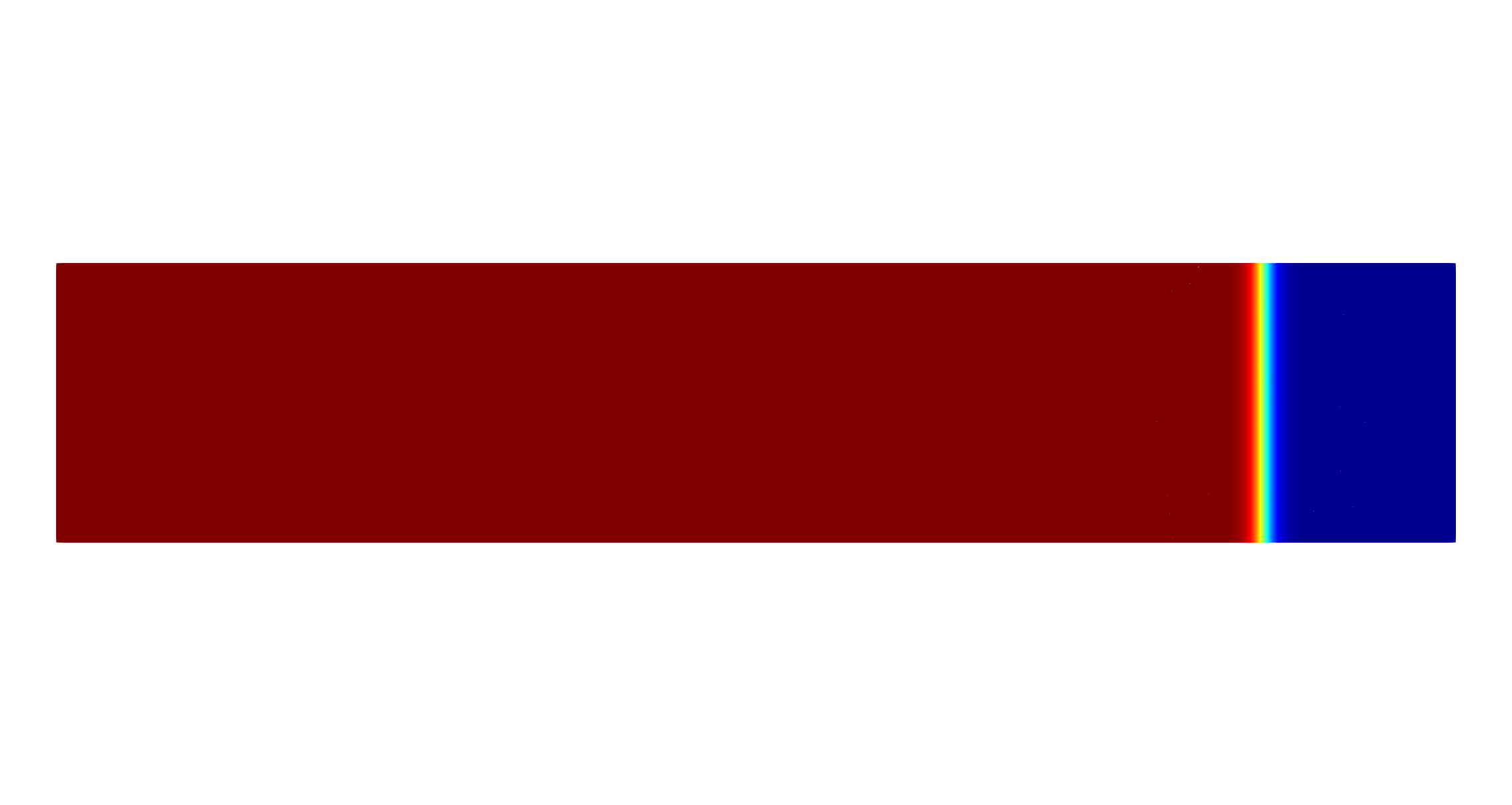} 
\includegraphics[trim={0 13cm 0 13cm},clip, scale=0.075]{singlewave_stationary/sol19.png}
\phantom{\includegraphics[trim={0 0cm 0 0.3cm},clip, scale=0.17]{singlewave_stationary/DGscalefinal.png}}
\includegraphics[trim={0 13cm 0 13cm},clip, scale=0.075]{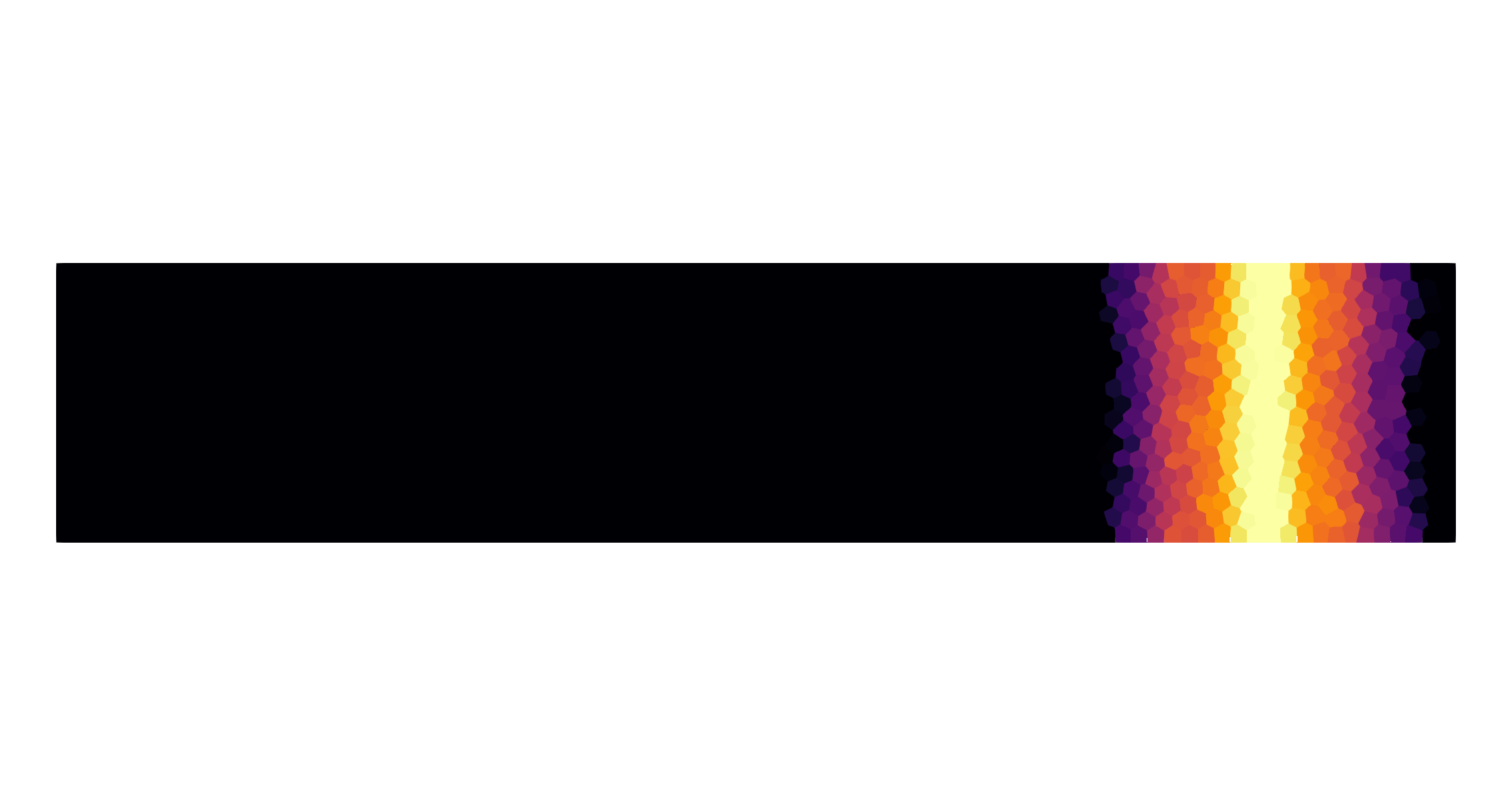}
\phantom{\includegraphics[trim={0 0cm 0 0.3cm},clip, scale=0.17]{singlewave_stationary/scalefinal.png}}
\includegraphics[trim={0 13cm 0 13cm},clip, scale=0.075]{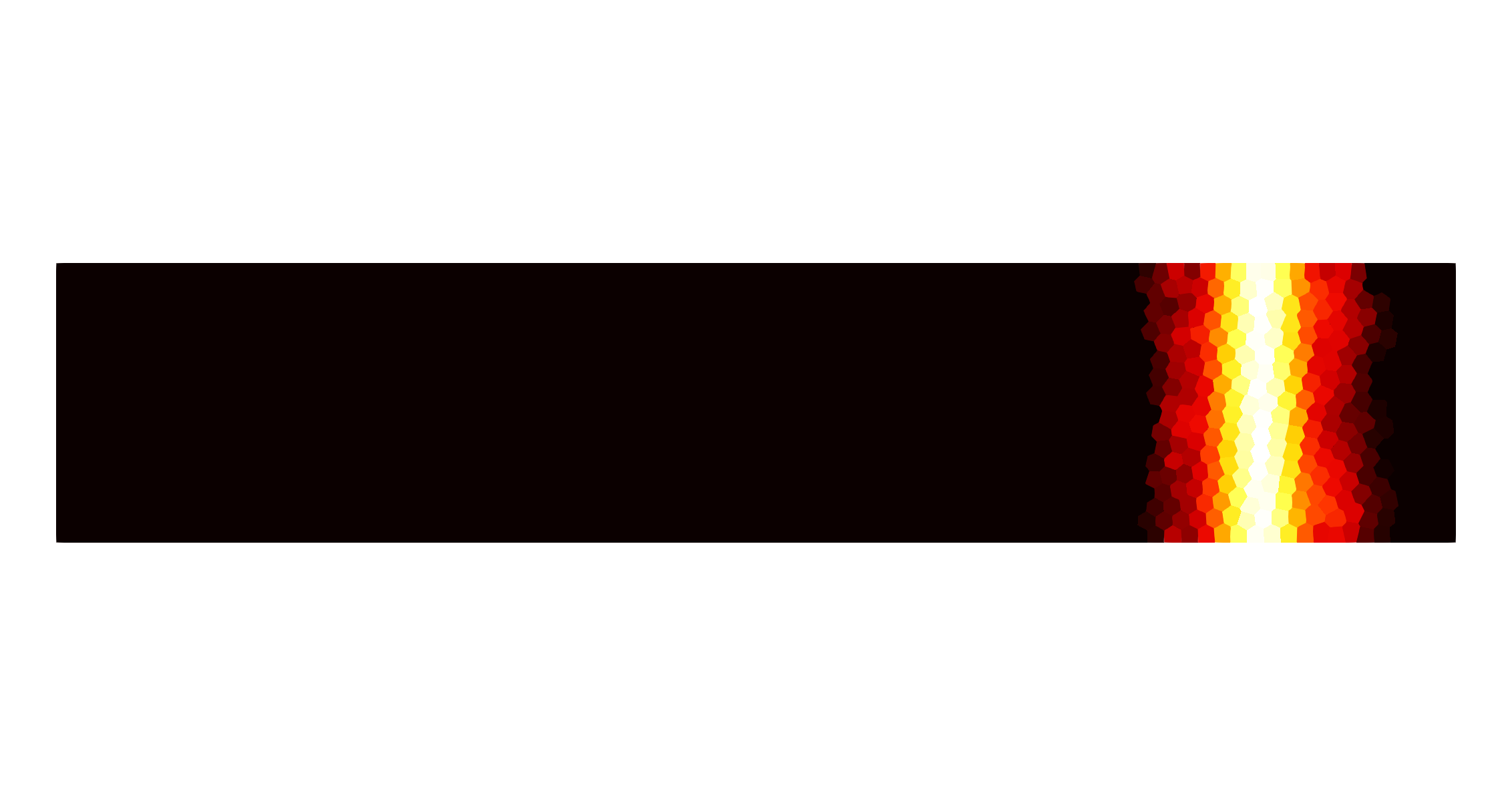}
\phantom{\includegraphics[trim={0 0cm 0 0.3cm},clip, scale=0.17]{singlewave_stationary/scale_p.png}}
\includegraphics[trim={0 13cm 0 13cm},clip, scale=0.075]{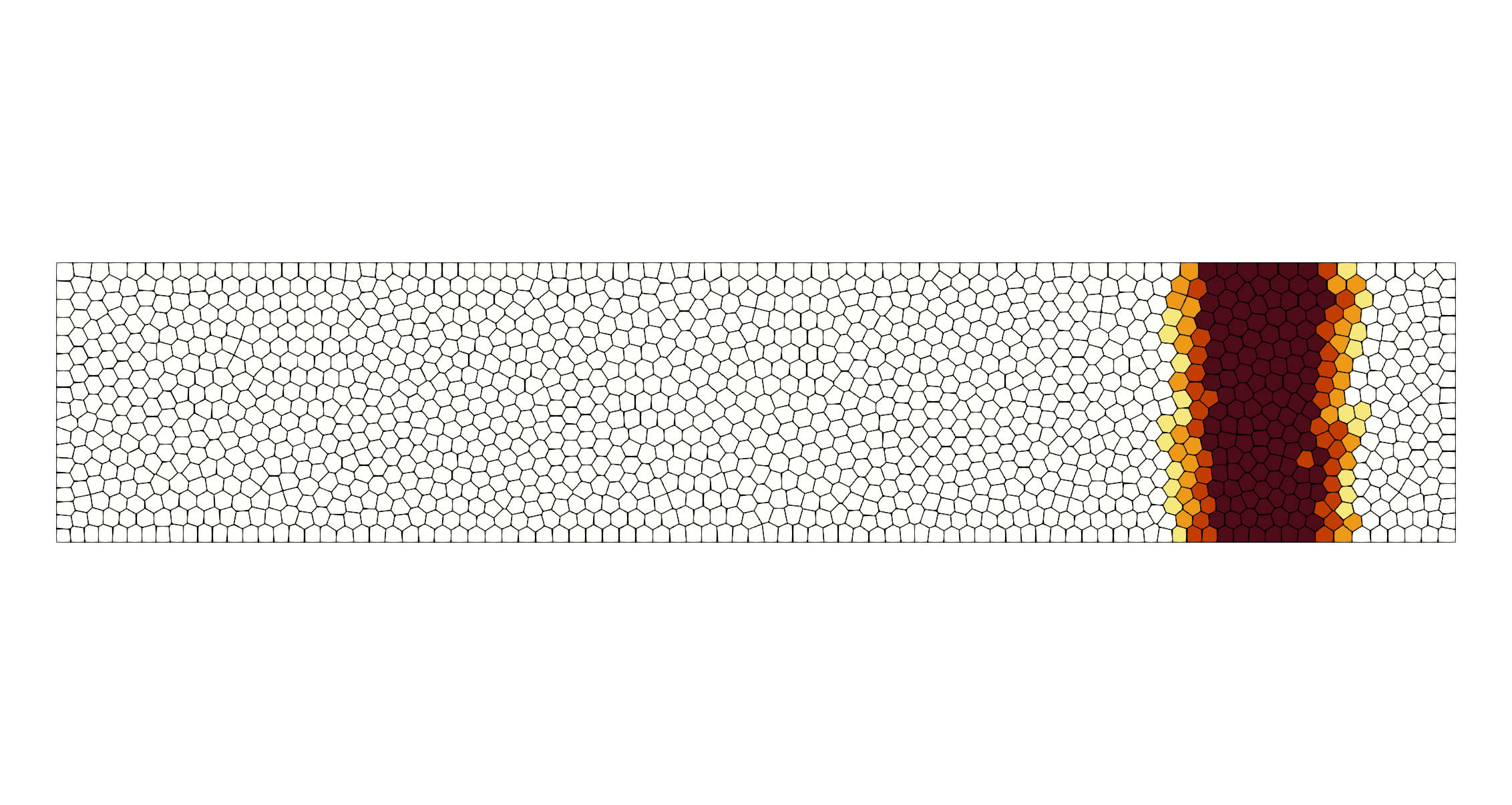}
\phantom{\includegraphics[trim={0 0.1cm 0 0.1cm},clip, scale=0.17]{singlewave_stationary/scale_t.png}}
\includegraphics[trim={0 13cm 0 13cm},clip, scale=0.075]{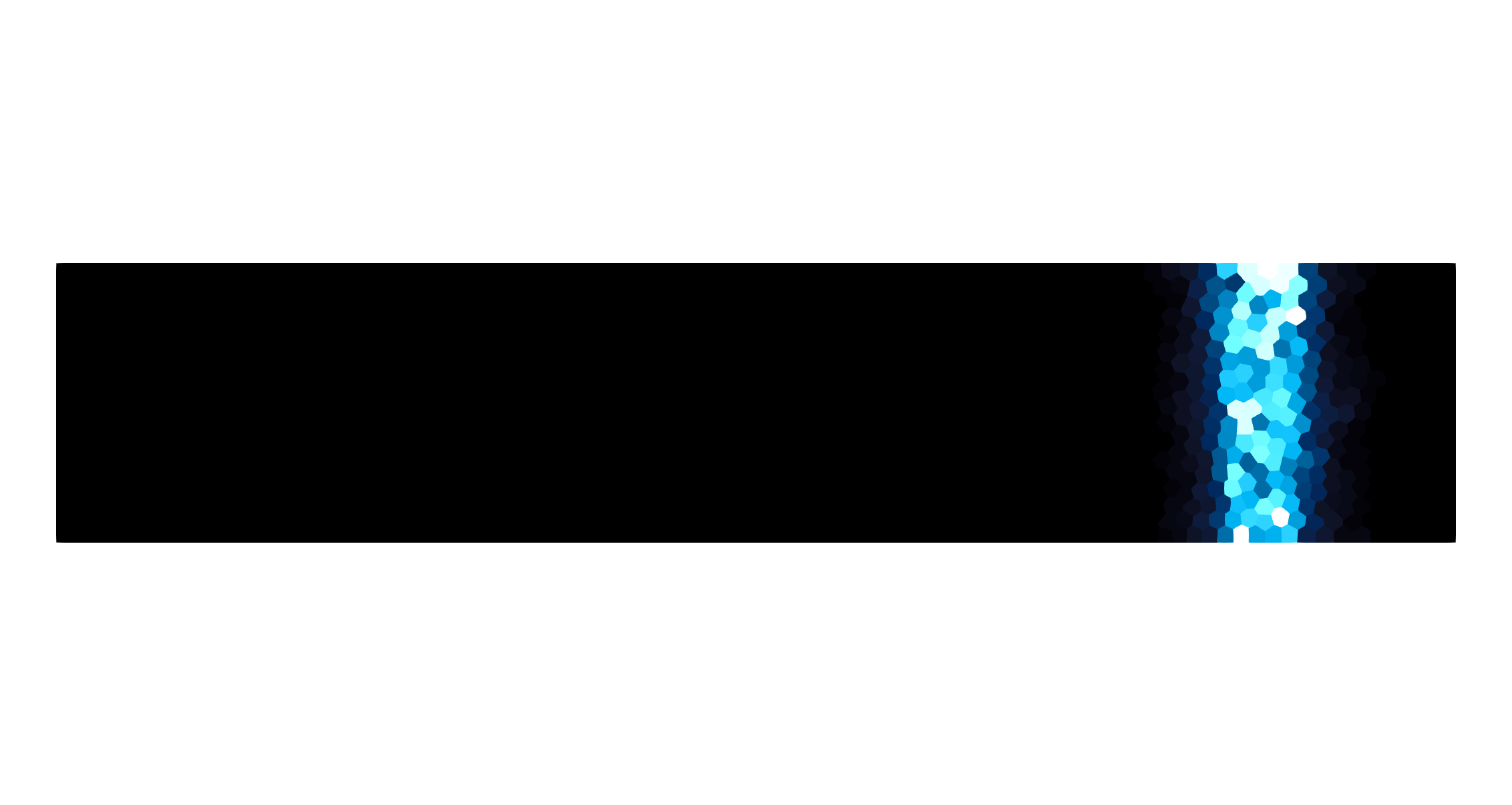} \caption{$t = 20.3 \, [\mathrm{ms}]$}  
  \end{subfigure}
\caption{Test case 1b.  Exact solution $u_{ex}$ (first row), numerical solution $u_h$ (second row), computed errors in the $DG$ and $L^2$ norms (third and fourth rows),  element-wise polynomial approximation degree (fith row) and local error indicator (sixth row) at two time instants $t=5\, [\mathrm{ms}]$ and $t=20.3\, [\mathrm{ms}]$.}
  \label{fig:: single_wave}
\end{figure}

In Figure~\ref{fig:: single_wave}, we report the results obtained at two time snapshots $t=5\, [\mathrm{ms}]$ and $t=20.3\, [\mathrm{ms}]$. 
The first two rows of Figure~\ref{fig:: single_wave} report the exact solution together with the corresponding numerical solution. In the third and fourth rows of Figure~\ref{fig:: single_wave}, we report the element-wise computed errors in the $DG$ and $L^2$ norms.
\begin{table}[!hp]
    \centering
    \caption{Test case 1b.  Computed errors in the $L2$  and $DG$  norms for different time snapshots  $t=5$ and $t=20.3  \, [\mathrm{ms}]$. Here $u_h^\text{ad}$ is the numerical solution computed with the $p$-adaptive scheme, $u_h^\text{un}$ is the one computed with a uniform polynomial degree $p_{\textrm{un}}=5$).}
    \label{tab:errorL2DGtable}
    \begin{tabular}{c|cc|cc}
    \hline
        &{$\|u_h^\text{un}(t) - u_\text{ex}(t)\|_{L^2(\Omega)}$} 
        & {$\|u_h^\text{un}(t) - u_\text{ex}(t)\|_{DG}$} &
        {$\|u_h^\text{ad}(t) - u_\text{ex}(t)_{{L^2}(\Omega)}$} 
        & {$\|u_h^\text{ad}(t) - u_\text{ex}(t)\|_{DG}$}   \\
        \hline
        $t=5 \, \textrm[ms]$ &$1.63e-3$&$1.15e-2$&$1.62e-3$&$1.14e-2$ \\
        $t=20.3\, \textrm[ms]$ &$5.43e-3$&$2.54e-2$&$5.50e-3$&$2.55e-2$\\
        \hline
    \end{tabular}
\end{table}
The element-wise polynomial degree adapted with our $p$-adaptive algorithm ($\tau_\text{threhsold} = 3.25e-3$) every 5 time steps of the Crank-Nicholson time advancing scheme is shown in Figure~\ref{fig:: single_wave} (fifth row) together with corresponding local error indicator defined as in Equation \eqref{eq::indicator} (Figure~\ref{fig:: single_wave}, sixth row). The results show that the proposed $p$-adaptive algorithm is able to track the wavefront correctly. The distribution of the local polynomial degree for each element of the mesh over time is consistent with the position of the wavefront. In particular, we observe that high-order polynomial degrees are automatically set by our $p$-adaptive algorithm to correctly capture the steep wavefront, while low order elements are employed where the solution is flat.
\begin{figure}[!htbp]
    \centering
    \begin{subfigure}[t]{0.3\textwidth}
\includegraphics[width=1\linewidth]{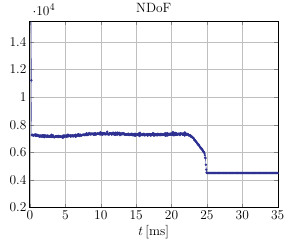}
\caption{\label{fig:dof_evolution_tv}}
    \end{subfigure}
        \begin{subfigure}[t]{0.3\textwidth}
\includegraphics[width=1\linewidth]{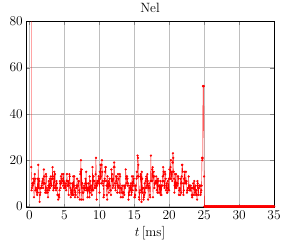}
\caption{\label{fig:comptau_singletv}} 
    \end{subfigure}
    \begin{subfigure}[t]{0.3\textwidth}
\includegraphics[width=1\linewidth]{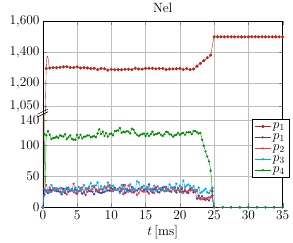}
\caption{\label{fig:elem_tv}} 
    \end{subfigure}
\caption{Test case 1b. Left: Evolution of the number of degrees of freedom ($\mathrm{NDoF}$) as a function of time, driven by our $p$-adaptive algorithm ($\tau_\text{threhsold} = 3.25e-3$). At the initial time $\mathrm{NDoF}=3.15e+4$. Center: Evolution of number of elements ($\mathrm{Nel}$) where the local polynomial degree is updated over time. 
Right: Evolution of the number of elements ($\mathrm{Nel}$) with local polynomial degrees equal to $1,2,3,4,5$ over time. As initial condition, all the elements are discretized with $p_\text{max}=5$; 
} \label{fig:NDoF_evolution_single}
\end{figure}

The evolution of the total $\mathrm{NDoF}$ as a function of time resulting from our $p$-adaptive algorithm is shown in Figure~\ref{fig:dof_evolution_tv}.
Starting at the initial time from a uniform polynomial degree distribution $p_\text{max}=5$, the $p$-adaptive algorithm is able to significantly reduce the total number of degrees of freedom while preserving the accuracy, resulting in a gain of 78\% compared to the uniform approximation. 
We also observe that, once the wave has passed, all elements are discretized with linear elements, leading to $4500$ degrees of freedom, resulting in a gain of 86\%.
Figure~\ref{fig:comptau_singletv} shows the evolution of the total number of modified elements ($\mathrm{Nel}$) over time, while Figure~\ref{fig:elem_tv} shows the evolution over time of the number of elements with local polynomial degree equal to $1,2,3,4,5$. As expected, most of the elements employ linear elements, except for a few elements in the proximity of the wavefront, where high-order polynomials are employed to correctly describe steep variations of the solution.
Exploiting the proposed $p$-adaptive algorithm, we observe a consistent reduction of the total number of degrees of freedom needed to discretize the system while maintaining accuracy in the approximation of the solution, as is evident from the results shown Table~\ref{tab:errorL2DGtable} where the computed errors in the $L2$  and $DG$  norms obtained with the $p$-adaptive scheme and with a uniform polynomial degree $p_{\textrm{un}}=5$ are compared at two different time snapshots  $t=5, 20.3 \, [\mathrm{ms}]$.

\FloatBarrier

\subsubsection{Test case 1c: comparison of different choices for the local error indicators}
\label{sec::comparison}
\begin{figure}[!htbp]
    \centering
    \begin{subfigure}[t]{0.45\textwidth}
    \centering
\includegraphics[width=0.8\linewidth]{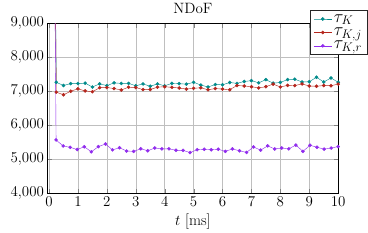}
\caption{\label{fig:dof_comparison}}
    \end{subfigure}
        \begin{subfigure}[t]{0.45\textwidth}
        \centering
\includegraphics[width=0.78\linewidth]{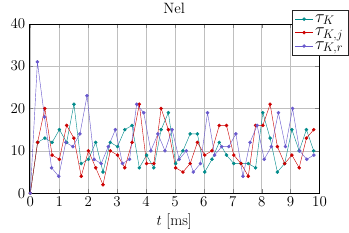}
\caption{\label{fig:modelem_comparison}}  
    \end{subfigure}
    \caption{Test case 1c. Time evolution of the total number of degrees of freedom (left) and of the total number of elements marked for $p$-refinement/de-refinement  (right) for the three different choices of the local error indicators driving the $p$-adaptive algorithm.}\label{fig:comparison}
\end{figure}

Next, we investigate different choices of local error indicators in driving our $p$-adaptive strategy. Specifically, we compare:
\begin{enumerate}
\item the "full" local error indicator $\tau_{K}$ as defined in Equation \eqref{eq::indicator},
\item the local error indicator relating only to the residual of the numerical solution, i.e., $\tau_{K,r}$,
\item the local error indicator related to the jump of the numerical solution, i.e., $\tau_{K,j}$.
\end{enumerate}
Figure~\ref{fig:comparison} shows the time evolution (with final time $T=10\, \mathrm{[s]}$) of the total number of degrees of freedom (left) and of the total number of elements marked for $p$-refinement/de-refinement when our adaptive algorithm is driven by the three different local error indicators (right). 
Our results indicate that the number of degrees of freedom resulting from using the residual indicator seems to be consistently smaller than the ones obtained by exploiting the other two indicators. 

\begin{figure}[!h]
\centering
\setlength{\tabcolsep}{2pt}
\renewcommand{\arraystretch}{1.0}
\begin{tabular}{c c c c}
 & 
\multicolumn{3}{c}{%
  \hspace{-4ex}\includegraphics[trim={0 0cm 0 0cm},clip, scale=0.16]{singlewave_stationary/scale_.png}} \\
\raisebox{7mm}{} &
\includegraphics[trim={2cm 13cm 2cm 13cm},clip, scale=0.07]{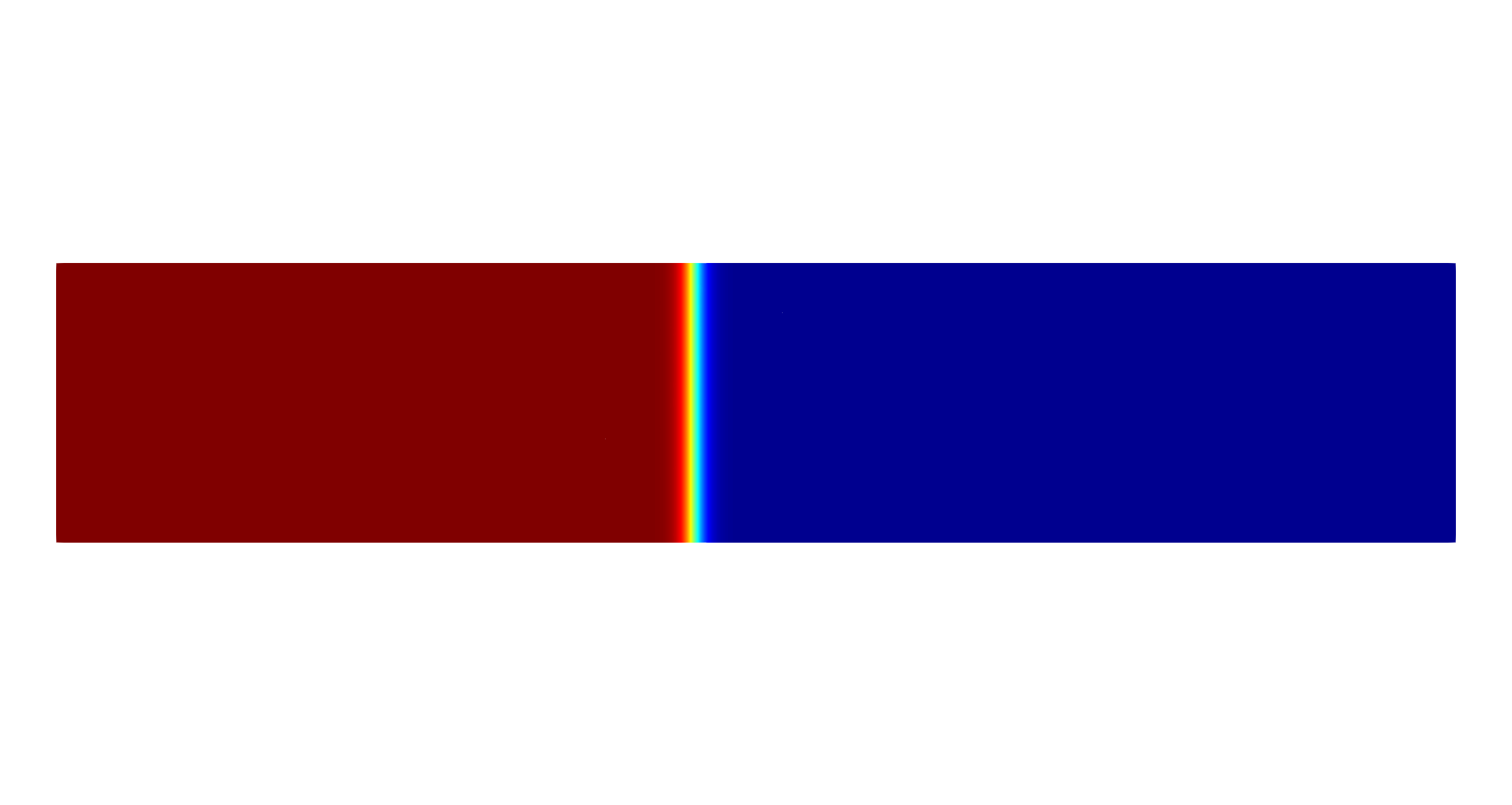} &
\includegraphics[trim={2cm 13cm 2cm 13cm},clip, scale=0.07]{comparisontau/sol.png} &
\includegraphics[trim={2cm 13cm 2cm 13cm},clip, scale=0.07]{comparisontau/sol.png} \\
 & 
\multicolumn{3}{c}{%
  \includegraphics[trim={0 0cm 0 0cm},clip, scale=0.16]{singlewave_stationary/scale_p.png}} \\
\raisebox{7mm}{} &
\includegraphics[trim={2cm 13cm 2cm 13cm},clip, scale=0.07]{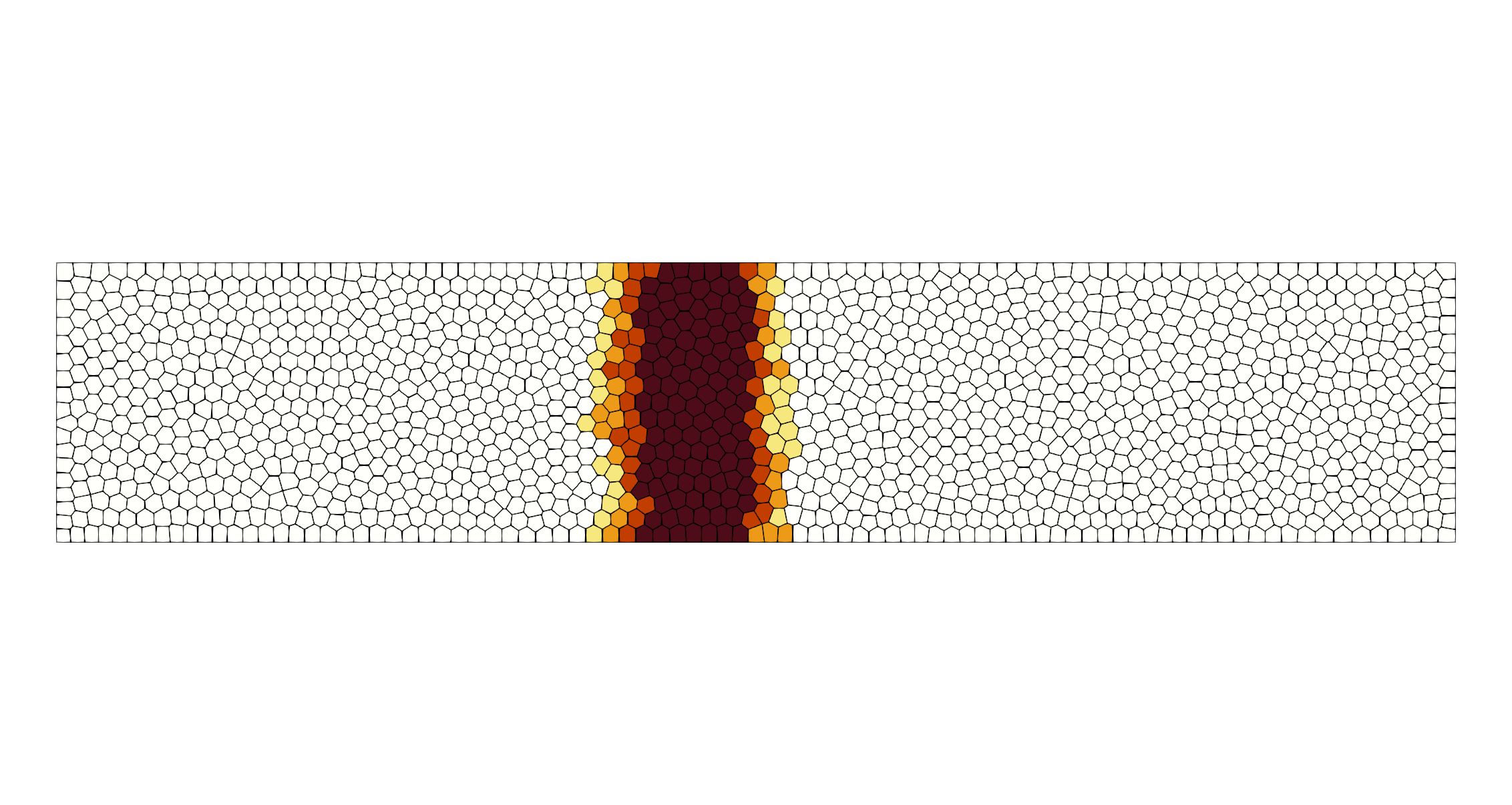} &
\includegraphics[trim={2cm 13cm 2cm 13cm},clip, scale=0.07]{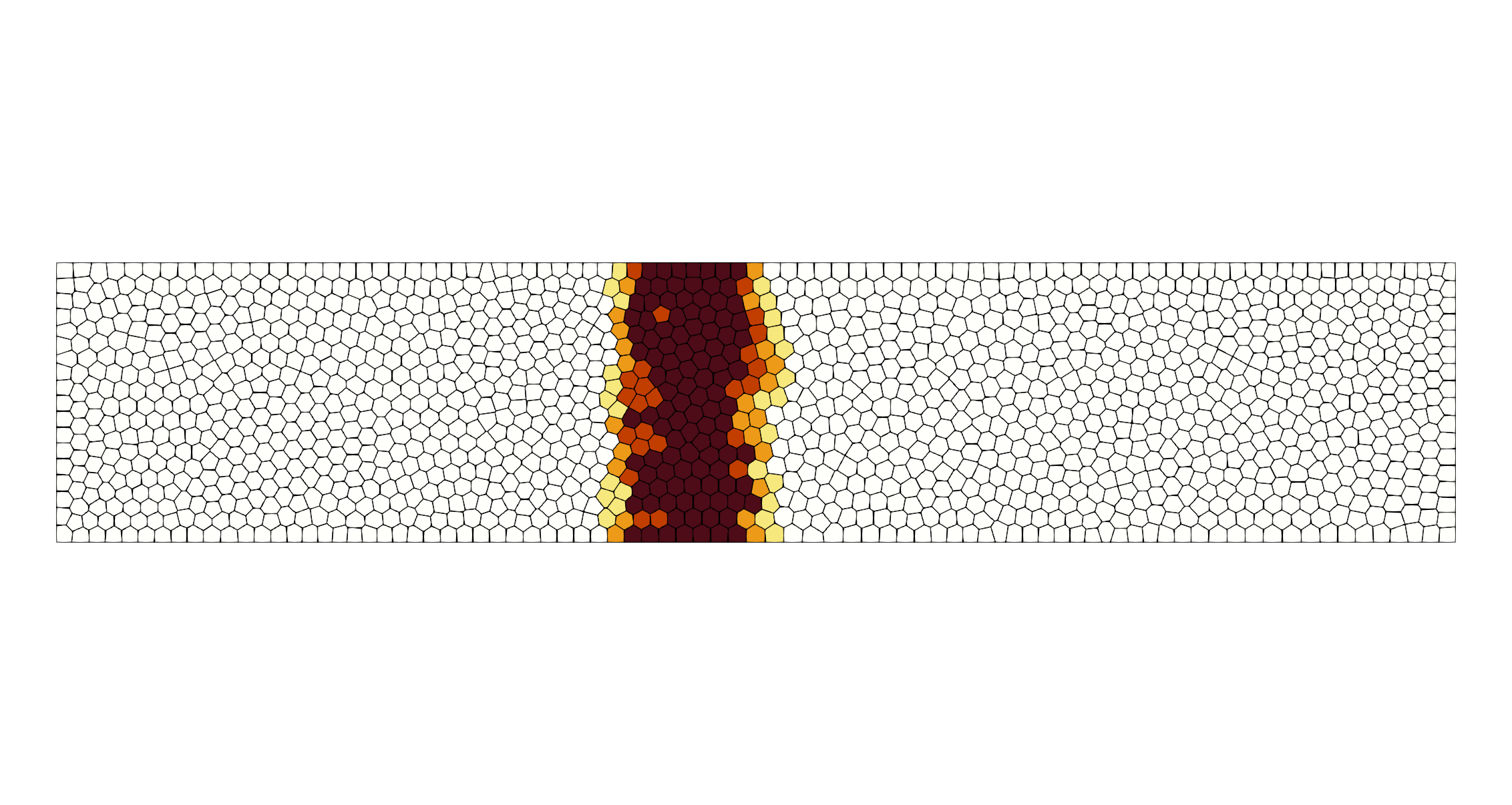} &
\includegraphics[trim={2cm 13cm 2cm 13cm},clip, scale=0.07]{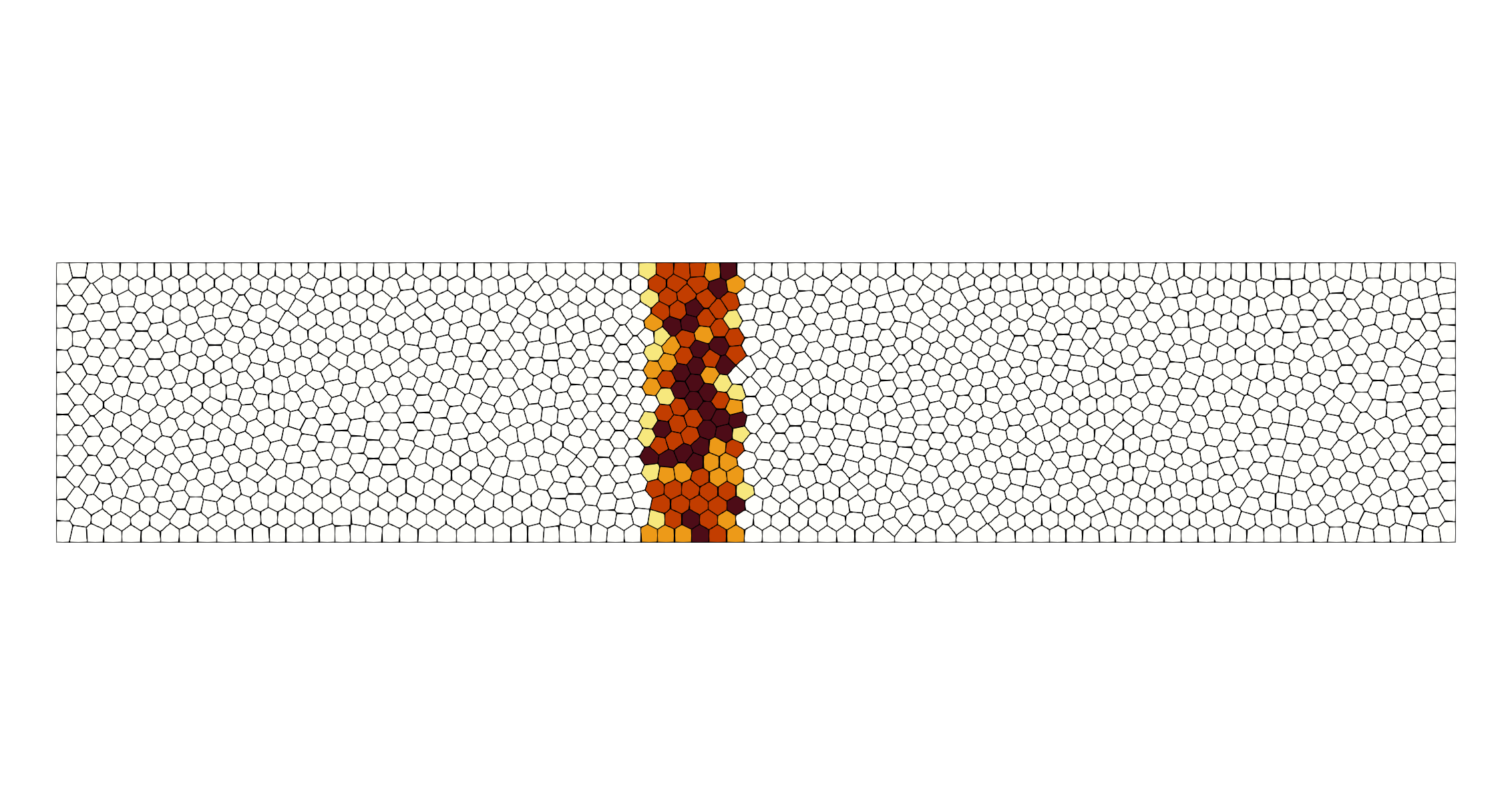} \\
&
\multicolumn{3}{c}{
  \includegraphics[trim={0 0cm 0 0cm},clip, scale=0.16]{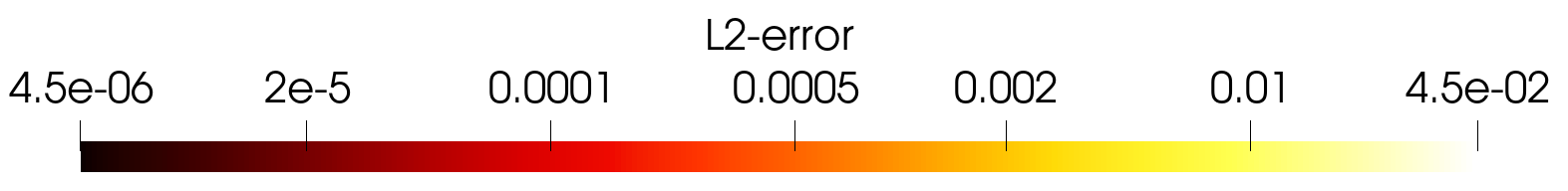}
}\\
\raisebox{7mm}{} &
\includegraphics[trim={2cm 13cm 2cm 13cm},clip, scale=0.07]{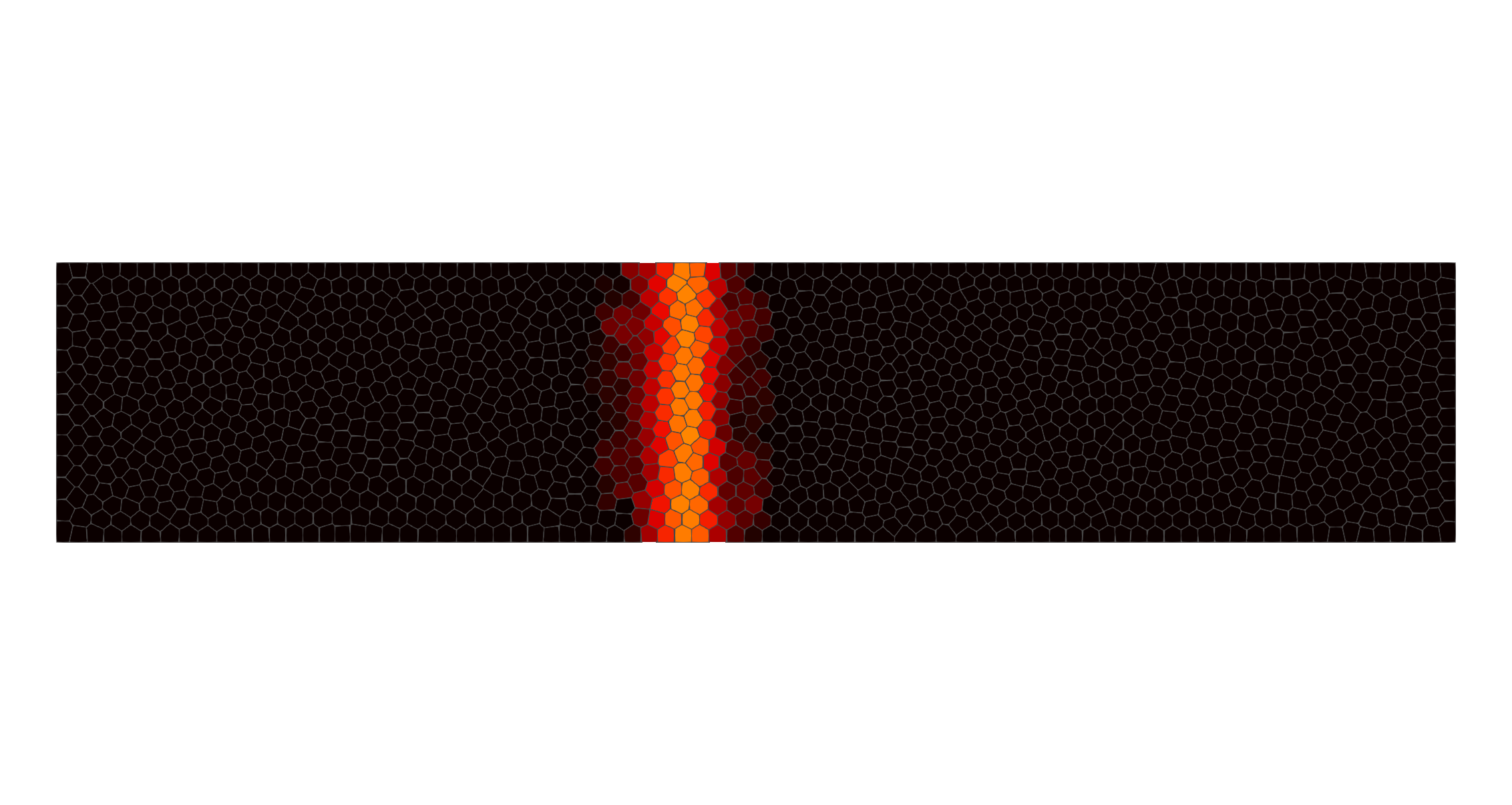} &
\includegraphics[trim={2cm 13cm 2cm 13cm},clip, scale=0.07]{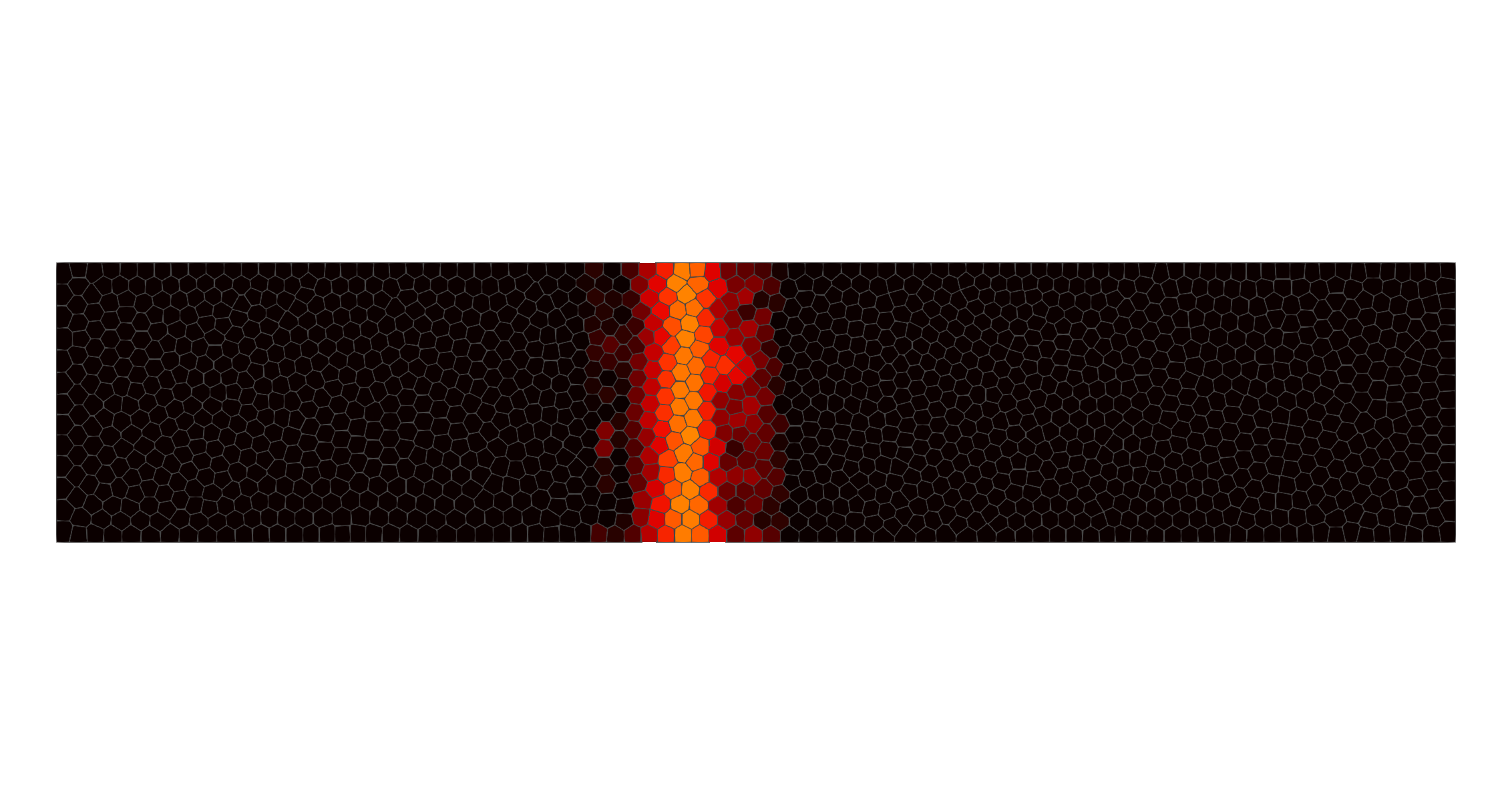} &
\includegraphics[trim={2cm 13cm 2cm 13cm},clip, scale=0.07]{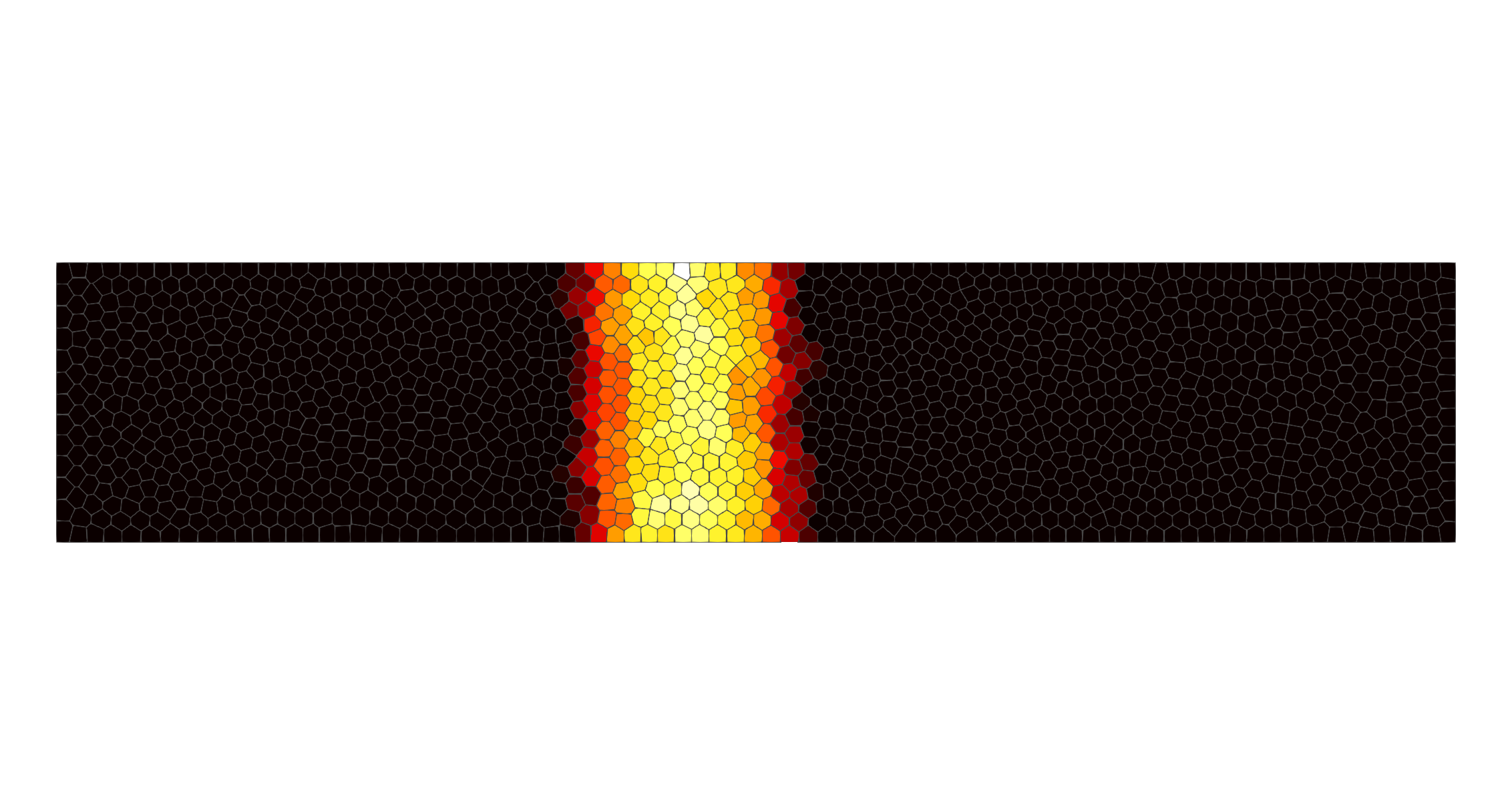} \\
&
\multicolumn{3}{c}{
     \includegraphics[trim={0 0cm 0 0cm},clip, scale=0.18]{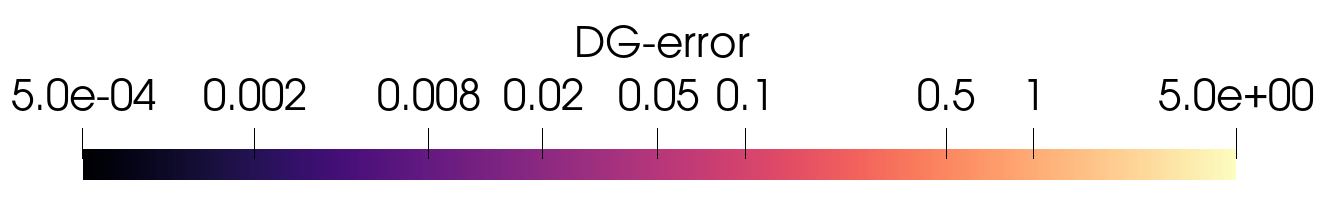}
} \\
\raisebox{7mm}{} &
\includegraphics[trim={2cm 13cm 2cm 13cm},clip, scale=0.07]{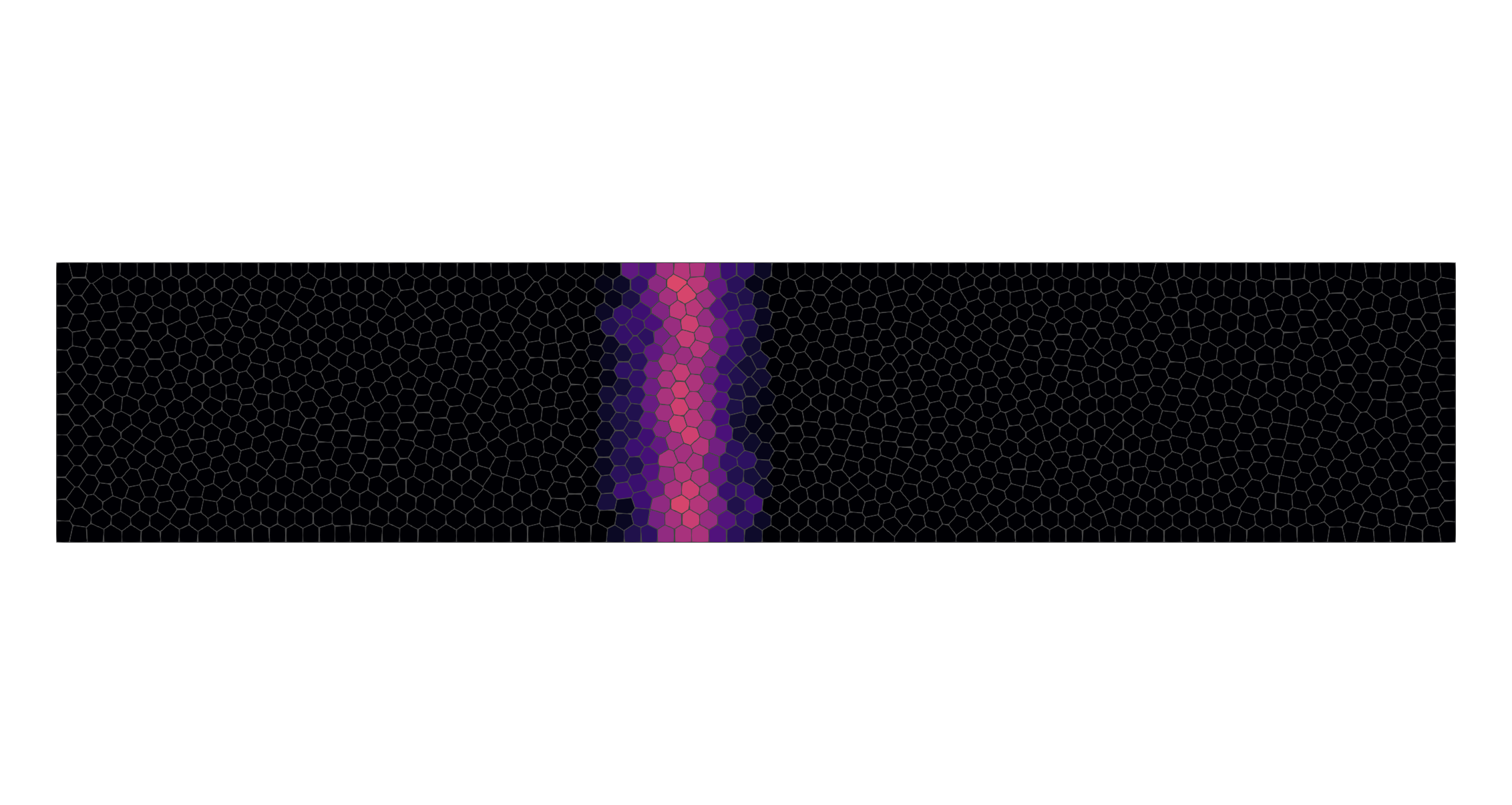} &
\includegraphics[trim={2cm 13cm 2cm 13cm},clip, scale=0.07]{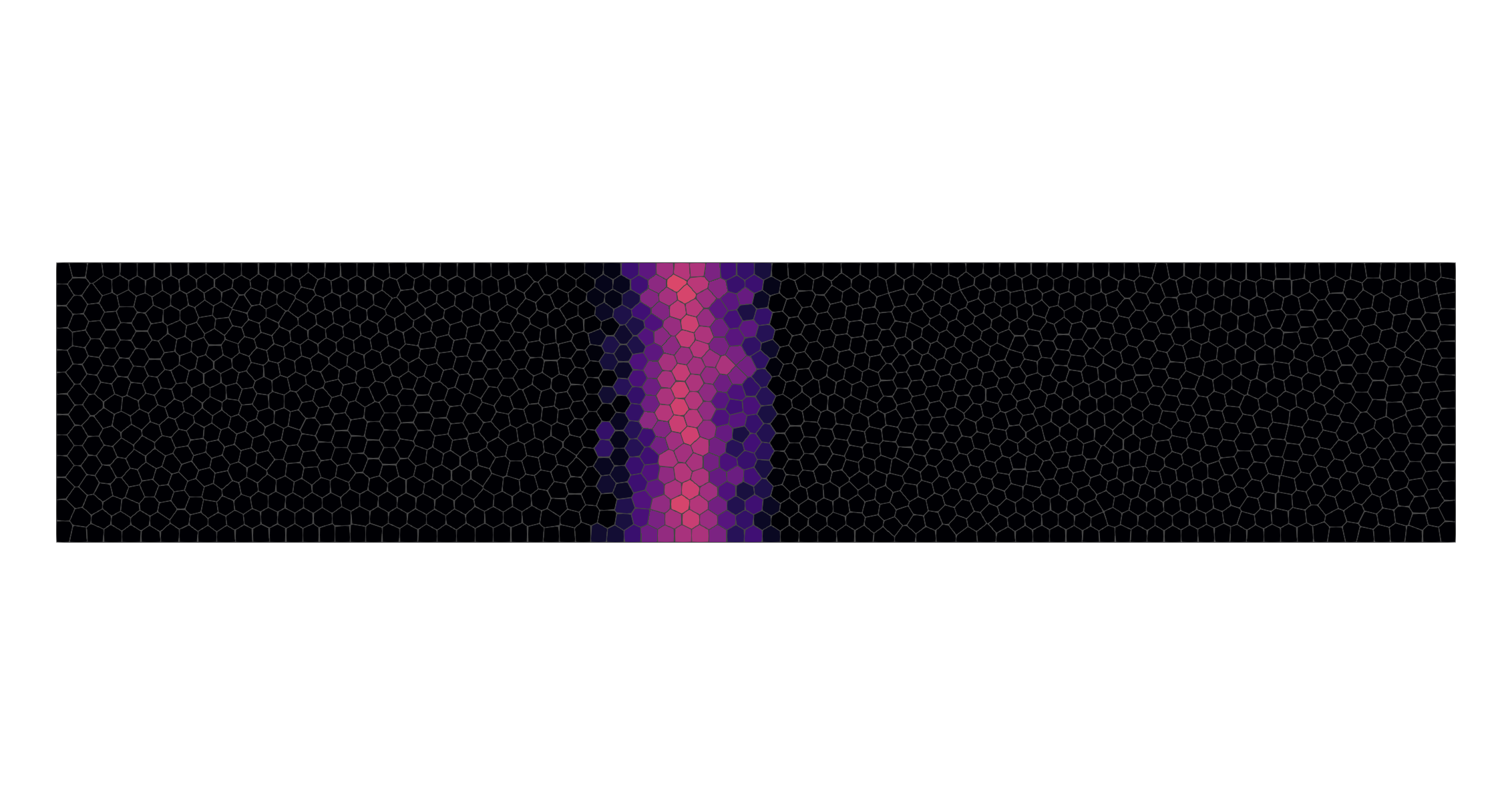} &
\includegraphics[trim={2cm 13cm 2cm 13cm},clip, scale=0.07]{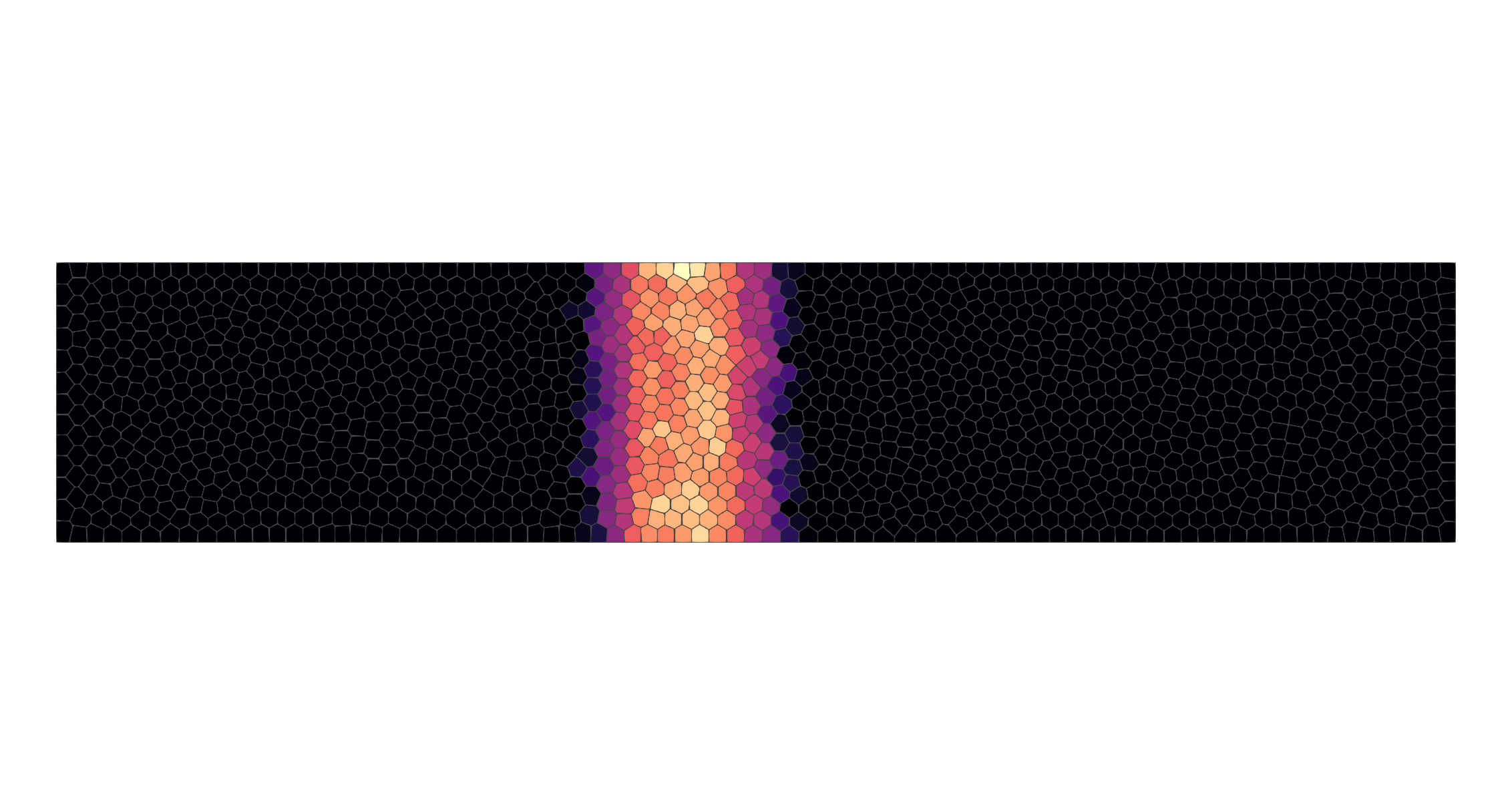} \\
& \includegraphics[trim={0 0cm 0 0cm},clip, scale=0.16]{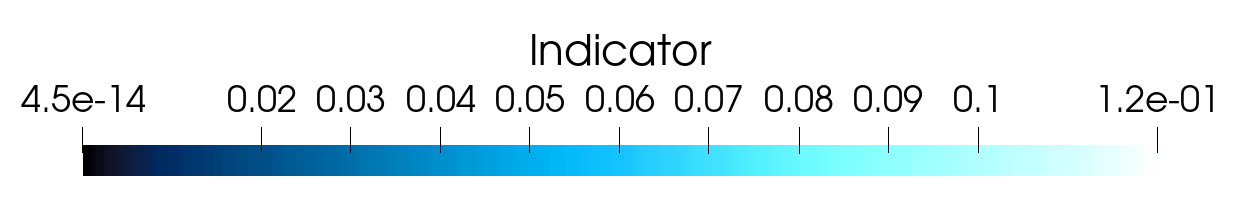}
 &  \includegraphics[trim={0 0cm 0 0cm},clip, scale=0.16]{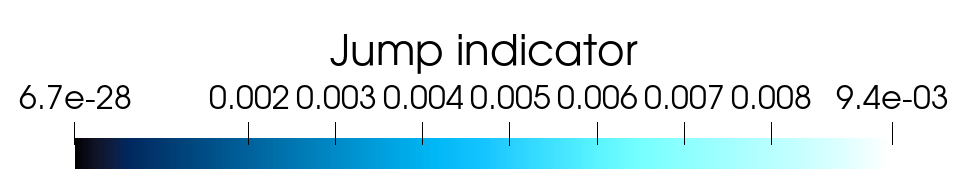}  &   \includegraphics[trim={0 0cm 0 0cm},clip, scale=0.16]{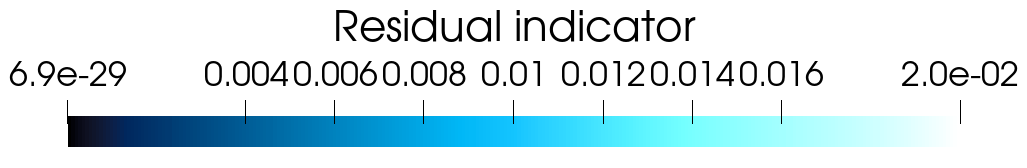} \\
\raisebox{7mm}{} &
  \includegraphics[trim={2cm 13cm 2cm 13cm},clip, scale=0.07]{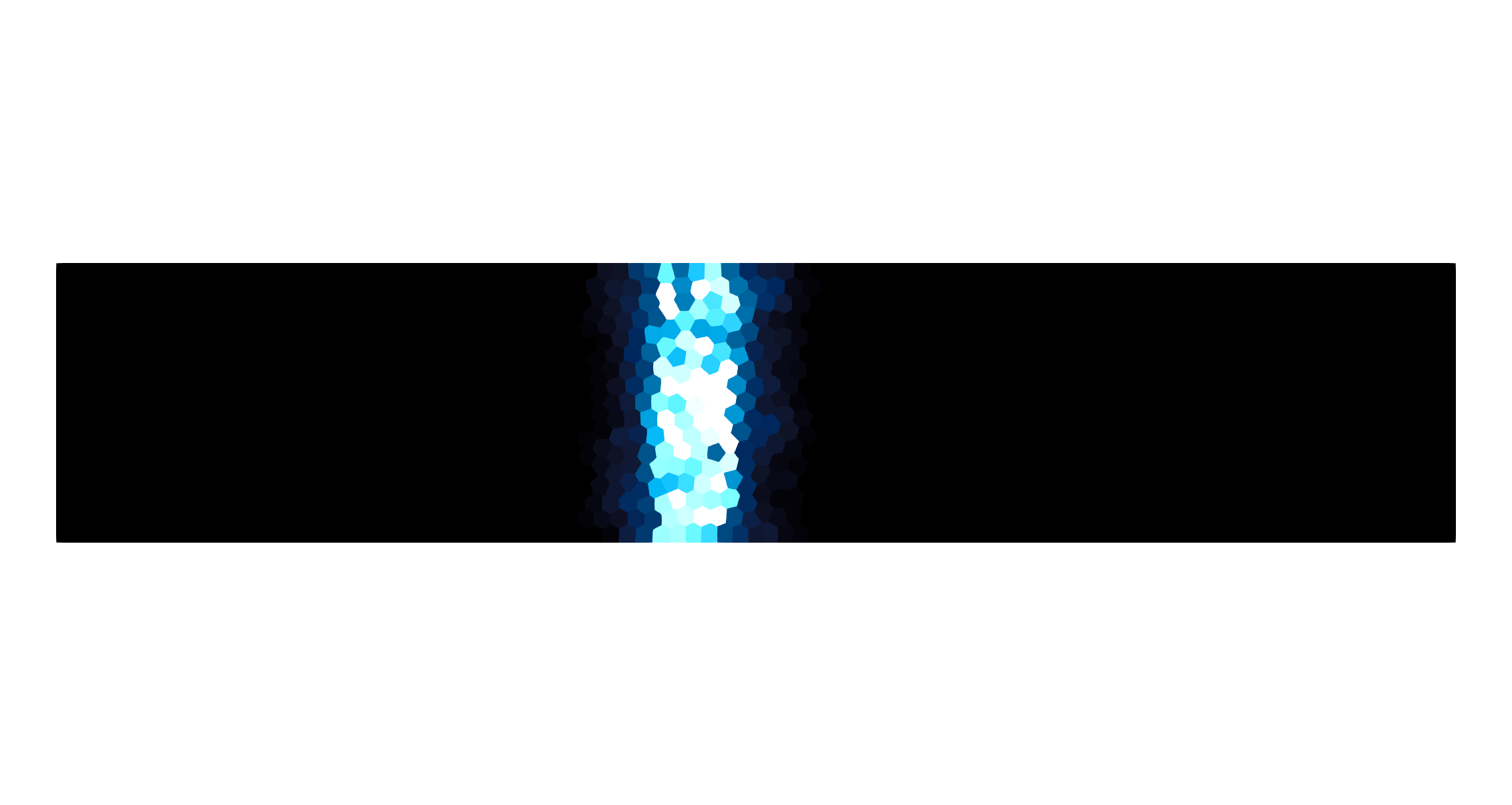}
&
  \includegraphics[trim={2cm 13cm 2cm 13cm},clip, scale=0.07]{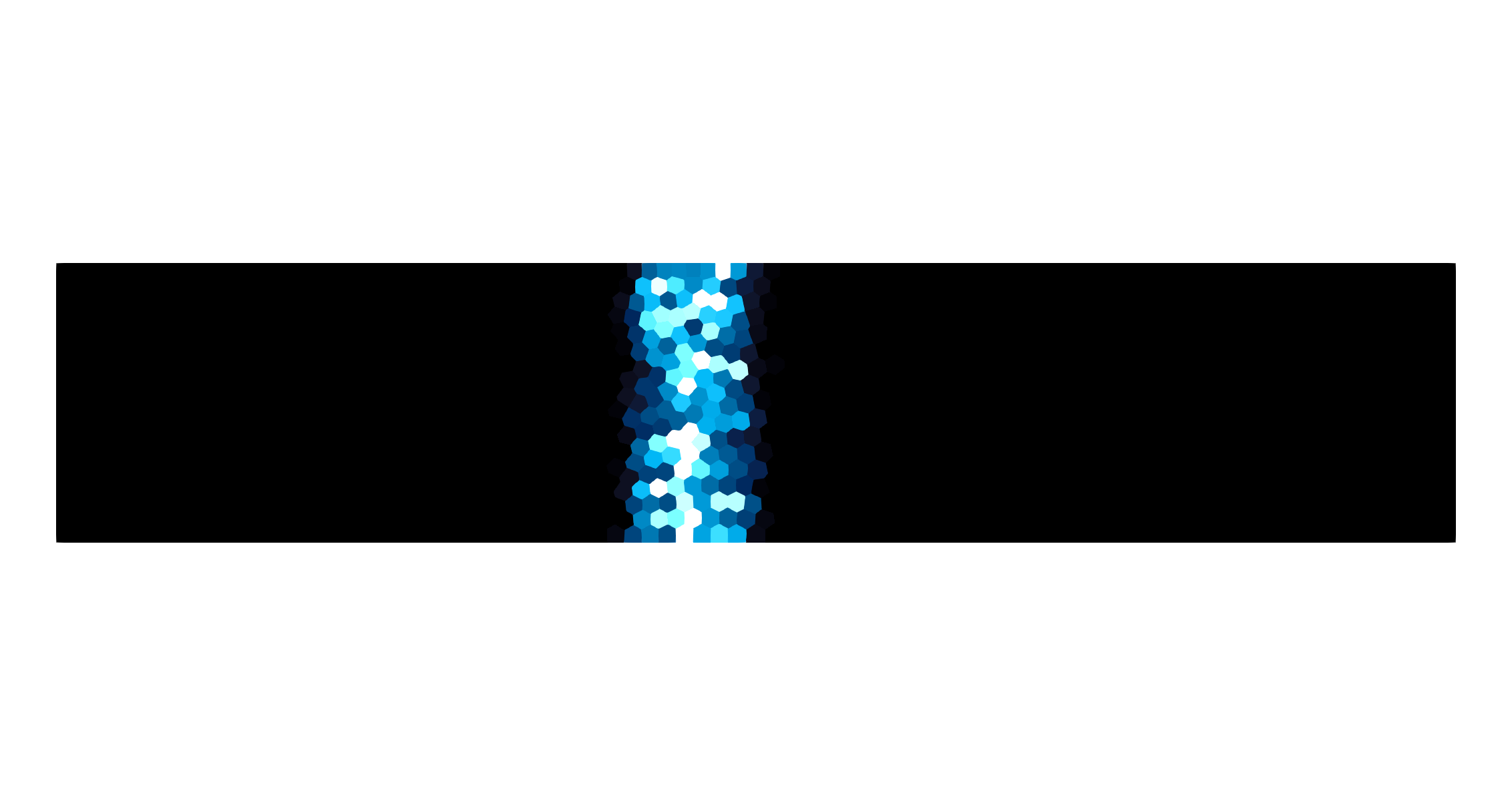}
&
  \includegraphics[trim={2cm 13cm 2cm 13cm},clip, scale=0.07]{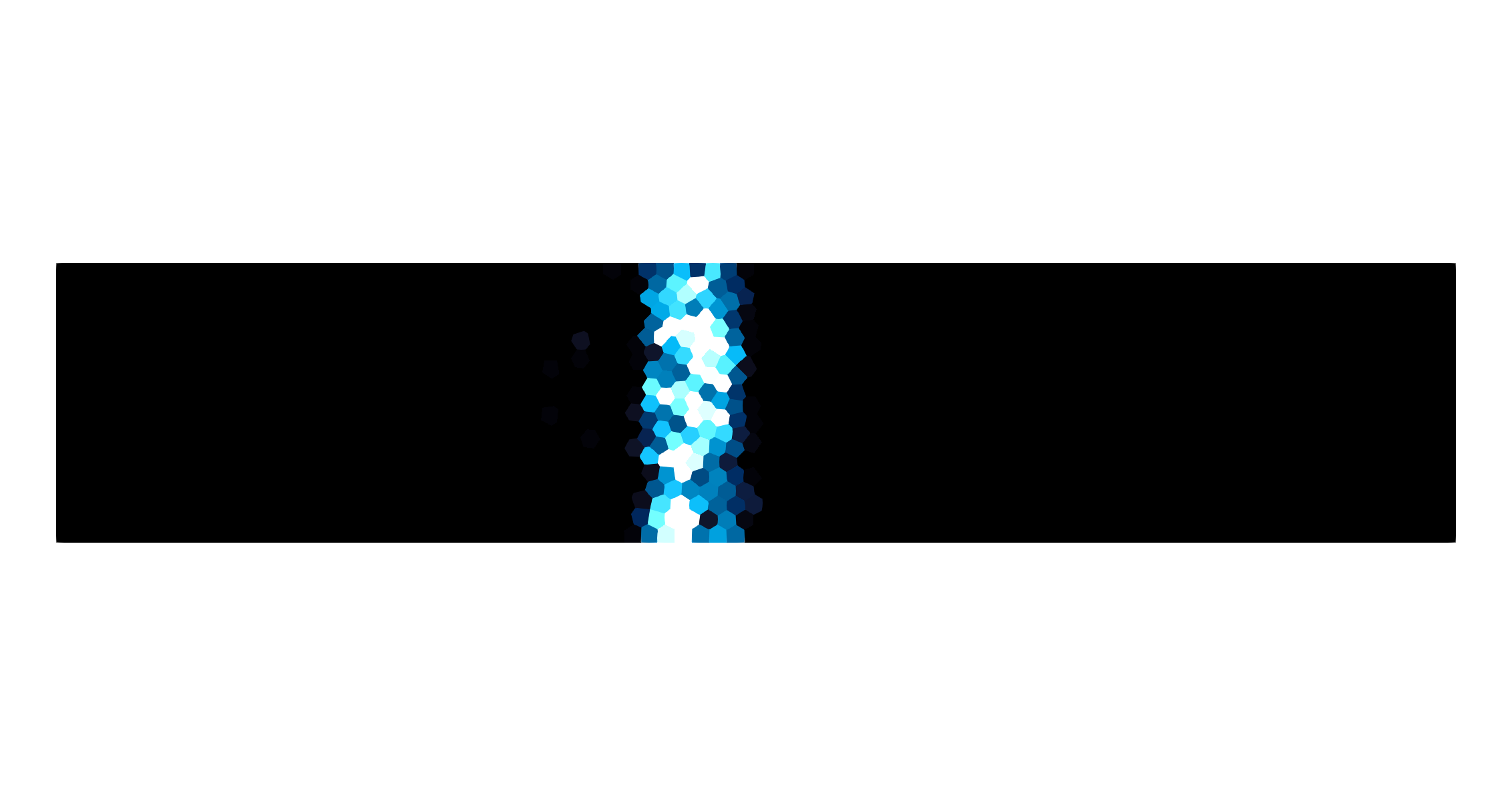}\\
 & (a) $\tau_K$ & (b) $\tau_{K,j}$ & (c) $\tau_{K,r}$
\end{tabular}
\caption{Test case 1c. From left to right: different local error indicators are employed to drive the $p$-adaptive algorithm. 
First row: numerical solution at $t = 2 \, [\mathrm{ms}]$; Second row: element-wise polynomial degree distribution; third and fourth rows: element-wise errors in the $L2$  and $DG$ norms; Fifth row: computed local error indicators.}
\label{fig::comparison_tau}
\end{figure}
\begin{figure}[htbp!]
    \centering
\begin{subfigure}[t]{0.3\textwidth}
    \centering
\includegraphics[width=1.05\linewidth]{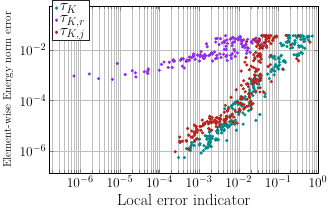}
\caption{$t=5 \, \textrm[ms]$\label{fig::tab1}}
\end{subfigure}\hfill
\begin{subfigure}[t]{0.3\textwidth}
    \centering
\includegraphics[width=1.05\linewidth]{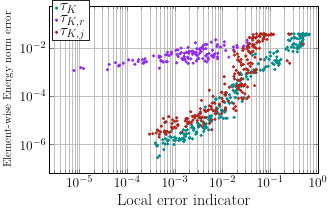}
\caption{$t=5 \, \textrm[ms]$\label{fig::tab2}}
\end{subfigure}\hfill
\begin{subfigure}[t]{0.3\textwidth}
    \centering
\includegraphics[width=1.05\linewidth]{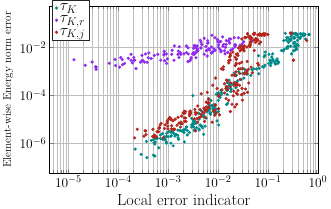}
\caption{$t=5 \, \textrm[ms]$\label{fig::tab3}}
\end{subfigure}
\caption{Test case 1c. Scatterplot of the element-wise error computed in energy norm defined in \eqref{eq::energy_norm} (y-axis) versus the local error indicators (x-axis) for three different time snapshots $t=1, 5, 10 \, \mathrm{[ms]}$. }
\label{tab:indicator_tab2}
\end{figure}

Figure~\ref{fig::comparison_tau} shows a comparison between the computed solution at time $t=2\, \mathrm{[ms]}$ (first row) as well as the element-wise polynomial approximation degree for the three choices of local error indicators (second row). In the same Figure (third and fourth rows) we report the computed errors in the $L2$ and $DG$ norms together with the element-wise values of the three error indicators (last row). From these results, we can deduce that the local polynomial degree distribution associated with the local residual indicator $\tau_{K,r}$ entails a lower number of degrees of freedom than the other two indicators. However, this choice leads to a less accurate numerical solution, as clear from the results reported in the third and fourth rows, where the errors computed based on employing $\tau_{K,r}$ are consistently larger. This suggests that the residual indicator is underestimating the error and, as a consequence, lowering the local polynomial degree even in at the wavefront, compromising the numerical scheme accuracy. 
The "full" and jump local error indicators outperform the residual indicator: indeed, we observe that they can correctly identify the wave, maintaining the required accuracy, as shown in from Figures \ref{tab:indicator_tab2}. 
We also observe that for the adaptive method based on the complete indicator $\tau_{K}$, the traveling wave is completely included in the band of elements where high-order approximations are employed, and the elements are then gradually scaled down to linear elements in correspondence with quiescence regions.  Furthermore, comparing the results obtained making use of the complete indicator ($\tau_K$) and the jump indicator ($\tau_{K,j}$), the time evolution of the total number of degrees of freedom is comparable and the error obtained by the adaptive algorithm considering $\tau_{K,j}$ as local error indicator seems to be slightly larger than the corresponding results obtained based on employing $\tau_K$. 

Figure~\ref{tab:indicator_tab2} illustrates the correlation between the error computed in the energy norm defined in Equation \eqref{eq::energy_norm} and three distinct a posteriori indicators.
More precisely, 
In Figure~\ref{tab:indicator_tab2} we report the element-wise scatterplot of the local error computed in energy norm (y-axis) when the $p$-adaptive algorithm is driven by the different local errors indicators versus the value of the three local error indicators (x-axis) for different time snapshots $t=1, 5, 10 \, \mathrm{[ms]}$. 
From the results reported in  Figure~\ref{tab:indicator_tab2}, the complete indicator $\tau_K$ proves to be the most effective indicator, as it exhibits a clear linear relation with the error computed in the energy norm, with values consistently bounded from below by the corresponding error.
This behavior suggests that $\tau_K$ provides a robust and reliable local error estimator and can be exploited within the $p$-adaptive algorithm. In fact, exploiting the complete indicator we obtain a $L2$ error of $3.82e-3$ and $\mathrm{DG}$-error of $4.12e-2$. Conversely, as expected, the residual indicator $\tau_{K, r}$ confirms its poor effectiveness: it is characterized by significant variance, lacks a clear correlation with the error, and frequently underestimates it, obtaining a $L2$ error of $1.05e-1$ and $\mathrm{DG}$-error of $1.32e+1$.  
Finally, the jump indicator $\tau_{K, j}$ exhibits an intermediate behavior, with a discernible positive correlation with the error and moderate variance; the values are generally of the same order of magnitude as $\tau_{K}$ but the indicator does not achieve the same level of accuracy and tends to underestimate the error, especially for lower error magnitudes. Here we obtain a $L2$ error of $3.85e-3$ and $\mathrm{DG}$-error of $4.74e-2$.
\subsection{Test case 2: "double" traveling wave solution}
In the second test case, we numerically investigate the ability of the $p-$adaptive algorithm to approximate multiple wavefronts. We consider a domain $\Omega = (-2,4) \times (-0.5,0.5)$ and two starting wavefronts defined through the following initial condition: 
\begin{equation}
    \begin{aligned}
      u^0(\boldsymbol{x}) =
      \frac{V_{\text{dep}}-V_{\text{rest}}}{2}\left(\tanh\left(\frac{\boldsymbol{x} + x_2}{\epsilon_2}\right) -  \tanh\left(\frac{\boldsymbol{x} + x_1}{\epsilon_1}\right) \right)  + V_\text{dep}.
\end{aligned}
\label{eq:initial_doublewave}
\end{equation}
We consider two configurations: a homogeneous case (see Section ~\ref{sec:HomogCase}), in which the wavefronts have the same velocity and steepness, and an heterogeneous case, in which the domain is split by assigning different conductivity values that cause wavefronts to propagate at different velocities (see Section ~\ref{sec:HeteroCase}).

\subsubsection{Test case 2a: homogeneous case} \label{sec:HomogCase}
\begin{figure}[h!]
    \centering
    \begin{subfigure}[t]{0.3\textwidth}
        \centering
\includegraphics[width=1\linewidth]{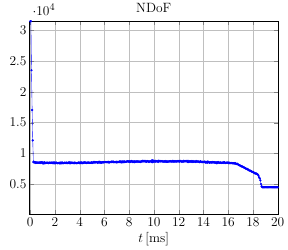}
\caption{\label{fig:dof_evolution_doubletv}} 
    \end{subfigure}
    \begin{subfigure}[t]{0.3\textwidth}
        \centering
\includegraphics[width=1\linewidth]{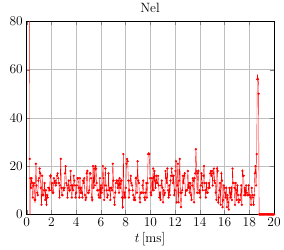}
\caption{\label{fig:comptau_doubletv}} 
    \end{subfigure}
\begin{subfigure}[t]{0.3\textwidth}
        \centering
\includegraphics[width=1.1\linewidth]{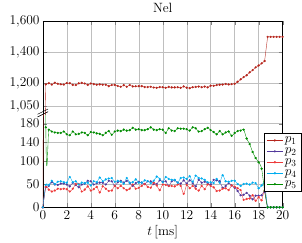}
\caption{\label{fig:elem_doubletv}} 
\end{subfigure}
\caption{Test case 2a. Left: Evolution of the number of degrees of freedom ($\mathrm{NDoF}$) as a function of time, driven by our $p$-adaptive algorithm ($\tau_\text{threshold} = 3.25e-3$). At the initial time $\mathrm{NDoF}=3.15e+4$. Center: Evolution of the number of elements ($\mathrm{Nel}$) where the local polynomial degree is updated over time.  Right: Time evolution of the number of elements ($\mathrm{Nel}$) with local polynomial degrees equal to $1,2,3,4,5$. At the initial time all the elements are discretized with $p_\text{max}=5$.} \label{fig:ndof_evolution}
\end{figure}

This test case consists of a simulation in which both waves move with the same velocity $c = 0.1212 \, [\mathrm{mm} \cdot \mathrm{ms}^{-1}]$.
The mesh is characterized by 1500 polytopal elements ($h = 0.1105$);  the initial condition is defined as in Equation \eqref{eq:initial_doublewave}, with $x_1=-3$ and $x_2=1.5$. 
We point out  that the two waves evolve in opposite directions towards the center of the domain at the same speed. 
We have repeated the same analysis of test case os Section ~\ref{dof} considering the complete indicator $\tau_K$ and $\tau_{\text{threshold}}=\tau_\text{threshold}^\text{min}$. Figure~\ref{fig:ndof_evolution} reports the time evolution of the number of degrees of freedom $\mathrm{NDoF}$ as a function of time. 
Starting from $3.15e+4$ $\mathrm{NDoF}$ (corresponding to the case where all the elements are discretized with $p_{\textrm{max}} = 5$) we reach around $8.6e+3$ $\mathrm{NDoF}$ (with a $73\%$ reduction).
The corresponding time evolution of the number of elements ($\mathrm{Nel}$) where the local polynomial degree is updated over time is shown in Figure~\ref{fig:ndof_evolution} (center) together with the time evolution of the number of elements with local polynomial degrees equal to $1, 2, 3, 4, 5$.  We also observe that after the two waves collide (approximately at $15.6 \,[\mathrm{ms}]$), a resting state is reached with a discretization of linear polynomials and $4.5e+3$ $\mathrm{NDoF}$ (corresponding to a $86\%$ reduction). In addition, the cost of computing the local error indicator is decreased significantly due to the localized computations of the indicator for the relevant elements only, as shown in Figure~\ref{fig:comptau_doubletv}. Also in this test case, we observe a significant reduction in the total number of degrees of freedom as a consequence of the $p$-adaptive scheme. 
\begin{table}[h!]
    \centering
    \caption{Test case 2a. Computed errors in the $L2$ and $DG$ norms for different time snapshots $t= 7.1 \, [\mathrm{s}]$ and $t = 16.6 [\mathrm{s}]$. As exact solution we considered the numerical solution computed with a uniform polynomial degree ($p_{\textrm{un}}=5$).}
    \label{tab:errorL2DG}
    \begin{tabular}{c|cc}
    \hline
         & $t= 7.1 \, [\mathrm{s}]$ &  $t = 16.6 \,[\mathrm{s}]$\\
        \hline
        $\|u_h^\text{ad}(t) - u_h^\text{un}(t)\|_{L^2(\Omega)}$ & $7.37e-5$ & $9.12e-5$\\
        $\|u_h^\text{ad}(t) - u_h^\text{un}(t)\|_\text{DG}$ & $1.80e-3$ & $2.04e-3$\\
        \hline
    \end{tabular}
\end{table}
 \begin{figure}[!h]
 \centering
\begin{subfigure}{0.45\textwidth}
\centering
\hspace*{40ex}\includegraphics[scale=0.14]{singlewave_stationary/scale_.png}
\includegraphics[trim={1cm 13cm 1cm 13cm},clip, scale=0.075]{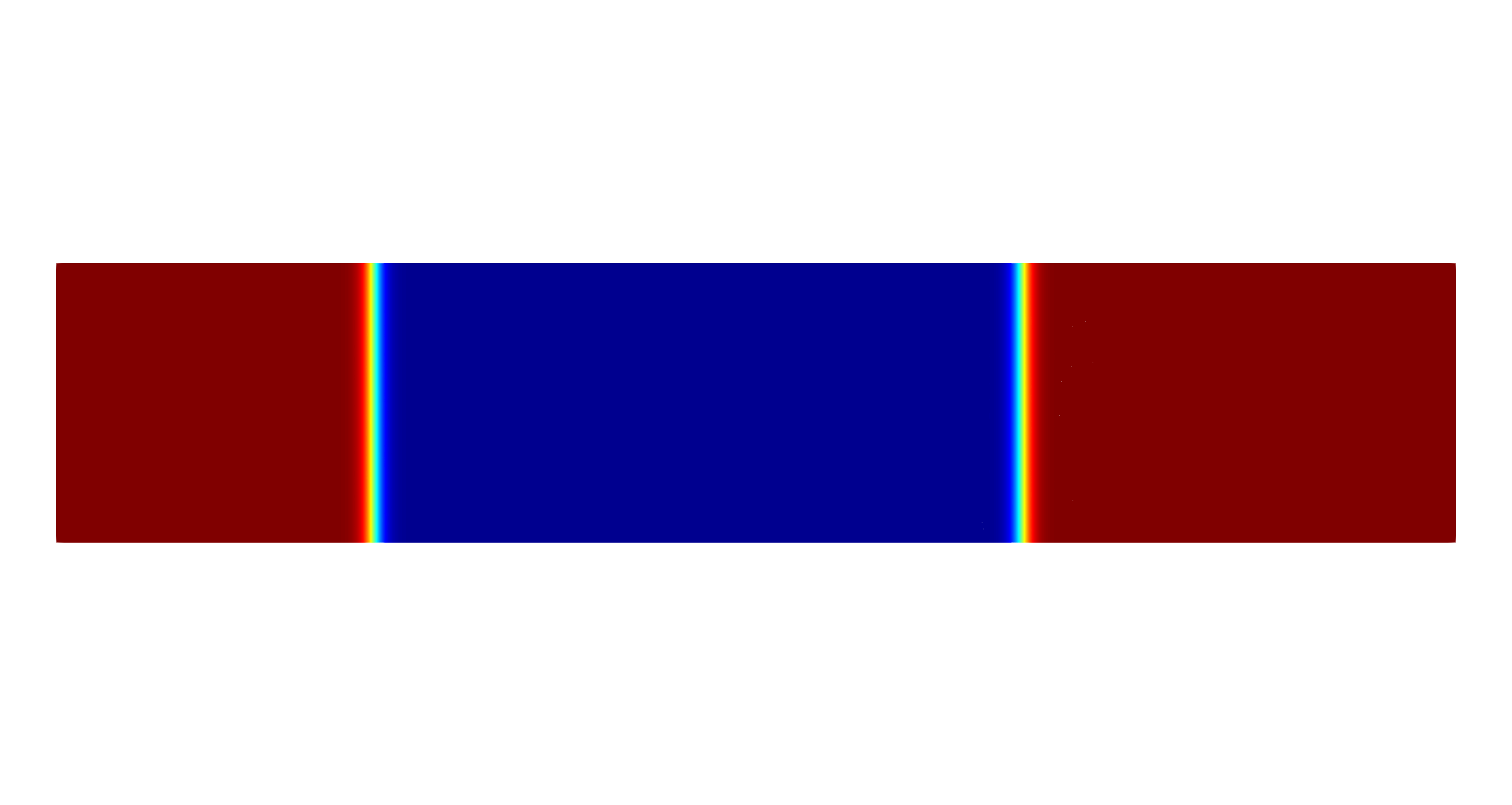}
\includegraphics[trim={1cm 13cm 1cm 13cm},clip, scale=0.075]{doublehomostationary/sol.70.png}
\hspace*{41ex}\includegraphics[trim={0 0cm 0 0cm},clip, scale=0.16]{singlewave_stationary/scale_p.png}
\includegraphics[trim={1cm 13cm 1cm 13cm},clip, scale=0.075]{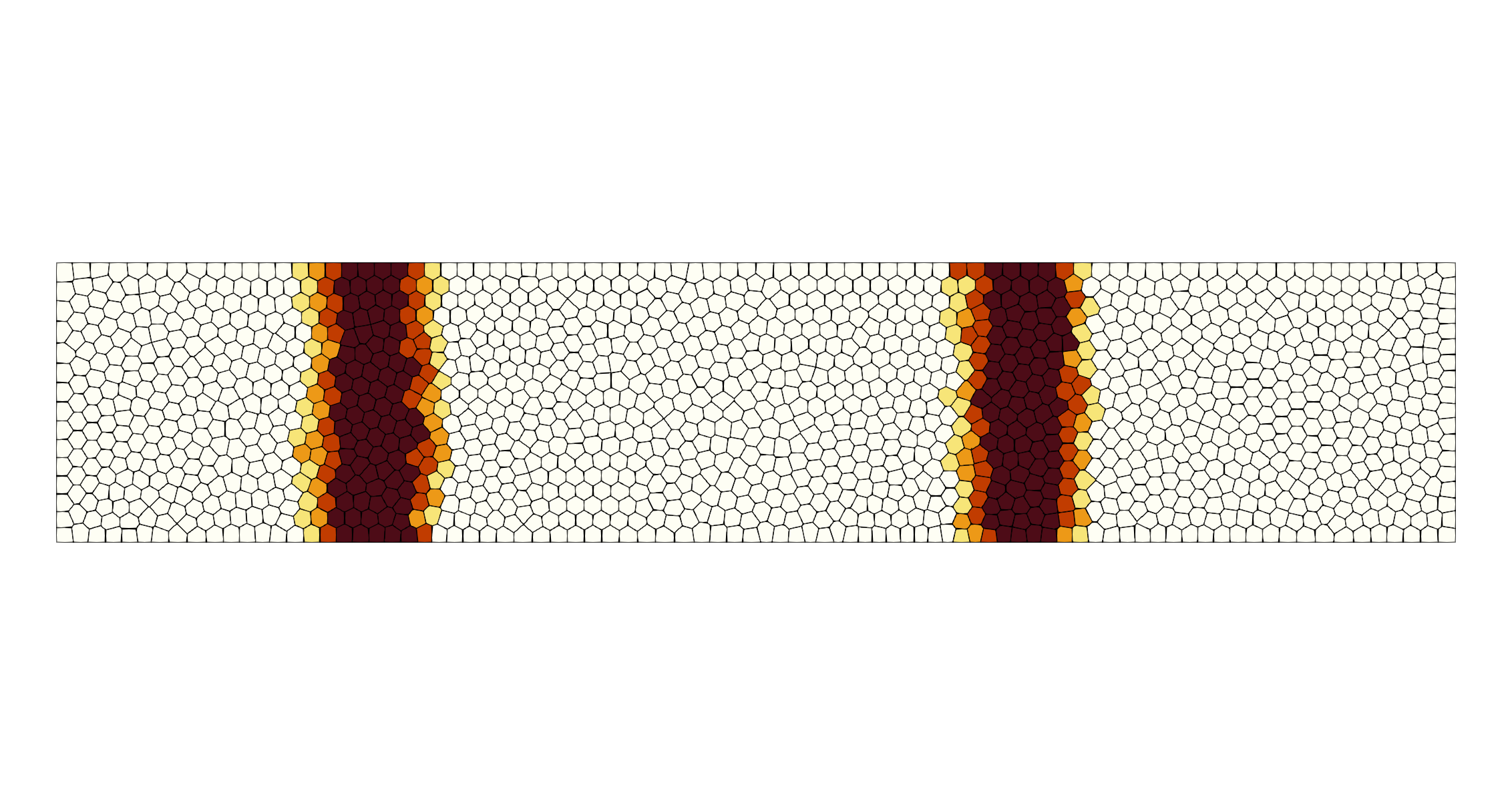}
\hspace*{35ex}\includegraphics[trim={0 0cm 0cm 0cm},clip, scale=0.16]{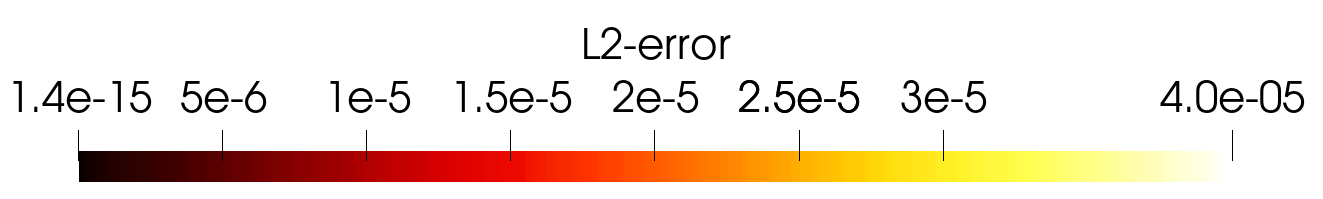}
\includegraphics[trim={1cm 13cm 1cm 13cm},clip, scale=0.075]{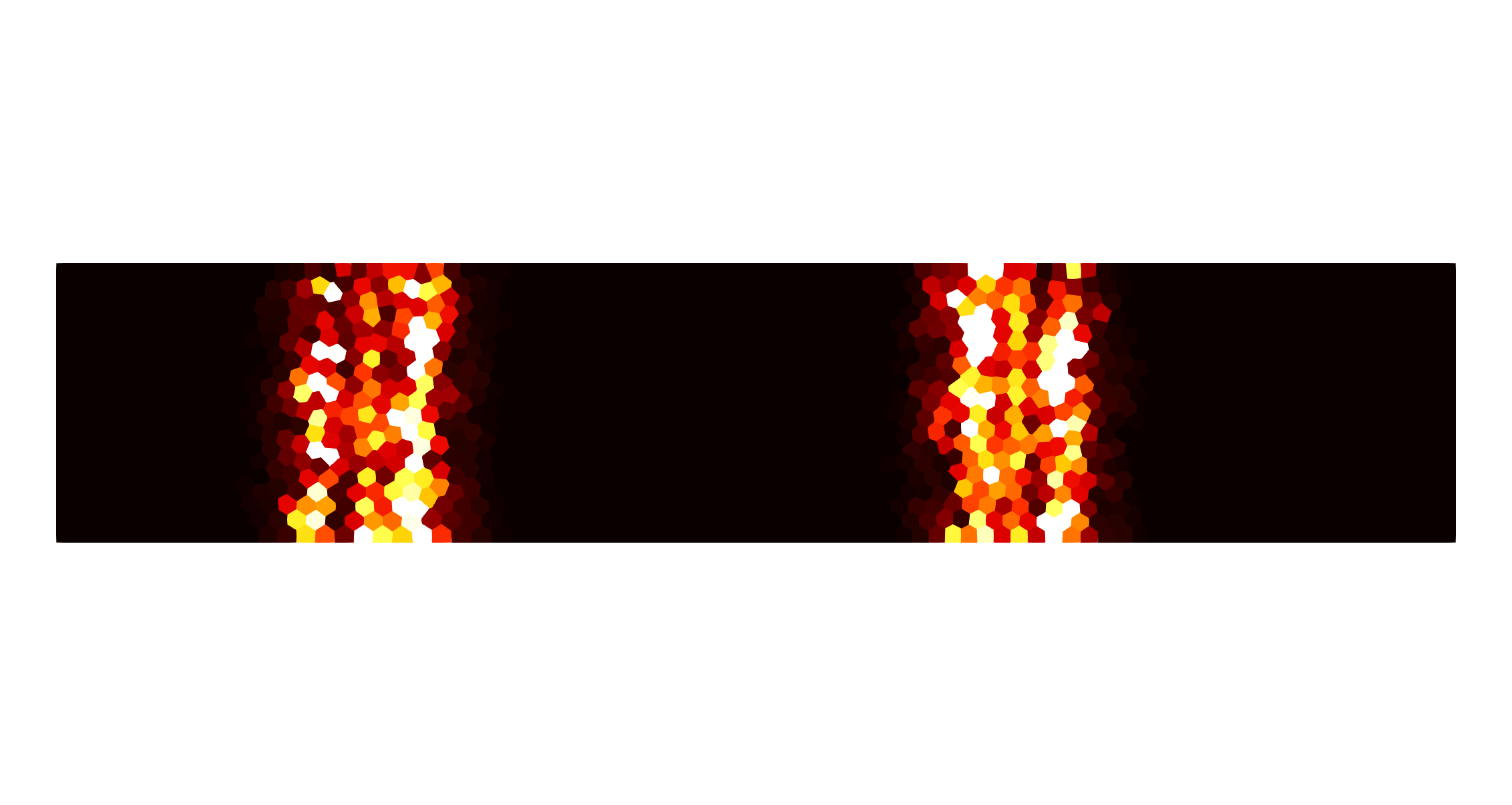} 
\hspace*{35ex}\includegraphics[trim={0 0cm 0 0cm},clip, scale=0.16]{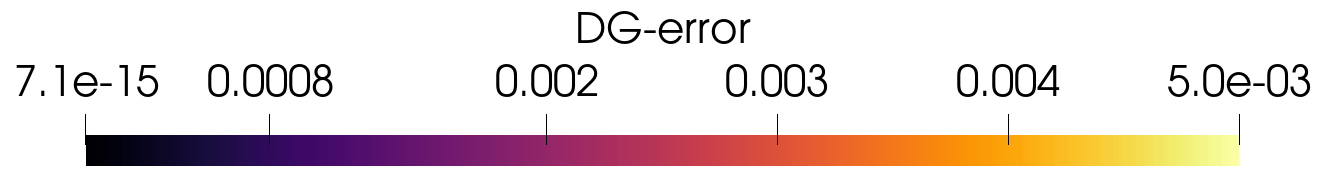}
\includegraphics[trim={1cm 13cm 1cm 13cm},clip, scale=0.075]{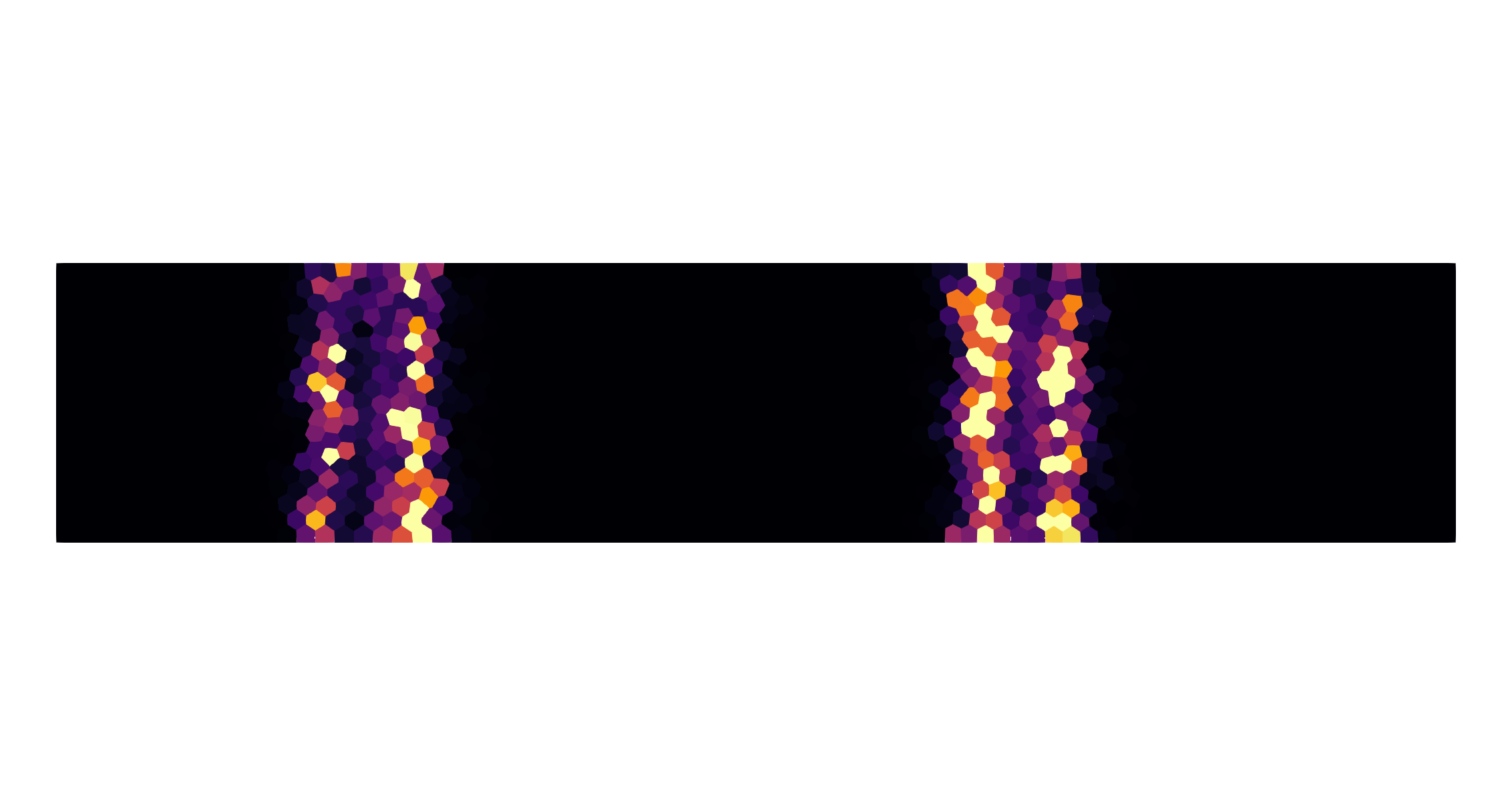}
\hspace*{30ex}\includegraphics[trim={0 0.1cm 0 0.1cm},clip, scale=0.16]{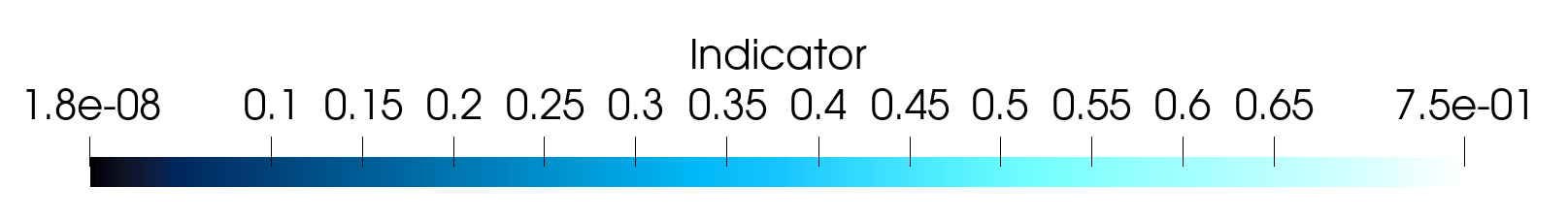}
\includegraphics[trim={1cm 13cm 1cm 13cm},clip, scale=0.075]{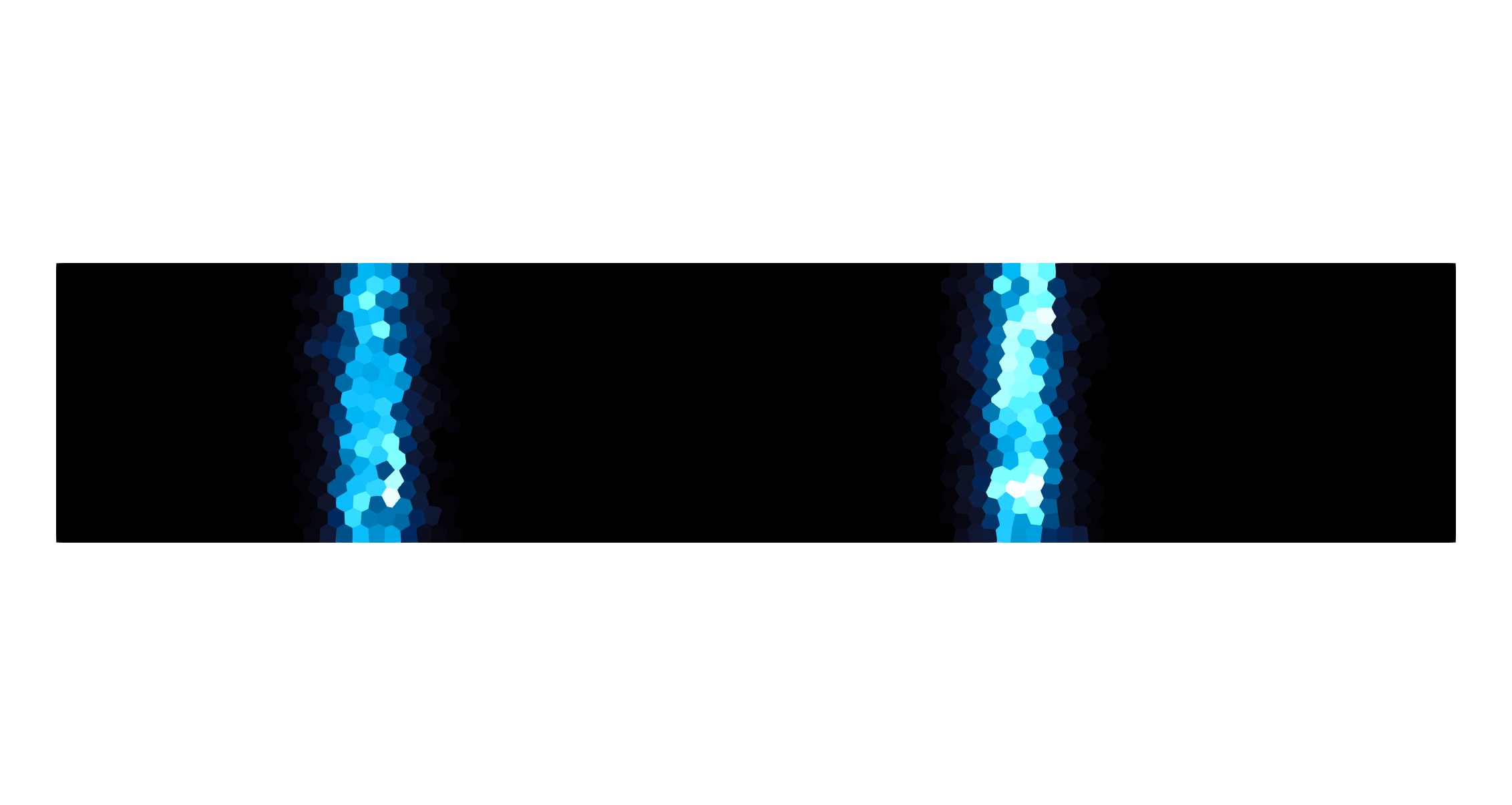} 
\subcaption{$t = 7.1 \; [\mathrm{ms}]$}  
\end{subfigure}
\begin{subfigure}{0.45\textwidth}
\centering
\phantom{\includegraphics[scale=0.14]{singlewave_stationary/scale_.png}}
\includegraphics[trim={1cm 13cm 1cm 13cm},clip, scale=0.075]{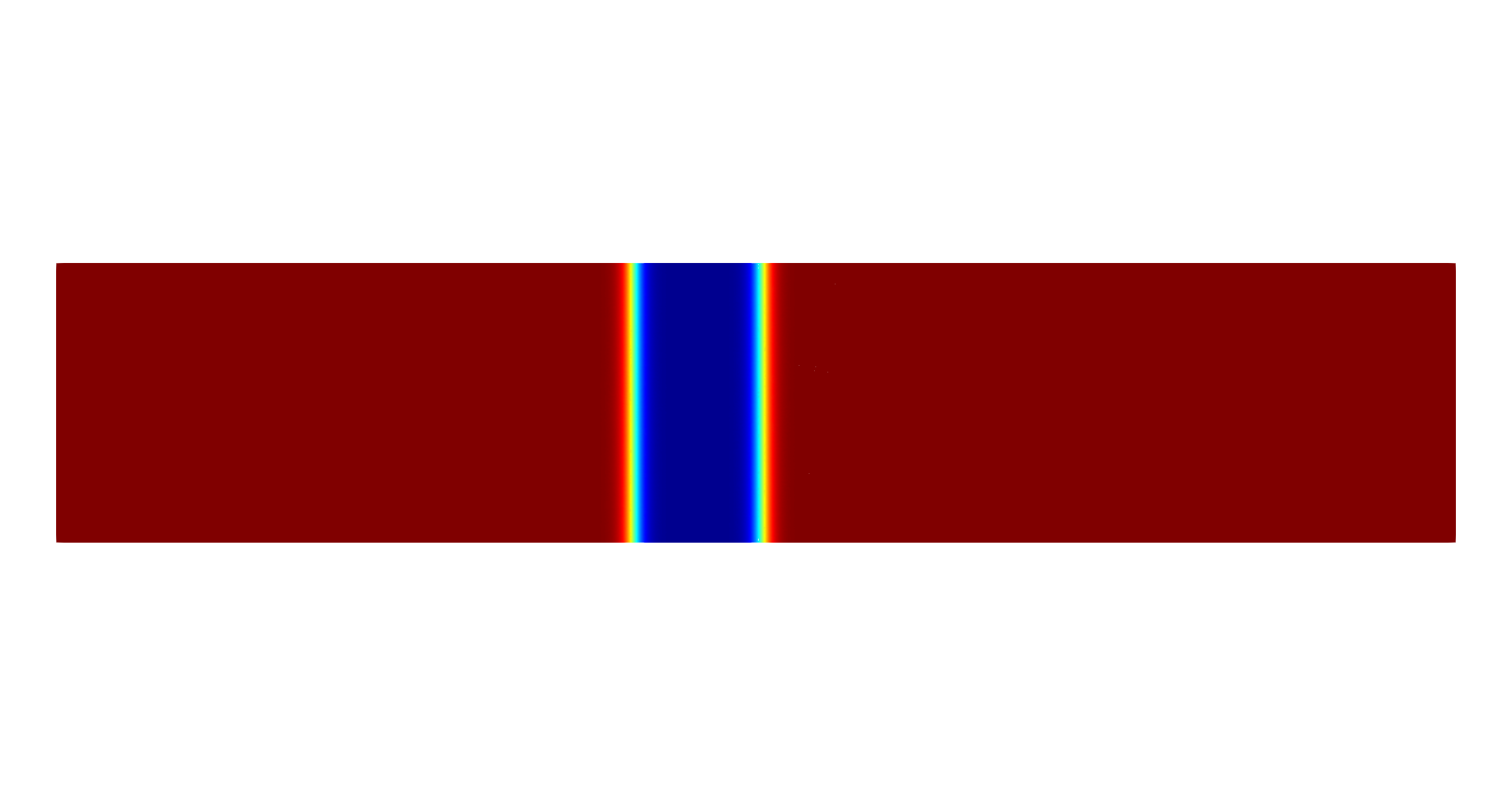}
\includegraphics[trim={1cm 13cm 1cm 13cm},clip, scale=0.075]{doublehomostationary/sol162.png}
\phantom{\includegraphics[trim={0 0cm 0 0cm},clip, scale=0.16]{singlewave_stationary/scale_p.png}}
\includegraphics[trim={1cm 13cm 1cm 13cm},clip, scale=0.075]{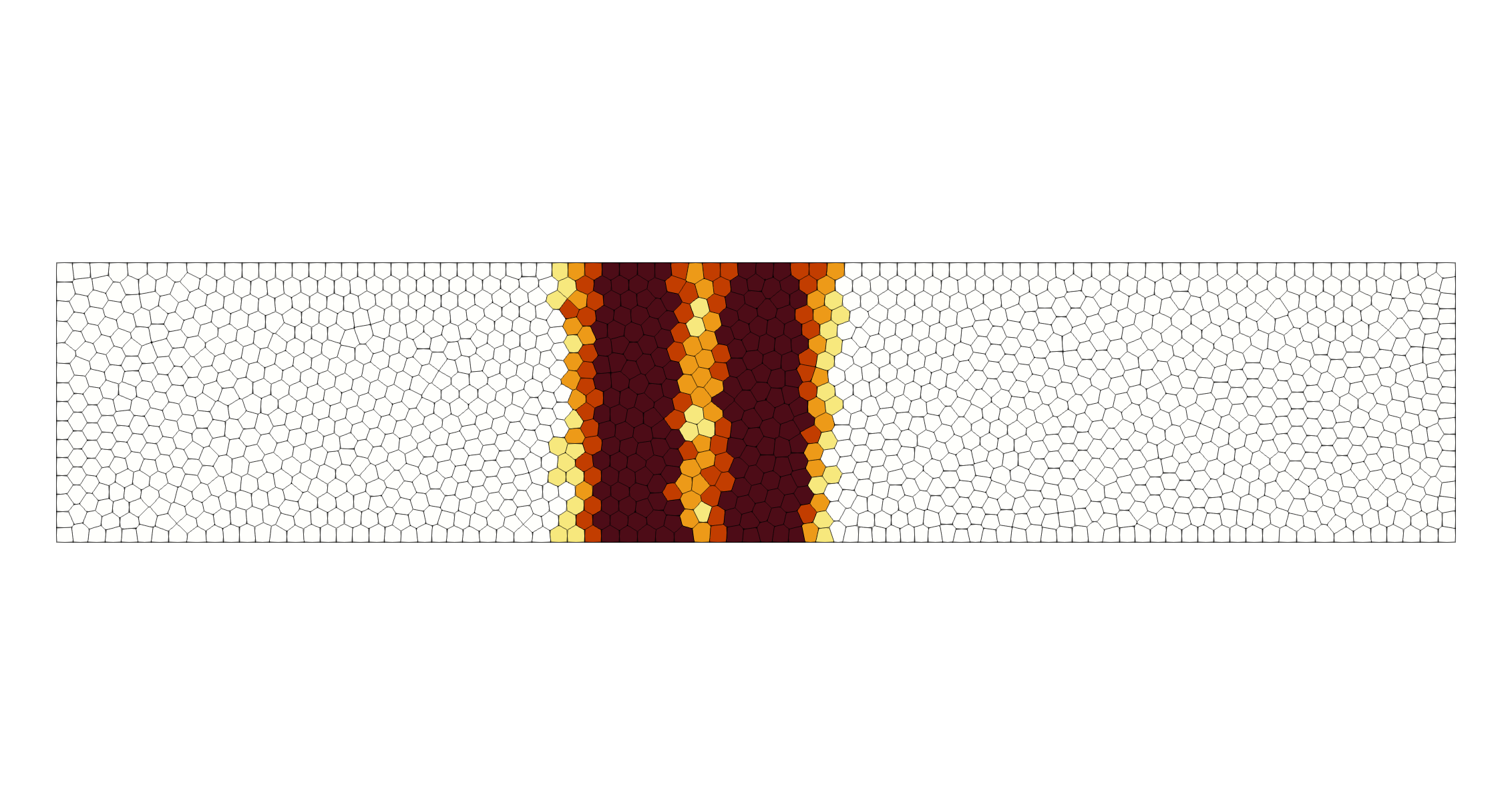}
\phantom{\includegraphics[trim={0 0cm 0cm 0cm},clip, scale=0.16]{doublehomostationary/l2scaleerror.png}}
\includegraphics[trim={1cm 13cm 1cm 13cm},clip, scale=0.075]{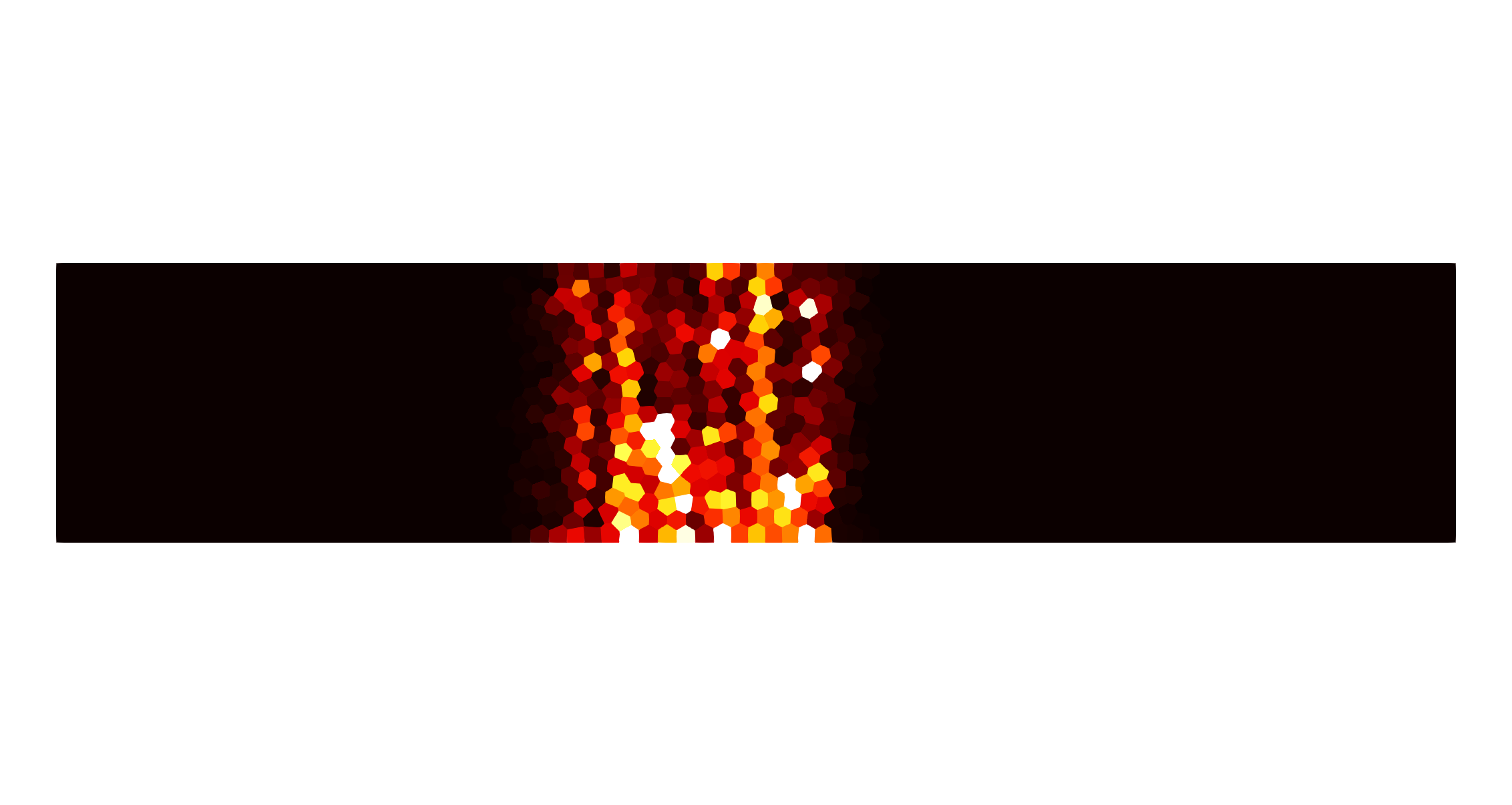}
\phantom{\includegraphics[trim={0 0cm 0 0cm},clip, scale=0.16]{doublehomostationary/DGscalenew.png}}
\includegraphics[trim={1cm 13cm 1cm 13cm},clip, scale=0.075]{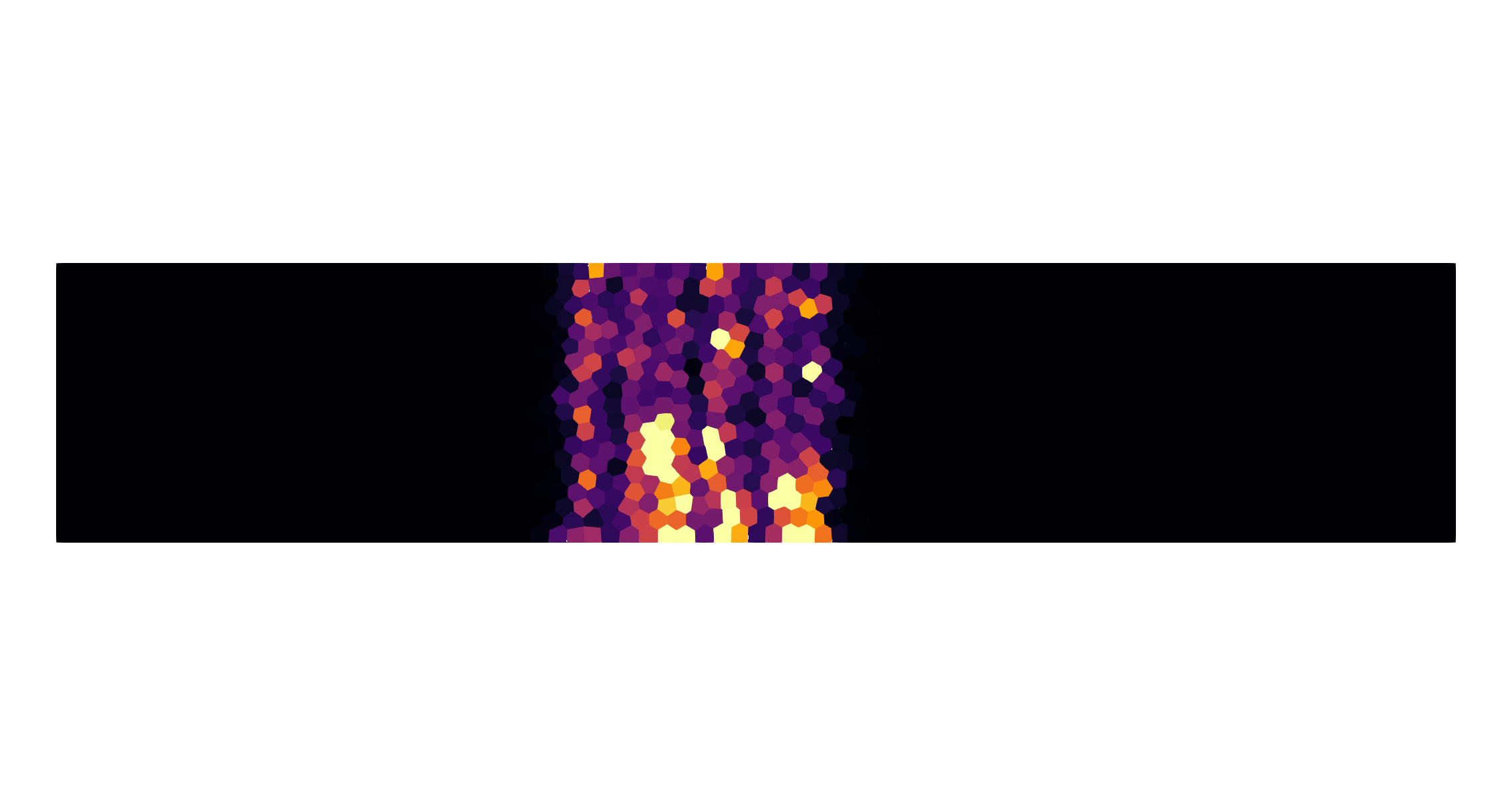}
\phantom{\includegraphics[trim={0 0.1cm 0 0.1cm},clip, scale=0.16]{doublehomostationary/scale_tau.png}}
\includegraphics[trim={1cm 13cm 1cm 13cm},clip, scale=0.075]{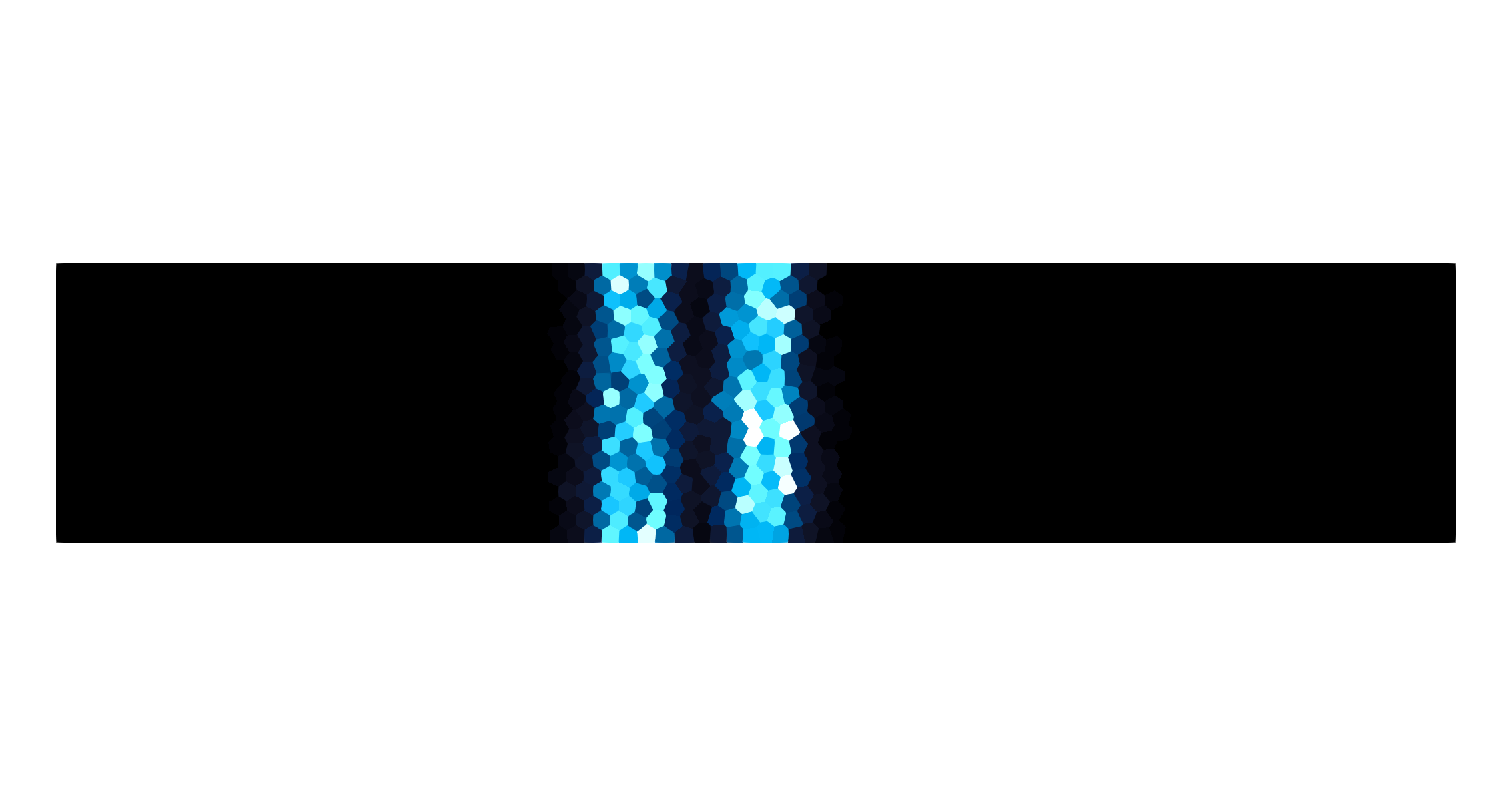} 
\subcaption{$t = 16.6 \; [\mathrm{ms}]$} 
 \end{subfigure}
\caption{Test case 2a. First and second rows: numerical solution computed based on employing a uniform polynomial degree ($p_{\textrm{un}}=5$) and the $p$-adaptive strategy (with $p_\text{max}=5$), respectively. Third row: polynomial approximation degree distribution. Third and fourth rows: 
errors in the $L2$ and $DG$ norms computed with the numerical solution obtained with a uniform polynomial degree ($p_{\textrm{un}}=5$). Fifth row: local error indicator.
The results refer to two different time snapshots $t= 7.1 \, [\mathrm{ms}]$ (left panel) and $t = 16.6 \, [\mathrm{ms}]$ (right panel).}
  \label{fig::ev_double_wave}
\end{figure}
In Figure~\ref{fig::ev_double_wave} (first and second rows) we report the numerical solution computed based on employing a uniform polynomial degree ($p_{\textrm{un}}=5$) and the $p$-adaptive strategy (with $p_\text{max}=5$), respectively, at two different time snapshots $t= 7.1 \, [\mathrm{ms}]$ (left panel) and $t = 16.6\, [\mathrm{ms}]$ (right panel). We observe that: i) the $p$-adaptive algorithm correctly tracks the two wavefronts, with a high-order polynomial degree employed in the proximity of the wavefronts; ii) the polynomial degree distribution is similar for both waves, cf. Figure~\ref{fig::ev_double_wave} (third row). This is consistent with what is expected because the two waves have the same slope and speed. Figure~\ref{fig::ev_double_wave} (third and fourth rows) also reports the errors in the $L2$ and $DG$ norms computed with the numerical solution obtained with a uniform polynomial degree ($p_{\textrm{un}}=5$). Comparing these results with the local error indicator reported in Figure~\ref{fig::ev_double_wave}, last row, we can observe that they exhibit the same behavior, indicating an accuracy of the proposed $p$-adaptive algorithm. This is also confirmed by the results reported in Table~\ref{tab:errorL2DG}, where we show the errors in the $L2$ and $DG$ norms computed for the test case presented, considering as reference solution the numerical solution computed using a uniform polynomial degree ($p_{\textrm{un}}=5$).

\subsubsection{Test case 2b: heterogeneous case} \label{sec:HeteroCase}
\begin{figure}[!htbp]
\centering
\begin{minipage}[b]{0.45\textwidth}
\centering
\includegraphics[trim={2cm 13cm 2cm 13cm},clip,scale=0.0755]{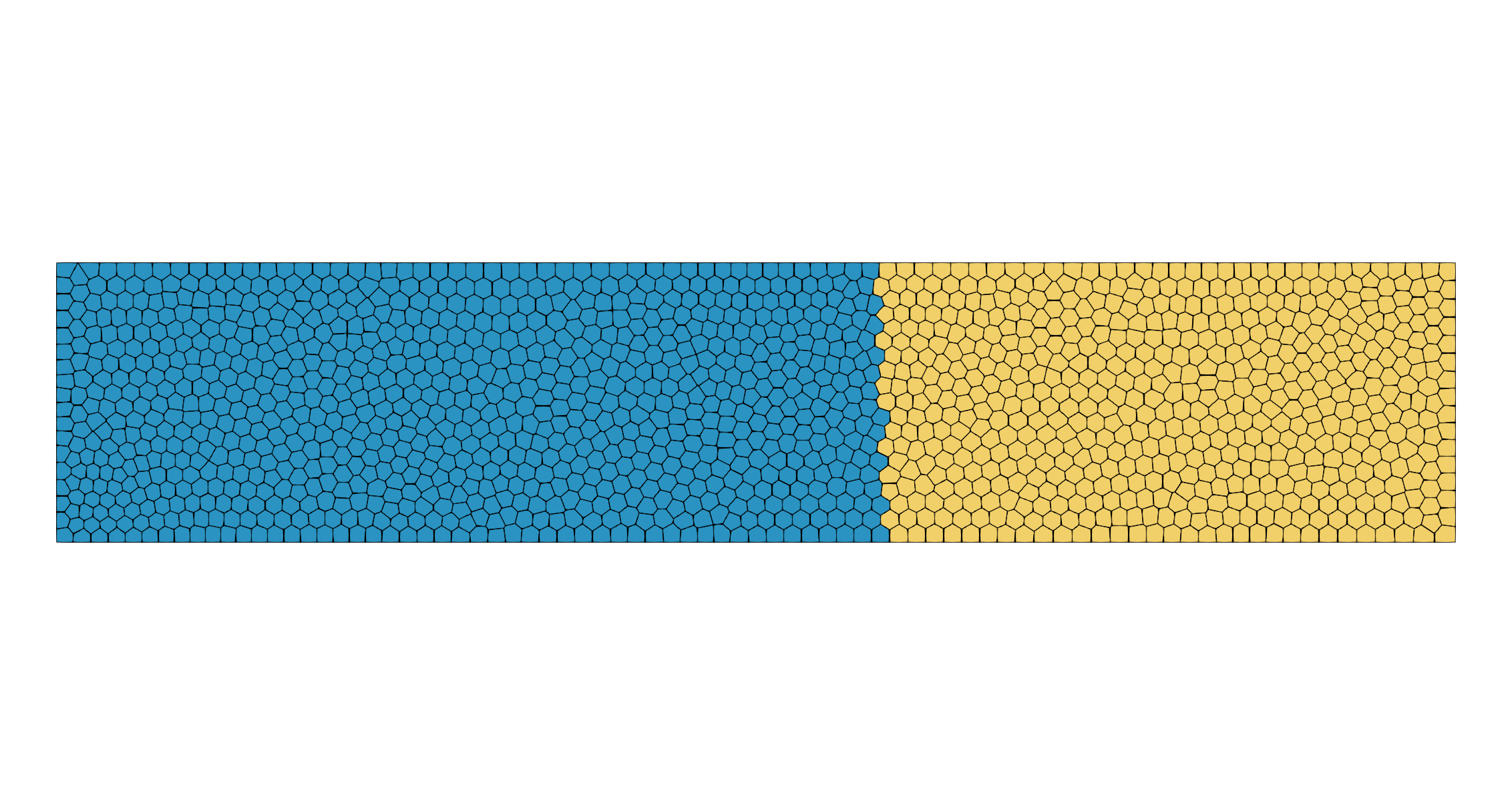}
\vspace{-10ex}
\begin{center}
\Large $\Omega_1$ \hspace{4em}  $\Omega_2$
\end{center}
\vspace{2ex}
\label{fig:hetero_mesh}
\vspace{-2ex}
\includegraphics[trim={0cm 0cm 0cm 0cm},clip,scale=0.19]{singlewave_stationary/scale_.png}
\hspace{-0.6em}\includegraphics[trim={2cm 13cm 2cm 13cm},clip, scale=0.075]{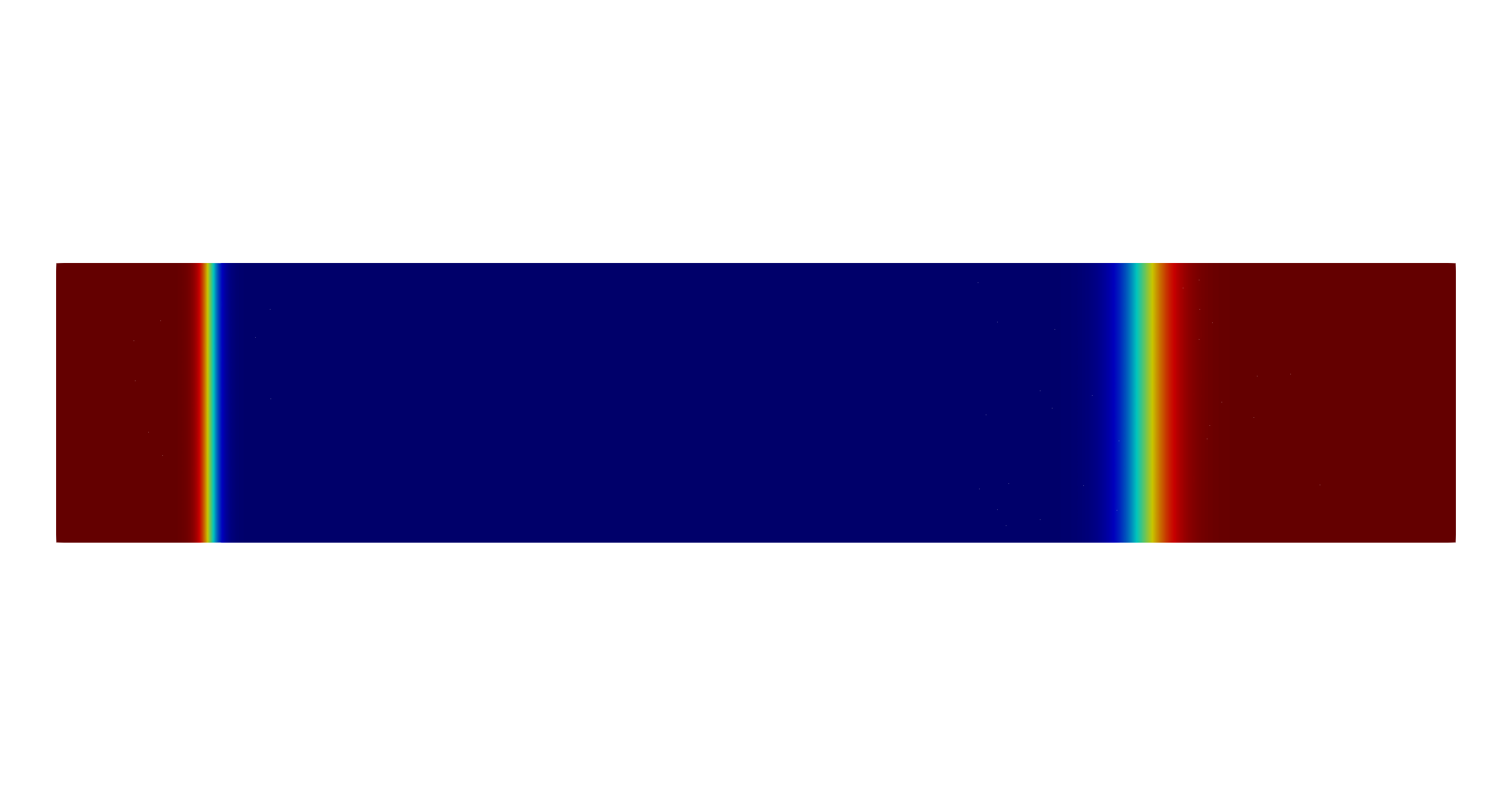}
\vspace{-0.1cm}
\subcaption{\label{fig:double_initial_hetero}}
\end{minipage}\hfill
\begin{minipage}[b]{0.45\textwidth}
\centering
\begin{tikzpicture}[scale=.5, transform shape]
\begin{axis}[
    width=4.6in,
    height=2.8in,
    domain=-3:5,
    samples=1000,
    xmin=-2, xmax=5,
    ymin=-89, ymax=34,
    xlabel={\footnotesize $x$},
    ylabel={\footnotesize $y$},
    ytick={-85, -60, -40, -20, 0, 10, 30},
    grid=both,
    minor tick num=1,
    axis lines=middle,
    enlargelimits=true,
    tick style={black},
    xmajorgrids, ymajorgrids,
    tick label style={font=\small},
]
\addplot[blue, thick, line width = 1.4pt] {115/2 * (tanh((x - 2.5)/0.4) - tanh((x + 1.5)/0.1)) + 30};
\end{axis}
\end{tikzpicture}
\vspace{0.1cm}
\subcaption{\label{fig:profile}}
\end{minipage}
\caption{Test case 2b. Computational mesh and initial condition for double traveling wavefront (left) and profile of the initial condition (right).}
\end{figure}

In this test case, we consider a heterogeneous conductivity tensor resulting in the two waves traveling at different speeds and with different slopes. This test was designed to assess the ability of the $p$-adaptive algorithm to capture two distinct waves with different dynamics.  For this test case, we consider a discontinuous conductivity tensor $\Sigma(\mathbf{x}) = \Sigma_1 \mathbbm{1}_{\Omega_1}(\mathbf{x}) + \Sigma_2 \mathbbm{1}_{\Omega_2}(\mathbf{x})$, where the subdomains $\Omega_1$ ad $\Omega_2$ are defined as: $\Omega_1 = (-2,1) \times (-0.5,0.5)$ and $\Omega_2 = (1,4) \times (-0.5,0.5)$ (cf. Figure~\ref{fig:double_initial_hetero}). Imposing the wave speed as $c_1 = 0.1212 \, [\mathrm{mm}\cdot \text{ms}^{-1}]$ and $c_2 = 0.3157 \, [\mathrm{mm}\cdot \text{ms}^{-1}]$, 
we obtain the conductivity tensors $\Sigma_1= 0.0081 \mathbbm{1} \, \mathrm{[mS \cdot mm^{-1}]}$ and $\Sigma_2= 0.0551 \mathbbm{1} \, \mathrm{[mS \cdot mm^{-1}]}$. This choice results in a faster wave evolving from the right to the left of the domain until it encounters a jump in conductivity at $x=1$; after that, the wave slows down and becomes steep exactly like the other one, evolving with the same speed.
The initial condition profile together with the computational grid is shown in Figure~\ref{fig:profile}.

\begin{figure}[h!]
    \centering
    \begin{subfigure}[t]{0.3\textwidth}
\centering
\includegraphics[width=1\linewidth]{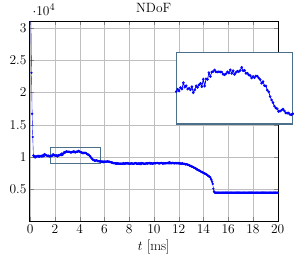}
\caption{\label{fig:dof_evolution_doubletvh}} 
    \end{subfigure}\hfill
    \begin{subfigure}[t]{0.3\textwidth}
\centering
\includegraphics[width=1\linewidth]{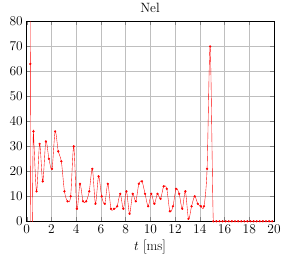}
\caption{\label{fig:comptau_doubletvhetero}} 
    \end{subfigure}\hfill
    \begin{subfigure}[t]{0.3\textwidth}
\centering
\includegraphics[width=1.05\linewidth]{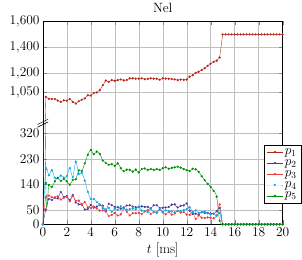}
\caption{\label{fig:elem_doubletvhetero}} 
\end{subfigure}
     \caption{Test case 2b. Left: Evolution of the number of degrees of freedom ($\mathrm{NDoF}$) as a function of time, driven by our $p$-adaptive algorithm.  Center: Evolution of number of elements ($\mathrm{Nel}$) where the local polynomial degree is updated over time.  Right: Evolution of the number of elements ($\mathrm{Nel}$) with local polynomial degrees $p_K=1,2,3,4,5$ over time.  At the initial time all the elements are discretized with $p_\text{max}=5$.} \label{fig:ndof_evolution_doublehetero}
\end{figure}

Also in this case we have repeated the same analysis of test case Section~\ref{dof}. Figure~\ref{fig:ndof_evolution_doublehetero} shows the time evolution of $\mathrm{NDoF}$ together with the evolution of the number of elements where the local polynomial degree is updated over time and the time evolution of the number of elements with local polynomial degrees equal to $1,2,3,4,5$.  
With respect to the homogeneous test case analyzed in Section~\ref{sec:HomogCase}, we observe an increase in the total number of degrees of freedom. The different wave speed induces the $p$-adaptive algorithm to increase the bandwidths where high-order polynomials are employed. Moreover, the presence of the discontinuity in the conductivity tensor leads to an additional increase in the total $\mathrm{NDoF}$ in the time period when the second wavefront passes through that specific portion of the domain (a zoom in this particular time window is reported in \ref{fig:dof_evolution_doubletvh}).
This increase in $\mathrm{NDoFs}$ lasts from $2 \, [\mathrm{ms}]$ up to $5\, [\mathrm{ms}]$. After that, the number of elements with each polynomial degree stabilizes, since the two waves have the same profile (see Figure~\ref{fig:elem_doubletvhetero}). From the numerical results, we can conclude that the proposed $p$-adaptive algorithm accurately identifies the two waves in the domain, using different polynomial degree distributions to effectively capture the fronts in a heterogeneous context. 
 \begin{figure}[!htbp]
 \centering
   \begin{subfigure}{\textwidth}
    \centering
    \hspace*{-.5em}\includegraphics[scale=0.18]{singlewave_stationary/scale_.png}
    \includegraphics[trim={0 0cm 0 0cm},clip, scale=0.18]{singlewave_stationary/scale_p.png}
    \includegraphics[trim={0 0.1cm 0 0.1cm},clip, scale=0.19]{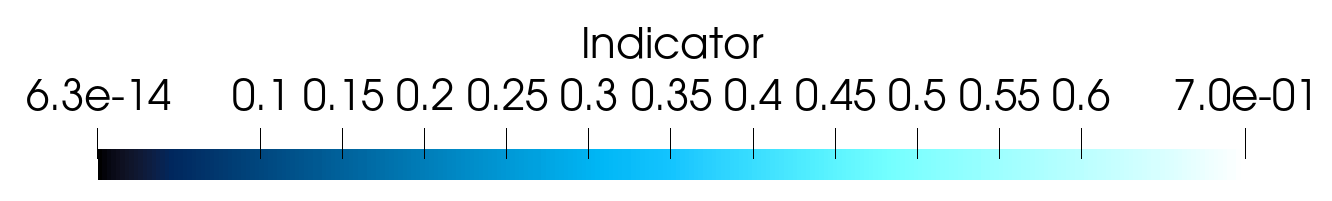}
    \end{subfigure}
 \begin{subfigure}{0.41\textwidth}
\centering
\includegraphics[trim={3cm 13cm 3cm 13cm},clip, scale=0.075]{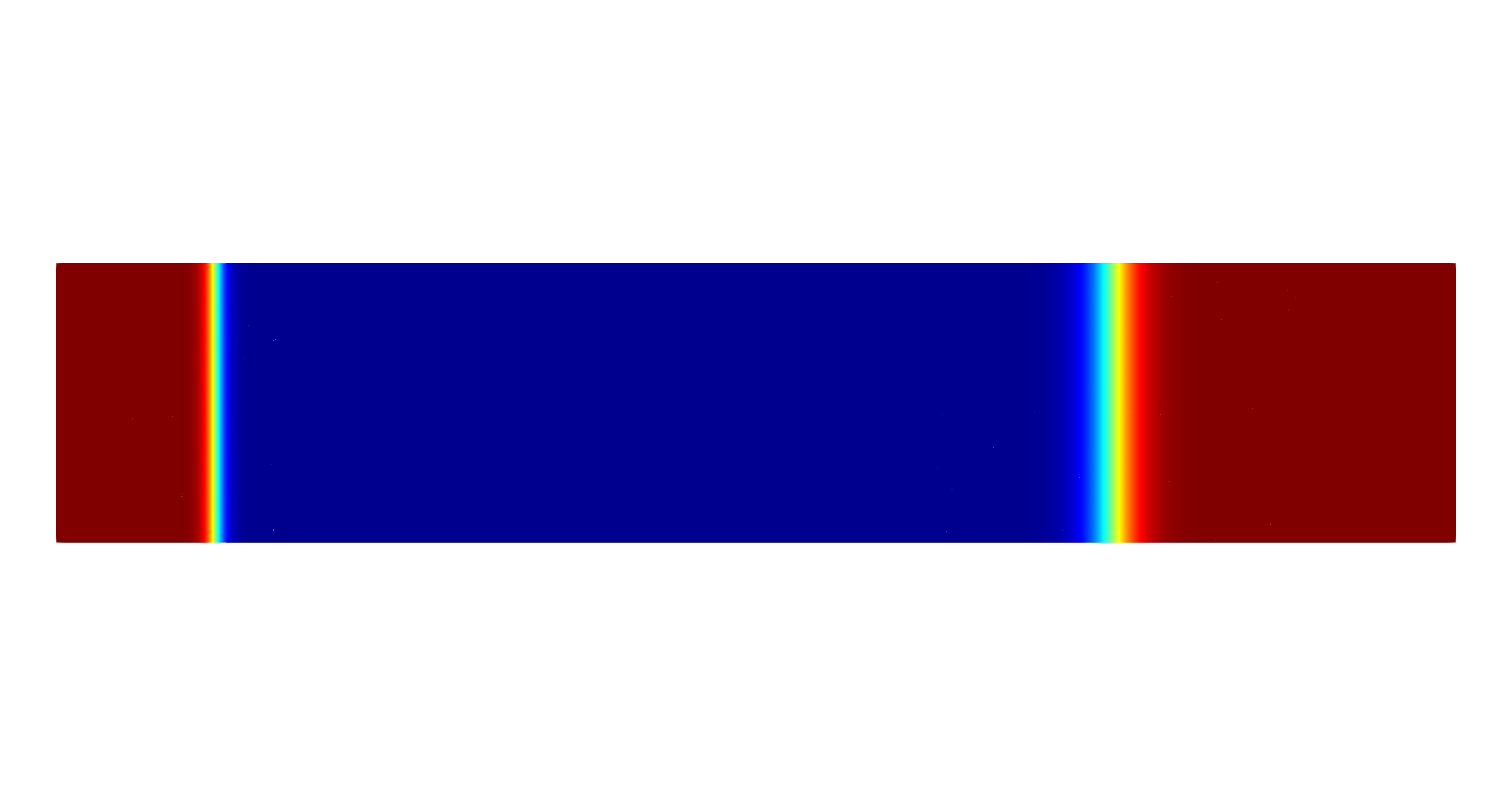}
\includegraphics[trim={3cm 13cm 3cm 13cm},clip, scale=0.075]{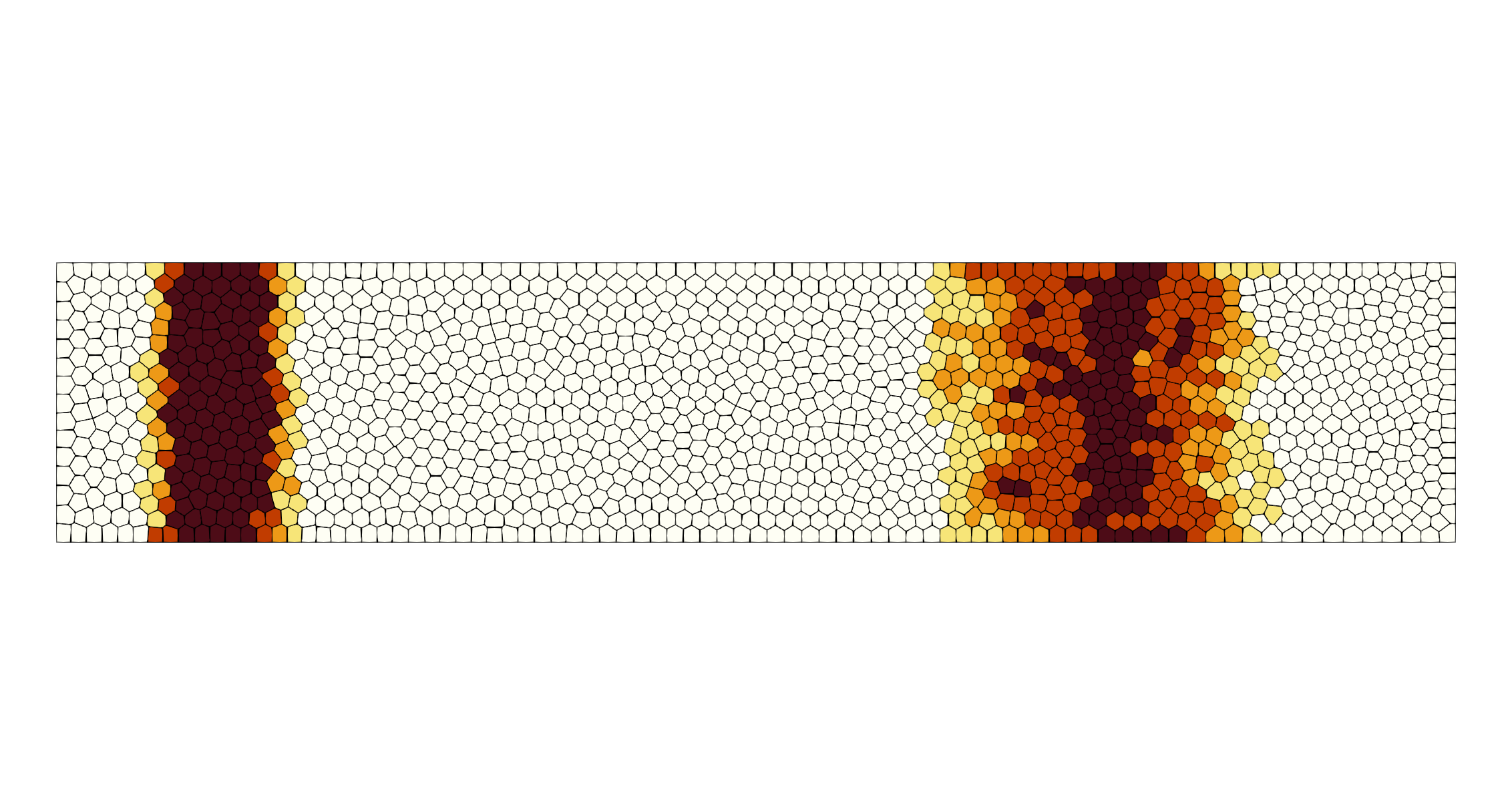} 
\includegraphics[trim={3cm 13cm 3cm 13cm},clip, scale=0.075]{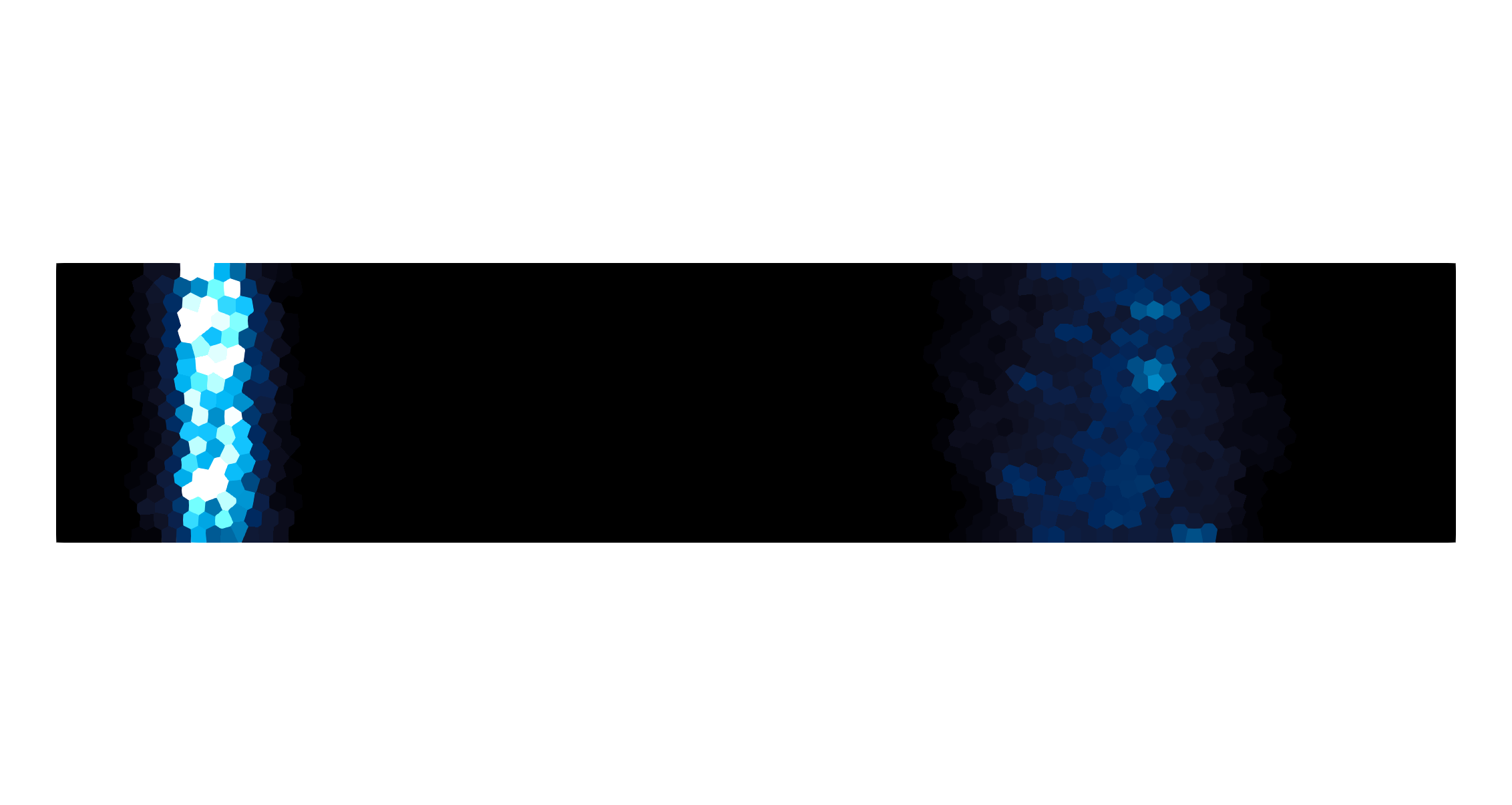}
\subcaption{$t = 1.5 \; \mathrm{ms}$}\label{fig:doublehetero01}  
 \end{subfigure}
 \begin{subfigure}{0.41\textwidth}
\centering
\includegraphics[trim={3cm 13cm 3cm 13cm},clip, scale=0.075]{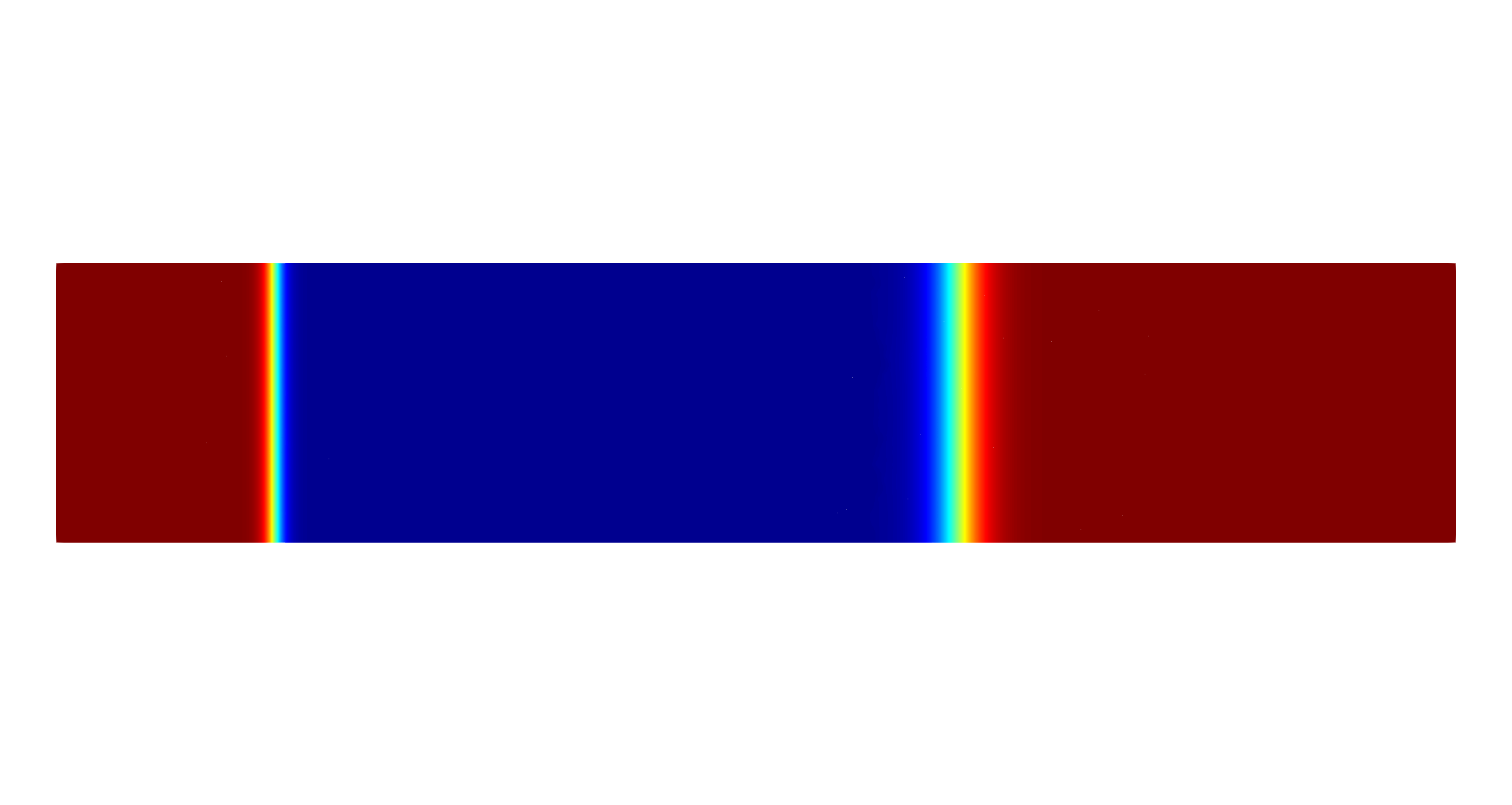}
\includegraphics[trim={3cm 13cm 3cm 13cm},clip, scale=0.075]{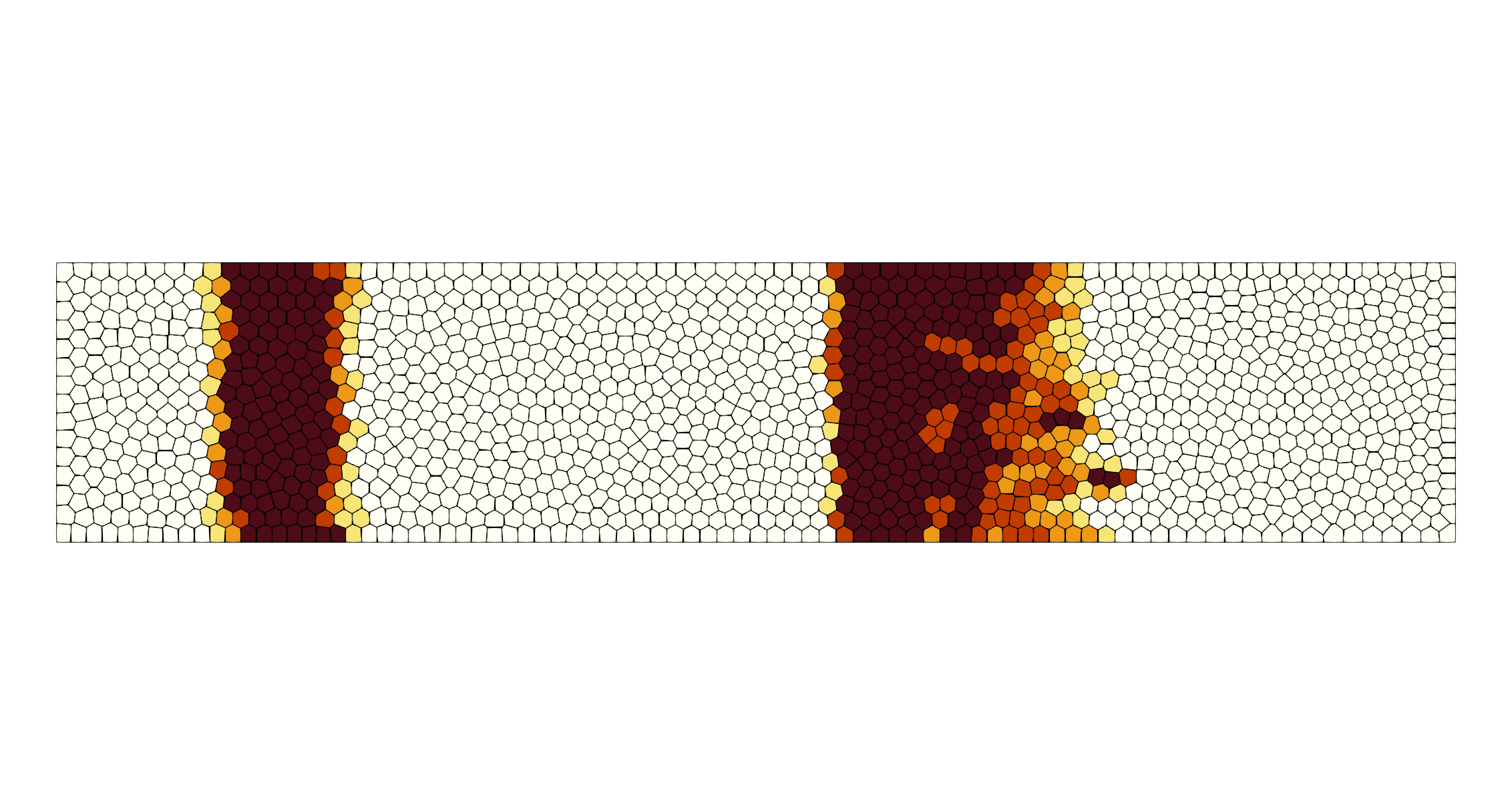}
\includegraphics[trim={3cm 13cm 3cm 13cm},clip, scale=0.075]{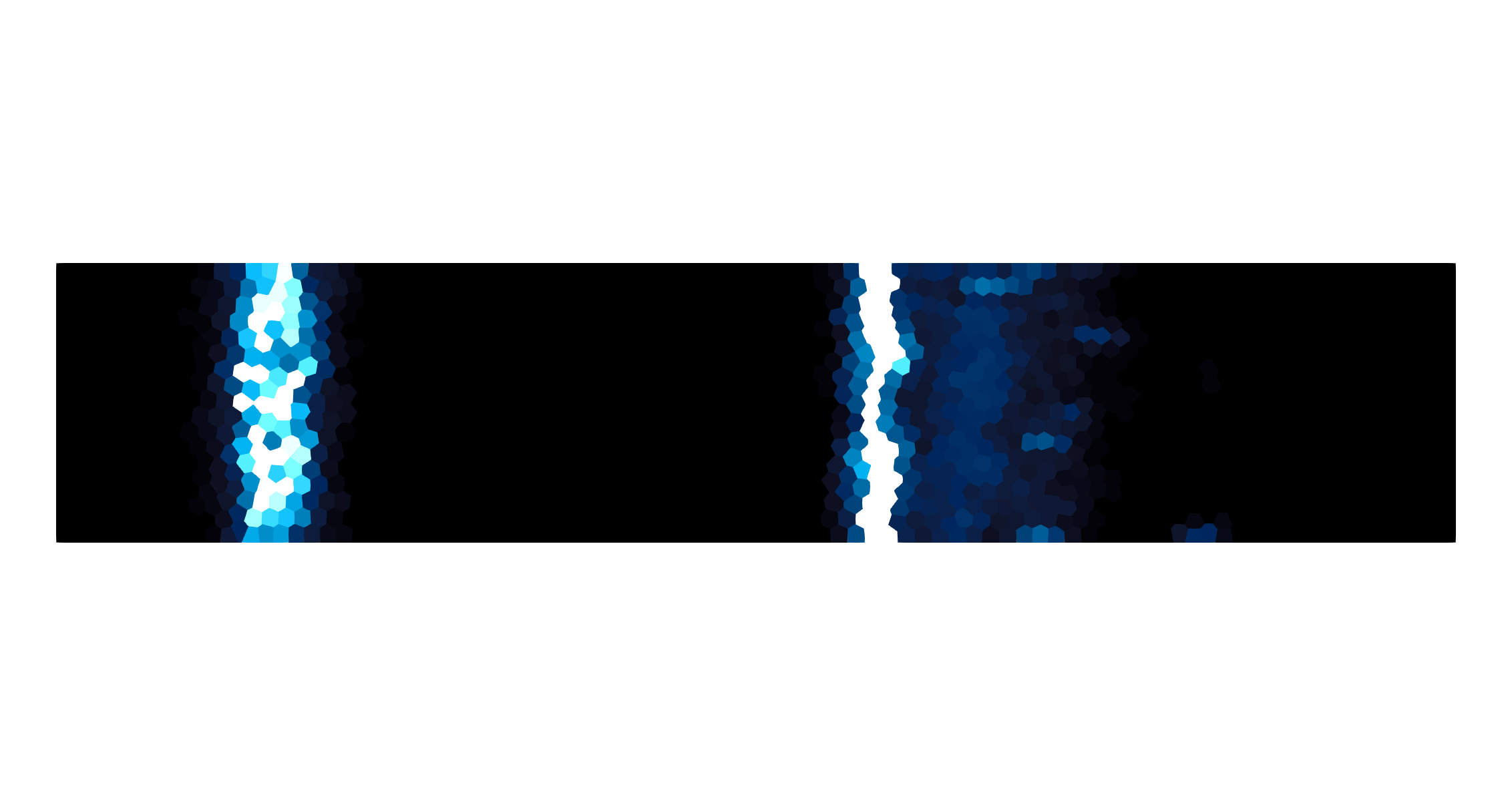} 
\subcaption{$t = 3.6 \; \mathrm{ms}$}  
\label{fig:doublehetero0}
 \end{subfigure}
    \begin{subfigure}{0.41\textwidth}
\centering
\includegraphics[trim={3cm 13cm 3cm 13cm},clip, scale=0.075]{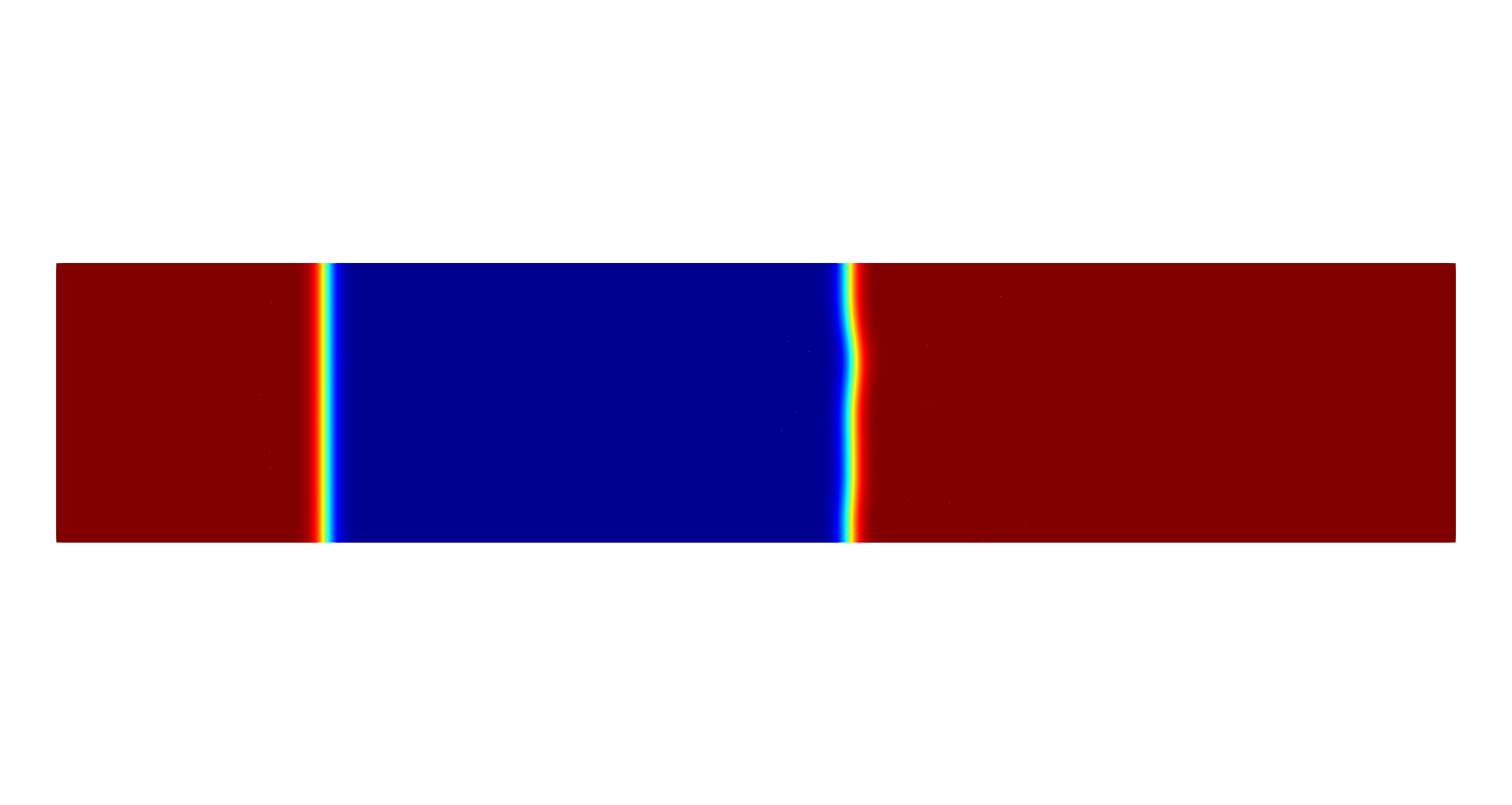}
\includegraphics[trim={3cm 13cm 3cm 13cm},clip, scale=0.075]{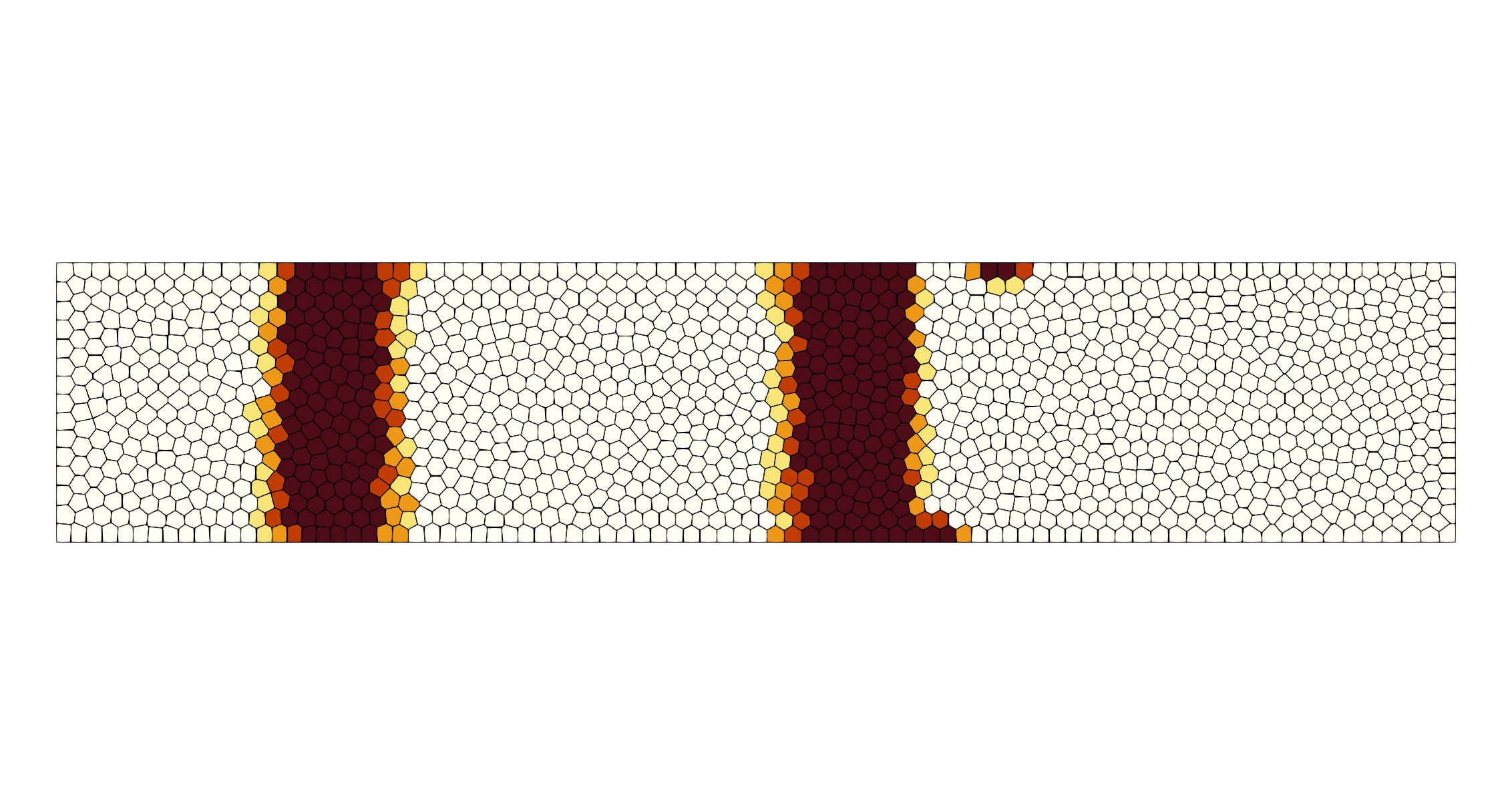}
\includegraphics[trim={3cm 13cm 3cm 13cm},clip, scale=0.075]{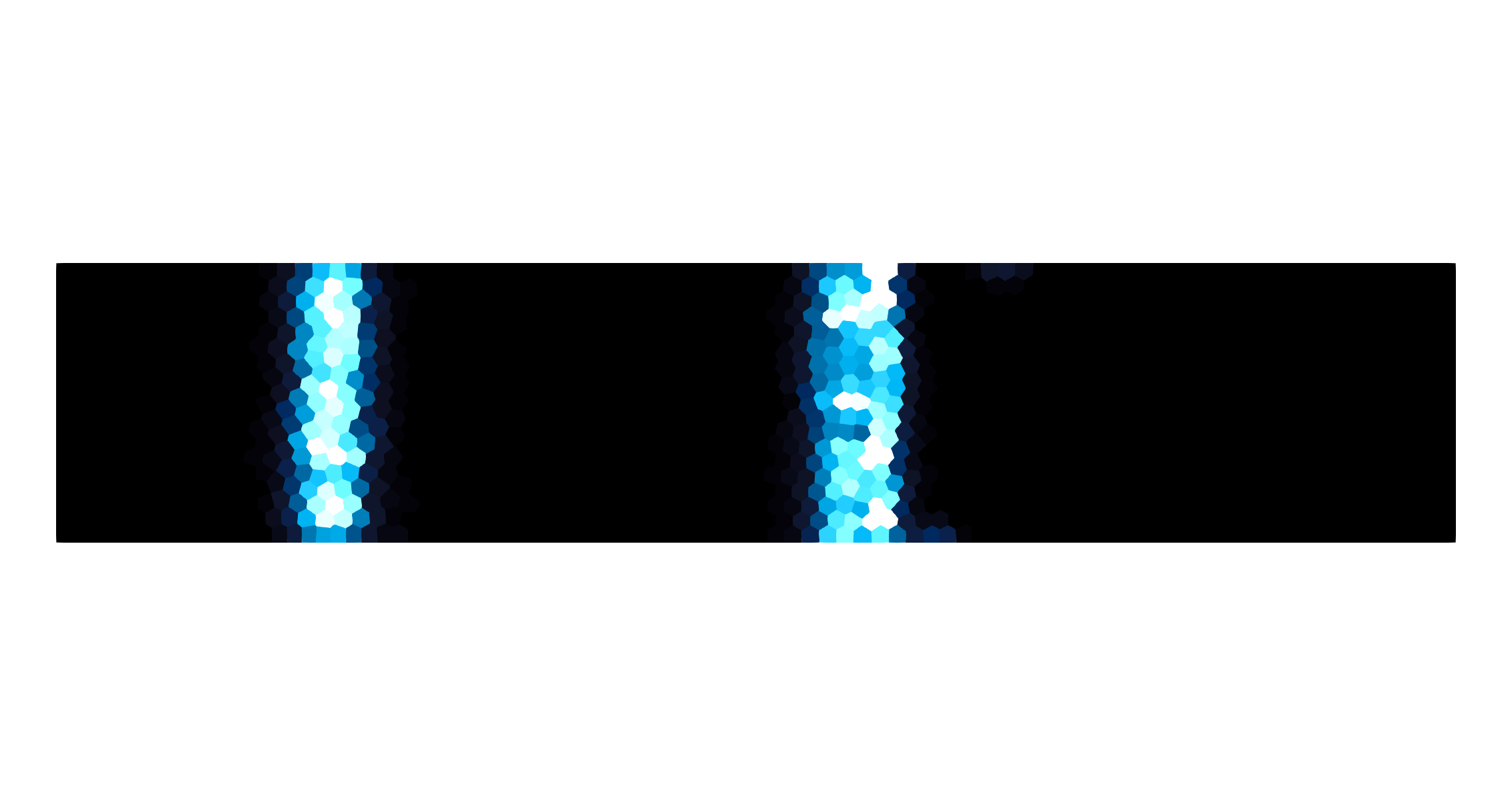} 
\subcaption{$t = 7.7 \; [\mathrm{ms}]$}  \label{fig:doublehetero0D}  
 \end{subfigure}
     \begin{subfigure}{0.41\textwidth}
\centering
\includegraphics[trim={3cm 13cm 3cm 13cm},clip, scale=0.075]{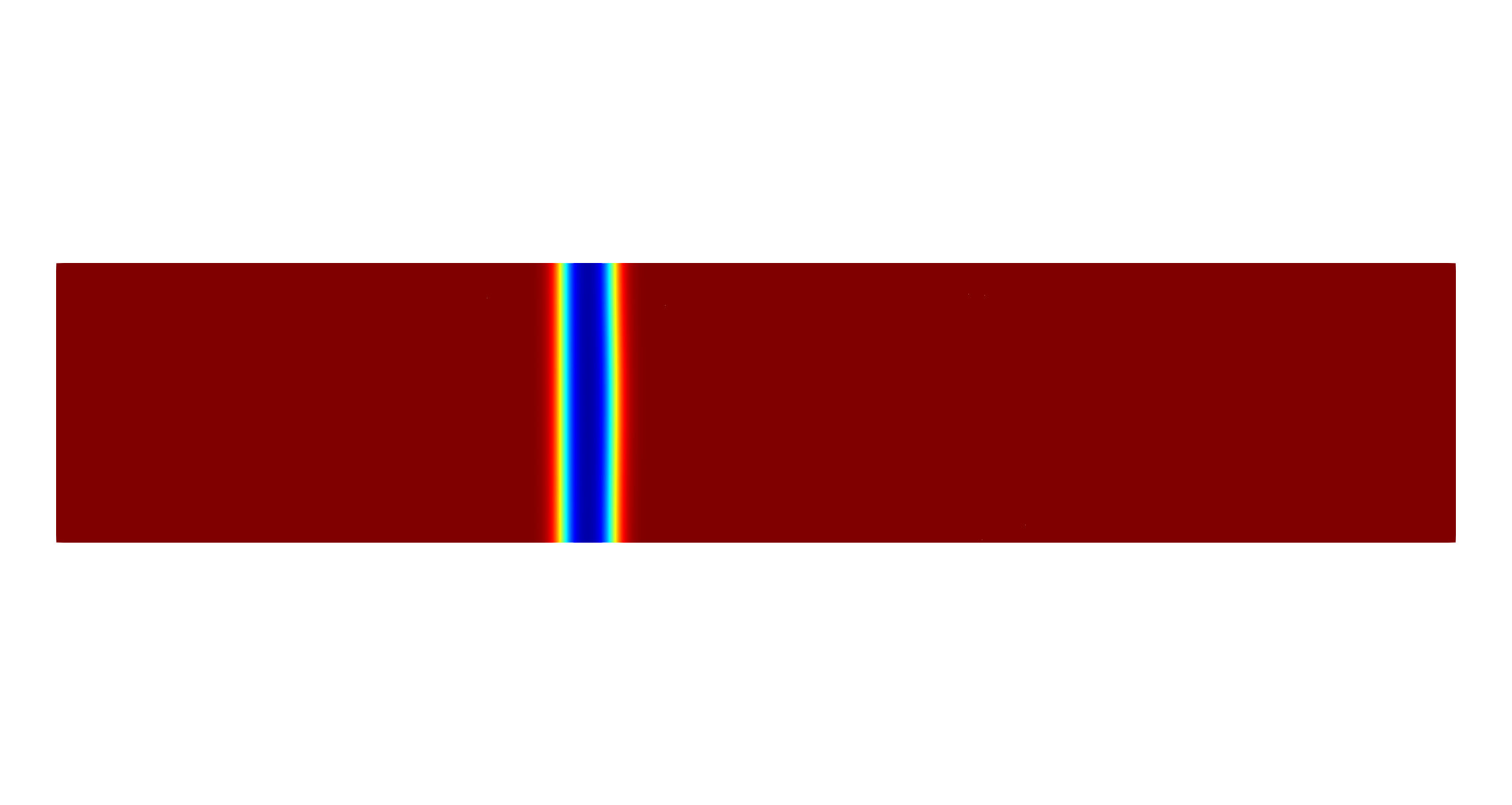}
\includegraphics[trim={3cm 13cm 3cm 13cm},clip, scale=0.075]{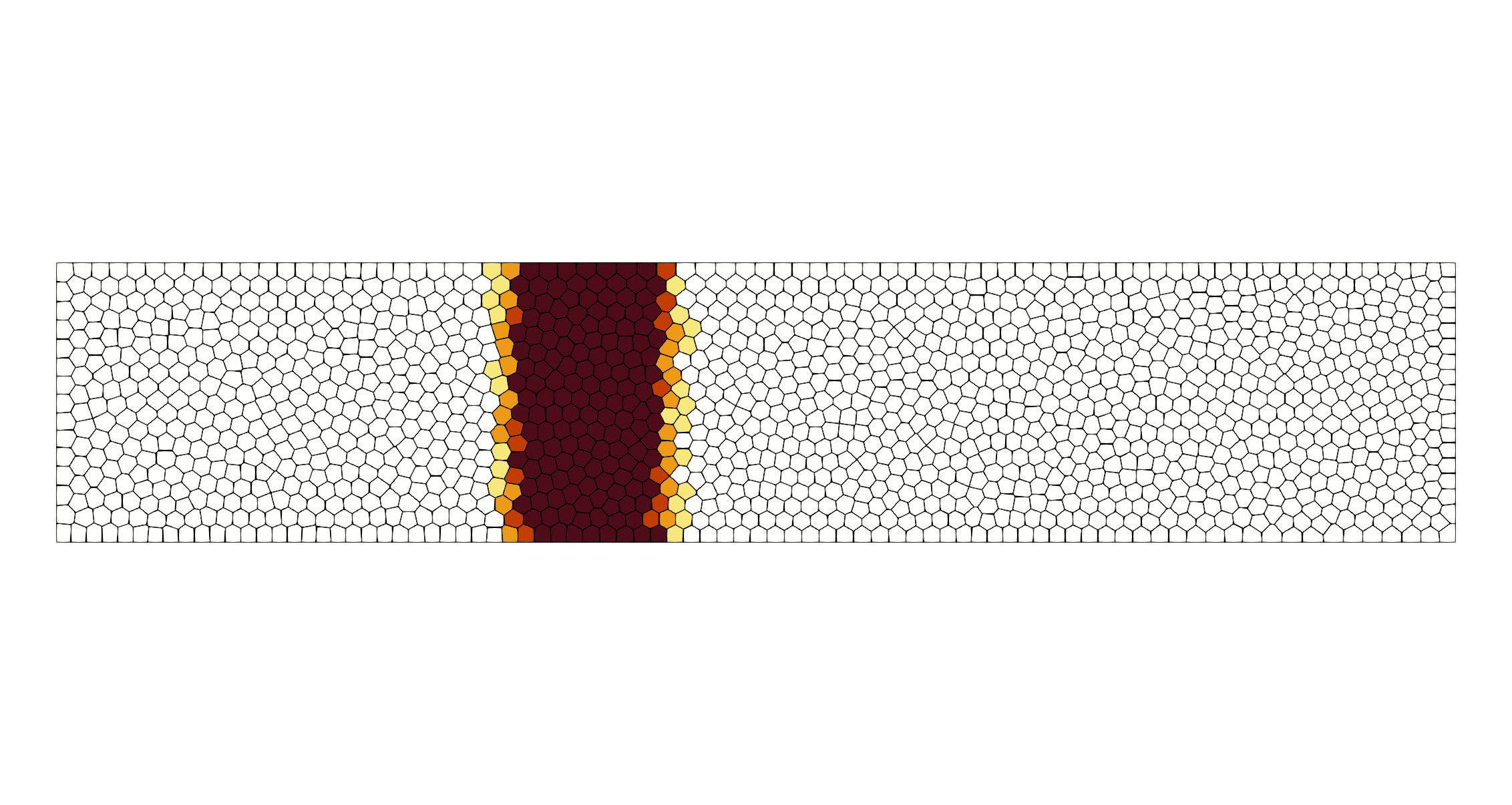}
\includegraphics[trim={3cm 13cm 3cm 13cm},clip, scale=0.075]{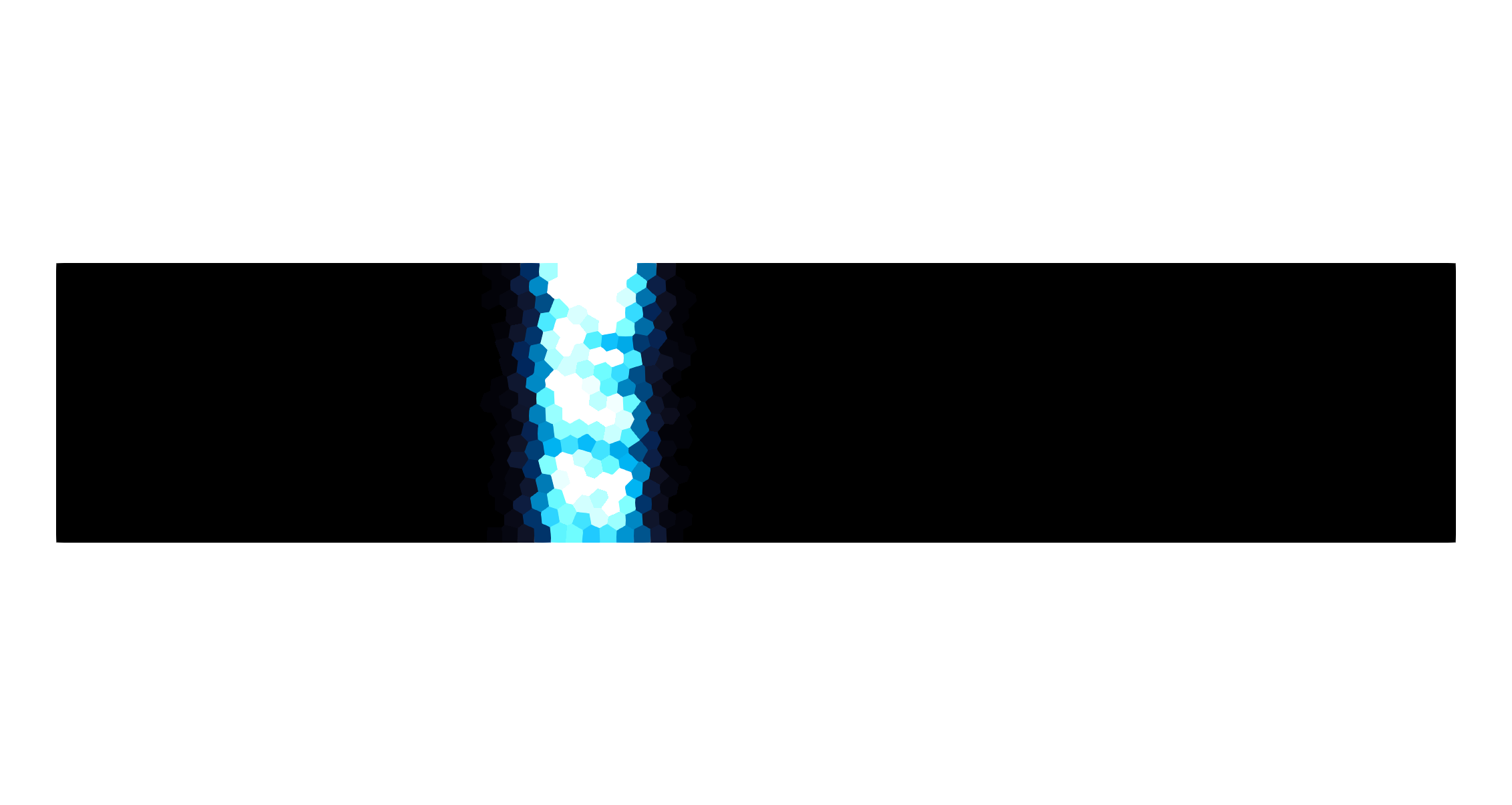} 
\subcaption{$t = 13.0 \; [\mathrm{ms]}$}  \label{fig:doublehetero0D}  
 \end{subfigure}
  \caption{Test case 2b. Computed solution,  element-wise polynomial approximation degree and local error indicator defined as in Equation \eqref{eq::indicator} for different time snapshots $t =1.5,3.6, 7.7, 13.0 \; [\mathrm{ms}]$.}
  \label{fig::ev_double_wavehetero}
\end{figure}
Figure~\ref{fig::ev_double_wavehetero} shows, for different time snapshots $t =1.5,3.6, 7.7, 13 \; [\mathrm{ms}]$, the computed solution,  the element-wise polynomial approximation degree, and the local error indicator. We observe, as expected, the distinct polynomial bandwidths visible for the two wavefronts. At the beginning of the simulation, the less steep, but faster, wave on the right-hand side of the domain shows a bandwidth where high-order polynomials are employed. Close to the jump in the conductivity, the wave speed and slope change, evolving similarly to the wave on the left, which is steeper. 
\FloatBarrier
\section{Application to neuronal electrophysiology}
\label{sec:application}
In this section, we present a test case demonstrating the capabilities of the proposed $p$-adaptive method in approximating physio-pathological scenarios of brain electrophysiology.
\subsection{Grey matter tissue}
\label{sec:test_case_3}
In this test case, we challenge the proposed $p$-adaptive algorithm in the context of neuronal electrophysiology, where we exploit the coupling of the monodomain equation with the Barreto-Cressman conductance-based ionic model \eqref{eq:monodomain}, which is capable of representing the dynamics of neuronal electrophysiology in brain tissues. We focus on the evolution of pathological neuronal impulses, representing dynamics generated by a potassium ion imbalance in the transmembrane space. This spiking behavior, typical of epileptic seizure, is investigated in portions of grey matter tissue modeled isotropically. 
\begin{figure}[!htbp]
\vspace{-0.2cm}
\centering
\begin{subfigure}[t]{0.45\textwidth}
\centering
\color{white}{test}
\end{subfigure}\hfill
\begin{subfigure}[t]{0.45\textwidth}
\centering
\hspace*{-10ex}\includegraphics[trim={0cm 0cm 0cm 0cm},clip,scale=0.17]{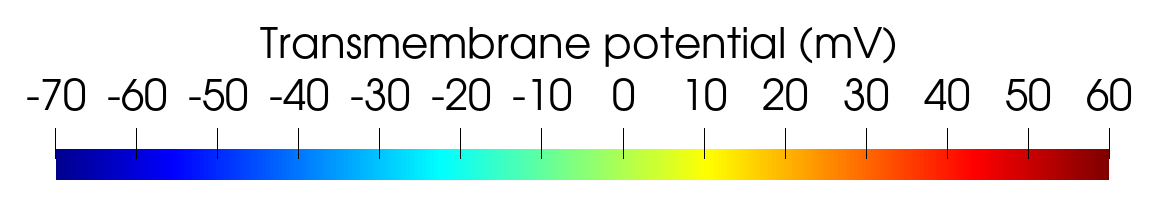}
\end{subfigure}\hfill
\begin{subfigure}[b]{0.45\textwidth}
\vspace*{-3.5ex}
    \centering
    \includegraphics[trim={2cm 7cm 2cm 7cm},clip,scale=0.07]{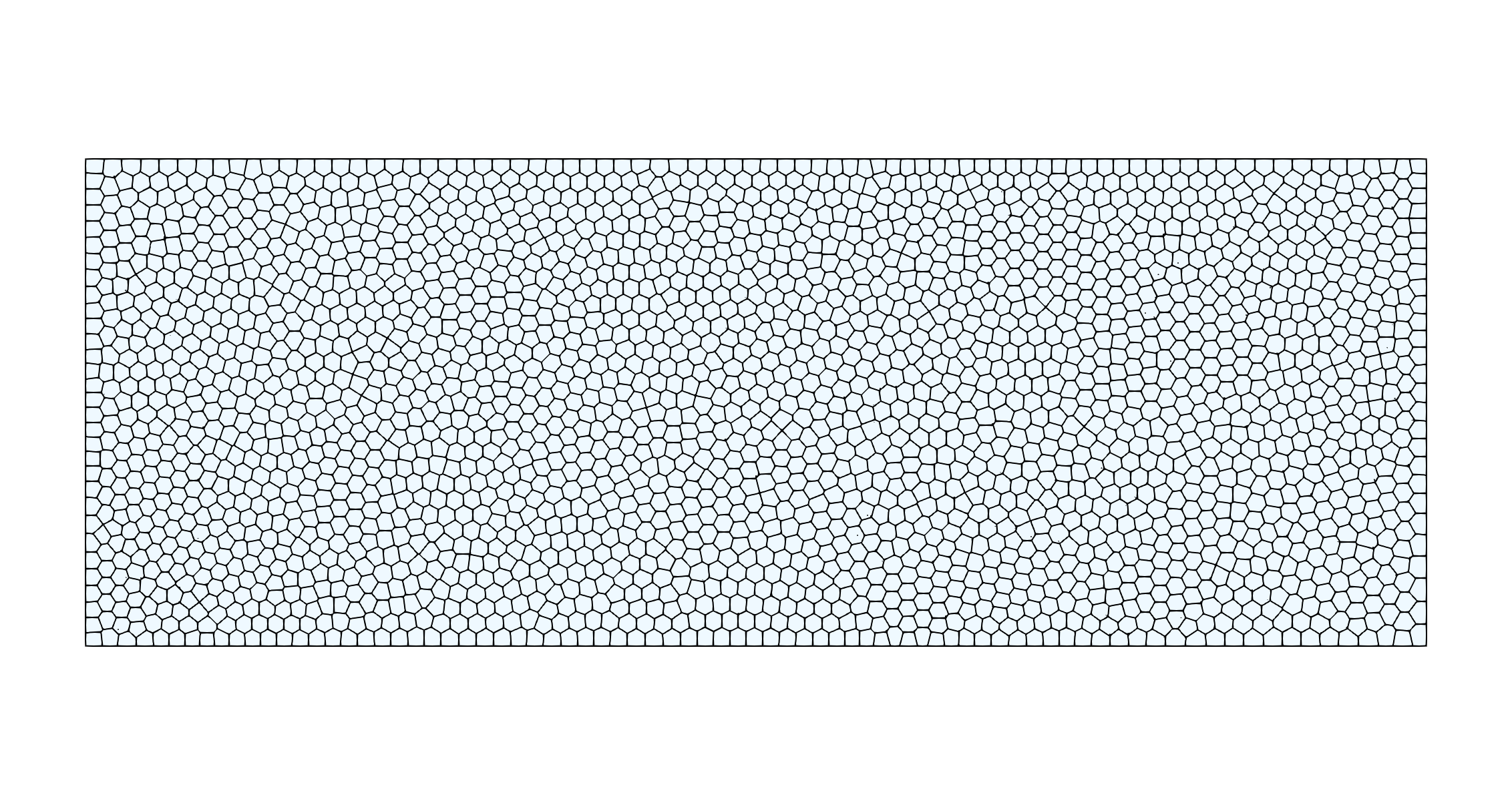}
     \caption{\label{fig::bc_mesh}}     
\end{subfigure}
\begin{subfigure}[b]{0.45\textwidth}
\vspace*{-3.5ex}
    \centering
\includegraphics[trim={2cm 2.2cm 2cm 2.2cm},clip,scale=0.16]{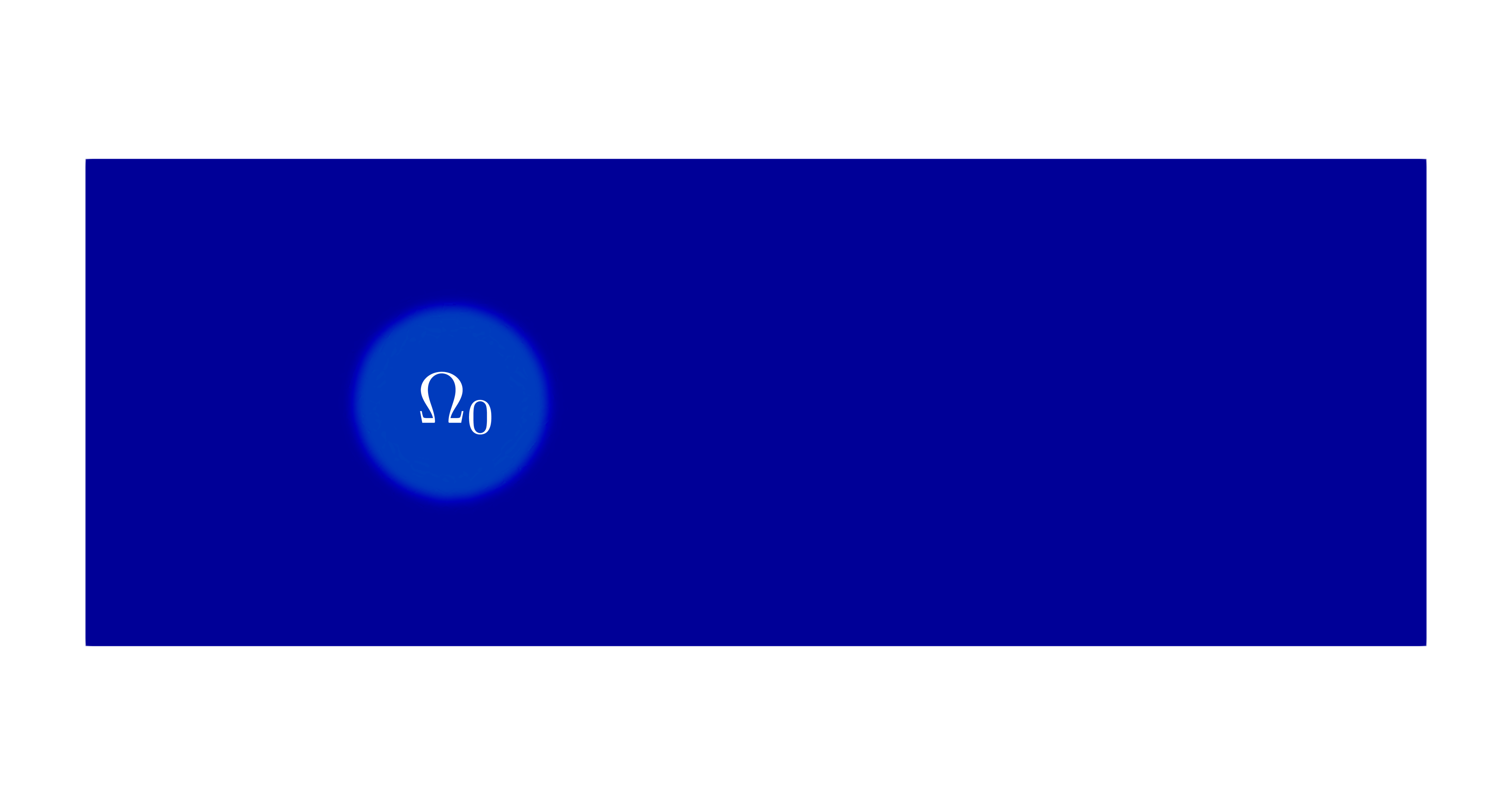} 
\caption{\label{fig::initial_bc}}
  \end{subfigure}
  \caption{Test case 3. Computational domain and polytopal grid with 2500 elements ($h = 0.1148$) (left) and initial condition for the transmembrane potential (right).}
    \label{fig::bc_initial_bc}
\end{figure}
The initial condition for the potential is defined by taking  $u^0 = -50\,[\mathrm{mV}]$ on a subregion  $\Omega_0$, while the remaining part of the domain has $u^0 = -67\,[\mathrm{mV}]$, as done in 
\cite{saglio2024high}. In this case, a pathological condition is simulated on a slab of grey matter tissue on a mesh of 2500 elements $(h=0.1148)$, cf. Figure~\ref{fig::bc_initial_bc}, and throughout the section we set $p_\text{max}=5$. 
\begin{table}[!htbp]
\small
\centering
    \caption{Test case 3. Parameters for the initial condition of the coupled monodomain Barreto-Cressman ionic model. Values taken from \cite{barreto2011ion}.}\label{table:BC_IC}
    \centering 
    \begin{tabular}{lll | lll}
    \toprule
   \textbf{Variable} &\textbf{Inital Value} & \textbf{Units} & \textbf{Variable} &\textbf{Inital Value} & \textbf{Units} \\
    \midrule
    $n^0$ &  $15.5$ & $\mathrm{mM}$ &  $c^0$ & $0$ & $\mathrm{mM}$     \\
    $k^0$ & $7.8$ & $\mathrm{mM}$ &  $g^{n,0}$ & $0.0936$ & -\\
    $g^{c,0}$  &  $0.08553$ & - &     $g^{k,0}$  &  $0.96859$ & -\\
    \bottomrule
    \end{tabular}
\end{table}
\begin{table}[!htbp]
\centering
 \footnotesize
    \caption{Test case 3. Parameters for Barreto-Cressman ionic model. Values taken from \cite{barreto2011ion}.}\label{table:G2}%
    \begin{tabular}{llll}  
    \toprule
    \textbf{Parameter} &  \textbf{Description} &  \textbf{Values} & \\
    \midrule
    $C_m$ &  \text{Membrane capacitance}& 10$^{-2}$ & $\mathrm{mF/cm^2}$   \\
    G$_{\text{AHP}}$ &  \text{Conductance of afterhyperpolarization current}& 0.01 &$\mathrm{mS/cm^2}$   \\
    G$_{\text{KL}}$ &\text{Conductance of potassium leak current} & $0.05$ &$\mathrm{mS/cm^2}$ \\
    G$_{\text{Na}}$ & Conductance of persistent sodium current & $100$ &$\mathrm{mS/cm^2}$  \\
    G$_{\text{CIL}}$ & \text{Conductance of chloride leak current}& $0.05$ &$\mathrm{mS/cm^2}$ \\
    G$_{\text{NaL}}$ & \text{Conductance of sodium leak current} & $0.0175$ &$\mathrm{mS/cm^2}$  \\
    G$_{\text{Ca}}$ & \text{Calcium conductance} & $0.1$ &$\mathrm{mS/cm^2}$ \\
    G$_{\text{K}}$ & \text{Conductance of potassium current} & $40.0$ &$\mathrm{mS/cm^2}$   \\
    G$_{\mathrm{glia}}$ & \text{Strength of glial uptake} & $66.66$ &$\mathrm{mM/s}$   \\
    K$_{\text{bath}}$ & \text{Conductance of potassium} & $8.0$ &$\mathrm{mM}$\\
    \bottomrule
\end{tabular}  
\end{table}
In Table~\ref{table:BC_IC} we report the initial conditions of all the variables of the ionic system, and in Table \ref{table:G2} we list the values of the parameters of the Barreto-Cressman model used for this simulation, cf. \cite{barreto2011ion}.
The conductivity values are taken as in \cite{schreiner2022simulating} for a portion of grey matter tissue, namely,  $\boldsymbol{\Sigma} = 0.7734 \mathds{1} \, \mathrm{[mS \cdot mm^{-1}]}$.
As extensively analyzed in \cite{saglio2024high,erhardt2020dynamics}, we exploit a pathological value of the parameter $K_\text{bath} = 8 \, \mathrm{[mM]}$ that allows the model to enter the Hopf-type bifurcation \cite{ullah2009influence}. This configuration results in a very fast and auto-induced bursting behavior. In this work, the pathological behavior is further accentuated by introducing an external forcing term, which increase the frequency of the spikes, accelerating the bursting behavior. The forcing term is defined as follows:
\begin{equation}
    \begin{aligned}
    I^\text{ext}(\boldsymbol{x},t) =  \frac{A}{1 + e^{ \sin(t)}}\mathds{1}_{ \Omega_0} (\boldsymbol{x}).
\end{aligned}
\label{eq:forcing_bc}
\end{equation}
Figure~\ref{fig:bc0d} depicts the acceleration of the pathological dynamics generated by this forcing: in the self-induced case, the action potential generates a second peak after $40\,[\mathrm{ms}]$ (Figure~\ref{fig:bc0d_1}); in contrast, the introduction of the external force anticipates the second peak at approximately $7\,[\mathrm{ms}]$ (Figure~\ref{fig:bc0d_0}). As before, we consider the error indicator for the adaptive choice of the polynomial degree defined in Equation \eqref{eq::indicator}.
\begin{figure}[!htbp]
\centering
\begin{subfigure}[t]{0.45\textwidth}
\centering
\resizebox{0.75\textwidth}{!}{
\begin{tikzpicture}{
      \begin{axis}[
        width=3.875in,
        height=2.26in,
        at={(2.6in,1.099in)},
        scale only axis,
        xmin=-0.05,
        xmax=100,
        xminorticks=true,
        xlabel = { $t$ [ms]},
        ylabel = { $u \text{ (mV)}$},
        ymin=-80,
        ymax=70,
        yminorticks=true,
        axis background/.style={fill=white},
        title={\Large \color{black} Transmembrane potential (mV) evolution for $A = 0 \;\frac{\text{mA}}{\text{mm}}$},
        xmajorgrids,
        xminorgrids,
        ymajorgrids,
        yminorgrids,
        legend style={legend cell align=left, align=left, draw=white!15!black}
        ]
        \addplot        table[x=t,y=v, col sep=comma,mark=none] {solution_adapt.csv};
      \end{axis}}
    \end{tikzpicture}}
     \caption{\label{fig:bc0d_1}} 
     \end{subfigure}
\begin{subfigure}[t]{0.45\textwidth}
\centering
\resizebox{0.75\textwidth}{!}{
\begin{tikzpicture}{
      \begin{axis}[
        width=3.875in,
        height=2.26in,
        at={(2.6in,1.099in)},
        scale only axis,
        xmin=-0.05,
        xmax=100,
        xminorticks=true,
        xlabel = { $t$ [ms]},
        ylabel = { $u \text{ (mV)}$},
        ymin=-80,
        ymax=70,
        yminorticks=true,
        axis background/.style={fill=white},
        title={\Large \color{black}  Transmembrane potential (mV) evolution for $A = 9 \;\frac{\text{mA}}{\text{mm}}$},
        xmajorgrids,
        xminorgrids,
        ymajorgrids,
        yminorgrids,
        legend style={legend cell align=left, align=left, draw=white!15!black}
        ]
        \addplot        table[x=t,y=v, col sep=comma,mark=none] {solution_adapt_forcing.csv};
      \end{axis}}
    \end{tikzpicture}}
     \caption{\label{fig:bc0d_0}} 
     \end{subfigure}
\caption{Test case 3. Figure \eqref{fig:bc0d_1} shows the dynamics of the auto-induced and bursting pathological behavior of the 0D-ionic model. In Figure \eqref{fig:bc0d_0} we shows the evolution of the transmembrane potential in the 0D-ionic model in presence of the forcing term defined in Equation \eqref{eq:forcing_bc} with amplitude $A = 9 \;\frac{\text{mA}.}{\text{mm}}$.} \label{fig:bc0d}
\end{figure}
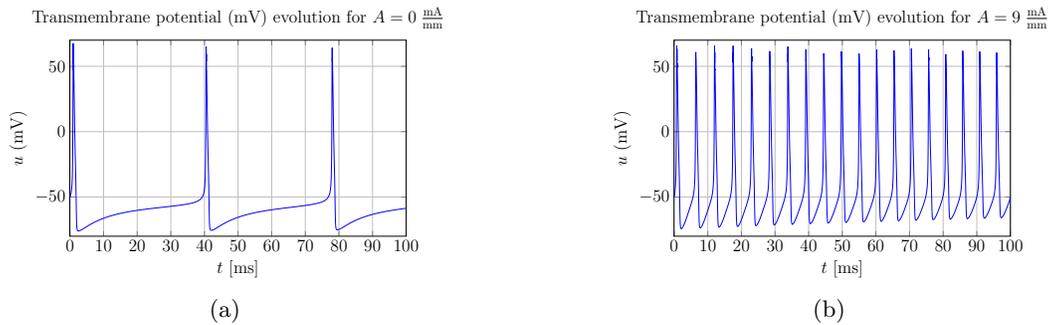
 \begin{figure}[h!]
 \vspace*{-2ex}
 \centering
  \begin{subfigure}[b]{\textwidth}
    \centering
    \includegraphics[scale=0.18]{photosbc/scale.png}
    \includegraphics[trim={0 0cm 0 0cm},clip, scale=0.18]{singlewave_stationary/scale_p.png}
    \includegraphics[trim={0 0.1cm 0 0.3cm},clip, scale=0.19]{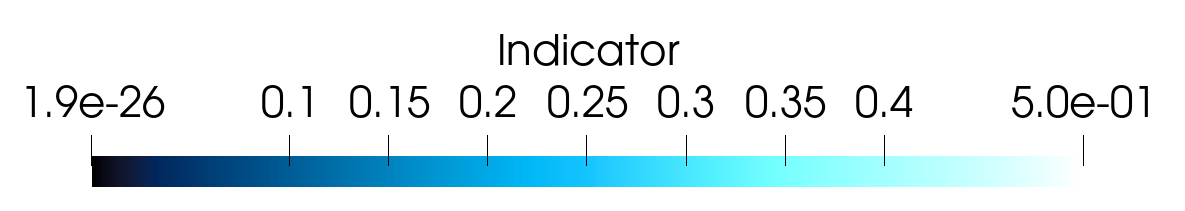}
    \end{subfigure}\hfill
    \begin{subfigure}[b]{0.04\textwidth}
\raisebox{3mm}{
\vspace{-2cm}
\hspace*{-4ex}
  \begin{tabular}{@{}l@{}}
$3.4\,[\mathrm{ms}]$\\[1.67cm]
$7.2\,[\mathrm{ms}]$\\[1.67cm]
$9.4\,[\mathrm{ms}]$\\[1.67cm]
$13.2\,[\mathrm{ms}]$
  \end{tabular}}
  \end{subfigure}
 \begin{minipage}{0.31\textwidth}
 \centering
 \hspace*{1.ex}\includegraphics[trim={5cm 7cm 5cm 7cm},clip, scale=0.07]{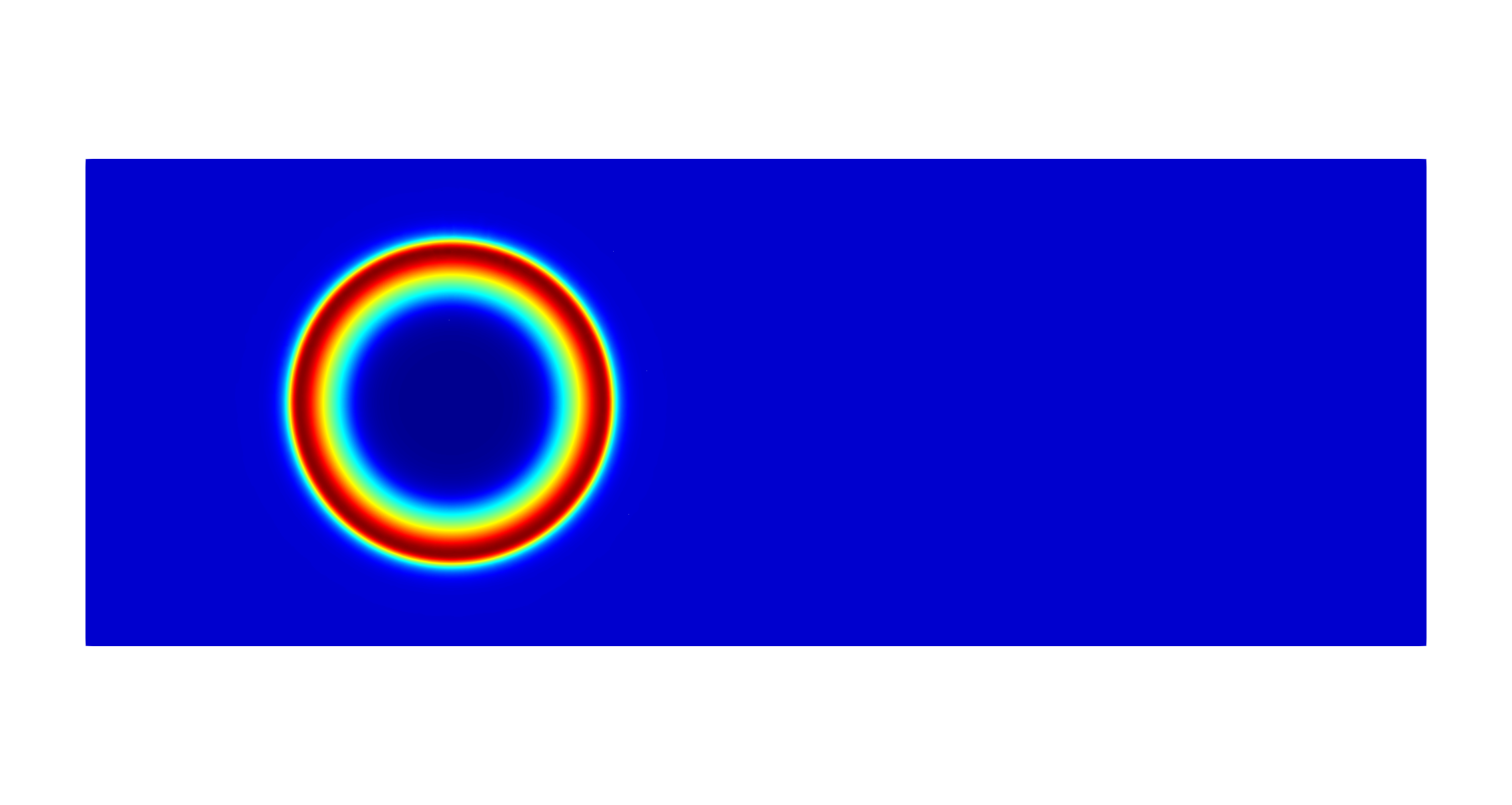}
\hspace*{1ex}\includegraphics[trim={5cm 7cm 5cm 7cm},clip, scale=0.07]{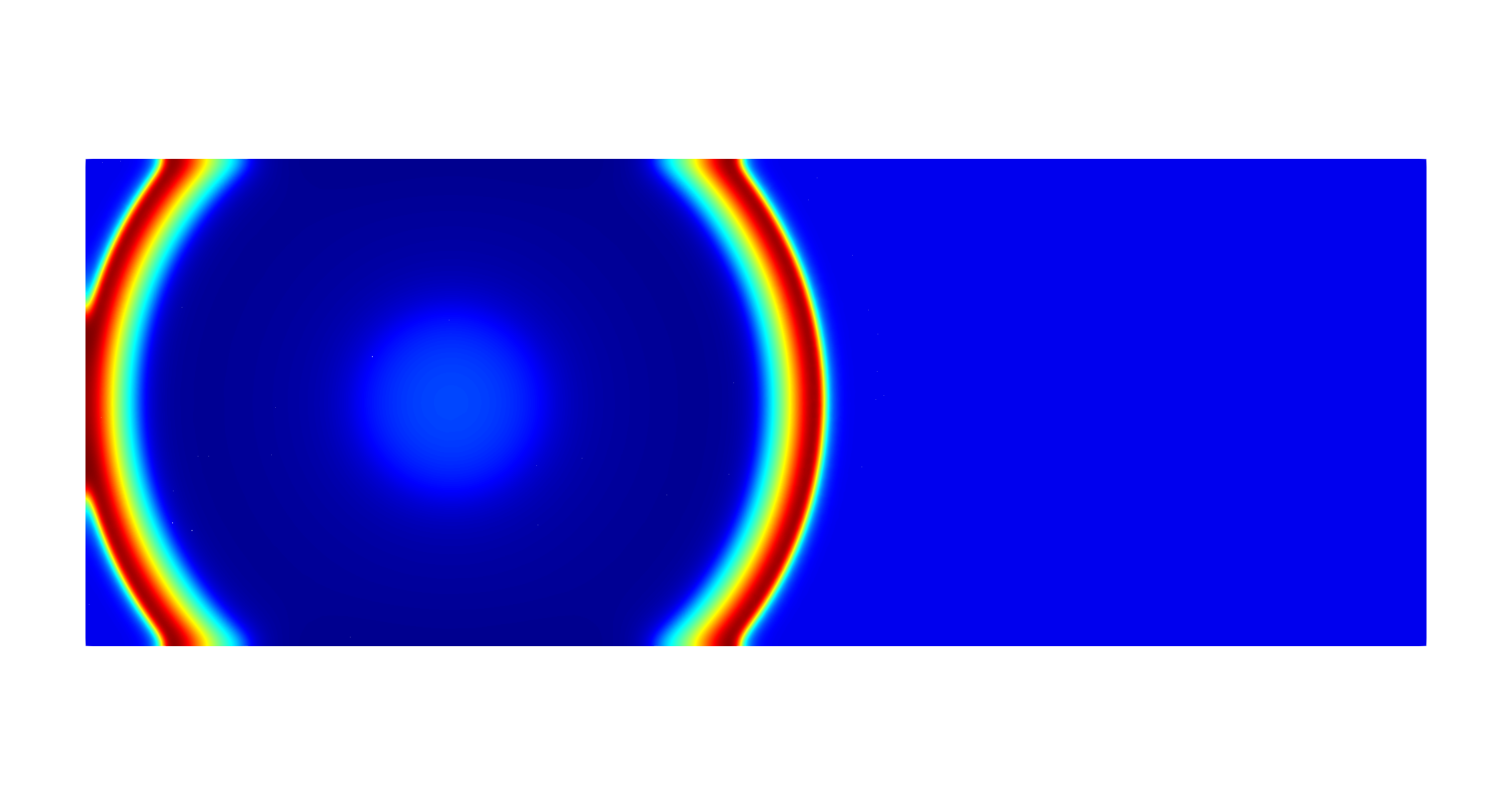}
\hspace*{1ex}\includegraphics[trim={5cm 7cm 5cm 7cm},clip, scale=0.07]{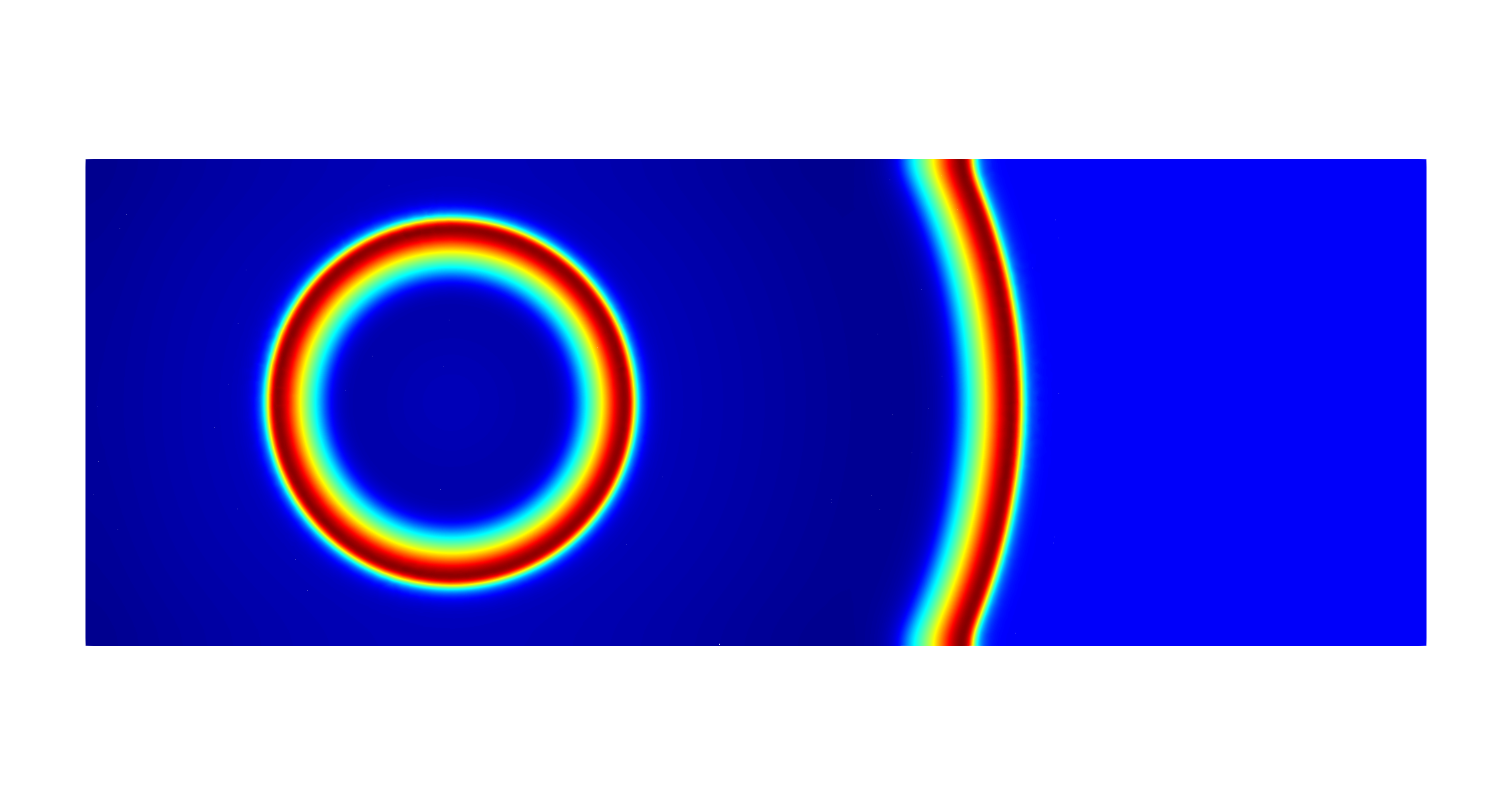}
\hspace*{1ex}\includegraphics[trim={5cm 7cm 5cm 7cm},clip, scale=0.07]{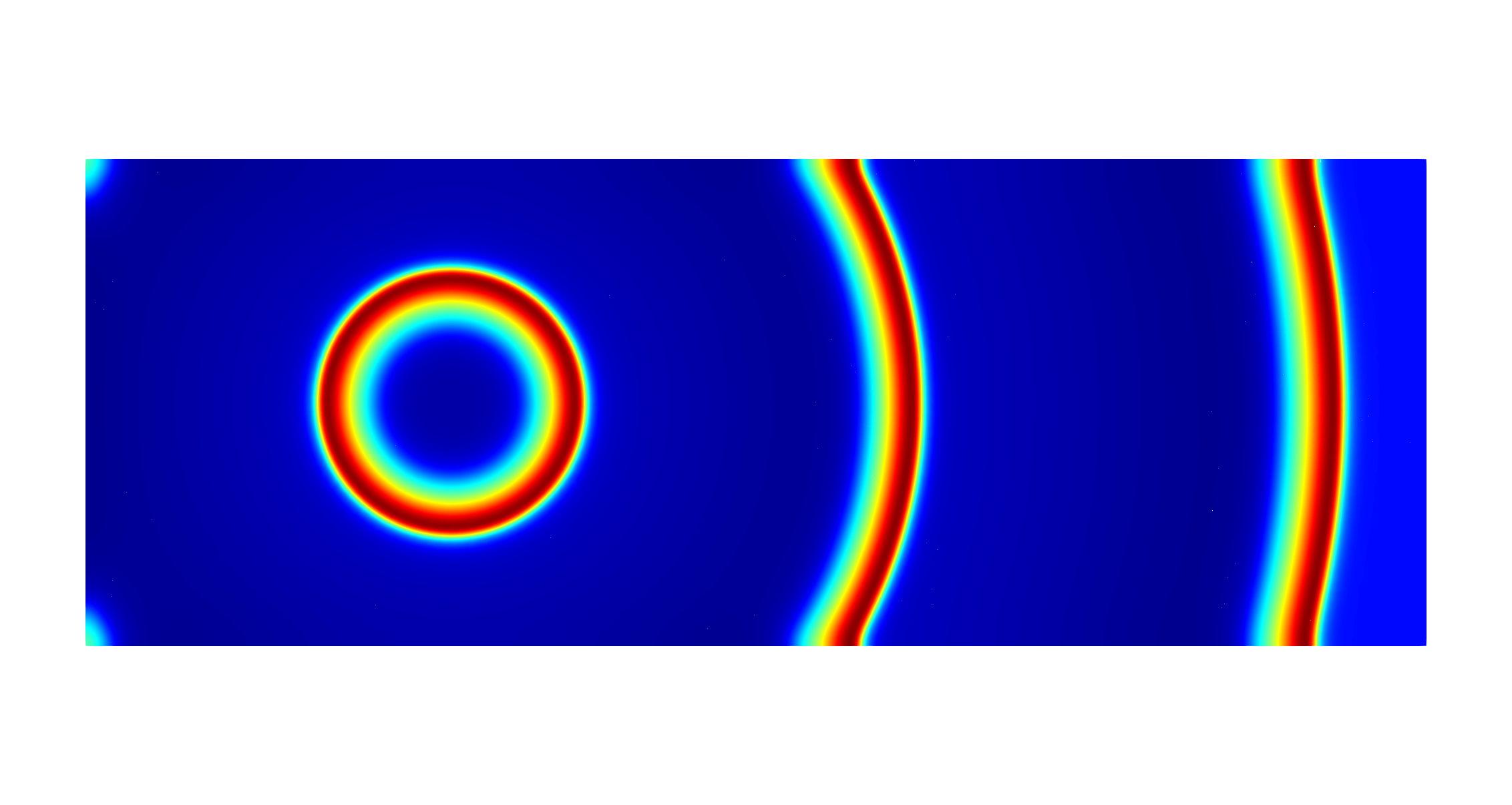}
 \subcaption{Transmembrane potential} \label{fig:bc_tranpot}
\end{minipage}
\begin{minipage}{0.31\textwidth}
\centering
\hspace*{-0.ex}\includegraphics[trim={5cm 7cm 5cm 7cm},clip, scale=0.07]{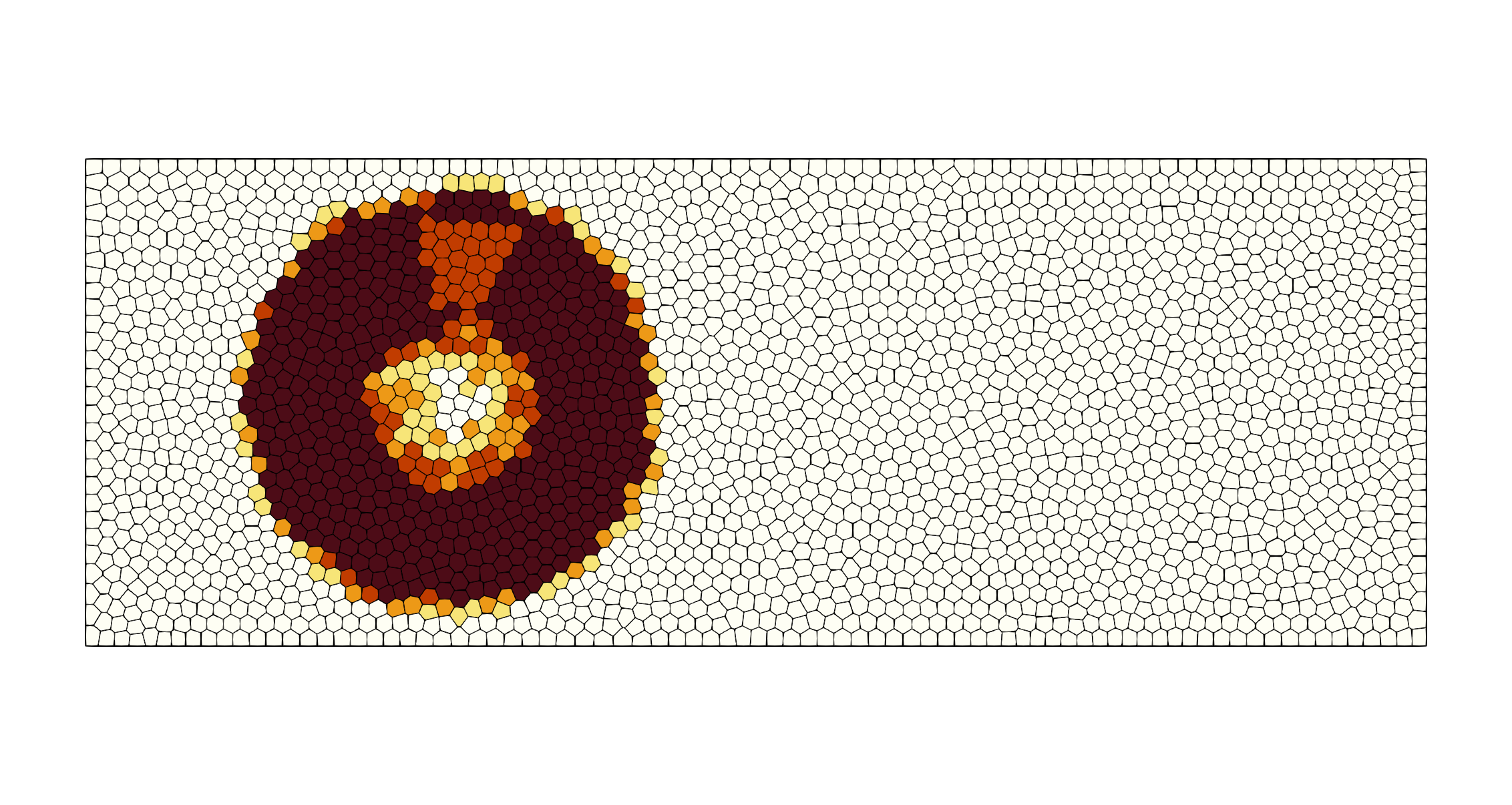}
\hspace*{-0.ex} \includegraphics[trim={5cm 7cm 5cm 7cm},clip, scale=0.07]{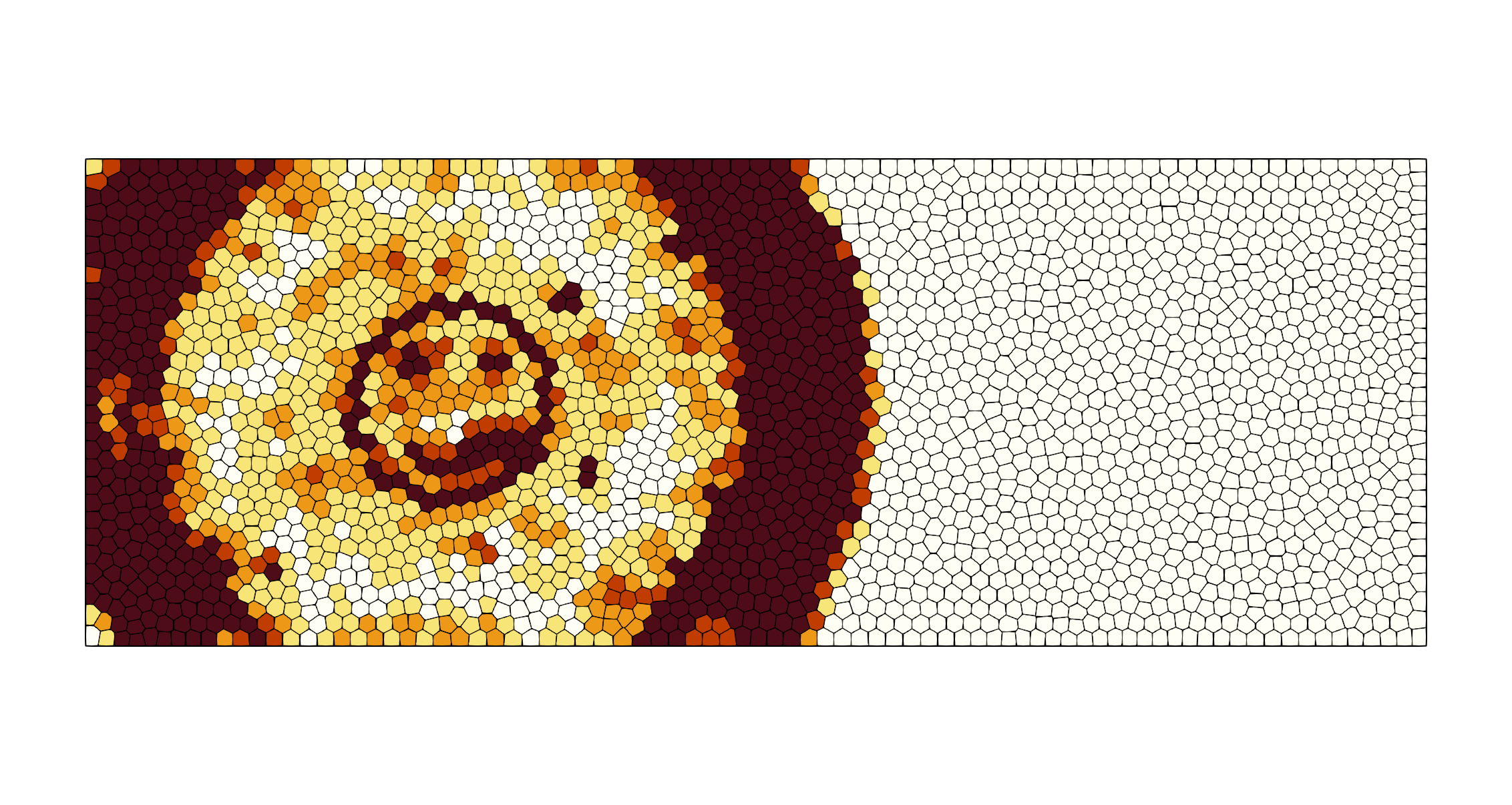}
\hspace*{-0.ex}\includegraphics[trim={5cm 7cm 5cm 7cm},clip, scale=0.07]{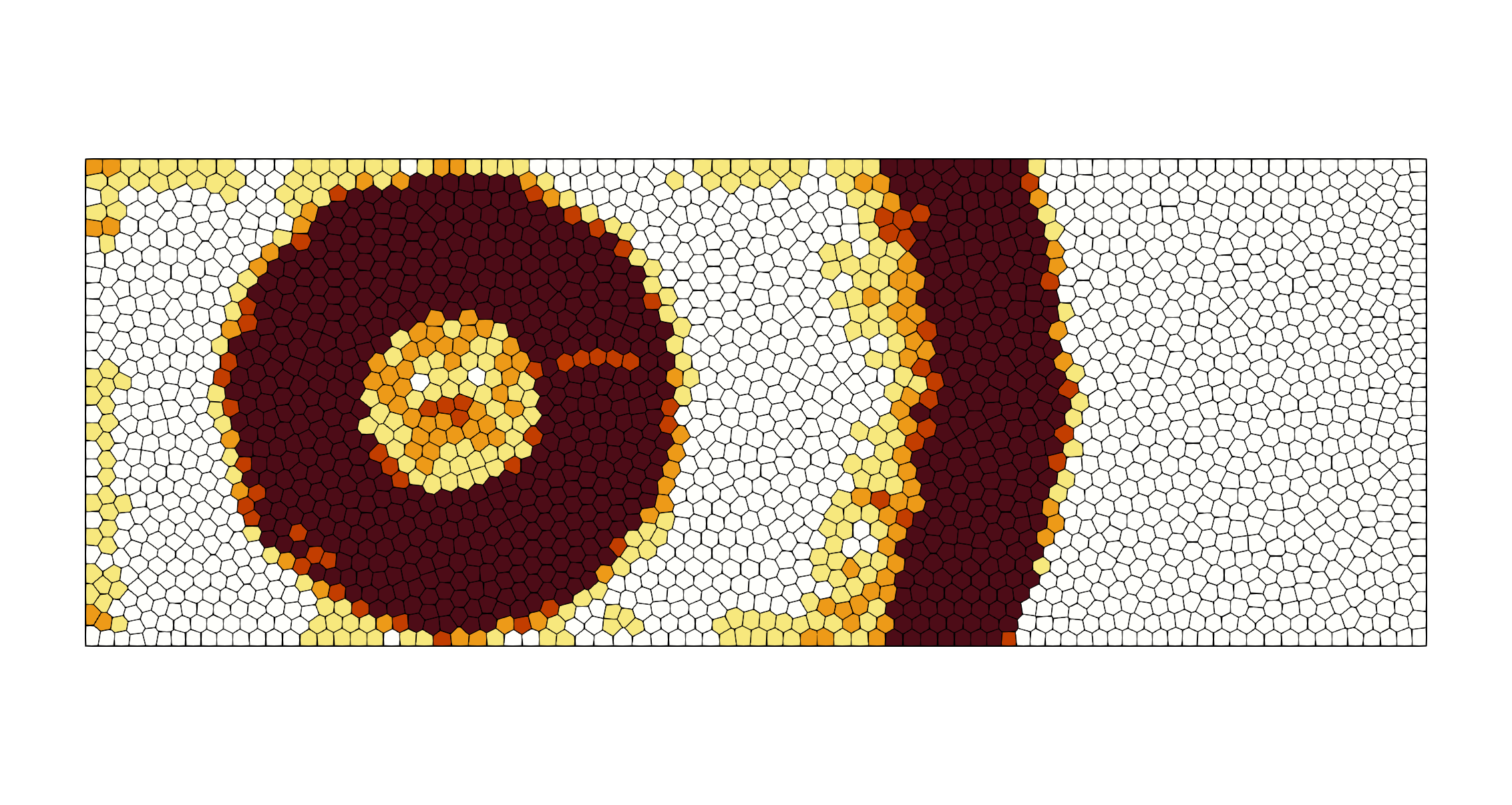}
\hspace*{-0.ex}\includegraphics[trim={5cm 7cm 5cm 7cm},clip, scale=0.07]{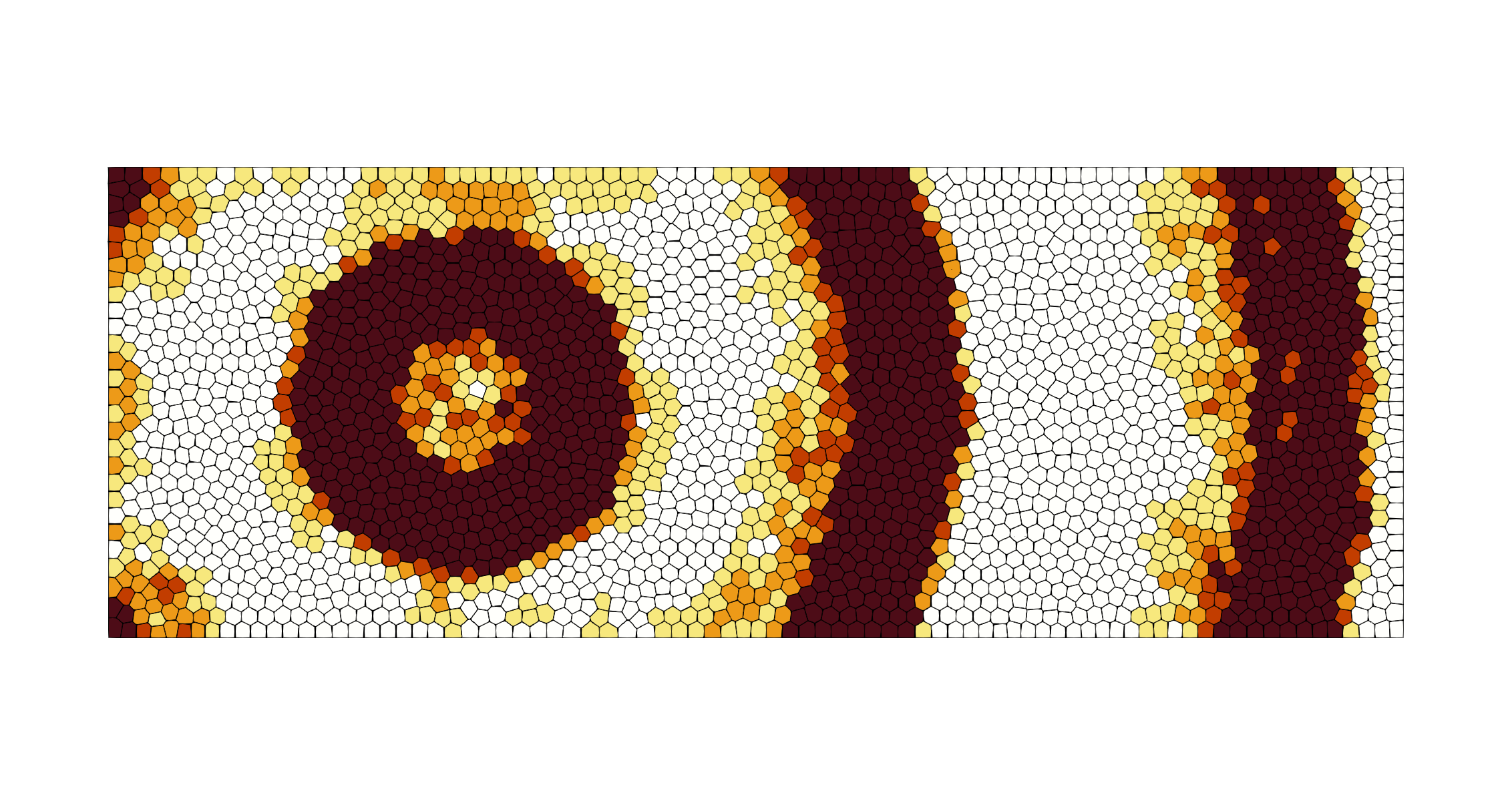}
      \subcaption{Polynomial approximation degree}  \label{fig:p_tranpot}
      \end{minipage}
\begin{minipage}{0.31\textwidth}
\centering
\hspace*{-0.ex}\includegraphics[trim={5cm 7cm 5cm 7cm},clip, scale=0.07]{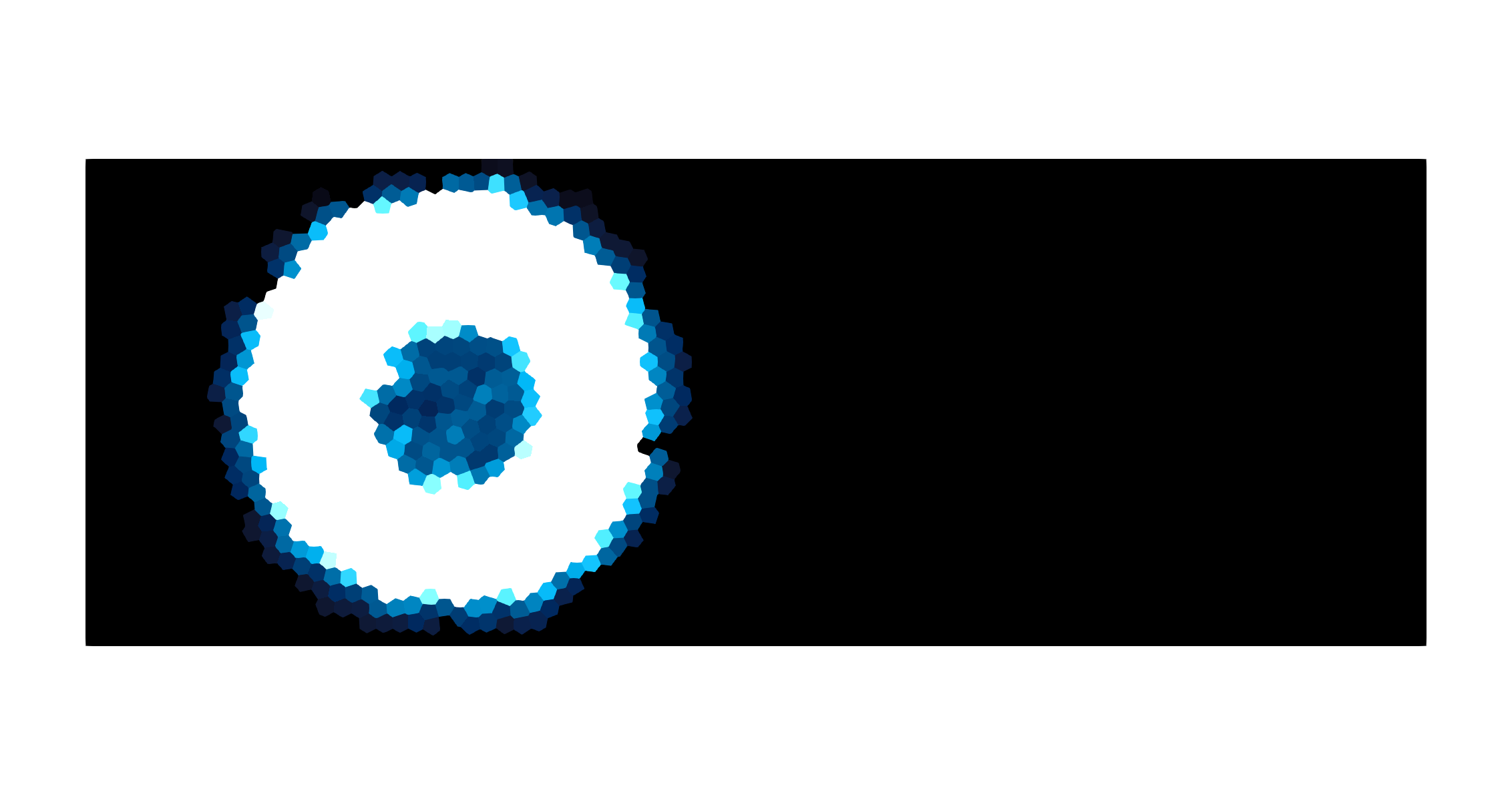}
\hspace*{-0.ex}\includegraphics[trim={5cm 7cm 5cm 7cm},clip, scale=0.07]{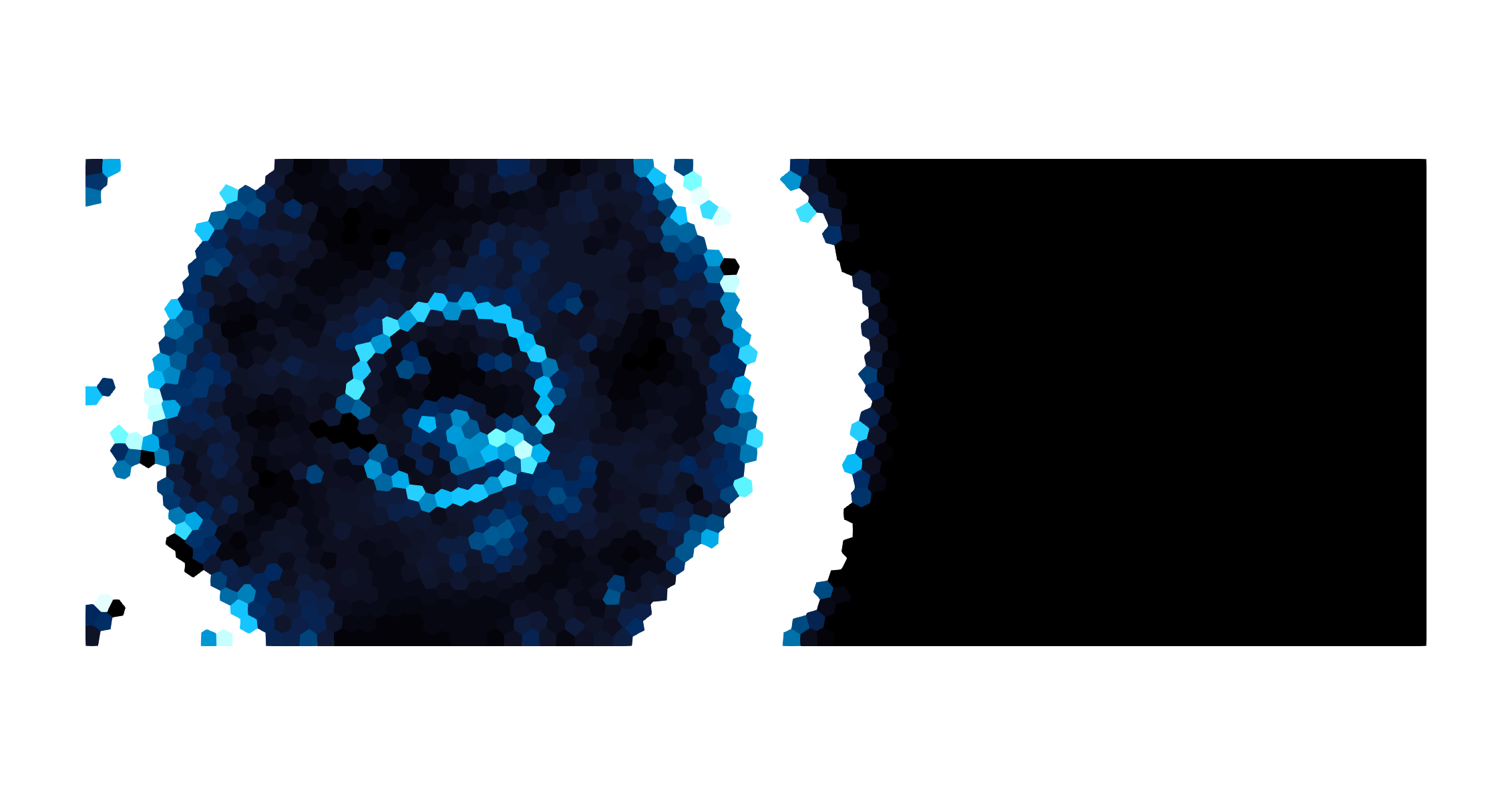}
\hspace*{-0.ex}\includegraphics[trim={5cm 7cm 5cm 7cm},clip, scale=0.07]{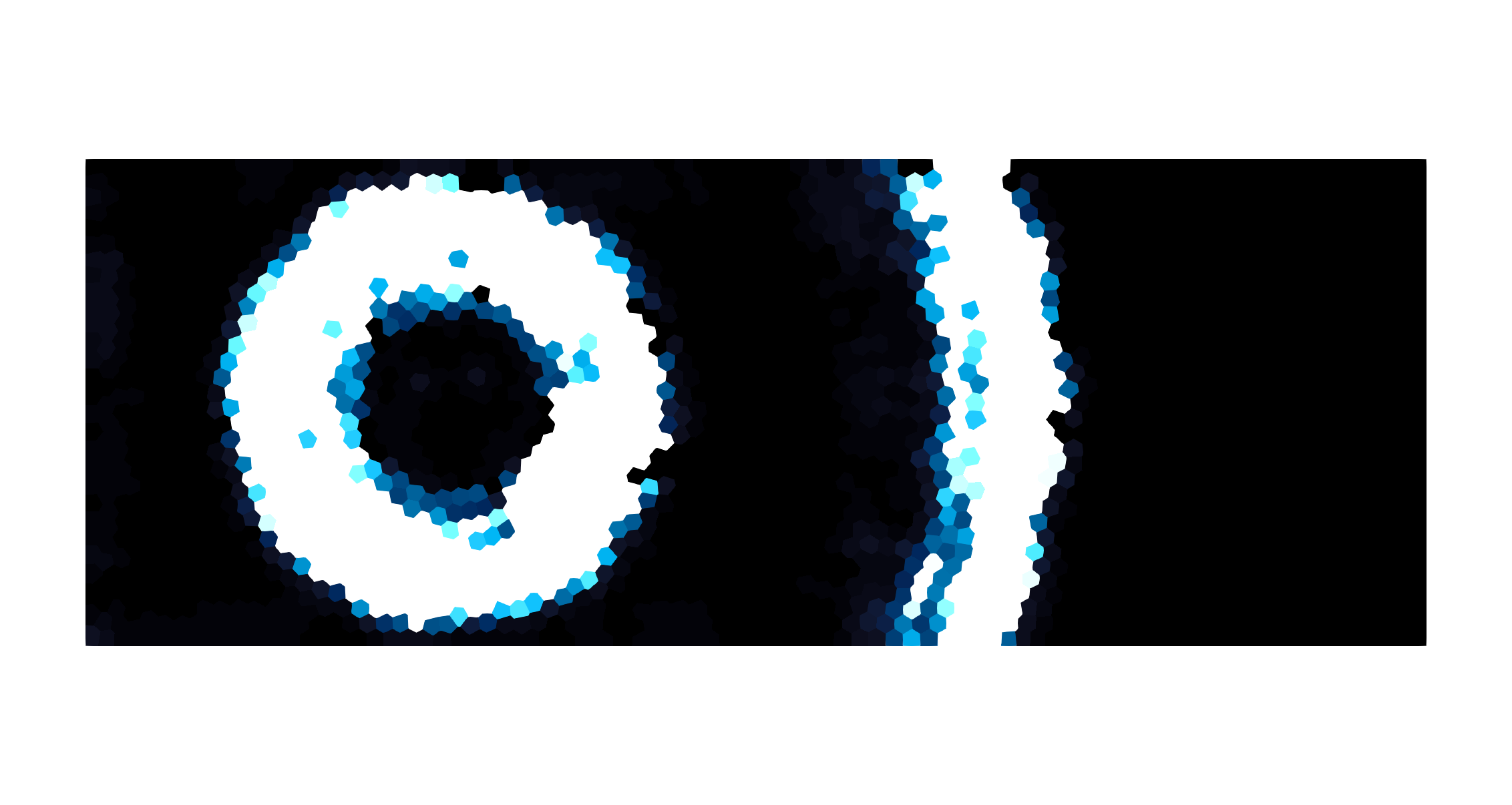}
\hspace*{-0.ex}\includegraphics[trim={5cm 7cm 5cm 7cm},clip, scale=0.07]{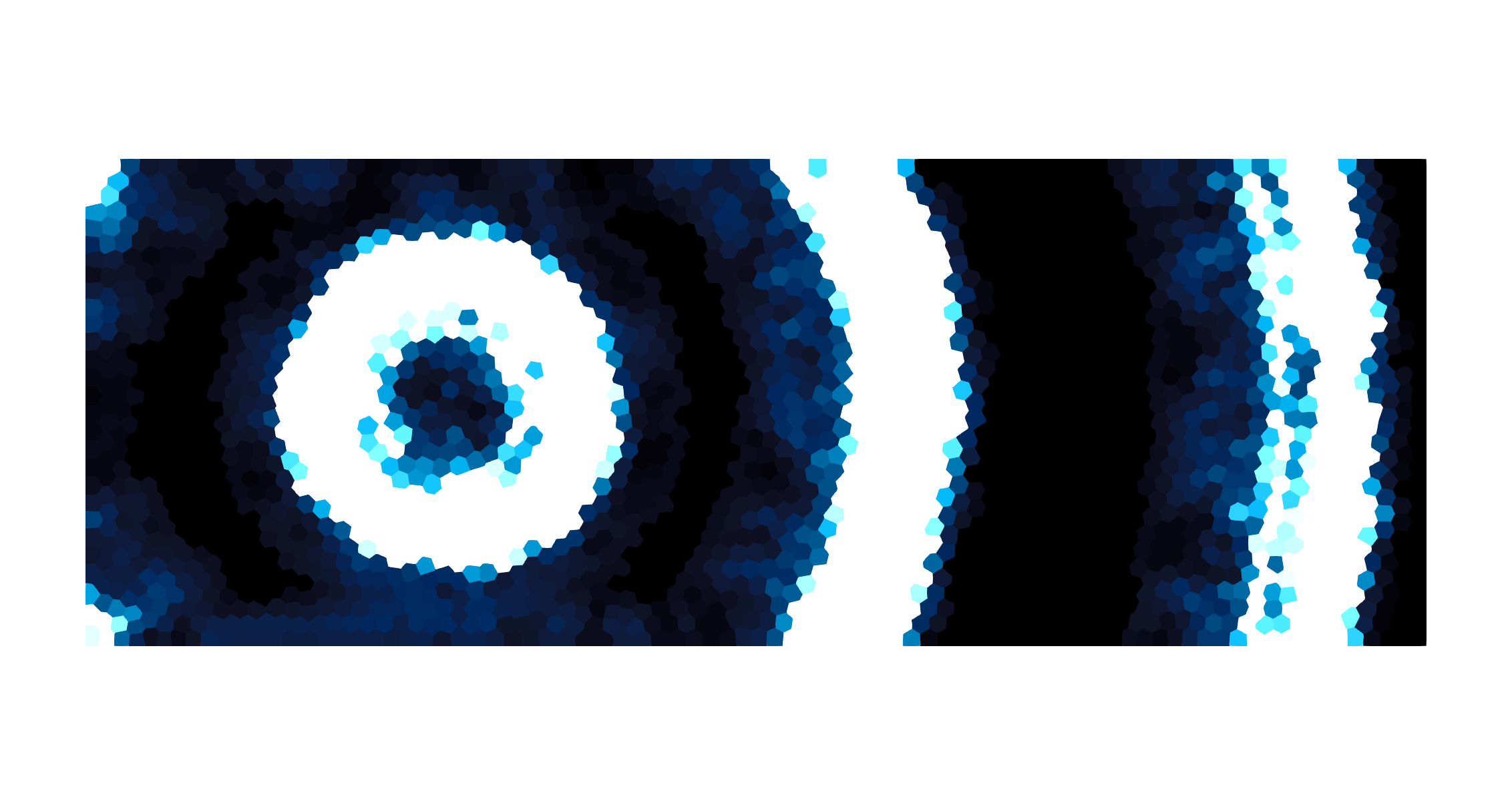}
      \subcaption{Local error indicator}  \label{fig:p_tranpot}
      \end{minipage}
  \caption{Test case 3.  Computed transmembrane potential in a section of grey matter tissue (left), element-wise polynomial approximation degree (middle), and local error indicator defined as in \eqref{eq::indicator} (right) at different time snapshots $t=3.4, 7.2, 8.4, 13.2 \, [\mathrm{ms}]$.}
  \label{fig::bc_wave}
  \vspace{2ex}
\end{figure}
Figure~\ref{fig::bc_wave} reports the computed solution, the element-wise polynomial approximation degree, and the computed total local error indicator defined as in \eqref{eq::indicator} at different time snapshots $t=3.4, 7.2, 8.4, 13.2 \, [\mathrm{ms}]$. As expected, larger values of the local error indicator and, as a consequence, higher-order polynomials are observed in correspondence with the depolarization and polarization wavefronts. As it can also be seen from the evolution of the 0D-model in Figure~\ref{fig:bc0d_1}, an equally rapid decay (polarization) follows the potential peak, which also requires high-order polynomial degrees for the accurate local approximation.
\begin{figure}[!htbp]
    \centering
    \begin{subfigure}[t]{0.3\textwidth}
        \centering
\includegraphics[width=0.98\linewidth]{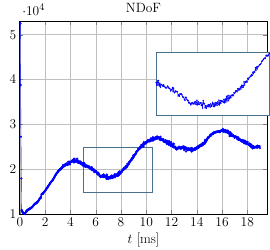}
\caption{\label{fig:dof_evolution_bc}}
    \end{subfigure}\hfill
    \begin{subfigure}[t]{0.3\textwidth}
        \centering
\includegraphics[width=1\linewidth]{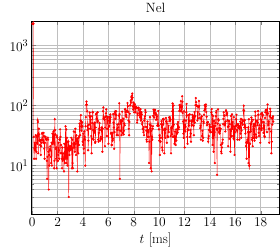}
   \caption{\label{fig:comptau_bc}}
    \end{subfigure}\hfill
     \begin{subfigure}[t]{0.3\textwidth}
        \centering
\includegraphics[width=1.03\linewidth]{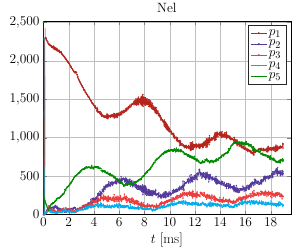}
     \caption{\label{fig:elem_bc}}
    \end{subfigure}
 \caption{Test case 3. Left: Evolution of the number of degrees of freedom ($\mathrm{NDoF}$) as a function of time, driven by our $p$-adaptive algorithm.  Center: Evolution of the number of elements ($\mathrm{Nel}$) where the local polynomial degree is updated over time.  Right: Time evolution of the number of elements ($\mathrm{Nel}$) with local polynomial degrees equal to $1,2,3,4,5$ over time.  At the initial time all the elements are discretized with $p_\text{max}=5$.} \label{fig:ndof_evolution_bc}
\end{figure}
Regions of quiescence correspond instead to linear polynomial approximation degrees, where the indicator values are also lower. A new pattern can be appreciated in the results reported in Figure~\ref{fig:ndof_evolution_bc}, where we show the time evolution of the $\mathrm{NDoFs}$. As expected, there is a significant increase in $\mathrm{NDoFs}$ approximately at $7\, [\mathrm{ms}]$ and $13.5\, [\mathrm{ms}]$ due to the generation of the new traveling waves. Figure~\ref{fig:ndof_evolution_bc} (center and right) also reports the time evolution of the total number of elements where refinement/de-refinement takes place and the total number of elements where local polynomial degrees of order $1,2,3,4,5$ are employed. The number of elements where linear polynomials are employed is the majority throughout the simulation, but slowly decreases in the first $14\, [\mathrm{ms}]$ of the simulation till reaching a dynamic plateau. This is consistent with the underlying physics, since the multiple neuronal waves are steep and localized. The computational costs for computing the indicator slightly increase when the second (and third) wave is introduced, as the number of elements to be checked increases accordingly. From the numerical results, it can be concluded that our $p$-adaptive algorithm demonstrates the ability to automatically adjust the polynomial degree where neuronal waves are located.

\subsection{Brain Stem}
\label{sec:test_brainstem}
\begin{figure}[!htbp]
\centering
\begin{subfigure}[b]{0.45\textwidth}
    \centering
    \includegraphics[trim={10cm 0.5cm 10cm 0.5cm},clip,scale=0.085]{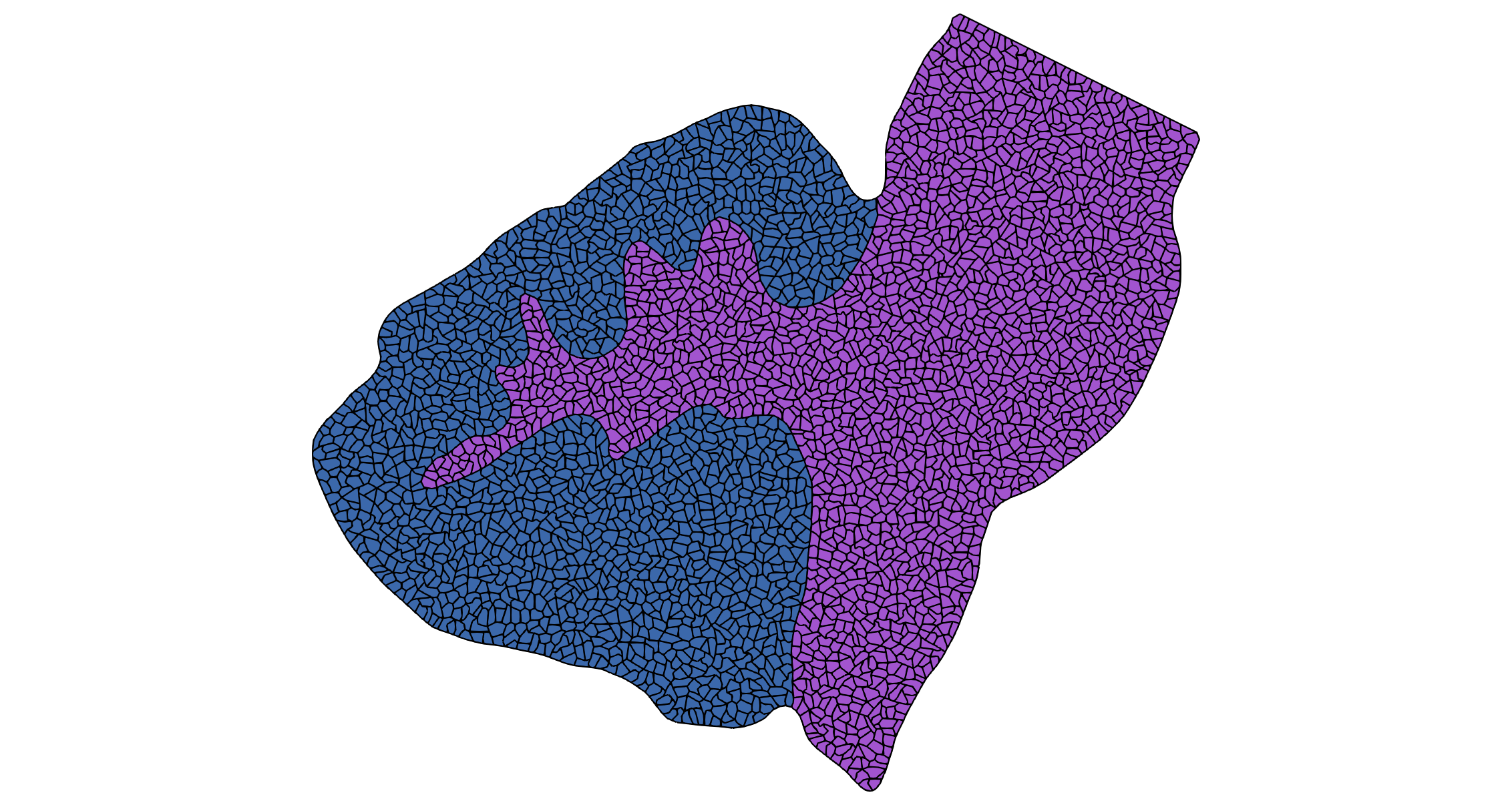}
     \caption{\label{fig::bc_meshStem}}     
\end{subfigure}
\begin{subfigure}[b]{0.45\textwidth}
    \centering
\hspace*{-2ex}\includegraphics[trim={10cm 0.5cm 10cm 0.5cm},clip,scale=0.085]{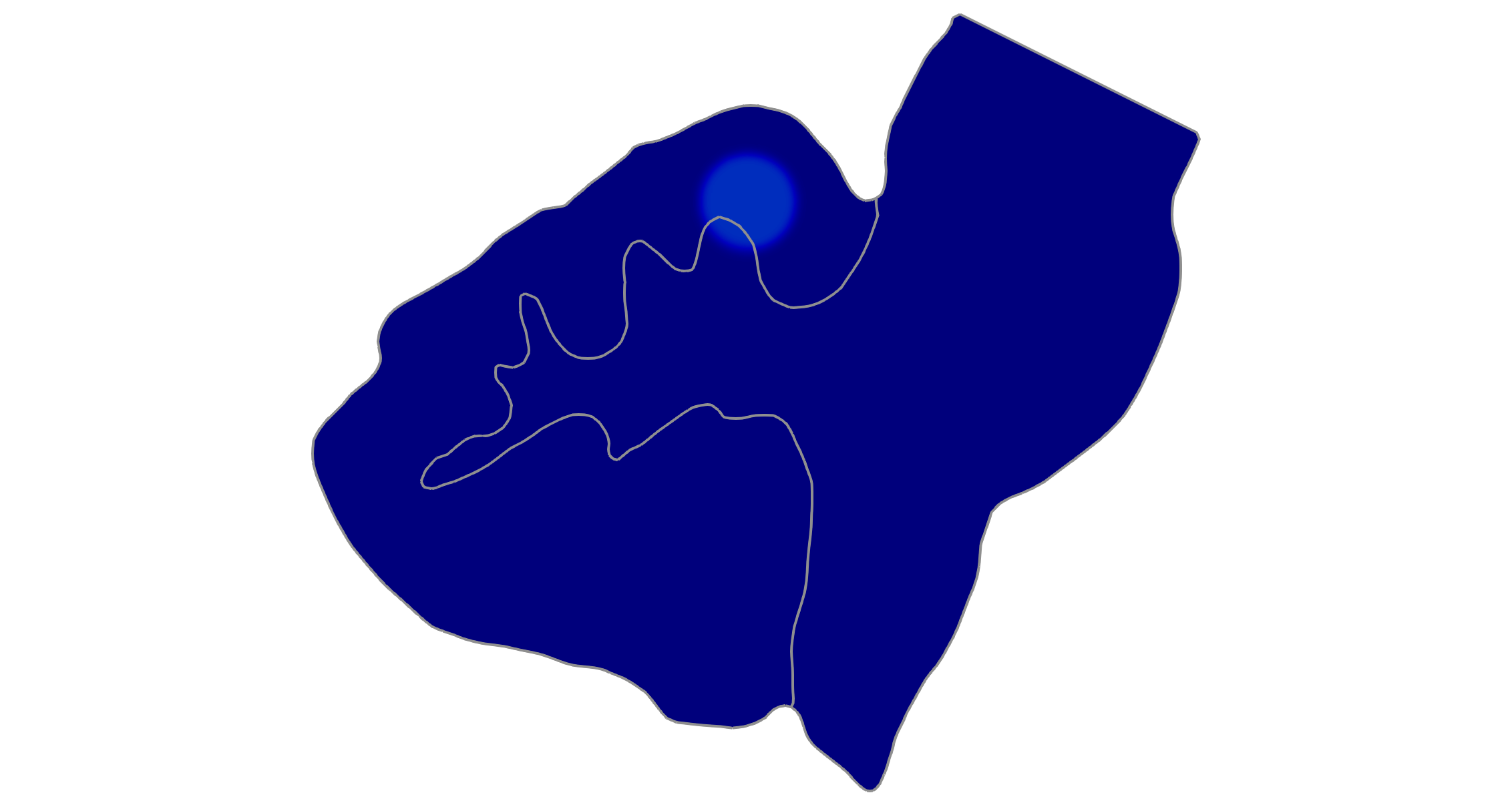} 
\vspace{1ex}
\hspace*{0ex}\includegraphics[trim={0cm 0cm 0cm 0cm},clip,scale=0.17]{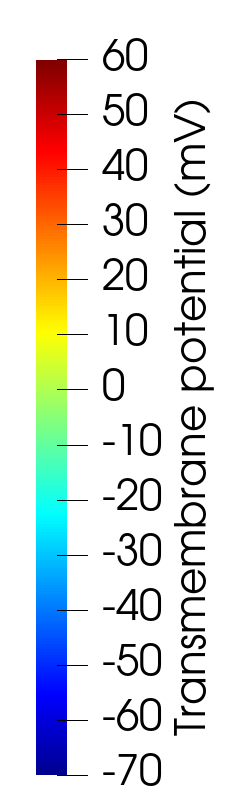} 
\caption{\label{fig::initial_bcStem}}
  \end{subfigure}
  \caption{Test case 4.  Agglomerated polytopal grid consisting of 3523 elements of the sagittal section of a human brain stem
  with differentiation of tissues as white matter (purple) and grey matter (blue) (left); pathological initial condition for the transmembrane potential ($u_0 = -50 \, [\mathrm{mV}]$ in $\Omega_0$).}
\label{fig::bc_initial_bcStem}
\end{figure}

\begin{figure}[h!]
    \centering
    \begin{subfigure}[t]{0.3\textwidth}
        \centering
    \includegraphics[width=0.95\linewidth]{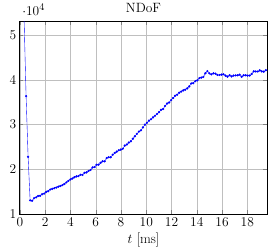}
\caption{\label{fig:dof_evolution_bcStem}}
    \end{subfigure}\hfill
    \begin{subfigure}[t]{0.3\textwidth}
    \includegraphics[width=1\linewidth]{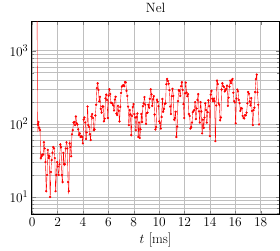}
    \caption{\label{fig:comptau_bcStem}} 
    \end{subfigure}\hfill
    \begin{subfigure}[t]{0.3\textwidth}
\includegraphics[width=1.05\linewidth]{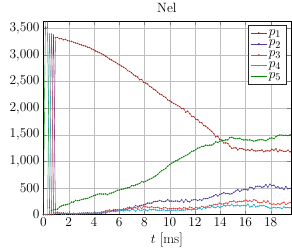}
        \caption{\label{fig:elem_bcStem}}
    \end{subfigure}
\caption{Test case 4. Left: Evolution of the number of degrees of freedom ($\mathrm{NDoF}$) as a function of time, driven by our $p$-adaptive algorithm.  Center: Evolution of the number of elements ($\mathrm{Nel}$) where the local polynomial degree is updated over time.  Right: Evolution of the number of elements ($\mathrm{Nel}$) with local polynomial degree equal to $1,2,3,4,5$ over time.  At the initial time, all the elements are discretized with $p_\text{max}=5$.}
\label{fig:ndof_evolution_bcStem}
\end{figure}

In this test, we evaluate the performance of the proposed $p$-adaptive PolyDG method in a geometrically complex domain representing a sagittal section of a human brain stem reconstructed starting from structural Magnetic Resonance Images (MRI) from the OASIS-3 database \cite{lamontagne2019oasis}. The considered brain sagittal section has a maximum longitudinal length of $7.5$ cm and a vertical length of $7 \, \mathrm{cm}$.
Here, we simulate pathological seizures dynamics in a brain section characterized by heterogeneous tissue, such as grey and white matter. As in the previous sections, we consider an external forcing term defined as in Figure \ref{fig:bc0d_0}.
An initial fine triangular grid has been agglomerated using the approach based on neural networks proposed in  \cite{antonietti2022agglomeration}, leading to a polytopal mesh of 3523 elements, cf. Figure \eqref{fig::bc_meshStem}, where we also show the subdivision into grey and white matter. 
Out of the $3523$ polygons constituting our agglomerated mesh, $2012$ are in the white matter and the remaining in the grey matter.
In Figure \eqref{fig::initial_bcStem}, we report the unstable initial condition. The behavior of white and grey matter in response to transmembrane ionic imbalances or external stimuli differs because of the preponderance of cell bodies in grey matter and axons in white matter. This makes the white matter anisotropic; in this simulation, we consider the following values for the conductivity tensors for grey and white matter
\begin{equation}
    \begin{aligned}
      \Sigma_\text{GM} = \begin{bmatrix}
          0.7354 & 0 \\ 0 & 0.7354
      \end{bmatrix} \; , \; \Sigma_\text{WM} = \begin{bmatrix}
          1.7354 & 0 \\ 0 & 1.2854
      \end{bmatrix},
\end{aligned}
\label{eq:sigma_Stem}
\end{equation}
respectively. For this simulation, we consider the same values of the parameters and initial condition as before, cf. Table~\ref{table:G2} and Table~\ref{table:BC_IC}. As before, the initial condition is discretized $p_\text{max}=5$ in all the mesh elements, corresponding to $7.3e+4$ degrees of freedom.
\begin{figure}[h!]
 \centering
   \begin{subfigure}[b]{\textwidth}
    \centering
    \includegraphics[trim={0 0cm 0 0cm},clip,scale=0.18]{photosbc/scale.png}
    \includegraphics[trim={0 0cm 0 0cm},clip, scale=0.18]{singlewave_stationary/scale_p.png}
    \includegraphics[trim={0 0.1cm 0 0.1cm},clip, scale=0.19]{photosbc/scale_taublue.png}
    \end{subfigure}\hfill
\begin{subfigure}[b]{0.03\textwidth}
\raisebox{10mm}{
\hspace*{-0.3cm}
  \begin{tabular}{@{}l@{}}
$t = 1[\mathrm{ms}]$\\[2.87cm]
$t = 9.4[\mathrm{ms}]$\\[2.87cm]
$t = 17.2[\mathrm{ms}]$
  \end{tabular}}
  \end{subfigure}
 \begin{minipage}{0.31\textwidth}
 \centering
\includegraphics[trim={12cm 1cm 15cm 1cm},clip, scale=0.085]{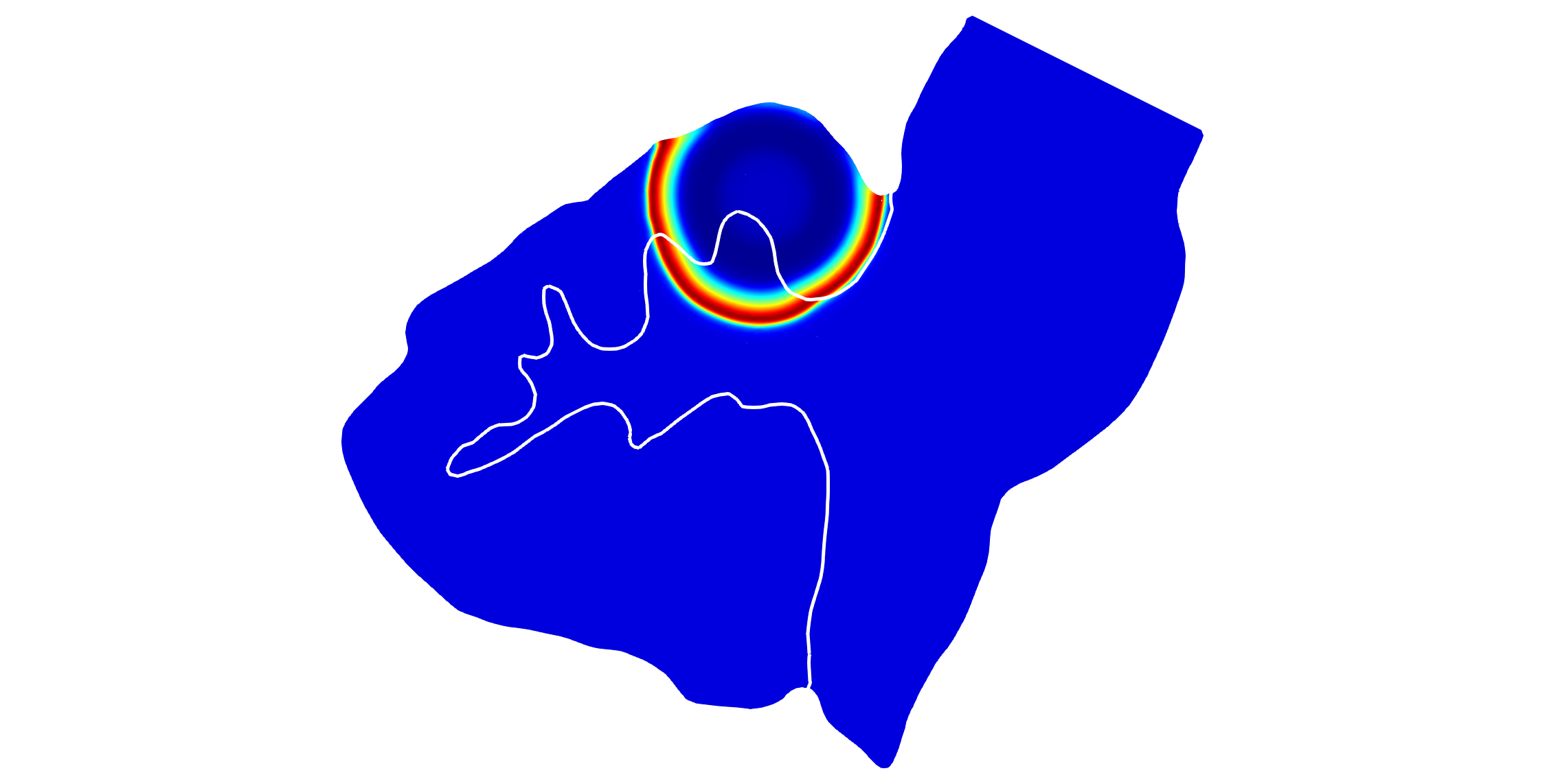} 
 \includegraphics[trim={12cm 1cm 15cm 1cm},clip, scale=0.085]{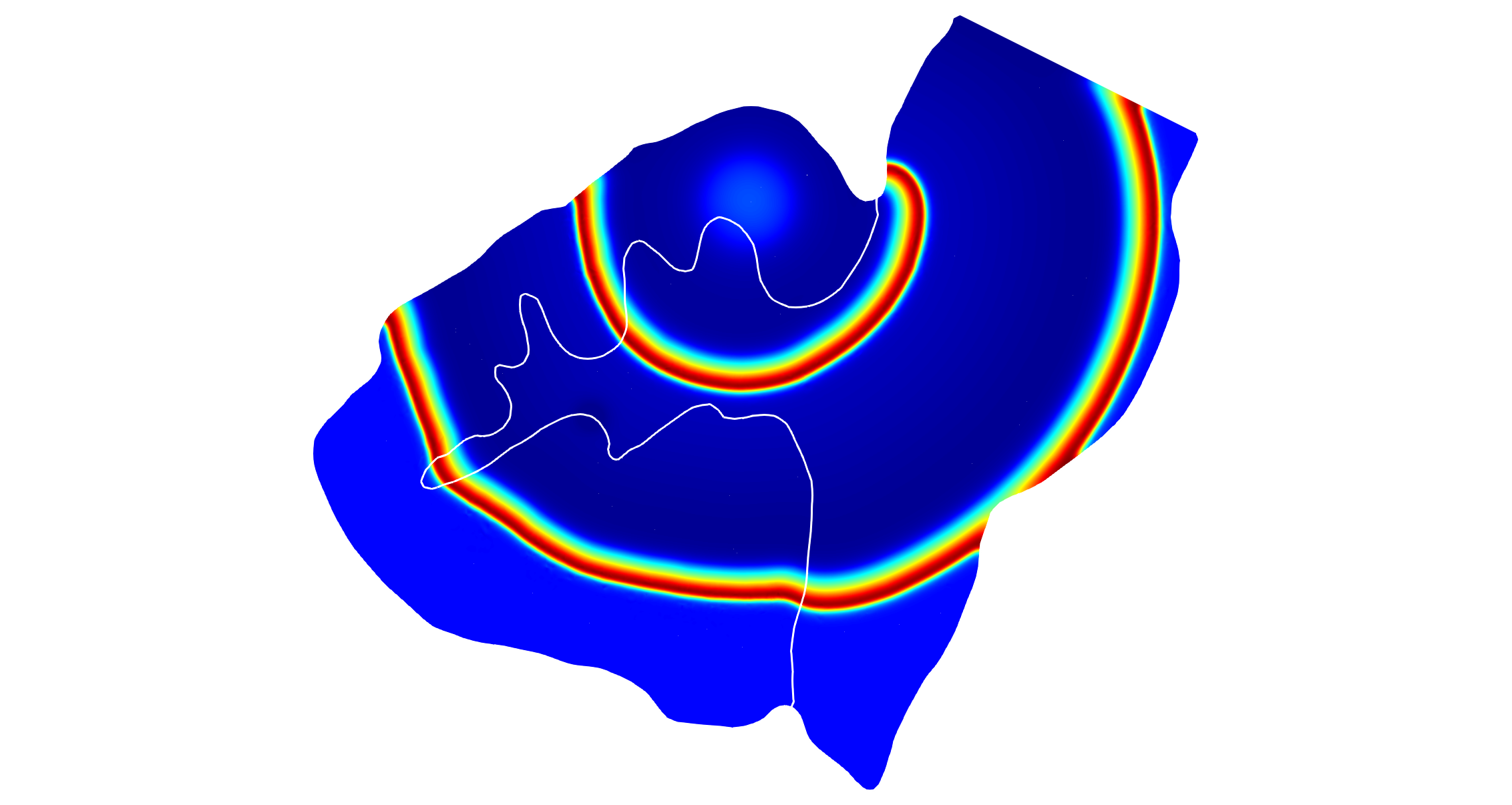}
 \includegraphics[trim={12cm 1cm 15cm 1cm},clip, scale=0.085]{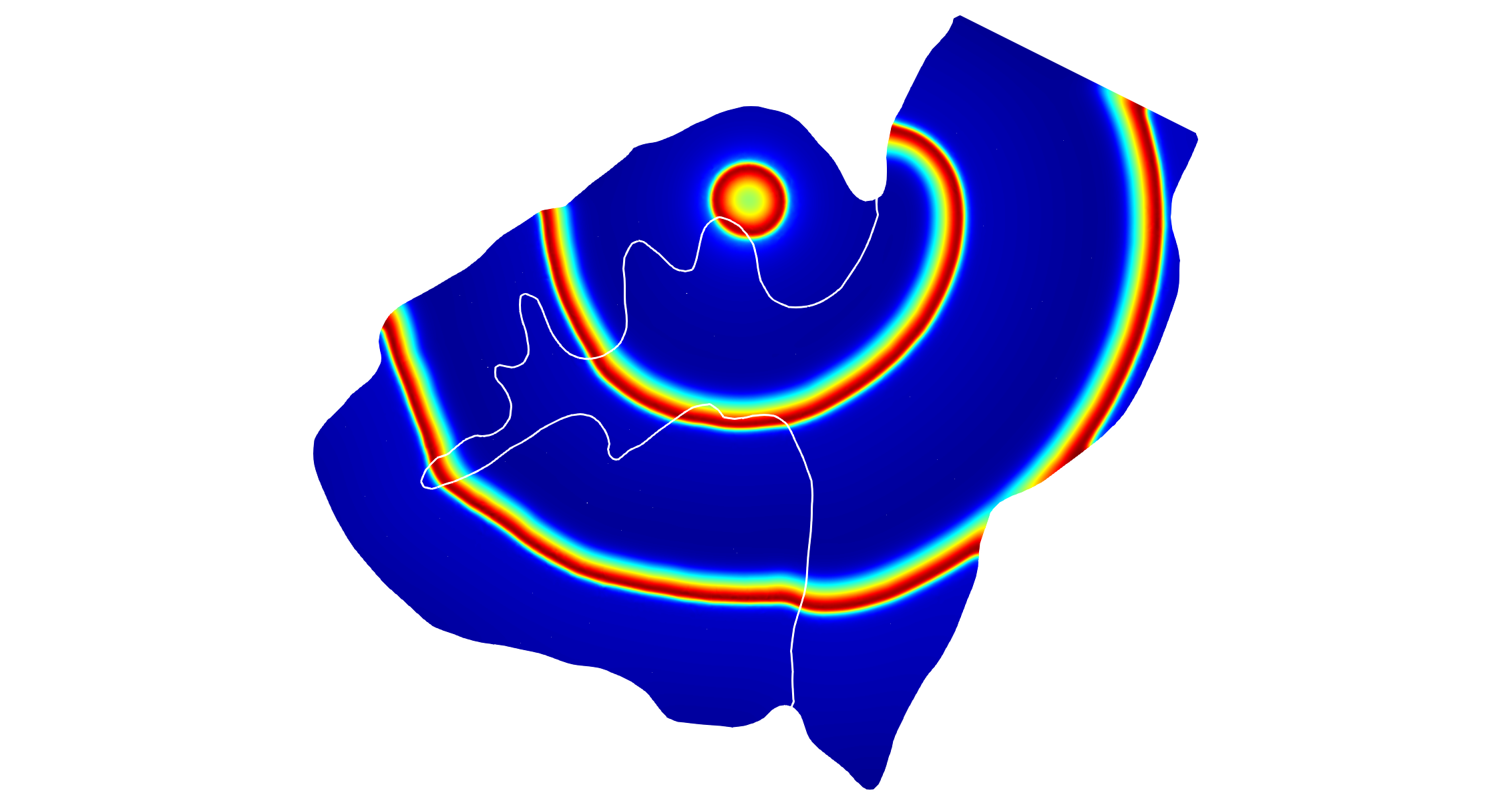}
 \subcaption{Transmembrane potential} \label{fig:bc_tranpot}
\end{minipage}
\begin{minipage}{0.31\textwidth}
\centering
\includegraphics[trim={12cm 1cm 15cm 1cm},clip, scale=0.085]{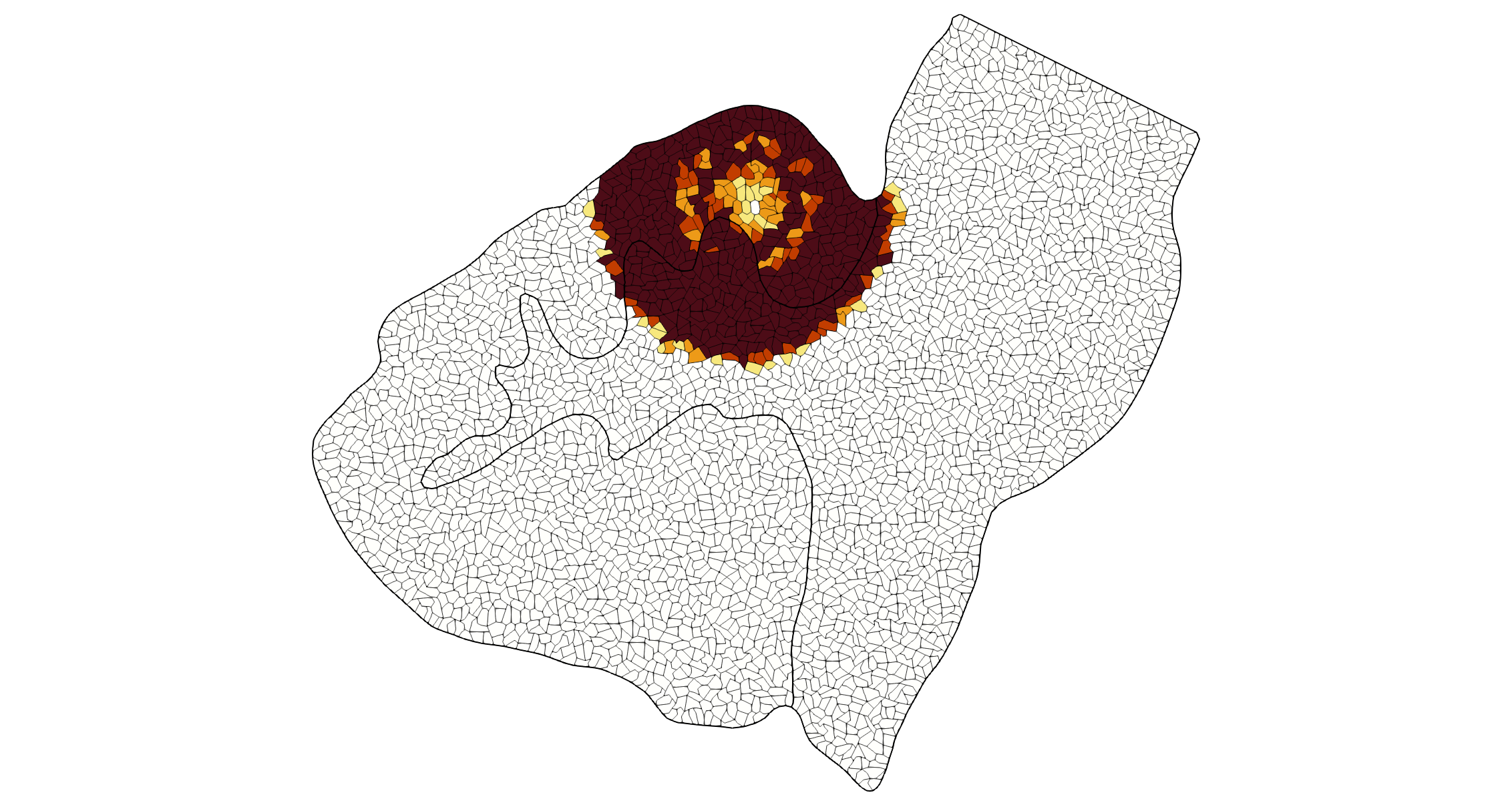}
\includegraphics[trim={12cm 1cm 15cm 1cm},clip, scale=0.085]{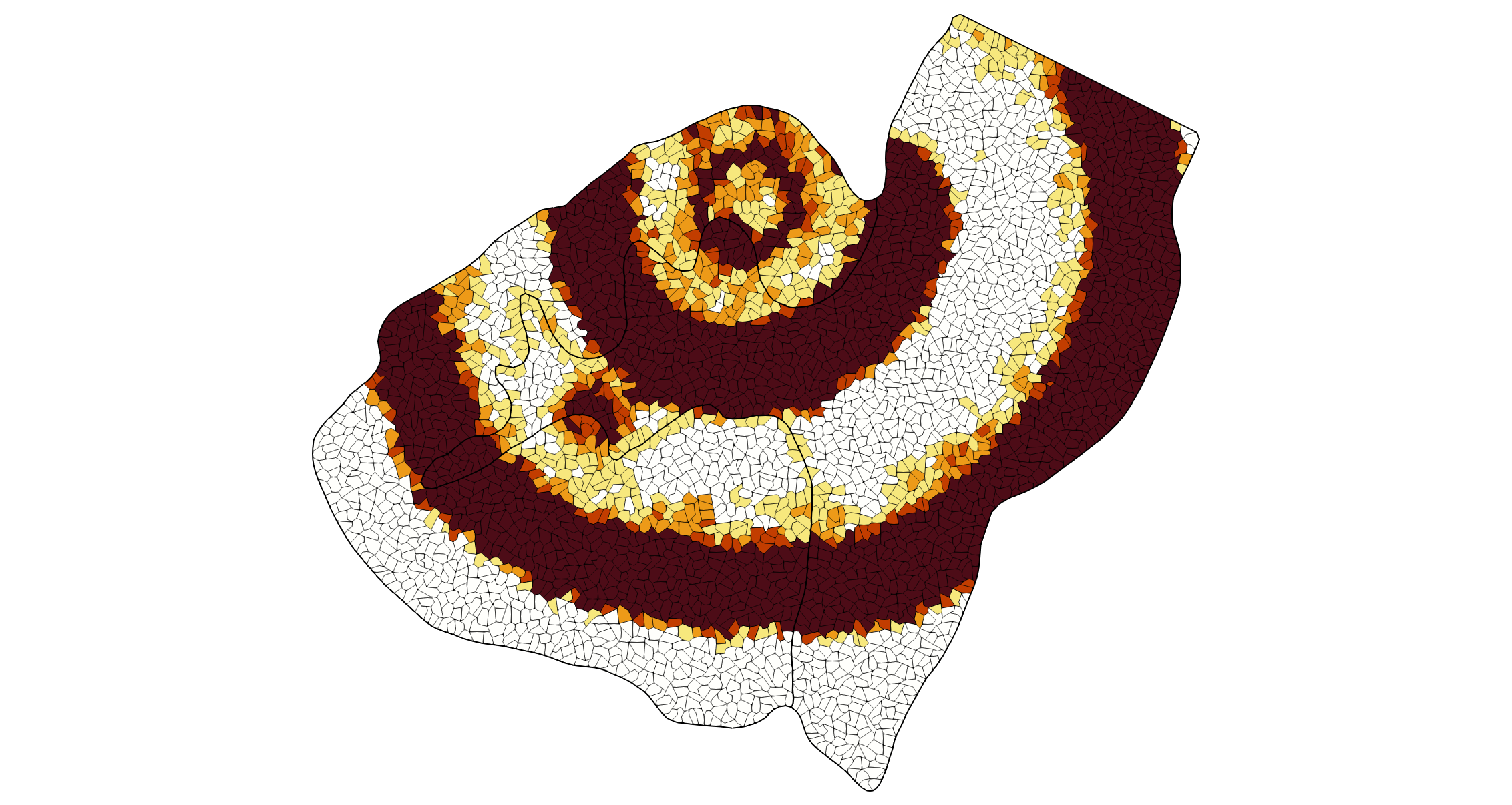}
\includegraphics[trim={12cm 1cm 15cm 1cm},clip, scale=0.085]{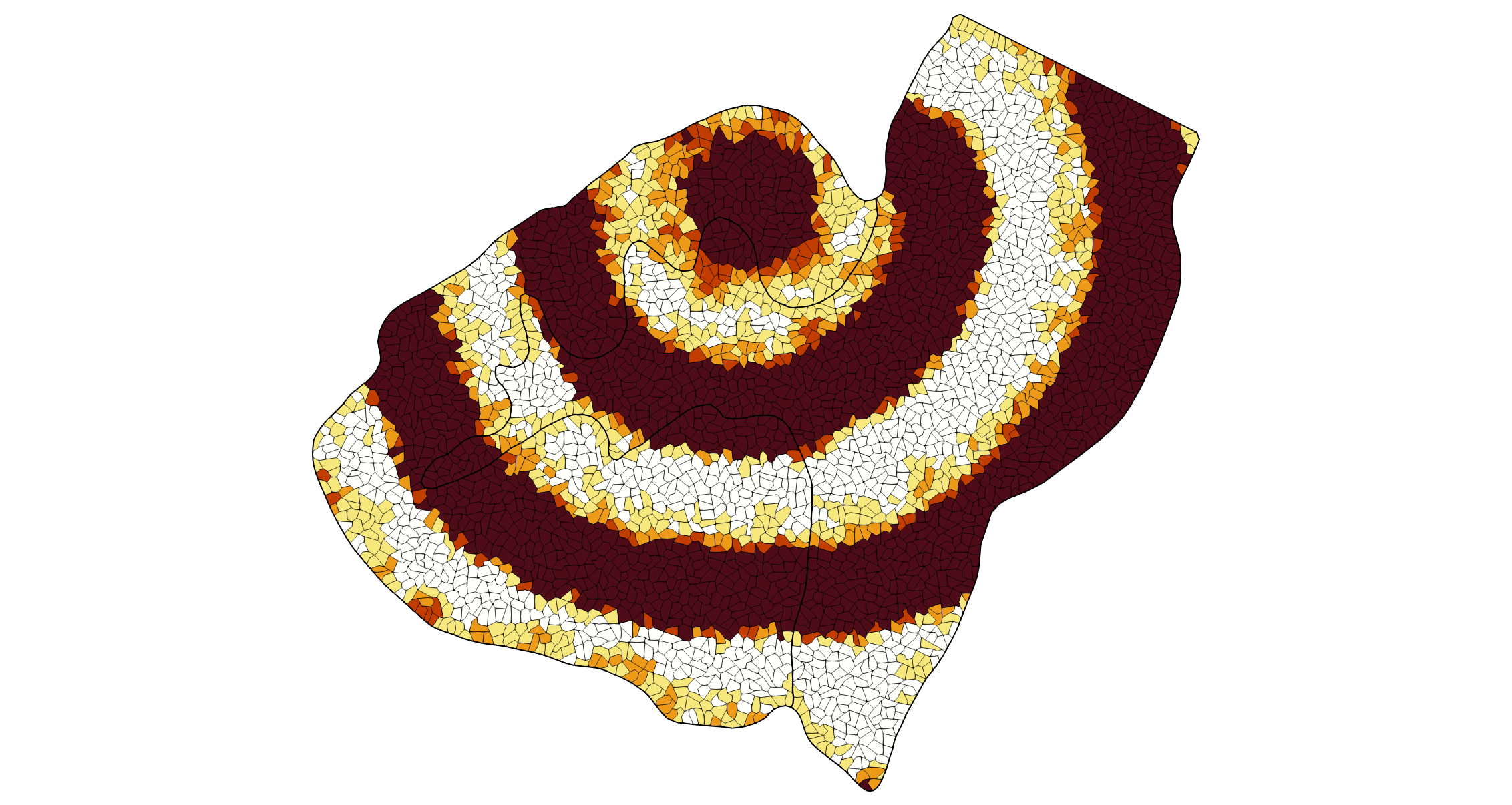}
\subcaption{Polynomial degree}  \label{fig:p_tranpot}
\end{minipage}
\begin{minipage}{0.31\textwidth}
\centering
\includegraphics[trim={11cm 1cm 13cm 1cm},clip, scale=0.085]{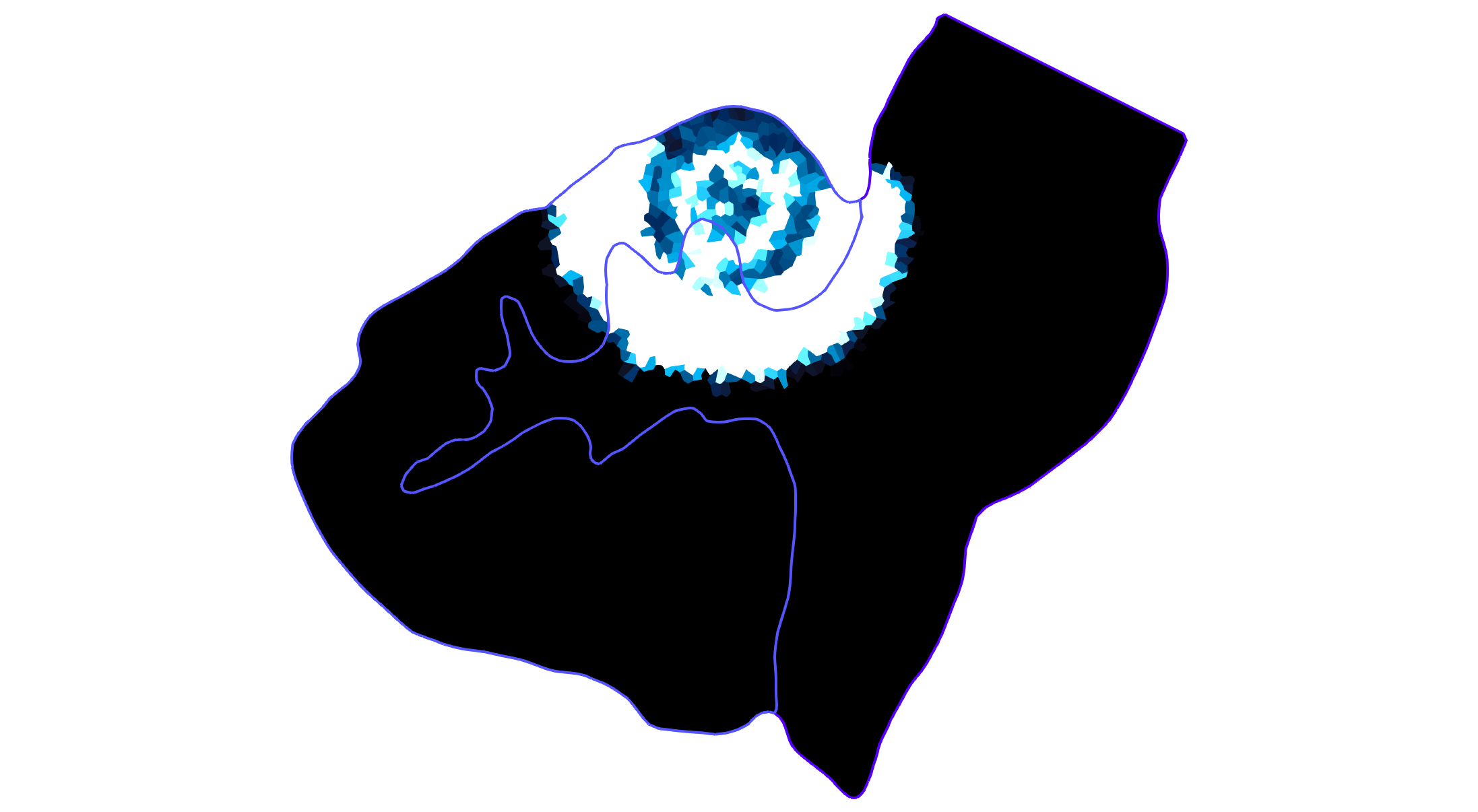}
\includegraphics[trim={11cm 1cm 13cm 1cm},clip, scale=0.085]{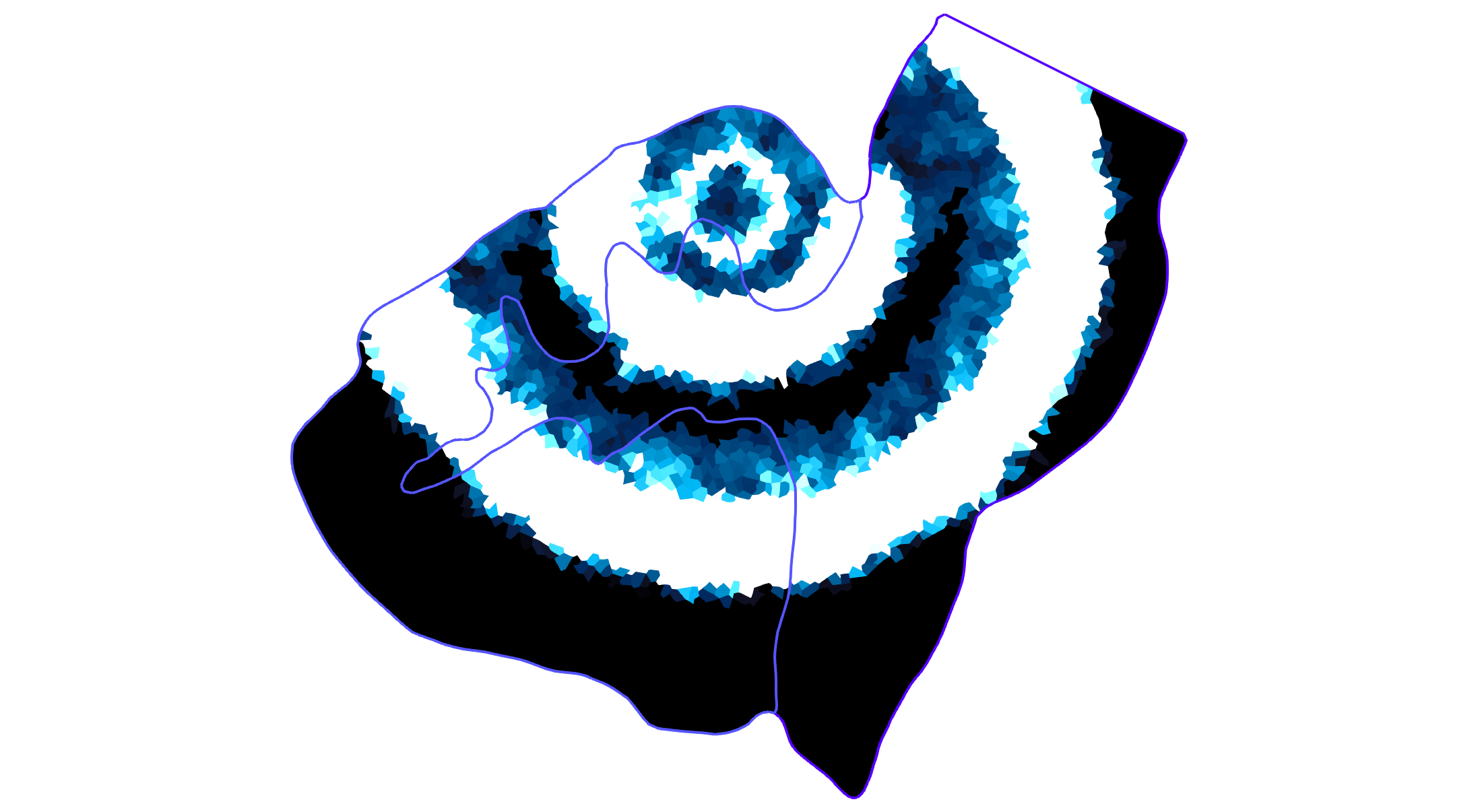}
\includegraphics[trim={12cm 1cm 13cm 1cm},clip, scale=0.085]{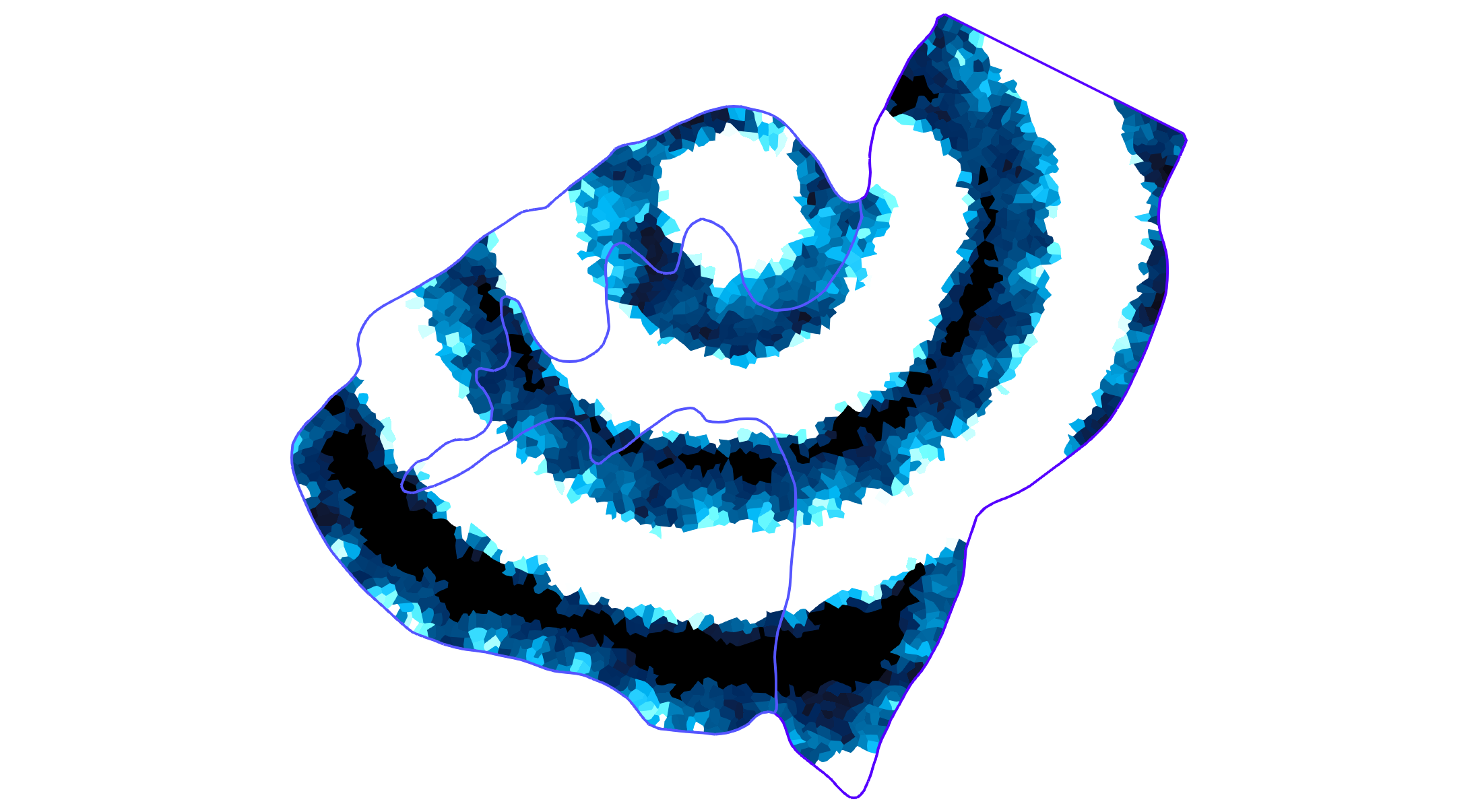}
      \subcaption{Local error indicator}  \label{fig:p_tranpot}
      \end{minipage}
  \caption{Test case 4. Computed transmembrane potential, element-wise polynomial approximation degree and local error indicator at three different time snapshots
$t=1.0, 9.4, 17.2 \, [\mathrm{ms}]$.}
  \label{fig::bc_waveStem}
\end{figure}
Figure \ref{fig::bc_waveStem} shows the computed transmembrane potential together with the element-wise polynomial approximation degree and the local error indicator at three different time snapshots
$t=1.0, 9.4, 17.2 \, [\mathrm{ms}]$. Also in this simulation, the proposed $p$-adaptive strategy is able to correctly capture the depolarization and polarization wavefronts and automatically and dynamically adapt the local polynomial degree to capture the steep variations of the transmembrane potential, while reducing the computational burden.
This is also confirmed by the results reported in Figure~\ref{fig:ndof_evolution_bcStem}, where we show the time-evolution of the $\mathrm{NDoF}$. We observe that the $p$-adaptive methods significantly reduced the $\mathrm{NDoF}$ compared to the corresponding scheme employing a uniform polynomial approximation degree (for which  $\mathrm{NDoF}=7.3e+4$).  After $14ms$, the $\mathrm{NDoF}$ drops to $4e+4$, reducing the total number of degrees of freedom by approximately $40\%$. Similar to the previous test case, we observe a gradual increase in the $\mathrm{NDoF}$ over time, driven by the appearance of multiple wavefronts that need to be properly tracked. The same trend can be observed in the evolution of the number of elements checked by the $p$-adaptive scheme and the number of mesh elements where the polynomial degree is equal to $1,2,3,4, \text{ and } 5$ (cf. Figure~\ref{fig:ndof_evolution_bcStem} center and right). Once the first wavefront has passed the whole domain, a stable distribution is reached both in terms of the total $\mathrm{NDoF}$ and the polynomial approximation degree distribution. 

\FloatBarrier
\section{Conclusion}
\label{sec:10}
\noindent This work addressed some numerical challenges in the numerical approximation of neuronal electrophysiology phenomena characterized by complex geometries, multiscale phenomena, and steep wavefronts. The considered mathematical model consists of the monodomain equation coupled with a specific neuronal ionic model to simulate the pathological behavior of the transmembrane potential. To overcome the associated computational burden, we proposed a novel $p$-adaptive algorithm based on the discontinuous Galerkin method on polytopal meshes coupled with Crank-Nicolson time marching scheme. Our automatic $p$-adaptive strategy is driven by a suitable a-posteriori error estimator to dynamically adjust the local polynomial approximation degree, enabling the dynamic and efficient capture of wavefronts while simultaneously reducing the total number of degrees of freedom. Thanks to a novel clustering algorithm, we automatically select a threshold value to efficiently identify elements for
adaptive updates, further improving efficiency. 
We tested the proposed $p$-adaptive algorithm on a wide set of benchmark test cases, demonstrating that the number of degrees of freedom can be significantly reduced while preserving the accuracy of the numerical solution.
We also demonstrated the capabilities of the proposed $p$-adaptive method in approximating physio-pathological scenarios of brain electrophysiology. The results highlighted the $p$-adaptive method's capability to capture the complex interplay of physiological and pathological phenomena at a reasonable computational cost.

\section*{Declaration of competing interests}
The authors declare that they have no known competing financial interests or personal relationships that could have appeared to influence the work reported in this article.

\section*{Acknowledgments}
The brain MRI images were provided by OASIS-3: Longitudinal Multimodal Neuroimaging: Principal Investigators: T. Benzinger, D. Marcus, J. Morris; NIH P30 AG066444, P50 AG00561, P30 NS09857781, P01 AG026276, P01 AG003991, R01 AG043434, UL1 TR000448, R01 EB009352. AV-45 doses were provided by Avid Radiopharmaceuticals, a wholly-owned subsidiary of Eli Lilly.
\bibliographystyle{hieeetr}
\bibliography{main.bib}

\end{document}